\documentclass[11pt]{article}
 \usepackage{amscd,amsfonts,amssymb,amsthm,amscd,amsmath,latexsym}
  \usepackage[hypertex]{hyperref}
\usepackage{makeidx}
\usepackage[all]{xy}
\makeatletter
\makeindex
\renewcommand{\subsubsection}{\@startsection
{subsubsection}
{1}
{0mm}
{0mm}
{0mm}
{\normalfont\normalsize\itshape}}

\makeatother
\begin{titlepage}

\title{Index theory, eta forms, and Deligne cohomology}
\author{Ulrich Bunke\thanks{Mathematisches Institut, Universit{\"a}t G{\"o}ttingen,
Bunsenstr. 3-5, 37073 G{\"o}ttingen, GERMANY,  bunke@uni-math.gwdg.de} }
\end{titlepage}
\begin{document}
\maketitle
\vspace{9cm}
MSC :\\
 58j28 (primary), \\
55S35 (secondary)

\newpage
\begin{abstract}
 \begin{enumerate}
\item The paper sets up a language to deal with Dirac operators on manifolds with corners of arbitrary codimension. In particular we develop a precise theory of boundary reductions.
\item We introduce the notion of a taming of a Dirac operator as an invertible  perturbation by a smoothing operator. Given a Dirac operator on a manifold with boundary faces
we use the tamings of its  boundary reductions in order to turn the operator into a Fredholm operator.
Its index is an obstruction against extending the taming from the boundary to the interior. 
In this way we develop an inductive procedure
to associate Fredholm operators to Dirac operators on manifolds with corners and develop the associated obstruction theory.
\item A central problem of index theory is to calculate the Chern character of the index of a family of Dirac operators. Local index theory uses the heat semigroup of an associated super-connection in order to produce differential forms representing this Chern character.  
In this paper we  develop a version of local index theory for families of Dirac operators on manifolds with corners. The resulting de Rham representative of the Chern character  is a sum of the local index form and $\eta$-form contributions from the boundary faces. If the index of the family vanishes and we have chosen a taming, then local index theory in addition gives a transgression form whose differential trivializes this Rham representative.
This transgression form plays an important role in the construction of secondary invariants.
\item  
Assume that the $K$-theoretic index of a family of Dirac operators (on a family of closed manifolds) vanishes on all $i-1$-dimensional subcomplexes of the parameter space.
The obstruction against increasing $i$ by one is an 
 $i$-dimensional integral cohomology class.
One of the main goals of this paper is to use the additional information given by local index theory in order to refine this obstruction class to a class in $i$-th integral Deligne cohomology.
As a byproduct we get a lift of the $i$-th Chern class of the index of a family of Dirac operators 
to Deligne cohomology. 
\item In low degrees $\le 3$ integral Deligne cohomology classifies well-known geometric objects
like $\mathbb{Z}$-valued functions, $U(1)$-valued smooth functions, hermitean line bundles with connections and geometric gerbes. Such objects have been previously associated to families of Dirac operators.
We verify that  these constructions are  compatible with our definitions. 
\end{enumerate}

\end{abstract}
\newpage

\newcommand{\Fr}{{\tt Fr}}
\newcommand{\Fred}{{\tt Fred}}
\newcommand{\Comp}{{\tt K}}
\newcommand{\vC}{\v{C}}

\newcommand{\DD}{\mathbf{v}}
\newcommand{\Z}{\mathbb{Z}}
\newcommand{\orient}{{\rm or}}
\newcommand{\cSet}{\mathcal{S}et}
\newcommand{\by}{{\bf y}}
\newcommand{\bE}{{\bf E}}
\newcommand{\Face}{{\tt Face}}
\newcommand{\cDelta}{\mathbf{\Delta}}
\newcommand{\LIM}{{\tt LIM}}
\newcommand{\diag}{{\tt diag}}
 \newcommand{\dist}{{\tt dist}}
\newcommand{\kaaa}{{\frak k}}
\newcommand{\paaa}{{\frak p}}
\newcommand{\vp}{{\varphi}}
\newcommand{\taaa}{{\frak t}}
\newcommand{\haaa}{{\frak h}}
\newcommand{\R}{\mathbb{R}}
\newcommand{\Hh}{{\bf H}}
\newcommand{\Rep}{{\tt Rep}}
\newcommand{\Hb}{mathbb{H}}
\newcommand{\Q}{\mathbb{Q}}
\newcommand{\str}{{\tt str}}
\newcommand{\Ind}{{\tt ind}}
\newcommand{\triv}{{\tt triv}}
\newcommand{\bD}{{\bf D}}
\newcommand{\bF}{{\bf F}}
\newcommand{\tX}{{\tt X}}
\newcommand{\Cliff}{{\tt Cl}}
\newcommand{\tY}{{\tt Y}}
\newcommand{\tZ}{{\tt Z}}
\newcommand{\tV}{{\tt V}}
\newcommand{\tR}{{\tt R}}
\newcommand{\Fam}{{\tt Fam}}
\newcommand{\Cusp}{{\tt Cusp}}
\newcommand{\bT}{{\bf T}}
\newcommand{\bK}{{\bf K}}
\newcommand{\bo}{{\bf o}}
\newcommand{\K}{\mathbb{K}}
\newcommand{\tH}{{\tt H}}
\newcommand{\bS}{\mathbf{S}}
\newcommand{\bB}{\mathbf{B}}
\newcommand{\tW}{{\tt W}}
\newcommand{\tF}{{\tt F}}
\newcommand{\bA}{\mathbf{A}}
\newcommand{\bL}{{\bf L}}
 \newcommand{\bom}{\mathbf{\Omega}}
\newcommand{\bundle}{{\tt bundle}}
\newcommand{\ch}{\mathbf{ch}}
\newcommand{\ve}{{\varepsilon}}
\newcommand{\C}{\mathbb{C}}
\newcommand{\gen}{{\tt gen}}
\newcommand{\cTop}{\mathcal{T}op}
\newcommand{\bP}{\mathbf{P}}
\newcommand{\Naaa}{\mathbf{N}}
\newcommand{\image}{{\tt image}}
\newcommand{\gaaa}{{\frak g}}
\newcommand{\zaaa}{{\frak z}}
\newcommand{\saaa}{{\frak s}}
\newcommand{\laaa}{{\frak l}}
\newcommand{\bN}{\mathbf{N}}
\newcommand{\stimes}{{\times\hspace{-1mm}\bf |}}
\newcommand{\ausg}{{\rm end}}
\newcommand{\bff}{{\bf f}}
\newcommand{\maaa}{{\frak m}}
\newcommand{\aaaa}{{\frak a}}
\newcommand{\naaa}{{\frak n}}
\newcommand{\brr}{{\bf r}}
\newcommand{\res}{{\tt res}}
\newcommand{\Aut}{{\tt Aut}}
\newcommand{\Pol}{{\tt Pol}}
\newcommand{\Tr}{{\tt Tr}}
\newcommand{\cT}{\mathcal{T}}
\newcommand{\dom}{{\tt dom}}
\newcommand{\Line}{{\tt Line}}
\newcommand{\db}{{\bar{\partial}}}
\newcommand{\Sf}{{\tt  Sf}}
\newcommand{\g}{{\gaaa}}
\newcommand{\cZ}{\mathcal{Z}}
\newcommand{\cH}{\mathcal{H}}
\newcommand{\cM}{\mathcal{M}}
\newcommand{\interi}{{\ttint}}
\newcommand{\singsupp}{{\tt singsupp}}
\newcommand{\cE}{\mathcal{E}}
\newcommand{\ccR}{\mathcal{R}}
\newcommand{\hol}{{\tt hol}}
\newcommand{\cV}{\mathcal{V}}
\newcommand{\cY}{\mathcal{Y}}
\newcommand{\cW}{\mathcal{W}}
\newcommand{\dR}{{\tt dR}}
\newcommand{\del}{{\tt del}}
\newcommand{\bdel}{\mathbf{del}}
\newcommand{\cI}{\mathcal{I}}
\newcommand{\cC}{\mathcal{C}}
\newcommand{\cK}{\mathcal{K}}
\newcommand{\cA}{\mathcal{A}}
\newcommand{\cU}{\mathcal{U}}
\newcommand{\Hom}{{\tt Hom}}
\newcommand{\vol}{{\tt vol}}
\newcommand{\cO}{\mathcal{O}}
\newcommand{\End}{{\tt End}}
\newcommand{\Ext}{{\tt Ext}}
\newcommand{\rk}{{\tt rank}}
\newcommand{\im}{{\tt im}}
\newcommand{\sign}{{\tt sign}}
\newcommand{\spann}{{\tt span}}
\newcommand{\symm}{{\tt symm}}
\newcommand{\cF}{\mathcal{F}}
\newcommand{\cD}{\mathcal{D}}
\newcommand{\bC}{\mathbf{C}}
\newcommand{\bbeta}{\mathbf{\eta}}
\newcommand{\bOmega}{\mathbf{\Omega}}
\newcommand{\bbz}{{\bf z}}
\newcommand{\bc}{\mathbf{c}}
\newcommand{\bb}{\mathbf{b}}
\newcommand{\bd}{\mathbf{d}}
\newcommand{\Ree}{{\tt Re }}
\newcommand{\Res}{{\tt Res}}
\newcommand{\Imm}{{\tt Im}}
\newcommand{\inter}{{\tt int}}
\newcommand{\clo}{{\tt clo}}
\newcommand{\tg}{{\tt tg}}
\newcommand{\ee}{{\tt e}}
\newcommand{\Li}{{\tt Li}}
\newcommand{\cN}{\mathcal{N}}
 \newcommand{\conv}{{\tt conv}}
\newcommand{\op}{{\tt Op}}
\newcommand{\tr}{{\tt tr}}
\newcommand{\ctg}{{\tt ctg}}
\newcommand{\degg}{{\tt deg}}
\newcommand{\Ad}{{\tt  Ad}}
\newcommand{\ad}{{\tt ad}}
\newcommand{\codim}{{\tt codim}}
\newcommand{\Gr}{{\tt Gr}}
\newcommand{\coker}{{\tt coker}}
\newcommand{\id}{{\tt id}}
\newcommand{\ord}{{\tt ord}}
\newcommand{\nat}{\mathbb{N}}
\newcommand{\supp}{{\tt supp}}
\newcommand{\sing}{{\tt sing}}
\newcommand{\spec}{{\tt spec}}
\newcommand{\Ann}{{\tt Ann}}
 \newcommand{\Or}{{\tt Or }}
\newcommand{\Diff}{\mathcal{D}iff}
\newcommand{\cB}{\mathcal{B}}
\def\imath{{i}}
\newcommand{\cR}{\mathcal{R}}
\def\hB{\hspace*{\fill}$\Box$ \newline\noindent}
\newcommand{\varho}{\varrho}
\newcommand{\ind}{{\tt index}}
\newcommand{\Indu}{{\tt Ind}}
\newcommand{\Fin}{{\tt Fin}}
\newcommand{\cS}{\mathcal{S}}
\newcommand{\orig}{\mathcal{O}}
\def\hB{\hspace*{\fill}$\Box$ \\[0.2cm]\noindent}
\newcommand{\cL}{\mathcal{L}}
 \newcommand{\cG}{\mathcal{G}}
\newcommand{\Mat}{{\tt Mat}}
\newcommand{\cP}{\mathcal{P}}
\newcommand{\bv}{\mathbf{v}}
\newcommand{\cQ}{\mathcal{Q}}
 \newcommand{\cX}{\mathcal{X}}
\newcommand{\bH}{\mathbf{H}}
\newcommand{\bW}{\mathbf{W}}
\newcommand{\pr}{{\tt pr}}
\newcommand{\bX}{\mathbf{X}}
\newcommand{\bY}{\mathbf{Y}}
\newcommand{\bZ}{\mathbf{Z}}
\newcommand{\bz}{\mathbf{z}}
\newcommand{\bkappa}{\mathbf{\kappa}}
\newcommand{\ev}{{\tt ev}}
\newcommand{\bV}{\mathbf{V}}
\newcommand{\Gerbe}{{\tt Gerbe}}
\newcommand{\gerbe}{{\tt gerbe}}
\newcommand{\hA}{\mathbf{\hat A}}
\newcommand{\Sets}{\mathcal{S}ets}

\newtheorem{prop}{Proposition}[section]
\newtheorem{lem}[prop]{Lemma}
\newtheorem{ddd}[prop]{Definition}
\newtheorem{theorem}[prop]{Theorem}
\newtheorem{kor}[prop]{Corollary}
\newtheorem{ass}[prop]{Assumption}
\newtheorem{con}[prop]{Conjecture}
\newtheorem{prob}[prop]{Problem}
\newtheorem{fact}[prop]{Fact}

\setcounter{tocdepth}{2}
 
 \tableofcontents
\parskip1ex
\newpage
\part{} 

\section{Introduction}

\subsection{The index of families of Dirac operators}

\subsubsection{}
Let $M$ be a smooth manifold.
A \index{generalized Dirac operator} generalized Dirac operator $D$ on $M$ is a first order differential
operator acting on sections of a complex vector bundle $V\rightarrow M$.
It is characterized amongst all first order differential operators
by the property that the symbol of its square
has the form
$$\sigma(D)(\xi)=g^M(\xi,\xi) \id + O(\xi)\ ,\quad \xi\in T^*M\ ,$$
where $g^M$ is a Riemannian metric on the underlying manifold.

\subsubsection{}
In the framework of index theory the operators have an additional
symmetry. We assume that the bundle $V$ has a hermitian metric.
Then we can define a $L^2$-scalar product between 
compactly supported sections of $V$. It is generally assumed that
$D$ is formally selfadjoint, i.e. it is symmetric on the space of
smooth sections with compact support in the interior of $M$.

If the dimension of $M$ is even, then in addition we require that
$V$ has a selfadjoint  involution $z\in \End(V)$ (i.e. a
$\Z/2\Z$-grading) which anti-commutes
with $D$. Then we can decompose $V=V^+\oplus V^-$
into the $\pm 1$-eigenspaces of $z$ and write
$$D=\left(\begin{array}{cc} 0&D^-\\D^+&0\end{array}\right)\ .$$

\subsubsection{}

Assume that $M$ is even-dimensional and closed.
Then $D^+:C^\infty(M,V^+)\rightarrow C^\infty(M,V^-)$
has a finite dimensional kernel and cokernel. By definition
$$\ind(D):=\dim(\ker(D^+))-\dim(\ker(D^-))\ .$$
This number  can also be written as \index{index!of a Dirac operator}
$$\ind(D)= \Tr_s P\ ,$$
where $P$ is the orthogonal projection onto the kernel of $D$ and
$\Tr_sA :=\Tr z A$.

The question of classical index theory is to compute
$\ind(D)\in\Z$ in terms of the symbol of $D$. It was solved by the
index theorem of Atiyah-Singer \cite{atiyahsinger63}.

\subsubsection{}

Let $B$ be some auxiliary compact topological space.
Let us consider a family $(D_b)_{b\in B}$ of Dirac operators
which is continuously parameterized by $B$. 

Assume that $M$ is compact and odd-dimensional.
Then we can form the family
$(F_b)_{b\in B}$ of selfadjoint Fredholm operators on $L^2(M,V)$, where
$F_b:=D_b(D_b^2+1)^{-1/2}$ is defined by functional calculus.
The family $(F_b)_{b\in B}$ is not continuous in the norm topology of bounded
operators. But for all $\psi\in L^2(M,V)$ the family
$(F_b\psi)_{b\in B}$ is a continuous family of vectors in the Hilbert space, and the family
$(1-F_b^2)_{b\in B}$ is a norm continuous family of 
compact operators.

Note that $F_b$ has infinite dimensional positive and negative eigenspaces.

\subsubsection{}

If $H$ is a separable Hilbert space, then we can consider the space
$\bK^1$ of all selfadjoint Fredholm operators $F$
such that $1-F^2$ is compact and $F$ has infinite dimensional positive
and negative eigenspaces. We equip this space with
the smallest topology such that for all $\psi\in H$ the families
$\bK^1\ni F\mapsto F\psi\in H$,  
and the family $\bK^1\ni F\rightarrow 1-F^2$ are norm continuous.

One can show that $\bK^1$ has the homotopy type of the classifying
space of the complex $K$-theory functor $K^1$.

\subsubsection{}

If $M$ is closed and odd-dimensional, then our family $(D_b)_{b\in B}$
gives rise to a continuous map $F:B\rightarrow \bK^1$
and therefore to a homotopy class
$$\ind((D_b)_{b\in B})=[F]\in [B,\bK^1]=K^1(B)\ .$$

\subsubsection{}
 
Let $H$ be a $\Z/2\Z$-graded separable Hilbert space.
We consider the space $\bK^0$ of all selfadjoint Fredholm operators $F$
which are  odd and   such that $F^2-1$ is
compact. In order to define the topology we consider $\bK^0$ as a
subset of $\bK^1$.  We then equip this subset with the induced topology.

Again one can show that $\bK^0$ has the homotopy type of the
classifying space of the complex $K$-theory functor $K^0$.

\subsubsection{}\label{uuh21}
If $M$ is even-dimensional and we set $F_b:=D_b(D_b^2+1)^{1/2}$ as before,
then $F_b\in \bK^0$.
The family $(F_b)_{b\in B}$ gives rise to a continuous map
$F:B\rightarrow \bK^0$ and therefore to a homotopy class
$$\ind((D_b)_{b\in B})=[F]\in [B,\bK^0]=K^0(B)\ .$$

\subsubsection{}

One issue which we  have suppressed here is that this definition
involves an unitary identification of $H$ with $L^2(M,V)$. Note that
by Kuiper's theorem the space of such unitary identifications is
contractible so that the construction above is independent of the
choice.

In fact, the scalar product on the Hilbert space $L^2(M,V)$ in general also
depends on $b\in B$ since the volume measure depends on the Riemannian
metric on $M$ which is determined by the symbol of $D$.
So what we must in fact choose is a trivialization of the bundle of
Hilbert spaces $(L^2(M,V,<.,.>_b))_{b\in B}$ which exists and is again
unique up to homotopy by Kuiper's theorem.

Arrived at this point we see that the construction above has the
following generalization. We consider the family of generalized Dirac operators
$(D_b)_{b\in B}$ as a family of fiber-wise differential operators
on the trivial fiber bundle $B\times M\rightarrow B$.
It is now straight forward to generalize the construction of the index
to the case of a family of fiber-wise generalized  Dirac operators on a merely
locally trivial bundle $E\rightarrow B$ with fiber $M$.

\subsubsection{}

Let $(D_b)_{b\in B}$  be family of fiber-wise generalized Dirac operators
on a fiber bundle $E\rightarrow B$ with even-dimensional closed fibers and
compact base $B$.
 After a perturbation of the family we can assume that
$\dim(\ker(D_b))$ is independent of $b\in B$.
In this case the family of vector spaces $(\ker(D_b))_{b\in B}$ forms
a $\Z/2\Z$-graded vector bundle $\ker(D)$ over the base $B$.

If one considers $K^0(B)$ as the Grothendieck group generated by
isomorphism
classes of vector bundles over $B$, then the class
$[\ker(D)]\in K^0(B)$ corresponds to
the index $\ind((D_{b})_{b\in B})\in K^0(B)$ as defined in \ref{uuh21}
under the usual identification of the two pictures
of the $K^0$-functor.

\subsubsection{}\label{taming22}

The most natural interpretation of the index is as an obstruction
class.
Let  $(D_b)_{b\in B}$  be a family of fiber-wise generalized Dirac operator
on a fiber bundle $E\rightarrow B$ with closed nonzero-dimensional fibers over a compact
base $B$.
Then it is a natural question if there exists a family
$(Q_b)_{b\in Q}$ of selfadjoint integral operators with smooth integral kernels
(which are odd with respect to the $\Z/2\Z$-grading in the
even-dimensional case) such that the perturbed family
$(D_b+Q_b)_{b\in B}$ is invertible for every $b$.
We call such a family $(Q_b)_{b\in Q}$ a taming.

Then we have the following assertion (see Lemma \ref{fambtt}): {\em  The family 
$(D_b)_{b\in B} $ admits a taming if and only if $\ind((D_{b})_{b\in
  B})=0$.}

\subsection{Local index theory for families}

\subsubsection{}
Let $D$ be a generalized Dirac operator on an even dimensional closed manifold.
By the McKean-Singer formula we can write
$$\ind(D)=\Tr_s \ee^{-t^2D^2}\ ,$$
where $t>0$.
The heat operator \index{heat operator} $\ee^{-t^2D^2}$ has a smooth integral kernel
$\ee^{-t^2D^2}(x,y)$ so that we can express the trace by an integral  
$$\Tr_s \ee^{-t^2D^2}=\int_M \tr_s \ee^{-t^2D^2}(x,x)\ .$$
Here $\ee^{-t^2D^2}(x,x)\in\End(V_x)\otimes \Lambda_x$, where $V_x$ and
$\Lambda_x$ are the fibers of $V$ and the density bundle
$\Lambda\rightarrow M$
over $x\in M$, and $\tr_s$ is the super trace on
$\End(V_x)$. 

It is well known  \index{local!trace}(see
\cite{berlinegetzlervergne92}, Ch. 2.) that the local trace of the heat operator admits an
asymptotic expansion
\begin{equation}\label{expan}\tr_s \ee^{-t^2D^2}(x,x)\stackrel{t\to 0}{\sim}\sum_{k\ge 0}
t^{2k-n} a_{2k-n}(x)\ ,\end{equation}
where $n=\dim(M)$ and the coefficients $a_{2k-n}(x)$ are locally determined by the
operator
$D$.
Thus we can write
$$\ind(D)=\int_M a_0(x)\ .$$
This reduces the computation of the index to the determination of
the coefficient $a_{0}$ in the local heat trace asymptotic.

\subsubsection{}\label{looi7}
 
\index{$a_0$}
The determination of $a_{0}$ is particularly easy for compatible Dirac
operators. These are Dirac operators which are associated to a Dirac
bundle structure $\cV$ on $V$ (see \ref{diracbbb} for a definition).
For compatible Dirac operators we have $a_{k}=0$ for $k<0$, and $a_0$
is given as follows.

Assume that $M$ is oriented so that the density bundle is trivialized
and $a_0$ is a highest degree form on $M$.
If $M$ admits a spin structure,
then the Dirac bundle $\cV$ is isomorphic to a twisted spinor bundle
$\cS(M)\otimes \bW$, where $\cS(M)$ is a spinor bundle of $M$ and $\bW=(W,h^W,\nabla^W)$
\index{$\hA(\nabla^{TM})$}
\index{$\nabla^W$}
is an auxiliary hermitian $\Z/2\Z$-graded hermitian vector bundle
with connection called the twisting bundle.
Under these assumptions we have the equality (see
\cite{berlinegetzlervergne92}, Ch. 4.)
$$a_0=[\hA(\nabla^{TM}) \ch(\nabla^{W})]_n\ ,$$
a formula which us usually called the {\em local index theorem}.
The forms on the right-hand side are the Chern-Weyl representatives of
the corresponding characteristic classes of $TM$ and $W$ associated to
the Levi-Civita connection $\nabla^{TM}$ and $\nabla^W$.

In the general case the assumptions above are satisfied locally on $M$
 so that $a_0$ can be determined locally be the same formula.

\subsubsection{}
Let us now consider a {\em smooth} family of fiber-wise Dirac operators
\index{$\ch$}\index{Chern!character}$(D_b)_{b\in B}$ on a smooth fiber bundle $\pi:E\rightarrow B$
with even-dimensional fibers.
Then the heat kernel method above can be generalized in order to
compute a de Rham representative of the Chern character of the index
of the family. Thus let $\ch:K^*(B)\rightarrow H^*(B,\Q)$ be the Chern character, 
and $\dR: H^*(B,\Q)\rightarrow H_{dR}(B)$ be the de Rham map.

The main idea is due to Quillen and known under the name {\em super-connection formalism} (our general reference for all that is \cite{berlinegetzlervergne92}). If we fix a horizontal distribution
$T^h\pi\subseteq TE$, i.e. a complement of the vertical bundle
$T^v\pi=\ker(d\pi)$, and a connection on $V$, then we obtain an
unitary connection $\nabla^u$ on the bundle of Hilbert spaces
$(L^2(E_b,V_{|E_b}))_{b\in B}$. \index{super connection}
We  define the family of super-connections
$$S_t:=tD+\nabla^u\ .$$ For $t>0$
the curvature $$S_t^2=t^2D^2 +  higher\:degree\: forms$$ is a differential form on $B$ with
values in the fiber-wise differential operators.
Its exponential is a form on the base with coefficients in the fibre-wise smoothing operators.
Thus $\Tr_s \ee^{-S^2_t}$ is a differential form on $B$.

The generalization of the McKean-Singer formula asserts now that for all $t>0$
\begin{eqnarray}
d \Tr_s \ee^{-S^2_t}&=&0 \label{trans32}\ ,\\
(2\pi i)^{-\deg/2}
[\Tr_s \ee^{-S^2_t}]&=&\dR(\ch(\ind((D_b)_{b\in B}))) \nonumber\ ,
\end{eqnarray}
where $\deg$ is the $\Z$-grading operator on differential forms on
$B$.

\subsubsection{}

The integral kernel of $\ee^{-S^2_t}$ again has an asymptotic
expansion of the form (\ref{expan}) with locally determined
coefficients $a_k$ which are now differential forms on $B$ with values
in fiber-wise densities.

Compared with the case of a single operator the situation is now more
complicated because of the following. The differential form
$\Tr_s \ee^{-S^2_t}$ depends on $t$, but we have a transgression
\index{transgression formula} formula which is formally a consequence of (\ref{trans32})
\begin{equation}\label{u761}\Tr_s \ee^{-S^2_t}-\Tr_s \ee^{-S^2_s}=-d \int_s^t \Tr_s
 \frac{\partial S_u}{\partial t} \ee^{-S_u^2} du\ .\end{equation} 
Thus in order to use the local asymptotic expansion of the heat kernel
in order to compute a de Rham representative of the Chern character of
the index we would like to require that the limit  $s\to 0$ of the integrals 
in (\ref{u761}) converges.
 
\subsubsection{}

\index{Bismut super-connection}As observed by Bismut this is the case for families of compatible
Dirac operators if one modifies the super-connection to the Bismut
super-connection
$$A_t:=tD+\nabla^u+\frac{1}{4t} c(T)\ ,$$
where $c(T)$ is the Clifford multiplication with the curvature of the
horizontal distribution which can be considered as a two form on $B$
with values in the vertical vector fields.

In this case the limit
$$\lim_{t\to 0} (2\pi i)^{-\deg/2}\Tr_s \ee^{-A^2_t}=:\Omega(\cE_{geom})$$
exists and defines the local index form \index{local!index form}
$\Omega(\cE_{geom})$.
Here $\cE_{geom}$ is our notation for a geometric family which is just
a shorthand for the collection of data needed to define the Bismut
super-connection (see Definition \ref{geo987}).
The following equality is the {\em local index theorem for families}
$$[\Omega(\cE_{geom})]=\dR(\ch(\ind(\cE_{geom})))\ ,$$
where $\ind(\cE_{geom})$ is our notation for the index of the family
of Dirac operators associated with $\cE_{geom}$.

\subsubsection{}
The local index form can be determined in a similar way as in
\index{$\Omega(\cE_{geom})$}\ref{looi7}. First assume that the vertical tangent bundle $T^v\pi$
has a spin structure. Then we can write the family of Dirac bundles as
twisted spinor bundle $\cV=\cS(T^v\pi)\otimes \bW$.
In this case we have
$$\Omega(\cE_{geom})=\int_{E/B} \hA(\nabla^{T^v\pi})\ch(\nabla^{W})\
,$$
where the connection $\nabla^{T^v\pi}$ is induced by the data of the
geometric family (see \cite{berlinegetzlervergne92}, Ch. 10). In the general case the spin assumption is satisfied
locally on $E$ so that we can obtain the integrand for the local index form by the same formula.

\subsubsection{}

If $\ind(\cE_{geom})=0$, then the form $\Omega(\cE_{geom})$ is exact.
The main idea of secondary index theory is to find a reason why
the index is trivial and to employ this reason in order to define a
form $\alpha$ such that $d\alpha=\Omega(\cE_{geom})$.

In the present paper the reason for   $\ind(\cE_{geom})=0$ is that
$\cE_{geom}$ admits a taming (see \ref{taming22}).
Let $\cE_t$ be a geometric family with a choice of a taming.
Using the taming $Q$  we further modify the super-connection
to 
$$A_t(\cE_t):=tD+\nabla^u+\frac{1}{4t} c(T) + t \chi(t) Q\ ,$$
where $ \chi$ is a cut-off function which vanishes for $t\le 1$ and is
equal to one for $t\ge 2$.
This modification has the following effects. For small $t$ the
modified
super-connection is the Bismut super-connection so that we can use the
knowledge about the small $t$-behavior.
For large $t$ the zero-form part is the invertible operator
$t(D+Q)$. This has the effect that $\Tr_s \ee^{-A^2_t(\cE_t)}$
and the integrand $\Tr_s
 \frac{\partial A_t(\cE_t)}{\partial t} \ee^{-A_t^2(\cE_t)}$ in (\ref{trans32})
vanish exponentially for large $t$.

\index{$\eta$-form}\index{$\eta(\cE_t)$}We can define the $\eta$-form
$$\eta(\cE_t):= (2\pi i)^{-\deg/2} \int_0^\infty \Tr_s
 \frac{\partial A_u(\cE_t)}{\partial u} \ee^{-A_u^2(\cE_t)} du\ .$$
Then we have
$$d\eta(\cE_{t})=\Omega(\cE_{geom})\ .$$

Note that the eta form depends on the taming. In fact, if
$\cE^\prime_t$ is given by a second taming, then $\eta(\cE_t)-\eta(\cE_t^\prime)$ is closed
and represents the Chern character of an element of $K^1(B)$
which measures the difference of the two tamings (see Corollary \ref{etajump}).
 
\subsubsection{}

The picture described above for the even-dimensional case has a
analogous odd-dimensional counterpart. 

\subsection{Absolute and relative secondary invariants}

\subsubsection{}
To any geometric family $\cE_{geom}$ we can associate its opposite
$\cE_{geom}^{op}$ by switching orientations and gradings. In
particular we have
$$\Omega(\cE_{geom}^{op})=-\Omega(\cE_{geom})\ ,\quad \ind(\cE_{geom}^{op})=-\ind(\cE_{geom})\ .$$

Let us fix $x\in K(B)$. Assume that $\cE_{geom,i}$, $i=0,1$, are
are two (nonzero-dimensional) geometric families such that $\ind(\cE_{geom,i})=x$.
Then the difference of local index forms
$\Omega(\cE_{geom,1})-\Omega(\cE_{geom,0})$ is exact.
The reason is that the fiber-wise sum $\cF_{geom}:=\cE_{geom,0}^{op}+\cE_{geom,1}$
has trivial index and thus admits a taming $\cF_t$.
We therefore can write
\begin{equation}\label{transss1}
\Omega(\cE_{geom,1})-\Omega(\cE_{geom,0})=d\eta(\cF_t)\ .
\end{equation}
Note that if we change the taming, then
$\eta(\cF_t)$ changes by a closed form with rational periods.

\subsubsection{}

The class $x\in K^*(B)$, $*\in \Z/2\Z$, determines a rational cohomology class
$\ch(x)\in H^{*}(B,\Q)$ \footnote{In the following we consider cohomology and
forms as $\Z/2\Z$-graded.}. If we represent $x$ as the index of
a fixed geometric family $\cE_{geom,0}$, then we obtain the additional
information of a de Rham representative $\Omega(\cE_{geom,0})$ of
$\ch(x)$. If we take another representative $\cE_{geom,1}$, then have
the secondary information $\eta(\cF_t)$ such that $(\ref{transss1})$
holds true. 

Thus after fixing $\cE_{geom,0}$, we have a relative invariant
defined on the set of all geometric families with index $x$ which
takes values in the quotient
$\cA_B^{*-1}(B)/\cA^{*-1}_B(B,d=0,\Q)$ of all differential forms on $B$ by closed
forms with rational periods. 

\subsubsection{}\label{gv21}

One may ask if one can turn this relative invariant into an absolute one.
\index{Deligne cohomology!rational}
\index{$H^*_{Del}(B,\Q)$}In fact there is a group\footnote{This in fact a ring valued functor, but the
  ring structure is not important in the present paper.} valued
functor $H^*_{Del}(B,\Q)$ called rational Deligne cohomology
(see Definition \ref{delcoh}) which may capture this kind of
information. Rational Deligne cohomology comes with two natural transformations
$$H^*_{Del}(B,\Q)\stackrel{R}{\rightarrow} \cA^*_B(B,d=0)\ ,\quad
H^*_{Del}(B,\Q)\stackrel{\DD}{\rightarrow} H^*(B,\Q)\ ,$$ which are
called the curvature and the characteristic \index{characteristic class}class,
such that $[R^u]=\dR(\DD(u))$ for all $u\in H^*_{Del}(B,\Q)$.
Furthermore, there is a natural transformation
$$\cA^{*-1}_B(B)\stackrel{a}{\rightarrow} H^*_{Del}(B,\Q)$$
such that $R^{a(\alpha)}=d\alpha$
and
$$0\rightarrow \cA^{*-1}(B)/\cA^{*-1}_{B}(B,d=0,\Q)\rightarrow 
 H^*_{Del}(B,\Q)\stackrel{\DD}{\rightarrow } H^*(B,\Q)\rightarrow 0$$
is exact, where the first map is induced by $a$.

If we fix an element $\ind_{Del,\Q}(\cE_{geom,0})\in  H^*_{Del}(B,\Q)$
such that $R^{\ind_{Del,\Q}(\cE_{geom,0})}=\Omega(\cE_{geom})$, then we
would obtain an invariant \index{$\ind_{Del,\Q}$}
$\ind_{Del,\Q}(\cE_{geom,1})\in  H^*_{Del}(B,\Q)$ for all geometric families
$\cE_{geom,1}$ with $\ind(\cE_{geom,1})=x$ 
by the prescription 
$$\ind_{Del,\Q}(\cE_{geom,1}):=\ind_{Del}(\cE_{geom,0})+a([\eta(\cF_t)])\
.$$

This definition has the drawback of being
not natural with respect to pull-back.
In fact, if  $f:B^\prime\rightarrow B$ is a smooth map, then in general 
$f^*\ind_{Del,\Q}(\cE_{geom,0})\not=\ind_{Del,\Q}(f^*\cE_{geom,0})$
since equality would be by an accidental choice.

It is one of the main objectives of the forthcoming paper \cite{bunkeschick031}
to show how this approach can be modified in order to define a
natural $\ind_{Del,\Q}(\cE_{geom})\in H_{Del}^*(B)$.

\subsection{Integral secondary invariants}

\subsubsection{}\label{gg77}

The group $K(B)$ admits a natural decreasing filtration 
$$\dots\subseteq  K_{n}(B)\subseteq K_{n-1}\subseteq \dots \subseteq K_{0}(B)=K(B)\
.$$\index{Atiyah-Hirzebruch!filtration}\index{$K_p^*(X)$}\index{filtration!of $K$-theory}
By definition $x\in K_n(B)$, if $f^*x=0$ for all $n-1$-dimensional
$CW$-complexes $A$ and continuous maps $f:A\rightarrow B$.
In the present paper we index the Chern classes by their degrees.
The odd-degree classes are related to the even degree ones by
suspension (see \ref{fgr1}).

Let $k,m\in\nat$ be such that $k=2m$ or $k=2m-1$.
If $x\in K_k(B)$, then
$c_l(x)=0$ for all $l<k$. This implies
that
$$c^\Q_k(x)=(-1)^{m-1} (m-1)! \ch_k(x)\ .$$
Thus if $x$ is the index of a geometric family $\cE_{geom,0}$,
then 
$$[(-1)^{m-1}(m-1)!\Omega^k(\cE_{geom,0})]= \dR(c^\Q_k(x))\ .$$
In particular, this multiple of the local index form has integral periods.
This leads to a relative secondary invariant as follows.
If $\cE_{geom,1}$ is another family with index $x$,
then we have
$$(-1)^{m-1}(m-1)!\Omega^k(\cE_{geom,1})-(-1)^{m-1}(m-1)!\Omega^k(\cE_{geom,0})=(-1)^{m-1}(m-1)!d\eta^{k-1}(\cF_t)\
.$$

\subsubsection{}\label{obbs}
Observe that for any $u\in K(B)$ the rational class
$(-1)^{m-1}(m-1)!\ch_{k-1}(u)$ has integral periods.
Therefore using Corollary \ref{etajump}, after fixing  $\cE_{geom,0}$, we can define
the relative invariant
with values in $\cA^{k-1}_B(B)/\cA^{k-1}_B(B,d=0,\Z)$
such that it associates to the family
$\cE_{geom,1}$ the class
$[(-1)^{m-1}(m-1)!\eta^{k-1}(\cF_t)]$.

\subsubsection{}

Let us turn this into an absolute invariant right now.
It takes values in the integral Deligne cohomology
$H^k_{Del}(B)$ which has similar structures as its 
rational counterpart \ref{gv21}.

If $\cE_{geom}$ is a geometric family with $\ind(\cE_{geom})\in
K_k(B)$,
then we want to define a natural class
$\hat c_k(\cE_{geom})\in H^k_{Del}(B)$ such that
$R^{\hat c_k(\cE_{geom})}=(-1)^{m-1}(m-1)!\Omega^k(\cE_{geom})$ and
$\DD(\hat c_k(\cE_{geom}))=c_k(\ind(\cE_{geom}))$.
We will use the fact that a class $c\in H^k_{Del}(B)$ 
is determined by its holonomy $H(c):Z_{k-1}(B)\rightarrow \R/\Z$,
where
$Z_{k-1}(B)$ denotes the group of smooth singular cycles on $B$.

Consider a cycle $z\in Z_{k-1}(B)$. The trace of $z$ is the union of
the images of the singular simplices belonging to $z$. This trace
admits a neighborhood $U$ which up to homotopy looks like an at most
$k-1$-dimensional $CW$-complex. By assumption $\ind(\cE_{geom})_{|U}=0$
so that $(\cE_{geom})_{|U}$ admits a taming
$(\cE_{geom})_{|U,t}$.
In Subsection \ref{next21} we show that
there exists a  unique class
$\hat c_k(\cE_{geom})\in H^k_{Del}(B)$ such that
$$H(\hat c_k(\cE_{geom}))(z)=[(-1)^{m-1}(m-1)!\int_{z}
\eta^{k-1}((\cE_{geom})_{|U,t})]_{\R/\Z}$$
for all $z\in Z_{k-1}(B)$.
Note that the right-hand side is independent of the choice of the
taming by the observation \ref{obbs}.

\subsubsection{}

If $x\in K_k(B)$, then in fact $\ch_k(x)$ itself has integral periods,
or equivalently, $c_k(x)$ is divisible by $(m-1)!$. However the reason
for this is more complicated.

Let us assume that $B$ has the structure of a $CW$-complex with
filtration \index{filtration!of $CW$-complex by skeletons} $B^{k-1}\subseteq B^{k}\dots \subseteq B$ by skeletons.
Then the restriction $x_{|B^{k-1}}$ is trivial. In fact what we need
here is a trivialization of a geometric object which represents $x$.  Then we could ask for an extension of this
trivialization to $B^k$. 

In the present paper we represent $x$ as the
index of a geometric family $\cE_{geom}$. Then the trivialization is a
taming of the restriction of this family to $B^{k-1}$ (or better some
open neighborhood). We ask for an extension of the taming
\index{obstruction!class}\index{$o^k$} to $B^k$. In this situation there is an obstruction theory 
with obstruction class $o^k\in H^k(B,\Z)$. This class depends on the
choice of the restriction of the taming to $B^{k-2}$, and if $o^k$ vanishes,  then there
exists and extension of the taming from $B^{k-2}$ to 
$B^k$.

By a result of Kervaire \cite{kervaire59} we know that the image of $o^k$ in rational
cohomology is equal to $(-1)^m\ch_k(x)$.
So the geometric reason for the integrality of the periods of $\ch_k(x)$
is the existence of the taming of the restriction  of $\cE_{geom}$ 
to $B^{k-1}$. Thus one could again ask for an absolute invariant \index{$\hat o^k$}
in $\hat o^k\in H^k_{Del}(B)$ such that $\DD(\hat o^k)=o^k$ and
$R^{\hat o^k}=(-1)^m\Omega^k(\cE_{geom})$.

\subsubsection{}

It seems not to be very natural to construct a
global object like $\hat o^k$ from data defined in a neighborhood of $B^{k-1}$.
In fact from the geometric point of view  even to fix a $CW$-structure
on $B$ is not natural.
Thus we will introduce another picture of a trivialization of 
the index associated to a geometric family. The main notion is that of
a tamed $k-1$-resolution. In the present paper we discuss two sorts of tamed
resolutions  \ref{tkres} and \ref{chdalter}.
To either sort of a tamed $k-1$-resolution we associate integral Deligne
cohomology classes in  $H^k_{Del}(B)$ with the expected properties and study their
dependence on the choices. Unfortunately it is to technical to state
the precise result here in the introduction. We refer to Section
\ref{klo98} for further details.

\subsubsection{}

It turns out that for $k=0,1,2,3$ the Deligne cohomology class
associated to a tamed $k-1$-resolution of $\cE_{geom}$ is (up to sign)
equal to the class $\hat c_k(\cE_{geom})$. Thus our analytic
obstruction theory provides an alternative construction of that class.

For these small degrees the integral Deligne cohomology classifies
well-known geometric objects, namely continuous $\Z$-valued functions,
smooth $\R/\Z$-valued functions, hermitian line bundles with
connection, and geometric gerbes.

In fact, such objects have been associated previously to geometric
families. For even-dimensional geometric families we have the
index as $\Z$-valued function and the determinant line bundle with
Quillen metric and Bismut-Freed connection.
To odd-dimensional families we associate the spectral asymmetry,
i.e. the $\eta$-invariant, and the index gerbe (recently constructed
by Lott \cite{lott01}).
In Section \ref{fc34} we show that $\hat c_k(\cE_{geom})$
indeed classifies the expected  object.

\subsection{Dirac operators, boundaries, and corners}

\subsubsection{}

A large part of the present paper is devoted to the local index theory
for geometric families with corners.
The main motivation for developing this theory was its use in the
definition and analysis of tamed resolutions. But in the form
presented in Part \ref{ggh43} it seems to be interesting and useful in
its own right.

Before we give a more detailed overview we shall express one warning right
now. A manifold with corners for us is mainly a combinatorial object.
It allows to talk about faces and product structures. But the
index theory is in fact a $L^2$-index theory over the extended
manifold. Here the extension is obtained by glueing cylinders to all faces in
order to obtain a manifold without singularities like boundaries or
corners. The operators will be extended in a translation invariant
manner over the cylinders.

\subsubsection{}

Consider a generalized Dirac operator on a manifold
with boundary. Assume that we have collar neighborhood
of the boundary where the operator is translation invariant.
We now glue a cylinder in order to complete the manifold and
extend the operator naturally.

It turns out that the operator is Fredholm
(in the sense that zero does not belong to the essential spectrum) if and only if 
the boundary reduction of the Dirac operator is invertible.
This condition is  of course not fulfilled in general.
But it is true that the index of the boundary reduction vanishes
by the bordism invariance of the index. Therefore we can find a taming of
the boundary reduction.
This taming can be lifted to the cylinder and then extended as a
bounded operator to the whole manifold using a cut-off function.

We use this lift of the taming in order to perturb the original
Dirac operator. We call this a boundary taming.
The resulting boundary tamed operator is now Fredholm and we can
consider its index. In this case our index theorem \ref{etaprop} is
just a version of the classical  Atiyah-Patodi-Singer
\cite{atiyahpatodisinger75} theorem.

In fact all this works for families and yields the a version of the
results of Bismut and Cheeger \cite{bismutcheeger90},
\cite{bismutcheeger901}.
In fact, our approach 
 is very close to that of Melrose-Piazza \cite{melrosepiazza97} who just consider a
slightly different way of lifting the taming in order to produce
the Fredholm operator.

\subsubsection{}

If the boundary has several components then there is a 
complication which already appears in the case of the Dirac operator
on the interval $[0,1]$. If we consider all boundary components  as one boundary
face, then we can proceed as before. But the indices of the boundary
reductions of the Dirac operator 
to the individual components do not necessarily vanish separately. 
So if we wish to consider the components separately, then we are stuck.
This problem becomes even more complicated in the presence of corners
of higher codimension. 

\subsubsection{}

Let us now assume that we have a corner of codimension two.
Then locally there are two boundary components which meet in this
corner. The corner now appears as a boundary of the boundary components.
In order to find a Fredholm perturbation we first must choose a taming
of the reduction of the Dirac operator to the corner. 
This induces  Fredholm perturbations of the
boundary reductions  as explained above.
Assuming that their indices vanish we can choose tamings of the
boundary reductions and finally get a Fredholm perturbation of the
original operator.

For corners of higher codimension we proceed in a similar inductive
manner. In each intermediate step where we encounter a Fredholm
perturbation whose index is an obstruction to proceed further.

\subsubsection{}
In order to get some control about this we develop an obstruction
theory. We introduce the notion of an admissible face decomposition.
In general a face in this decomposition  may have several connected
components, which we call atoms. Admissibility first of all requires
that 
each face must be embedded. So e.g. the one-eck, a two dimensional
manifold with one corner point and one boundary component, is excluded.
Furthermore is required that if a face of higher codimension meets
a face of lower codimension, then it is in fact contained in the
latter.

The combinatorial datum of the inclusion relation between faces is
used to define the face complex, a chain complex over $\Z$.
Its cohomology receives the obstructions against continuing
the inductive process of choosing tamings. For details we refer to Subsection
\ref{obst}.

\subsubsection{}

What we said above extends to families.
For a boundary tamed family of compatible Dirac operators
on a family of manifolds with corners we prove a local index theorem
\ref{famind}.

The boundary faces of a boundary tamed family are tamed. In Subsection \ref{ttrr4411} we extend the construction of
the $\eta$-forms to tamed geometric families with corners.
Then the index theorem has the expected form. {\em The de Rham
representative of the Chern character of the index is the
sum of the contributions of the local index form and the sum of
$\eta$-forms of the boundary faces.}

\subsubsection{}

Above we used the notion of the boundary reduction in a sloppy way.
In fact, we must form iterated boundary reductions in order to reduce to
corners of higher codimension.
The difficulty is that geometric families may have non-trivial
automorphisms. 

A boundary face of a geometric family is a well-defined
isomorphism class of geometric families. Since in our
framework a manifold with corners comes with a distinguished collar
there is an essentially canonical way of restricting Dirac bundles.
If we take this canonical restriction of the Dirac bundle we call the
resulting geometric family a canonical model of the boundary face.
Recall that we want to lift a taming of the boundary reduction to a
boundary taming of the original operator. In order to do this we must
be precise with the identifications.

But consider now a corner of codimension two. Then we have two local
boundary components meeting in this corner. Therefore we have two ways
of reducing the Dirac operator to the corner.
Already the induced orientations of the corner are opposite.
In fact it turns out that the isomorphism classes of geometric
families obtained in these two ways are opposite to each other.
It is an important observation that if we take canonical models in
these two ways, then
there is a preferred (orientation and grading reversing) isomorphism
between them (see Lemma \ref{hiphop1}).

This will be important for the follows reason.
Let us consider the two boundary faces separately. Then we want to fix
tamings of their boundaries. If possible we then extend these boundary
tamings (of the boundary faces) to tamings.
We want to say that if the chosen  tamings on 
the boundaries of the boundary faces coincide,
then all choices together induce a boundary taming of the original
family. In order to make this comparison we need the preferred
 isomorphism.

In Section
\ref{jjhh7} we develop a precise language in order to deal with this
kind of problems.

\subsubsection{}

Note that the problem already appears for the simplest kind of Dirac
operators. Consider a Riemannian spin manifold. Then we have a
well-defined {\em isomorphism class} of the spinor bundle.
If one takes into account that the spinor  bundle has a real structure,
then the isomorphisms between different representatives of the spinor bundle are determined
up to sign.

Assume that the manifold is odd-dimensional. 
If the manifold has a boundary, then we have an
induced spin structure on the boundary. Furthermore, the restriction
of the spinor bundle to the boundary is isomorphic to the spinor
bundle of the boundary.
But if we want to consider a boundary value problem
(e.g. fix a spinor on the boundary and look for harmonic spinors on
the interior
which have this prescribed boundary restriction)
we must fix this isomorphism. There is no canonical choice.

\section{Further remarks}

\subsection{Orientation conventions}\label{hhhas}



\subsubsection{}

In various places we must take a boundary of an oriented manifold
with the induced orientation. Here is our {\em orientation convention}.
The  model of a manifold with boundary is the upper half space
$\R^n_+:=\{x_n\ge 0\}$. We define the induced boundary orientation
such that in this example standard orientations are preserved.

\subsubsection{}

We define the de Rham differential \index{de Rham differential} as usual such that
$$d(fd\omega)=df\wedge d\omega$$
for a function $f$ and a form $\omega$.

If $\omega$ is a compactly supported  $n$-form on $\R^n$,
then we define
$$\int_{\R^n}\omega :=\int_{\R^n} \omega(e_1,\dots,e_n)(x) dx\ ,$$
where $(e_i)_{i=1}^n$ is the standard basis and $dx$ is the Lebesgue measure.

\subsubsection{}\label{stokess}

 These conventions determine the sign in Stoke's theorem \index{sign in Stoke's theorem}
$$\int_{\R^n_+} d\omega =(-1)^n \int_{\R^{n-1}} \omega_{|\R^{n-1}}\ .$$

\subsubsection{}

A spinor module \index{spinor!module} of the Clifford algebra $\Cliff(\R^1)$
is a one-dimensional complex vector space $\Delta^1$ \index{$\Delta^n$} such that the generator
$e_1\in\R^1$ acts as multiplication by $i$.
This determines the isomorphism classes of the spinor modules
$\Delta^n$ of $\Cliff(\R^n)$ for all $n\in\nat$ by
the rule $\Delta^{n+m}\cong \Delta^n*\Delta^m$, where the
$*$-operation is defined in Definition \ref{startr}.

We consider \index{Spin group} \index{$Spin(n)$}$Spin(n)\subset \Cliff(\R^n)$. In this way  we have fixed the
isomorphisms classes of the spinor representations as well.

\subsection{About the contents}

\subsubsection{}

In Section \ref{ggh43} we consider the analysis of Dirac operators on
manifolds with corners (or rather their extensions). We make precise
the notion of a boundary reduction and a taming.
Besides introducing language the main result of this section is the
obstruction
theory in Subsection \ref{obst}.

It is only the language part which is needed in subsequent Sections.

\subsubsection{}

In Section \ref{uuzz66} we introduce geometric families and generalize
the notion of boundary faces and tamings to the family case.
The main results are the construction of the $\eta$-form for a tamed
family with corners and the computation of its 
de Rham differential in  
Theorem  \ref{etaprop}. It turns out that the differential
of the $\eta$-form is the sum of the expected contribution of the
local index form and a boundary correction which is again given by
$\eta$-forms.

We use this result in order to deduce a local index theorem for
boundary tamed families \ref{famind} which works in the even as well
in the odd-dimensional case. Already in the  case without boundary
the proof seems to be simpler than the usual proofs since there is no
complicated discussion of the $t\to\infty$ limit as e.g. in Ch. 9 
of \cite{berlinegetzlervergne92}.

\subsubsection{}

In Section \ref{hhgg1} 
we collect some material about $K$-theory.
In particular we introduce the Atiyah-Hirzebruch filtration and
various pictures of the associated obstruction theory.
This section contains no new results.

\subsubsection{}

In Section \ref{d555} we introduce the notion of geometric and tamed
chains. We develop an obstruction theory which governs the prolongation
of tamings of chains. The main results
are Lemmas \ref{olem} and \ref{ww}.
The theory developed here is
more general than needed later. But we feel that it may have
application elsewhere. It is an open problem to understand the groups
of bordism classes of geometric and tamed $k$-chains
$G^k(B)$ and $G^k_t(B)$ in detail.

\subsubsection{}

In Section \ref{kj65}
we consider a special sort of chains, namely resolutions of a
geometric family. In this special case we relate the problem of finding tamed lifts
with the obstruction theory associated to the Atiyah-Hirzebruch
filtration of $K$-theory.
The main result is Theorem \ref{compara} in which we show that the
corresponding two obstruction sets are equal.

\subsubsection{}
In Section \ref{dd212} we introduce  Deligne cohomology.
Then we associate Deligne cohomology classes to tamed resolutions.
The main result is the construction of Deligne-cohomology valued lifts
of the obstruction classes \ref{delcoh1} and Chern classes
\ref{chern3}.

\subsubsection{}

In the last Section \ref{fc34} we discuss in detail the small degree cases.

\subsection{Acknowledgment}

\subsubsection{}
The author started to work on this project after he has received
  the first version of Lott's paper on the index gerbe \cite{lott01}
  in June 2001. The construction of higher-dimensional Deligne
  cohomology classes was not contained in  this first version. After
  an E-mail conversation with J. Lott in September  2001 it became clear that we had the same project, but different approaches. 

\subsubsection{}

The topic of the present paper was the starting point of a greater
project. In this project we consider smooth versions of generalized
cohomology theories in the same sense as Deligne cohomology is a
smooth version of integral cohomology. The main idea is to look for
lifts  of cohomological identities (like index theorems) involving
geometric objects to the smooth extensions using secondary information
(i.e. the information encoded in the geometric proofs of the identities,
if such proofs exist). We refer to
\cite{bunkeschick031} and \cite{bunkeschick022}
for further information.

One of purposes of the writing and rewriting of the present paper is
that its results can be used in this related work without any further adaption.

\subsubsection{}

I thank Th. Schick for helpful remarks and corrections.
Furthermore, I thank the referees for their criticism
which  forced a considerable improvement in precision and presentation of
the paper.

\part{Index theory for families with corners}\label{ggh43}

\section{Dirac operators on manifolds with corners}\label{jjhh7}

\subsection{Dirac bundles}
\newcommand{\Cl}{{\tt Cl}}
\subsubsection{}
In this subsection we fix some basic conventions about Dirac
bundles. Roughly speaking a Dirac bundle is a complex vector bundle  together with all
those structures needed to define the Dirac operator. For the purpose
of the present paper the reduction of Dirac bundles over products 
is of particular importance. In order to provide an effective notation
for the investigation Dirac operators on manifolds with corner
singularities of higher codimension we develop a reduction calculus.

\subsubsection{}
Let $(M,g^M)$ be a Riemannian manifold. By $\Cl(TM)$  \index{$\Cl(TM)$} we denote the bundle of complex Clifford algebras
associated to the bundle of euclidian vector spaces
$(TM,g^M)$.  
\begin{ddd}\label{diracbbb}
Let $M$ be even-dimensional.
A Dirac bundle \index{Dirac bundle} over $(M,g^M)$ is a tuple
\index{$\cV$} $\cV=(V,h^V,\nabla^V,c,z)$, \index{($V,h^V,\nabla^V,c,z)$} where
\begin{enumerate}
\item
$V$ is a complex vector bundle over $M$,
\item $h^V$ is a hermitian metric on $V$, \index{$h^V$}
\item $\nabla^V$ is a connection on $V$ which is compatible with $h^V$,
\item $c:TM\rightarrow \End(V)$ is a bundle homomorphism \index{$c$}
which is parallel and extends to a $*$-homomorphism
$c:\Cl(TM)\rightarrow \End(V)$, i.e. it satisfies
\begin{enumerate}
\item $c(X)^*=-c(X)$ for all $X\in TM,$
\item $c(X)^2=-\|X\|^2_{g^M}$ for all $X\in TM$,
\item $\nabla^V_Yc(X)-c(X)\nabla^V_Y=c(\nabla^{TM}_YX)$ for
  $X\in C^\infty(M,TM)$, $Y\in TM$, where $\nabla^{TM}$ is the     
Levi-Civita connection on $TM$.
\end{enumerate}
\item $z$ is a $\Z/2\Z$-grading of $V$ which is parallel, i.e.
$[\nabla^V_X,z]=0$ for all $X\in TM$, and which satisfies
 $\{c(X),z\}=0$\footnote{Here $\{X,Y\}:=XY+YX$  denotes the anticommutator.} for all $X\in TM$.
\end{enumerate}
If $M$ is odd-dimensional, then a Dirac bundle
is a tuple
$\cV=(V,h^V,\nabla^V,c)$ of objects as above, but without grading. 
 \end{ddd}
\index{grading}If $M$ is even-dimensional and oriented, then the Clifford multiplication by the volume form induces a \index{grading!by volume form}canonical grading which can be written in terms of an oriented orthogonal frame as
$i^{\frac{\dim(M)}{2}}c(X_1)\dots c(X_n)$.  But we shall need other gradings as well. They arrise e.g. if we twist by auxiliary graded hermitian bundles, see \ref{twisis}.

\subsubsection{}
\index{Dirac operator!associated to a Dirac bundle}\index{Dirac bundle!associated Dirac operator}\index{$D(\cV)$} To a Dirac bundle $\cV$  we associate the Dirac operator
$D(\cV):C^\infty(M,V)\rightarrow C^\infty(M,V)$ which is the first-order elliptic formally selfadjoint differential operator given by the composition
$$C^\infty(M,V)\stackrel{\nabla^V}{\rightarrow} C^\infty(M,T^*M\otimes
V)\stackrel{g^M}{\rightarrow} C^\infty(M,TM\otimes
V)\stackrel{c}{\rightarrow} C^\infty(M,V) .$$
\subsubsection{}
\index{spinor!bundle}\index{$\cS(M)$}A typical example is the Dirac bundle structure $\cS(M)$
on the spinor bundle $S(M)$ of a Riemannian
spin manifold. Let $M$ be $n$-dimensional. A spin structure on $M$ is
a $Spin(n)$-principal bundle $Q\rightarrow M$ together with an
isomorphism of $Q\times_{Spin(n)} SO(n)\to M$ with the $SO(n)$-principal bundle of
orthonormal oriented frames. Recall that we consider $Spin(n)\subset
\Cl(\R^n)$.
If we choose a spinor module $\Delta^n$ of $\Cl(\R^n)$, then
the spinor bundle is represented by $S(M):=Q\times_{Spin(n)}\Delta^n$.
It carries a Dirac bundle structure in a natural way.

Note that $\Delta^n$ and thus
$S(M)$ are only well-defined as an
isomorphism class. The group of automorphisms of $\Delta^n$ as a hermitean 
$\Cl(\R^n)$-module is $U(1)$. Similarly, the group of automorphisms of the Dirac bundle
$S(M)$ consists of locally constant functions with values in $U(1)$.

A finer consideration (compare
e.g. \cite{broeckertomdieck85}, Ch. 6.6) shows that $\Delta^n$ and $S(M)$
 come with an additional symmetry, namely a
real or quaternionic structure depending on the class of $n$ in
$\Z/8\Z$. We refer to \cite{bunkeschick022}, Sec. 2.2, where we have
written out these structures explicitly. In the present paper we will not need the explicit formulas, but only their existence. In fact, we will use these symmetries in order to reduce the group of automorphisms.
The group of automorphisms of $\cS(M)$ which preserve this additional
structure is the group $\Z/2\Z$ acting through scalar multiplication
by $\pm 1$. Therefore if we require an isomorphism between two
models of the spinor bundles of $M$ to preserve the additional symmetries, then it is uniqely determined up to a sign ambiguity.

\subsubsection{}
\index{$S(\R)$} In order to rigidify the
definition of boundary reductions of Dirac bundles later we will fix once and for
all a spinor bundle $S(\R)$, which induces our choices of the spinor bundles
on all open subsets of $\R$. We could take $S(\R)=\R\times \C$ as hermitean vector bundle with connection such that the standard basis vector $e_1\in
\R$ (i.e. $e_1=1$) acts as multiplication by $i$.
Since $\dim(\R)=1$ is odd, there is no grading to be fixed.

\subsubsection{}\label{twisis}
Given a Dirac bundle $\cV$ on $M$ (e.g. the spinor bundle),
 then other Dirac bundles can be  obtained  by twisting. \index{twisting}
By $\bV$ we denote the 
hermitian vector bundle with connection underlying $\cV$.
If  
\index{$\bW$}\index{$(W,h^W,\nabla^W,z_W)$}$\bW:=(W,h^W,\nabla^W,z_W)$ is an auxiliary hermitian vector bundle
with metric connection on $M$ which is $\Z/2\Z$-graded by $z_W$, then we form the twisted Dirac bundle
$\cV\otimes\bW$. The Dirac bundle structure on the underlying
hermitian vector bundle with connection $\bV\otimes \bW$ is given 
by the Clifford multiplication $c(X)\otimes z_W$ (and the $\Z/2\Z$-grading
$z\otimes z_W$, if $M$ is even-dimensional).

\subsubsection{}
\begin{ddd}
\index{$\cV^{op}$}\index{opposite!Dirac bundle}By $\cV^{op}$ we denote the opposite Dirac bundle given
by $$\cV^{op}:=(V,h^V,\nabla^V,-c,-z)$$ in the even-dimensional, and by
$$\cV^{op}:=(V,h^V,\nabla^V,-c)$$ in the odd-dimensional case.
\end{ddd}
Note that in the even-dimensional case the grading operator $z$ induces an isomorphism
$$(V,h^V,\nabla^V,c,z)\cong (V,h^V,\nabla^V,-c,z)\ .$$
\index{$\bW^{op}$}\index{opposite!graded vector bundle}In general, if  we set $\bW^{op}:=(W,h^W,\nabla^W,-z_W)$, then we have
\begin{equation}\label{euidewdew}
\cV\otimes\bW^{op}\cong (\cV\otimes\bW)^{op}\cong \cV^{op}\otimes
\bW\ .
\end{equation}

\subsubsection{}
Assume that $M$ is a Riemannian spin-manifold. 
Let $P$ be the $SO(n)$-principal bundle of oriented orthonormal
frames. 
Let $Q\rightarrow M$ and $u:Q\times_{Spin(n)}SO(n)\cong P$,  represent the
spin structure. \index{spin structure}
Let $P^{op}\rightarrow M$ be the  $SO(n)$-principal bundle of oriented
orthonormal frames of $M^{op}$, where $M^{op}$ denotes the Riemannian manifold $M$ with the reversed orientation.
If $n$ is odd, then we define an isomorphism of $SO(n)$-principal
bundles
$\theta:P\rightarrow P^{op}$ which maps the frame $(X_1,\dots X_n)$ to
\index{$P^{op}$}$(-X_1,\dots,-X_n)$. In this case the opposite spin structure is
represented by
$Q\rightarrow M$ and $Q\times_{Spin(n)}SO(n)\stackrel{u}{\cong}
P\stackrel{\theta}{\cong} P^{op}$.
 
Assume now that $n$ is even. We consider the element
$e_1\in \R^n\subset Pin(n)\subset Cl(\R^n)$. Its image $E_1\in O(n)$ acts as
$(x_1,\dots,x_n)\mapsto (x_1,-x_2,\dots ,-x_n)$.
We define the map
$\theta:P\rightarrow P^{op}$ by $(X_1,\dots,X_n)\mapsto
(X_1,-X_2,\dots,-X_n)$. It becomes an isomorphism of $SO(n)$-principal
bundles if we twist the action of $SO(n)$ on $P$ by $E_1$.
We write $\tilde P$ for the bundle $P$ with this twisted action.

Let $Q^{op}\rightarrow M$ be the $Spin(n)$-principal bundle
which is given by $Q\rightarrow M$ with the $Spin(n)$-action twisted by
$e_1$.
Then we have an isomorphism of $SO(n)$-principal bundles
$Q^{op}\times_{Spin(n)} SO(n)\stackrel{\tilde u}{\cong} \tilde P\stackrel{\theta}{\cong}
P^{op}$, where $\tilde u[q,s]:=u[q,E_1sE_1]$.
In this way we represent the opposite spin structure.
\index{opposite!spin structure}

The following assertion is now easy to check.
Let  $-M$ denote the manifold with the opposite orientation and spin
structure. 
\begin{lem} \label{euifweiufwe}
We have $\cS(-M)\cong \cS(M)^{op}$. 
\end{lem}
The point of this discussion is that $\cS(-M)$ is defined as an isomorphism class of Dirac bundles only, while
$\cS(M)^{op}$ gives a canonical representative of this class once we have fixed
a representative $\cS(M)$.

\subsubsection{}
If $M$ is a Riemannian spin manifold, then every Dirac bundle
on $M$ is of the form $\cS(M)\otimes \bW$, where $\bW$ is uniquely determined
up to isomorphism and called the twisting bundle.
In fact, if $M$ is even-dimensional, then we have $W\cong \Hom_{\Cl(TM)}(S(M),V)$ with induced metric,
connection and $\Z/2\Z$-grading.
If $M$ is an odd-dimensional, then
$W:=\Hom_{\Cl(TM)}(S(M)\oplus S(M)^{op},V)$,
where the $\Z/2\Z$-grading is induced from the grading $\diag(1,-1)$
of $S(M)\oplus S(M)^{op}$.


\subsubsection{}
If $f:M\rightarrow N$ is a local  isometry of Riemannian manifolds, and
$\cV$ is a Dirac bundle over $N$, then we have a pull-back Dirac bundle
$f^*\cV$ over $M$. The underlying $\Z/2\Z$-graded vector bundle is
$f^*\bV$ which is defined as usual. We define the Clifford
multiplication as follows. Let $m\in M$. Then we have a natural
isomorphism $F_m:(f^*V)_m\stackrel{\sim}{\rightarrow} V_{f(n)}$. If $X\in
T_mM$, then we define $c_{f^*\cV}(X)$ such that
$F_m(c_{f^*\cV}(X)(v))=c_{\cV}(df(X))F_m(v)$ for all $v\in (f^*V)_m$.

\subsection{Operations with Dirac bundles}\label{operations}

\subsubsection{}
In the theory of boundary value problems for Dirac operators
we have the standard simplifying assumption of a product structure on a collar neighborhood of the boundary. This product structure allows us to write the Dirac operator on the collar in a simple form in terms of its boundary reduction.
In the case of higher-codimensional singularities
like corners we need a generalization of the notions of a product
structure and a boundary reduction.

\subsubsection{}
Let $(H,h^H)$ be a connected Riemannian spin manifold with spinor
bundle $S:=S(H)$. Assume, that $M=N\times H$ is a Riemannian product. 
\begin{ddd}
\index{locally of product type}We say that the Dirac bundle $\cV$ on $M$ is locally of product type, if $R^V(X,Y)=0$ for all
$X\in TN$ and $Y\in TH$, where $R^V$ denotes the curvature of $\nabla^V$.
\end{ddd}

\subsubsection{}
\index{$\cV//H$} If $\cV$ is locally of product type, then we define
an isomorphism class of Dirac bundles 
$\cV//H:=\cW$ over $N$, the
reduction of $\cV$
along $H$, by the following construction.

Let us write $\cW=(W,h^W,\nabla^W,c_W,z_W)$ if $N$ is even-dimensional,
and  $\cW=(W,h^W,\nabla^W,c_W)$ if $N$ is odd-dimensional.
We fix some point $h\in H$ and let $S_h$ denote the fiber of $S$ over
$h$. 
The bundle $W$ is given by
\begin{equation}\label{deudwedewdewdw}
\begin{array}{|c|c|c|}\hline
W:=\dots&\dim(H)\equiv 0(2)&\dim(H)\equiv 1(2)\\
\hline
\dim(N)\equiv 0(2)&\Hom_{\Cl(T_hH)}(S_h,V_{|N\times\{h\}})&\Hom_{\Cl(T_hH)}(S_h\oplus S_h^{op},V_{|N\times\{h\}})\\\hline
\dim(N)\equiv
1(2)&\Hom_{\Cl(T_hH)}(S_h,V_{|N\times\{h\}})&\Hom_{\Cl(T_hH)}(S_h,V_{|N\times\{h\}})\\
\hline
\end{array}\ 
\end{equation}
We let $\nabla^W$ and $h^W$ be the induced connection and metric.

Let $X\in TN\cong T(N\times\{h\}) \subset TM$.
The Clifford multiplication $c_W(X)$ is given by
$$\begin{array}{|c|c|c|}\hline
c_W(X)\phi:=&\dim(H)\equiv 0(2)&\dim(H)\equiv 1(2)\\
\hline
\dim(N)\equiv 0(2)& c(X)\circ \phi\circ z_S&c(X)\circ \phi \circ \left(\begin{array}{cc}0&1\\1&0\end{array}\right)\\\hline
\dim(N)\equiv
1(2)& c(X)\circ \phi\circ z_S&i z\circ c(X)\circ \phi\\
\hline
\end{array}\ .$$
Finally, the grading is given by
$$\begin{array}{|c|c|c|}\hline
z_W\phi:=&\dim(H)\equiv 0(2)&\dim(H)\equiv 1(2)\\
\hline
\dim(N)\equiv 0(2)&z\circ \phi\circ  z_S&\phi\circ \left(\begin{array}{cc}1&0\\0&-1\end{array}\right)\\
 \hline
\dim(N)\equiv
1(2)&  & \\
\hline
\end{array}\ .$$

\subsubsection{}
The construction of $\cV//H$ depends on the choice of the base point $h\in H$ and the
choice of the spinor bundle $S(H)$. Let us write $\cV//(H,h,S(H))$ for the moment in order to indicate this
dependence. 
If $h_0,h_1\in H$, then we choose a path $\gamma$ from $h_0$ to $h_1$.
Parallel transport along this path induces isomorphisms
$\Cl(T_{h_0}H)\cong \Cl(T_{h_1}H)$ and $S_{h_0}\cong S_{h_1}$ which
preserve the Clifford algebra and module structures (and grading). Moreover we have 
isomorphisms  $\Cl(TM)_{{|N\times\{h_0\}}}\cong
\Cl(TM)_{{|N\times\{h_1\}}}$ and
$V_{|N\times\{h_0\}}\cong V_{|N\times\{h_1\}}$
 again preserving the Clifford algebra and module structures (and grading).
By our curvature assumption the last isomorphism also preserves the 
connection. Therefore these isomorphisms induce an isomorphism
$\Phi(\gamma):\cV//(H,h_0,S(H))\rightarrow \cV//(H,h_1,S(H))$ of Dirac bundles.
Note that this isomorphism may depend on the choice of the path
$\gamma$. Later we will impose further restrictions which imply
independence of the path.
If $S^\prime(H)$ is another choice for the spinor bundle of $H$, then
we have an isomorphism $\alpha: S(H)\rightarrow S^\prime(H)$ which induces an
isomorphism $\cV//\alpha:\cV//(H,h_0,S(H))\rightarrow
\cV//(H,h_1,S^\prime(H))$. For example we could take
$S^\prime(H)=S(H)$ and $\alpha=-1$. Then $\cV//\alpha=-1$.
If $\beta:\cV\rightarrow \cV^\prime$ is an isomorphism of Dirac
bundles,
then we have an induced isomorphism
$\beta//(H,h,S(H)):\cV//(H,h,S(H))\rightarrow \cV^\prime//(H,h,S(H))$.
Let $g:N^\prime\rightarrow N$ an isometry. Then we have a canonical
isomorphism
$g^\sharp:g^*(\cV//(H,h,S(H)))\rightarrow (g\times \id_H)^*\cV//(H,h,S(H))$.

\begin{ddd}
We define the reduction $\cV//H$ of $\cV$ along $H$ by the construction above.
\end{ddd}

The boundary reduction of a Dirac bundle is a special case of this construction
where $\dim(H)=1$.

\subsubsection{}\label{ccoa}
We now discuss the opposite process.
We start with a Dirac bundle $\cW$ on $(N,g^N)$. Furthermore, let
\index{$\cW*H$}$(H,g^H)$ be a Riemannian spin manifold with spinor bundle $S=S(H)$.
Then we define a Dirac bundle $\cV:=\cW*H$ (which is locally of product type)
 on the Riemannian product $M:=N\times H$ as follows. 

We define $V$ by 
\begin{equation}\label{ewoifcewc}
\begin{array}{|c|c|c|}\hline
V:=&\dim(H)\equiv 0(2)&\dim(H)\equiv 1(2)\\
\hline
\dim(N)\equiv 0(2)& W\otimes  S&W\otimes S\\\hline
\dim(N)\equiv
1(2)& W\otimes  S&W\otimes S\otimes \C^2\\
\hline
\end{array}
\end{equation}
(formally we should write here $\pr_H^*S$ etc, but we omit the
projections in order to simplify the notation).
This bundle comes with an induced metric and connection.
 Furthermore, we define the Clifford multiplication with $X\in TN$ by
$$\begin{array}{|c|c|c|}\hline
c(X):=&\dim(H)\equiv 0(2)&\dim(H)\equiv 1(2)\\
\hline
\dim(N)\equiv 0(2)&c_W(X)\otimes z_S& c_W(X)\otimes 1\\\hline
\dim(N)\equiv
1(2)& c_W(X)\otimes z_S&c_W(X)\otimes 1\otimes \left(
\begin{array}{cc}0&i\\-i&0\end{array}\right)\\
\hline
\end{array}\ ,$$
and with $Y\in TH$ by
$$\begin{array}{|c|c|c|}\hline
c(Y):=&\dim(H)\equiv 0(2)&\dim(H)\equiv 1(2)\\
\hline
\dim(N)\equiv 0(2)&1\otimes c_S(Y)&z_W\otimes c_S(Y)\\\hline
\dim(N)\equiv
1(2)& 1\otimes c_S(Y)&1\otimes c_S(Y)\otimes \left(
\begin{array}{cc}1&0\\0&-1\end{array}\right)\\
\hline
\end{array}\ .$$
Finally the grading is given by 
$$\begin{array}{|c|c|c|}\hline
z:=&\dim(H)\equiv 0(2)&\dim(H)\equiv 1(2)\\
\hline
\dim(N)\equiv 0(2)&z_W\otimes z_{S}& \\\hline
\dim(N)\equiv
1(2)& &1\otimes 1\otimes \left(
\begin{array}{cc}0&1\\1&0\end{array}\right)\\
\hline
\end{array}\ .$$

\begin{ddd}\label{startr}
We define the extension $\cW*H$ of $\cW$ by $H$ by the construction above.
\end{ddd}

\subsubsection{}
Note that $\cW*H$ depends on the model of the spinor bundle $S(H)$. However,
the isomorphism class of $\cW*H$ only depends on the Riemannian spin
manifold $H$.
 Let us write $\cW*(H,S(H))$ for the moment in order to
indicate this dependence on $S(H)$.
If $S^\prime(H)$ is another choice and $\alpha:S(H)\rightarrow
S^\prime(H)$ is an isomorphism, then we have an isomorphism
$\cW*\alpha:\cW*(H,S(H))\rightarrow \cW*(H,S^\prime(H))$.
If $\beta:\cW\rightarrow \cW^\prime$ is an isomorphism of Dirac
bundles, then we have an isomorphism
$\beta*(H,S(H)):\cW*(H,S(H))\rightarrow \cW^\prime* (H,S(H))$.
Let $g:N^\prime\rightarrow N$ be an isometry. Then we have a natural 
isomorphism
$g_\sharp:(g\times \id_H)^*(\cW*(H,S(H)))\rightarrow g^*\cW*(H,S(H))$.

\subsubsection{}
The following Lemma is a simple consequence of the structure of 
modules over complex Clifford algebras.
\begin{lem}\label{hinher}
\begin{enumerate}
\item
There is an isomorphism $(\cW*H)//H\cong \cW$.
In fact, we have a canonical isomorphism
$\cW*(H,S(H))//(H,h,S(H))\cong \cW$.
\item
If $H=H_1\times H_2$ is a product of Riemannian spin manifolds,  
then there is an isomorphism 
$\cW*(H_1\times H_2)\cong (\cW*H_1)*H_2$.
\item
Let $-H$ denote $H$ with the opposite orientation and spin structure.
Then we have isomorphisms $\cW*(-H)\cong (\cW*H)^{op}\cong \cW^{op}*H$.
\item If $N$ is a spin manifold, then we have an isomorphism
$\cS(N)*H\cong \cS(M)$, where $M$ has the product spin structure.
We also have $\cS(M)//H\cong \cS(N)$.
\end{enumerate}
\end{lem}
\proof
We give the canonical isomorphism claimed in 1.
All the isomorphisms in the following  are induced by canonical
isomorphisms of tensor algebra.
Let first  $H$ be even-dimensional.
Then this isomorphism is given by
$$\Hom_{\Cl(T_hH)}(S_h,W\otimes S_h)\cong 
\Hom_{\Cl(T_hH)}(S_h,S_h)\otimes W\cong \C\otimes W\cong W\ .$$
Let now $H$ be odd-dimensional and $N$ be even-dimensional.
Then
\begin{eqnarray*}
\Hom_{\Cl(T_hH)}(S_h\oplus S_h^{op},W\otimes S_h)&\cong&
\Hom_{\Cl(T_hH)}(S_h\oplus S_h^{op},W^+\otimes S_h\oplus W^-\otimes
S_h^{op})\\
&\cong& \Hom_{\Cl(T_hH)}(S_h,S_h)\otimes W^+\oplus
\Hom_{\Cl(T_hH)}(S_h^{op},S_h^{op})\otimes W^-\\&\cong& \C\otimes
W^+\oplus \C\otimes W^-\\&\cong& W^+\oplus W^-\\&\cong& W \ .\end{eqnarray*}
If $N$ and $H$ are odd-dimensional, then
we  have
\begin{eqnarray*}
\Hom_{\Cl(T_hH)}(S_h,W\otimes S_h\otimes \C^2)&\cong&
\Hom_{\Cl(T_hH)}(S_h,S_h\oplus S_h^{op})\otimes W\\
&\cong&\Hom_{\Cl(T_hH)}(S_h,S_h)\otimes W\\
&\cong&\C\otimes W\\
&\cong& W\ .
\end{eqnarray*}
We leave it to the interested reader to check that these isomorphisms
are compatible with the remaining Dirac bundle structures.

Next we discuss 4. Let $n_i:=\dim(H_i)$ and consider 
$n_i$-dimensional oriented  euclidian vector spaces $V_i$. We choose 
models $\Delta_i$ of the spinor modules of
$\Cl(V_i)$. Then a model of the spinor module $\Delta$ of 
$\Cl(V_1\oplus V_2)$ is given by
$\Delta_1*\Delta_2$ which is defined by analogous formulas as
the extension of spinor bundles.
A formal way to fix the details is as follows.
We choose spinor bundles $S(V_i)$ and
take 
$\Delta_i:=S(V_i)_0$. Then we form
$\cS(V_0\times V_1):=\cS(V_0)*(V_1,S(V_1))$ and define
$\Delta:= S(V_1\times V_2)_0$. 

Let $P_i\rightarrow H_i$ be the $SO(n_i)$-principal bundles given by
the oriented orthonormal frame bundles of $H_i$.
Let $Q_i\stackrel{\Z/2\Z}{\rightarrow} P_i\to H_i$ be the $Spin(n_i)$-principal
bundles
representing the spin structures of $H_i$. Then we get the model
$S(H_i):=Q_i\times_{Spin(n_i)}\Delta_i$ of the spinor bundle of $H_i$.
We have an embedding $Spin(n_1)\times Spin(n_2)/\Z/2\Z\rightarrow Spin(n_1+n_2)$.
The product spin structure of $H_1\times H_2$ is given by the
extension of structure groups from $Spin(n_1)\times Spin(n_2)/\Z/2\Z$ to
$Spin(n_1+n_2)$ of $Q_1\times Q_2/\Z/2\Z\rightarrow P_1\times
P_2\rightarrow H_1\times H_2$.
Therefore a model of $S(H_1\times H_2)$ is given by
$Q_1\times Q_2/\Z/2\Z\times_{Spin(n_1)\times Spin(n_2)/\Z/2\Z}\Delta$.
Writing out this explicitely gives the isomorphism
$$\cS(H_1)*(H_2,S(H_2))\cong \cS(H_1\times H_2)\ .$$
The second assertion of 4. follows from 1.

Note that this argument shows more. Namely, once we have fixed the
models for $\Delta_i$ and the representatives of spin structures $Q_i$
there are canonical models of $S(H_i)$ and $S(H_1\times H_2)$
such that the isomorphism $\cS(H_1)*(H_2,S(H_2))\cong \cS(H_1\times
H_2)$ is canonical. 
This will be employed in the proof of 2.
Assume that we have chosen $\Delta_i$ and $Q_i$ as in 4.
If $N$ is a spin manifold, then in addition we choose a representative
of the spin structure of $N$ and a model $\Delta_0$ of the
corresponding spinor module.
Furthermore, we must choose an isomorphism
$\Delta_0*(\Delta_1*\Delta_2)\cong (\Delta_0*\Delta_1)*\Delta_2$.
Then after  fixing these choices we have a sequence of canonical isomorphisms
 \begin{eqnarray*}
(\cS(N)*(H_1,S(H_1)))*(H_2,S(H_2))&\cong& \cS(N\times H_1)*(H_2,S(H_2))
\\&\cong& \cS((N\times H_1)\times
H_2))\\&\stackrel{*}{\cong}& 
\cS(N\times (H_1\times H_2))\\&\cong&
 \cS(N)*(H_1\times H_2,S(H_1\times H_2))\ ,
\end{eqnarray*}
where $*$ depends on the associativity isomorphism above.

This chain of isomorphisms extends two twisted spinor bundles.
If $\cV=\cS(N)\otimes \bW$, then
$$(\cV*(H_1,S(H_1)))*(H_2,S(H_2))\cong \cV*(H_1\times H_2,S(H_1\times H_2))\ .$$
In general, we can assume that $\cV$ is locally on $N$ a twisted Dirac bundle.
Thus we  obtain the required isomorphism locally on $N$. Since it is
canonical after the choices made on the models, these local
isomorphisms
glue and provide a global isomorphism.

Assertion 3. follows from  (\ref{euidewdew}) and Lemma \ref{euifweiufwe}.
\hB 

\subsubsection{}
Let $H$ be a Riemannian spin manifold with spinor bundle $S(H)$.
\begin{ddd}
\index{product structure}A product structure on a Dirac bundle $\cV$ over $N\times H$ is given
by a Dirac bundle $\cW$ over $N$ and an isomorphism $\cV\cong \cW*H$.
\end{ddd}
This notion of a product structure extends the usual notion 
of a product structure on a Dirac bundle near a boundary. In this case 
$H:=[0,1)$. Since $0\in H$ is the canonical base point  we conclude by  Lemma \ref{hinher} 1., that  the Dirac bundle $\cW$ on $N$ is
uniquely determined up to canonical isomorphism by the product structure. 
In the next paragraph we will generalize this observation to higher-dimensional manifolds $H$.

\subsubsection{}
Assume that $M:=N\times H$ is a product of Riemannian manifolds, where
$H$ is spin, flat and simply\index{product structure!existence} connected with spinor bundle $S(H)$. 
\begin{lem}\label{loglo}
A Dirac bundle $\cV$ over $M$ which is locally of product type and
satisfies $R^V(X,Y)=0$ for all $X,Y\in TH$ 
admits a product structure. In fact there is a Dirac bundle $\cW$
on $N$ which is well-defined up to a canonical isomorphism and a
canonical isomorphism $\cW*(H,S(H))\cong \cV$.

\end{lem}
\proof
For $h\in H$ we define $\cW(h):=\cV//(H,h,S(H))$. 
We fix $h_0\in H$ and set $\cW:=\cW(h_0)$.
If $h\in H$ and $\gamma$ is a path in $H$ from $h_0$ to $h$,
then we have an isomorphism $\Phi(\gamma):\cW(h_0)\rightarrow \cW(h)$.
By our additional assumptions $\Phi(\gamma)$ only depends on the
endpoints of the path so that we can write
$\Phi(\gamma)=\Phi(h,h_0)$. Therefore, $\cW$ is independent of the
choice of $h_0$ up to a canonical isomorphism.

Let $W(h)*H$ be the bundle underlying $\cW(h)*H$.
We define an isomorphism
$\Psi(h):W(h)*H_{|N\times H\{h\}}\rightarrow V_{|N\times\{h\}}$
by the following canonical maps of tensor algebra.
If $H$ is even-dimensional, then
\begin{eqnarray*}
W(h)*H_{|N\times \{h\}}&=&\Hom_{\Cl(T_hH)}(S_h,V_{|N\times \{h\}})\otimes S_h\\
&\cong&V_{|N\times \{h\}}\ .
\end{eqnarray*}
If $H$ is odd-dimensional and $N$ is even-dimensional, then 
$V_{|N\times \{h\}}$ as a module of $\Cl(T_hH)$ is the sum  of the
iso-typic components $V_{|N\times \{h\}}^{\pm}$ of type $S_h$ and $S_h^{op}$. Therefore again
\begin{eqnarray*}
W(h)*H_{|N\times \{h\}}&=&\Hom_{\Cl(T_hH)}(S_h\oplus
S^{op}_h,V_{|N\times \{h\}})\otimes S_h\\
&\cong&V_{|N\times \{h\}}^+\oplus V_{|N\times \{h\}}^-\\
&\cong&V_{|N\times \{h\}}\ .
\end{eqnarray*}
If $H$ and $N$ are odd-dimensional, then
the grading $z$ of $\cV$ induces an isomorphism $V_{|N\times \{h\}}^+\cong
V_{|N\times \{h\}}^-$ (here $V_{|N\times \{h\}}^\pm\subset V_{|N\times \{h\}}$ still denote the isotypic components of $\Cl(T_hH)$-modules of type $S_h$ and $S_h^{op}$.). Therefore we get
 \begin{eqnarray*}
W(h)*H_{|N\times \{h\}}&=&
\Hom_{\Cl(T_hH)}(S_h,V_{|N\times \{h\}})\otimes S_h\otimes \C^2\\
&\cong&\Hom_{\Cl(T_hH)}(S_h,V_{|N\times \{h\}}^+)\otimes S_h\otimes
\C^2\\
&\cong&V_{|N\times \{h\}}^+ \otimes
\C^2\\
&\cong&V_{|N\times \{h\}}^+\oplus V_{|N\times \{h\}}^+\\
&\stackrel{\id\oplus z}{\cong}&V_{|N\times \{h\}}^+\oplus V_{|N\times \{h\}}^-\\
&\cong&V_{|N\times \{h\}}\ .
\end{eqnarray*}
We now define the isomorphism
$\Psi:W*H\rightarrow V$ by
$$\Psi_{|N\times\{h\}}:=\Psi(h)\circ (\Phi(h,h_0)*H)\ ,$$
where $(\Phi(h,h_0)*H):W*H\rightarrow W(h)*H$ is induced by
$\Phi(h,h_0)$.
It is easy to check that this isomorphism preserves the Dirac bundle
structures.

The isomorphism $\Psi$ is canonical in the following sense.
If $\cW^\prime$ is defined using $h_0^\prime$, and $\Psi^\prime$ is
the corresponding isomorphism, then we have
$$\Psi^\prime\circ (\Phi(h_0^\prime,h_0)*H)=\Psi\ .$$
\hB

\subsubsection{}\label{hefhwefwe}

We now discuss lifts of operators. Let $\cW$ be a Dirac bundle over
the Riemannian manifold $N$ and $H$ be a Riemannian spin manifold with
spinor bundle $S(H)$.
We consider the Dirac bundle $\cV:=\cW*H$ over $M:=N\times H$.
If $Q$ is an operator on $C^\infty(N,W)$, then
we want to define the operator $L_N^M(Q)$ on $C^\infty(M,V)$.
E.g, if $N$ is even-dimensional and  we take for $Q$ the Dirac
operator $D(\cW)$, then the lift $L_N^M(D(\cW))$ should be the part of
$D(\cV)$, which differentiates in the $N$-direction.
If $X\in C^\infty(N,TN)$ is a vector field and $Q:=c_W(X)$,
then we want $L_N^M(c_W(X))=c(\tilde X)$, where $\tilde X\in C^\infty(M,TM)$
is the field induced by $X$. We will actually define $L_N^M$ such that it is an algebra homomorphism (Lemma \ref{eiufduiewfw}).

Similar, but slightly different relations are required in the case where $N$ is odd-dimensional, see \ref{ioioewf} for details.

\subsubsection{}\label{ioioewf}

First we assume that $H$ and  $N$ are even-dimensional.
If $Q$ is an operator on $C^\infty(N,W)$, then it splits
into an even and an odd part $Q=Q^++Q^-$. We
define the operator $L_N^M(Q)$ on
$C^\infty(M,V)$ as follows.
If $f\in C^\infty(M,V)$ is of the form
$\phi\otimes s$ with $\phi\in C^\infty(N,W)$ and $s\in C^\infty(H,S)$,
then  we set $L_N^M(Q)f:=Q^+\phi\otimes s+Q^-\phi\otimes z_Ss$.

If $H$ is odd-dimensional and $N$ is even-dimensional, then we define
$L_N^M(Q)f:=Q\phi\otimes s$. 

If $N$ is odd-dimensional,
then in fact we define a lift of operators of the form $Q=Q_1\otimes 1 + Q_2\otimes \sigma$,
where $Q_i$ are operators on $C^\infty(N,W)$, and
$\sigma$ is the generator of $\Cl^1$
(the Clifford algebra of $(\R,-x^2)$) with $\sigma^2=1$.
The role of this additional Clifford generator $\sigma$ will become clear later in the context of local index theory.
One should think of $\cW$ as coming from a Dirac bundle on $\R\times N$ by restriction to $\{0\}\times N$, and that $\sigma$ is the Clifford multiplication by $i n $, where $n\in T_0\R$  is the unit vector in positive direction.

Let first $H$ be odd-dimensional.
If  $f\in C^\infty(M,V)$ is of the form
$\phi\otimes s\otimes v$ with $\phi\in C^\infty(N,W)$, $s\in C^\infty(H,S)$,
and $v\in \C^2$,
then  $$L_N^M(Q)f:=Q_1\phi\otimes s\otimes v + Q_2\phi\otimes s\otimes  \left(
\begin{array}{cc}0&i\\-i&0\end{array}\right)v\ .$$
If $H$ is even-dimensional, the we set\index{$L_N^M(Q)$}
$$L_N^M(Q)f:=Q_1 \phi\otimes s +Q_2\phi\otimes z_S s\ .$$

\subsubsection{}
The following Lemma is easy to check.
\begin{lem}\label{eiufduiewfw}
We have $L_N^M(Q\circ R)\cong L_N^M(Q)\circ L_N^M(R)$.
\end{lem}
\subsubsection{}
We introduce the following notation.
\begin{ddd}\label{Ldef}
We set 
\index{$\bL_N^M(Q)$}$\bL_N^M(Q):=L_N^M(Q)$ if $N$ is even-dimensional, and
$\bL_N^M(Q):=L_N^M(Q\otimes \sigma)$, if $N$ is  odd-dimensional.
\end{ddd}
\index{${}_H\bL^M_N$}Sometimes we write ${}_H\bL^M_N$ in order to  indicate the significance
of $H$. Note that   $\bL_N^M$ is not multiplicative for odd-dimensional $N$.

If $X\in C^\infty(N,TN)$, then we have
$\bL_N^M(c_W(X))=c(X)$ and  $L_N^M(\nabla_X^W)=\nabla^V_X$.
In particular, $\bL_N^M(D(\cW))$ is the part of the Dirac operator $D(\cV)$ which
differentiates in the $N$-direction.

%
%
%
%
%

\subsubsection{}
We now consider Riemannian spin manifolds $H_1, H_2$ with spinor bundles
$S(H_i)$. Furthermore let $H_1\times H_2$ be the product of Riemannian
spin manifolds with spinor bundle $S:=S(H_1\times H_2)$. 
Let $\cW_0,\cW_1$ be Dirac bundles on $N$ such that there is an isomorphism
$\Phi:(\cW_1*H_1)*H_2\rightarrow  \cW_0*(H_1\times H_2)$. Note that this implies
by Lemma \ref{hinher}
that
$\cW_0\cong \cW_1$.
\begin{lem}\label{ssignn}
Given an operator $Q_0$ on $C^\infty(N,W_0)$ (which is odd if $\dim(N)$ is even), 
 there is a unique
operator $Q_1$ on $C^\infty(N,W_1)$ such that
$${}_{H_1\times H_2}\bL_{N}^{M}(Q_0)\circ \Phi=\Phi\circ {}_{H_2}\bL_{N\times
  H_1}^{M}({}_{H_1}\bL_N^{N\times H_1}(Q_1))\ .$$  
\end{lem}
\proof
We construct $Q_1$ as follows. We fix points $h_i\in H_i$, $i=1,2$.
Then we obtain a canonical isomorphism
\begin{equation}\label{bnlo}\Big(\big((\cW_1*H_1)*H_2\big)//\big(H_2,h_2,S(H_2)\big)\Big)//\big(H_1,h_1,S(H_1)\big)\cong \cW_1\
.\end{equation}
The operator
$$A:=\Phi^{-1}\circ {}_{H_1\times H_2}\bL_{N}^{M}(Q_0)\circ \Phi$$
commutes with the multiplication by functions of the form
$\pr^*f$, $f\in C^\infty(H_1\times H_2)$, $\pr:N\times H_1\times
H_2\rightarrow H_1\times H_2$.
Therefore we can restrict $A$ to the sub-manifold $N\times
\{h_1\}\times \{h_2\}$. We denote this restriction by
$A_{(h_1,h_2)}$. This operator anticommutes with the Clifford multiplication by elements
of $T_{(h_1,h_2)}(H_1\times H_2)$ and is odd if $N$ is even-dimensional.

%
%
%
Using  the tensor structure (\ref{ewoifcewc}) we see through a case by case discussion that there exists a unique operator $Q_1\in \End(C^\infty(N,W_1))$ (which is odd if $N$ is even-dimensional) such that
$$({}_{H_2}\bL^M_{N\times H_1}({}_{H_1}\bL_N^{N\times H_1}(Q_1)))_{(h_1,h_2)}=A_{(h_1,h_2)}\ .$$  Furthermore, one checks that $Q_1$ is independent of the choice of the points $h_i\in H_i$.
 
Let us explain in detail the case where $N$ and $H_i$ are even-dimensional.
In this case we form $(\Phi^{-1}\circ (1\otimes z_S)\circ \Phi)\circ A_{(h_1,h_2)}$.
This operator commutes with the Clifford action of $T_{(h_1,h_2)}(H_1\times H_2)$  an therefore induces an operator
$Q_1$ in view of the formula
$W_1=\Hom_{\Cl(T_{(h_1,h_2)}(H_1\times H_2))}(S_{(h_1,h_2)},V_{|N\times \{(h_1,h_2)\}})$ (see (\ref{deudwedewdewdw})).
\hB

%
%
%
%
%
%
%
%
%
%
%
%
%
%
%
%

The following Lemma follows immediately from the definitions.
\begin{lem}
The isomorphism $\cW*(-H)\cong \cW^{op}*H$ in Lemma \ref{hinher}, 3., is such that ${}_{(-H)}\bL_N^M=-({}_H\bL_N^M)$.
\end{lem}

\subsection{Manifolds with corners}

\subsubsection{}
For our constructions we need a category of manifolds 
in which boundaries are allowed, and in which we can form  products.
This naturally leads to manifolds with corners of arbitrary
codimension. As a special case a manifold with boundary is a manifold
with  corners of codimension one.

For the analysis of Dirac operators the presence of product structures
leads to considerable simplifications. This motivates our very rigid
notion of a manifold with corners where product structures are part of
the data. E.g. in the special case of a manifold with a boundary this means
that a  collar neighborhood is part of the structure. Later we will consider
Riemannian metrics and Dirac bundles which are compatible with
the product structures.
One price to pay for our rigid notion of a manifold with corners is
that  it requires some work to show
that simple manifolds like  the $n$-dimensional simplex  carries such a structure (see Lemma  \ref{symy1}).

\subsubsection{}
\index{$N(k,U)$}\index{manifold with corners!corner model}A manifold with corners is locally modeled on spaces of the form 
$N(k,U):=U\times [0,1)^k$, $k\ge 0$, where $U\subseteq
\R^{m}$ is an open subset. A point in this set is usually denoted by
$(x,r_1,\dots,r_k)$, $x\in U$, $r_i\in [0,1)$, or by $(x,r)$ with $r=(r_1,\dots,r_k)$.
We call these spaces corner models, and the $r_i$ the normal coordinates.
An isomorphism $N(k,U)\rightarrow N(k,U^\prime)$
of corner models is a bijection of  the form
$(x,r_1,\dots,r_k)\mapsto
(x^\prime,r_{\sigma(1)},\dots,r_{\sigma(k)})$ for some diffeomorphism
$x\mapsto x^\prime$ and permutation $\sigma\in \Sigma^k$.

\subsubsection{}
The space $N(k,U)$ has a filtration \index{filtration!of corner models}
$$\emptyset=N^{-1}\subset N^0\subset N^1\subset N^2\subset \dots \subset N(k,U)\ ,$$
where for $i\ge 0$ we have $(x,r_1,\dots,r_k)\in N^i$ iff $\sharp\{l|r_l=0\}\le i$.
Let $(x,r)\in N^i\setminus N^{i-1}$. For simplicity we assume that
$r_{k-i+1}=\dots=r_k=0$. Then the neighborhood
$U\times (0,1)^{k-i}\times [0,1)^i\subset N(k,U)$ of $x$ is a corner model
$N(i,U\times (0,1)^{k-i})$ in a natural way
by considering $U\times (0,1)^{k-i}\subset \R^{m}\times \R^{k-i}\cong \R^{m+k-i}$.
Further, if  $W\subseteq U\times (0,1)^{k-i}$ is an open subset,
then $W\times [0,1)^i\subset N(k,U)$ is a corner model naturally
isomorphic to
$N(i,W)$.

\subsubsection{}
Let $i\in \{0,\dots,k\}$. By $I_i(k)$ we denote the set of $i$-element subsets of $\{1,\dots,k\}$.
For $j\in I_i(k)$ let $$\partial_j N(k,U) :=\{(x,r_1,\dots,r_k)\in N(k,U)\:|\: \forall l\in j\: :\:r_l=0\}$$
be the corresponding face of codimension $i$. Note that $\partial_j
N(k,U)$ can be considered as a corner model $N(k-i,U)$.
We define the interior of the face $\partial_j N(k,U)$ by
$$\partial_j N(k,U)^\circ :=\{(x,r_1,\dots,r_k)\in N(k,U)|  l\in j\:
\Leftrightarrow\:r_l=0\}\ .$$

\subsubsection{}
We now define the notion of a manifold with corners.
\begin{ddd}
\index{manifold with corners}A manifold with corners is a metrizable  space which is
locally homeomorphic to corner models such that the transition maps
are isomorphisms of corner models. 
\end{ddd}
\index{manifold with corners!chart}The precise meaning of this is the following. The structure of a
manifold with corners is given by an atlas, where a chart $(\phi,V)$ is a
homeomorphism
$\phi:V\rightarrow N(k,U)$ of $V\subseteq M$ with a corner model
$N(k,U)$. If $(\phi^\prime,V^\prime)$  is a second chart in the atlas, then the
compatibility condition is the following. We require that
(after permuting the normal coordinates in $\phi$ and $\phi^\prime$)
$\phi(V\cap V^\prime)=N(i,W)=W\times [0,1)^i\subset N(k,U)$
such that $W\subset U\times (0,1)^{k-i}$, $\phi^\prime(V\cap V^\prime)=N(i,W^\prime)=W^\prime\times [0,1)^i\subset N(k^\prime,U^\prime)$
such that $W^\prime\subset U^\prime\times (0,1)^{k^\prime-i}$, and
the transition map $\phi^\prime\circ \phi^{-1}:N(i,W)\rightarrow
N(i,W^\prime)$ is an isomorphism of corner models.

In this paper a manifold (with or without corners) may consist of several connected components of different dimensions. Such examples will play an important role in the local index theory later, see e.g.  \ref{ufiwfwef}.

\subsubsection{}

\index{manifold with corners!morphism}The local model of a morphism of manifolds with corners
is a morphism of corner models $N(k,U)\rightarrow
N(k^\prime,U^\prime)$ which comes as composition of maps of the
following two forms. The first is a map
of the form $(x,r)\mapsto
(x^\prime,\sigma(p(r),0))$ for some smooth map $x\mapsto x^\prime$ and
$\sigma\in \Sigma^{k^\prime}$, where
$p:[0,1)^k\rightarrow [0,1)^{k^{\prime\prime}}$ is the projection onto
a subset of coordinates and $(p(r),0)\in[0,1)^{k^\prime}$
is the extension of the coordinates of $p(r)$ by zeros in the last
$k^\prime-k^{\prime\prime}$ entries, where $k^\prime,k\ge
k^{\prime\prime}$. The first possibility does not mix the $x$ and
$r$-coordinates. But this is the case
for the second possibility which is the inclusion of
the form $N(i,W)\hookrightarrow N(k,U)$,
where $W\subset U\times (0,1)^{k-i}$ is open. 


\begin{ddd}\label{hdiqwdqwd}
A morphism of manifolds with corners is a continuous map
which is locally modeled by morphisms of corner models.
\end{ddd}

\subsubsection{}

\index{manifold with corners!product}Note that the category of manifolds with corners has a product.
The interval $[0,1]$ is a manifold with corners.
The reflection, the inclusion of the endpoints,  and the projection to a point are morphisms
of manifolds with corners.
This implies that the cube $[0,1]^k$ has a natural structure of a manifold with corners.
The inclusion of faces, the projection onto faces, permutations and
reflections of coordinates are morphisms of manifolds with corners. 

\subsubsection{}

\index{simplex as manifold with corners}It is an important fact that the $n$-simplex $\Delta^n$ also admits the structure
of a manifold with corners such that the symmetric group
$\Sigma^{n+1}$ acts by morphisms of manifolds with corners.
Let $\Delta^n\subset \R^{n+1}$ be the standard $n$-simplex consisting of all
points $x=(x_0,\dots,x_n)$ with $x_i\in [0,1]$ and $\sum_{i=0}^n x_i=1$.  
The standard basis of $\R^{n+1}$ coincides with the set of vertices of
$\Delta^n$. The permutation group $\Sigma^{n+1}$ acts on $\R^{n+1}$ by
permutation of coordinates. This action restricts to $\Delta^n$.
In this picture $\Delta^n$ has no natural structure of a manifold with
corners. In the induced euclidian geometry the faces do not meet in
right angles. Therefore in the following Lemma we use a different model.

\begin{lem}\label{symy1}
$\Delta^n$ admits a structure of a manifold with corners such that
$\Sigma^{n+1}$ acts by morphisms of manifolds with corners.
\end{lem}
\proof
We first define $\Delta^{n}$ as a $\Sigma^{n+1}$-invariant
$C^1$-sub-manifold of $\R^{n+1}$. Let $0\le k\le n$.
The part of $\Delta^{n}$ in the quadrant with
$x_i\ge 0$, $i\le k$ and $x_i< 0$, $i >k$ is defined by the equation
$x_0^2+\dots +x_k^2=1$, and $x_i\in [-1,0]$ for $i>k$. This prescription
induces equations in the remaining quadrants by
$\Sigma^{n+1}$-invariance.
We now introduce local coordinates.
Let $x^0=(x^0_0,\dots,x^0_{n})\in\Delta^n$ be a boundary point such that
$x^0_i\ge 0$ for $i\le k$ and $x^0_{i}<0$ for $i>k$. 
Assume further, that $j\le k\le l< n$, $x^0_i\not=1$ for $1\le i\le j-1$,    
$x^0_{j}\not=0$ and $x^0_{j+1}=\dots=x^0_k=0$, $-1<x^0_{i}<0$ for $i=k+1,\dots,l$,
$x^0_{l+1}=\dots=x^0_{n}=-1$.  Then coordinates around $x^0$ are
$x\mapsto
((x_1,\dots,x_{j-1}),(x_{j+1},\dots,x_k),(x_{k+1},\dots x_l),(1-x_{l+1},\dots,1-x_{n}))\in
(0,1)^{j-1}\times (-1,1)^{k-j}\times (-1,0)^{l-k}\times [0,\infty)^{n-l}$.
Using the action of  $\Sigma^{n+1}$ we get coordinates
around all boundary points. The coordinate transitions preserve normal
coordinates up to permutation.
In the interior we a have a $C^1$-smooth structure.
We obtain a $C^\infty$-structure by choosing an appropriate sub-atlas.
\hB


\subsubsection{}\label{dsuidsdsd}
A manifold with corners $M$ has a filtration
 $$\emptyset=M^{-1}\subset M^0\subset M^1\subset \dots \subset M^k\subset \dots \subset M$$
\index{manifold with corners!filtration} \index{filtration!of a manifold with corners}such  $x\in M^i\setminus M^{i-1}$ if it admits a pointed
neighborhood which is modeled by  a pointed corner model
$(x_0,0)\in N(i,U)$. 

We now discuss the decomposition of a manifold with corners into
faces. Using the notion of faces we then will impose further restrictions on the global structure of
$M$. 

Let $A^\circ \subseteq M^{i}\setminus M^{i-1}$ be a connected component.
We construct a completion $A$ of $A^\circ$, which is again a manifold
with corners. Let $\bar A$ be the closure of $A^\circ$ in $M$.
In general it is not a manifold with corners. Consider for example
a two-dimensional disc $M$ with filtration $M^0\subset M^0\sqcup A^\circ\subset M$ (the one-eck). The boundary $\bar A$ of $M$  meets
itself in one corner, it is a broken circle. 
 But we want to complete  $A^\circ$ to a manifold $A$ with corners isomorphic to the interval, i.e. with two boundary points. The formal construction goes as follows.

If $x\in \bar A\setminus A^\circ$, then a neighborhood of $x$
has a local model $N(l,U)$ for some $l> i$ and connected $U$. There is a subset
$Q\subseteq I_{i}(l)$ such that $A^\circ \cap N(l,U)=\sqcup_{q\in
  Q}
\partial_q N(l,U)^\circ$.
Locally near $x$, the completion $A$ of $A^\circ$ is given by
$$A\cap N(l,U):=A^\circ\cap N(l,U)\cup_{\sqcup_{q\in Q}\partial_q N(l,U)^\circ}\sqcup_{q\in
  Q}\partial_q N(l,U)\ .$$
The local completions glue nicely. There is a surjective map
$A\rightarrow \bar A$, which is an isomorphism over $A^\circ$.
\begin{ddd}
\index{face!atom}An atom of faces of codimension $i$ is a completion of 
a connected component of $M^i\setminus M^{i-1}$ as constructed above. \index{$I^{atoms}_i(M)$}
By $I_i^{atoms}(M)$ we denote the set of atoms of faces of codimension
\index{face}$i$.
A
face of codimension $i$ is a subset of $I^{atoms}_i(M)$.
 \end{ddd}
Note that an atom of faces has an induced structure of a manifold
with corners. 

\subsubsection{}
The interval $[0,1]$ has two atoms of faces of codimension one, namely
the points $\{0\}$ and $\{1\}$. We can form the faces
$\{0\}$, $\{1\}$, $(\{0\}, \{1\})$.
The square $[0,1]^2$ has four atoms of faces of codimension one.
We could form e.g. the faces $\{0\}\times [0,1]$ or
$(\{0\}\times [0,1] , [0,1]\times \{1\})$.
The second face will not be allowed in an admissible face
decomposition later on, because it ``self-intersects'' at the 
corner point $(0,1)$ of codimension $2$.

\subsubsection{}
\begin{ddd}
\index{$\partial_mM$}\index{face!boundary face}If $m\subseteq I_i^{atoms}(M)$, then let $\partial_mM$ denote the
isomorphism class of pairs $(N,f)$, where $N$ is a manifold with
corners and $f:N\rightarrow \sqcup_{A \in m} A$ is a
\index{model of $\partial_mM$}diffeomorphism. A model of $\partial_m M$ is a representative of the isomorphism class.
\end{ddd}
If $(N,f)$ is a model of $\partial_jM$, $j\in I_k(M)$, then we have natural
maps $I_l^{atoms}(N)\rightarrow I_{l+k}^{atoms}(M)$, $l\ge 0$.

\subsubsection{}
Our reason to distinguish between faces and atoms of faces is the
following. Later we introduce perturbations of Dirac operators by
integral operators with smooth integral kernel. While Dirac operators
are local, integral operators are not. We want to consider integral
operators which are localized on faces, but not necessarily on atoms
of faces.

\subsubsection{}\label{mottuw4}
At the present stage the distinction between the isomorphism class
$\partial_jM$ and its models looks unnecessary complicated.
In fact, two models of $\partial_jM$ are isomorphic by a
unique isomorphism. But later we will consider faces with additional structures
like Dirac bundles
which allow for non-trivial automorphisms. Then this distinction will be unavoidable. For more motivation look at  \ref{zzuwwedd} where will introduce the notion of a distinguished model of a boundary face of a geometric manifold.
The main point of the discussion there is how to find a canonical lift of an isomorphism of the underlying manifolds-with-corner models of the boundary face 
to an isomorphism of the induced Dirac bundles. This lift is crucial if one  wants to formulate boundary value problems  for the Dirac equation properly. 
For further motivation we mention Lemma \ref{hiphop} which becomes important if one wants to understand  the compatibility of boundary conditions along faces
of codimension one at corner points, i.e. along faces of codimension two.
In the present paper the concept of a (pre-)taming (see \ref{tamtam1234}) is a replacement of the concept of boundary conditions, and the motivating remarks above apply to tamings as well.

\subsubsection{}
In the local model of $M$ near a point in $M^i\setminus M^{i-1}$
there are exactly $\binom{i}{j}$ atoms of faces of codimension $j$ of
the model
which meet in this point. Our analysis near corner points is based on
separation of variables. Here it is crucial that these $\binom{i}{j}$
local atoms belong to different faces of $M$. 

\subsubsection{}
We now introduce the notion of an admissible face decomposition.
\begin{ddd}
\index{face!face decomposition}A face decomposition of $M$ is given by partitions
$I_i(M)$ of $I_i^{atoms}(M)$ for all $i\ge 0$. \index{$I_i(M)$}
\end{ddd}
An element of $I_i(M)$ is thus a face of $M$.
We consider the following additional conditions. 
\begin{enumerate}
\item
The first condition is that for all $i\ge 0$ and $j\in I_i(M)$ the natural map
$f:N\rightarrow M$ for each model $(N,f)$ of $\partial_jM$
is an inclusion.
\item Under the first condition and if $(N,f)$ is a model of
  $\partial_jM$, $j\in I_i(M)$,  the natural map
$I^{atoms}_{m}(N)\rightarrow I^{atoms}_{m+i}(M)$ is an inclusion.
Let $I_m^{atoms}(\partial_j M)\subseteq I_{m+i}^{atoms}(M)$ denote the image which is independent
of the choice of the model.
The second condition is that for all $i\ge 0$, $j\in I_i(M)$, and 
$m\ge 0$ the subset $I^{atoms}_{m}(\partial_jM) \subseteq I_{m+i}^{atoms}(M)$
is a union of elements of $I_{m+i}(M)$.
\end{enumerate}
\begin{ddd}
\index{face decomposition!admissible}We call a face decomposition admissible if it satisfies the conditions
1. and 2. formulated above.
\end{ddd}
\subsubsection{}
The interval $[0,1]$ has two admissible face decompositions.
On the one hand we have the atomic face decomposition. On the other
hand we have the face decomposition with only one face of codimension
one.
\subsubsection{}
The existence of an admissible face decomposition 
 imposes global
restrictions on $M$ since in particular it implies, that
at each point in $M^{i}\setminus M^{i-1}$ meet exactly
$\binom{i}{j}$ atoms of faces of codimension $j$.
This is equivalent to the condition that for all atoms $A$ of faces the
canonical map $A\rightarrow M$ is an inclusion.
Consider e.g. a two-manifold $M$ (homeomorphic to a disk) with one boundary component, which
self-intersects at a corner point of codimension two.
This manifold can be equipped with the structure of a manifold with
corners (the one-eck). It has one atom $A$ of faces of codimension one which is
isomorphic to the interval. In this case the canonical map
$A\rightarrow M$ is not an inclusion since it identifies the endpoints
of the interval. The one-eck does not admit any admissible face decomposition.

\subsubsection{}\label{delk}
\begin{ddd}
We call a face decomposition of $M$ reduced if $I_0(M)$ consists of
one element. A face decomposition is called atomic if all faces are
atoms. A face decomposition $I^\prime_*(M)$ is called a refinement of
$I_*(M)$ if for all $i\ge 0$ the partition $I_i^\prime(M)$ refines
$I_i(M)$.  \index{face decomposition!atomic}
\index{face decomposition!refinement}
\end{ddd}

Note that the atomic face decomposition is the finest face decomposition. 
If $M$ admits an admissible face decomposition, then the atomic face
decomposition is admissible, too. 

Given any  face decomposition for $M$  we can form the associated \index{face decomposition!reduced}reduced face
decomposition by replacing the partition $I_0(M)$ by the trivial partition
consisting of one element. If the given face decomposition was
admissible, then so is its reduction. 

\subsubsection{}

If $M$ and $N$ are manifolds with corners and given admissible face
decompositions, then the union $M\sqcup N$ is a manifold with corners
which has an induced admissible face decomposition. Here $I_i(M\sqcup
N)$ is the partition of $I^{atoms}_i(M\sqcup N)\cong
I^{atoms}_i(M)\sqcup
I^{atoms}_i(N)$ induced by the partitions $I_i(M)$ and $I_i(N)$.

\subsubsection{}\label{sumdelk}
Assume that $M$ comes with an admissible face decomposition.
Let $i\ge 0$ and $j\in I_i(M)$.
Each model $(N,f)$ of
the manifold with corners $\partial_jM$ has an induced reduced face
decomposition. We identify $I^{atoms}_m(N)$ with the subset
$I^{atoms}_{m}(\partial_jM)\subseteq I^{atoms}_{m+i}(M)$.
Then the partition 
$I_{m}(\partial_j N)$ is the partition of 
$I^{atoms}_m(\partial_j M)$ which is induced from the partition of
$I_{m+i}(M)$. The induced face decomposition of $N$
is again admissible. 
We let $I_m(\partial_jM)\subseteq I_{m+i}(M)$ be the image of
$I_m(N)$ under the natural inclusion $I_m(N)\rightarrow
I_{m+i}(M)$ which is independent of the
choice of the model.

For $l\in I_m(\partial_j M)$ we have an equality
$\partial_l \partial_j M=\partial_l M$
which is given by a canonical isomorphism of models
together with an equality of the induced face decompositions.

\subsubsection{}
The consideration of faces of codimension one in faces of codimension
one leads to the notion of adjacent faces. Let  $i\in I_1(M)$
 and $j\in I_1(\partial_i M)$. Then we can view $j\in I_2(M)$ in a
 natural way. 
Let $f_i:N_i\rightarrow M$ and $f_j:N_j\rightarrow M$ be models of $i$ and $j$.
There is a unique $i^\prime\in
I_1(M)$, $i\not=i^\prime$, represented by a model
$f_{i^\prime}:N_{i^\prime}\rightarrow M$ such that $f_j(N_j)\subseteq 
f_i(N_i)\cap f_{i^\prime} (N_{i^\prime})$.

\begin{ddd}\label{conjface}
\index{face!adjacent}We call $i^\prime$ the adjacent face to $i$ with respect to $j$.
\end{ddd}

\subsubsection{}
If $M$ has an admissible face decomposition, then
there is a natural partial order on the set of faces.
Let $k$ and $j$ be faces of $M$. Let $(N,f)$ and $(N^\prime,f^\prime)$ be
models of these faces.
\begin{ddd}
We say that $k\le j$ if $f(N)\subset f^\prime(N^\prime)$ 
\end{ddd}

\subsubsection{}

The total boundary of $M$ is denoted by $\partial M$. It is
the  disjoint union $\bigsqcup_{i\in
  I_1(M)}\partial_i M$ with its induced face decomposition.
Note that $I_0(\partial M)=I_1(M)$. 
\subsubsection{}
If $M$ and $N$ are manifolds with corners, then $M\times N$ has a
natural structure of a manifold with corners. 
The atoms of faces of
codimension $i$ of $M\times N$ are given by $A\times B$, where
$A$ (and $B$) are atoms of faces of codimension $m$ of $M$ (resp. $n$ of $N$)
such that $m+n=i$.
Thus $I^{atoms}_i(M\times N)=\sqcup_{m+n=i} I^{atoms}_m(M)\times
I^{atoms}_n(N)$. Face decompositions of $M$ and $N$ induce a face
decomposition of $M\times N$ in a natural way.
If the given face decompositions of $M$ and $N$ were admissible, then
so is the induced face decomposition of $M\times N$.



\subsection{Orientations, metrics, and Dirac bundles on manifolds with corners}

\subsubsection{}\label{setupor}
Let $GL(n,\R)^0\subset GL(n,\R)$ denote the connected component of the
identity.
Let $V\rightarrow X$ be a real $n$-dimensional vector bundle over some
space
\index{frame bundle}\index{$\Fr(V)$}$X$. The frame bundle $\Fr(V)\rightarrow X$ is a principal bundle with
\index{orientation}structure group $GL(n,\R)$.
An orientation of $V$ is a reduction of the structure group of $\Fr(V)$ to
 $GL(n,\R)^0$, i.e.   an isomorphism class of pairs
$(P,\Phi)$, where $P\rightarrow X$ is a $GL(n,\R)^0$-principal bundle 
 and $\Phi$ is an  isomorphism of principal bundles
 $\Phi:P\times_{GL(n,\R)^0}GL(n,\R)\stackrel{\sim}{\rightarrow}
 \Fr(V)$. Note that via $\Phi$ we can and will consider $P$ as the
 sub-bundle $\Fr^+(V)\subset \Fr(V)$ of oriented frames.

If $f:X^\prime\rightarrow X$ is a continuous map,
then an orientation $(P,\Phi)$ of $V$ induces an orientation
$(f^*P,f^*\Phi)$ of $f^*V$ in a natural way.
 
\subsubsection{}\label{orired}
Assume that $V$   fits into an exact sequence of vector bundles
$$0\rightarrow W\rightarrow V\rightarrow X\times \R\rightarrow 0\ .$$
Let $s:X\times \R\rightarrow V$ be a split.  Consider
 $GL(n-1,\R)\subset GL(n,\R)$ as  the subgroup which fixes
the last basis vector of the standard basis of  $\R^n$ and the subspace spanned by the first $n-1$ basis vectors.
Then the split $s$ induces a reduction of the structure group of $\Fr(V)$
to $GL(n-1,\R)$ which is represented by the sub-bundle 
$\Fr(W)\subset \Fr(V)$. This inclusion maps a frame 
$(w_1,\dots,w_{n-1})$ of $W$ to the frame $(w_1,\dots,w_{n-1},s(1))$
of $V$.
We define the $GL(n-1,\R)^0$-principal bundle $\Fr(W)^+:=\Fr(V)^+\cap
\Fr(W)$. The sub-bundle $\Fr(W)^+\subset \Fr(W)$ 
represents the orientation of $W$ induced by the split $s$.

 
\subsubsection{}
An orientation of a manifold with corners $M$ is by definition an
orientation of $TM$. We choose a model $(R,i)$ of $\partial M$.
The pull-back $i^* TM$
sits in a natural exact sequence
$$0\rightarrow TR\rightarrow i^*TM\rightarrow N\rightarrow 0\ ,$$
where $N$ is the normal bundle.
We choose a section of $i^*TM$ 
consisting of inward pointing tangent vectors.
This section induces a trivialization $N\cong \partial M\times \R$
and a split $s:N\rightarrow i^*TM$.
\index{orientation!induced}The construction \ref{orired} produces the induced orientation of $TR$
and hence of $R$.

\subsubsection{}
This convention is compatible with the upper half space model.
Indeed, if we equip the upper half space $\R^n_+:=\R^{n-1}\times
[0,\infty)\subset \R^n$
with its standard orientation, then the orientation induced on
$\R^{n-1}$ is again the standard orientation. Note that with this
convention we have Stoke's theorem in the following form: For
$\omega\in C^\infty_c(\R^n_+,\Lambda^{n-1}T^*\R^n_+)$ we have
$$\int_{\R^n_+} d\omega= (-1)^n\int_{\R^{n-1}} \omega_{|\R^{n-1}}$$
(see also \ref{hhhas}).

\subsubsection{}

Next we introduce a class of Riemannian metrics on a manifold with
corners which is compatible with the local product structures.
\index{Riemannian metric!admissible}A Riemannian metric on a corner model $N(k,U)$ is admissible, if it is
equal to the  product
metric $g^U\oplus g^{[0,1)^k}$ in a neighborhood of the singular
stratum $U\times \{0\}$, where $g^{[0,1)^k}$
is the standard metric on $[0,1)^k$ and $g^U$ is some metric on $U$.
Note that isomorphisms of corner models preserve this structure of
metrics.
\begin{ddd}\label{udiwqdwqdq}
A Riemannian metric on a manifold with corners is admissible if it
locally induces admissible metrics on corner models.
\end{ddd}
\subsubsection{}
Note that any manifold with corners admits admissible Riemannian
metrics. In order to construct such a metric we glue given admissible
metrics on local corner models using an adapted partition of unity. 
Let us explain this in greater detail. We fix a locally finite
covering of $M$ by corner models $V_i\cong N(k_i,U_i)$. On each corner
model $V_i$ we choose a smooth cut-off function $\chi_i$ of compact
support such that $\chi_i(x,r)$ is independent of $r$ for small $r$.
We can assume that $\sum_i\chi_i>0$. We define the adapted partition of
unity  by $\rho_i:=(\sum_j\chi_j)^{-1}\chi_i$.
Using this partition of unity we can glue admissible Riemannian metrics
on the $V_i$ in order to obtain an admissible Riemannian metric
on $M$.

An admissible Riemannian metric on $M$ induces an admissible
Riemannian metric on all atoms of faces of $M$ by 
restriction.  

\subsubsection{}
\index{$\bar M$}\index{extension!of a manifold with corners}If $M$ is a manifold with corners, then we can form the extension $\bar M$
which is a smooth manifold without singularities containing $M$.
We first introduce the extension $\bar N(k,U):=U\times (-\infty,1)^k$
of the corner model $N(k,U)$. The inclusion $ [0,1)^k\rightarrow
(-\infty,1)^k$ induces an inclusion $N(k,U)\rightarrow \bar N(k,U)$.
Note that an isomorphism of corner models extends naturally to a
diffeomorphism of extensions. If $x\in N(k,U)^i\setminus
N(k,U)^{i-1}$
and $N(i,U^\prime)\subset N(k,U)$ is a neighborhood of $x$ which is
itself  a corner model, then there is a natural inclusion of
extensions $\bar N(i,U^\prime)\rightarrow \bar N(k,U)$.
Thus the transition maps between local charts of a manifold with
corners have natural extensions to the extended  chart domains. 
\begin{ddd}
We define the extension of $\bar M$
by replacing in all charts the corner models by their extensions
which are glued using the natural extensions of transition maps.
\end{ddd}
In more detail we consider a covering of $M$ by charts 
$V_i\cong N(k_i,U_i)$  which are isomorphic to corner models.
Then we have $M\cong \sqcup_i N(k_i,U_i)/\sim$, where the equivalence
relation relates points connected by a transition map.
Then we have $\bar M\cong \sqcup_i \bar N(k_i,U_i)/\sim$, where the
equivalence relation now relates points which are connected by a
canonical extension of a transition map.
Note that $\bar M$ is (non-canonically) diffeomorphic to $M^0$.

\subsubsection{}
The extension of the interval $[0,1]$ is the real line $\R$.
More general, the extension of the cube $[0,1]^k$ is $\R^k$.
If $M$ is a manifold with boundary $\partial M$ (equipped with the structure
of a  manifold with
corners), then its extension is 
$M\cup_{\partial M\cong \partial M\times\{0\}}\partial M\times (-\infty,0]$.
If $A$ is an atom of faces of $M$, then the canonical map
$A\rightarrow M$ extends to the extensions $\bar A\rightarrow \bar M$.

\subsubsection{}

In order to define a regularized trace on a certain class of integral
operators on $\bar M$
we need the following exhaustion $({}_r\bar M)_{r\ge 0}$ of $\bar M$
by compact subsets. 
For a corner model we set ${ }_r\bar N(k,U):=U\times [-r,1)$.
\begin{ddd}
\index{${}_r\bar M$}For $r\ge 0$ we define ${}_r\bar M\subset \bar M$ as the subset 
which is characterized by the property that in each chart  the inclusion ${}_r\bar M\rightarrow \bar M$
corresponds to the inclusion  ${}_r\bar N(k,U)\rightarrow \bar N(k,U)
$.
\end{ddd}

\subsubsection{}
E.g, we have ${}_r\overline{[0,1]}=[-r,r+1]\subset
\overline{[0,1]}\cong \R$.
If $A$ is an atom of faces of $M$ such that the canonical map
$A\rightarrow M$ is an inclusion, then we have
${}_r\bar A=\bar A\cap {}_r \bar M$ under the inclusion $\bar
A\rightarrow \bar M$ (note that $\bar A$ is the extension of $A$ and not the closure of $A^\circ$ in $M$ as in \ref{dsuidsdsd}).

\subsubsection{}
Assume that  $M$ has an admissible face decomposition.
Let
$k\ge 0$ and $j\in I_k(M)$ be a face of codimension $k$ of $M$.
Let $(N,f)$ be a model of $\partial_jM$. The embedding
$f:N\rightarrow M$ extends naturally to an embedding $\bar f:\bar
N\rightarrow \bar M$. If 
$N\times [0,1)^k\rightarrow M$ is a collar neighborhood
induced by the structure of a manifold with corners, then it extends
 naturally to an open embedding $\bar N \times
 (-\infty,1)^k\rightarrow \bar M$. In fact, these embeddings are
\index{distinguished!embedding} determined up to a permutation of the coordinates of $(-\infty,1)^k$.
We will call an embedding of this form a distinguished embedding.
We define the open set $U_j\subset \bar M$ to be the image of $\bar
N\times (-\infty,0)^k$. We call $U_j$ the cylinder over the face $j$
in $\bar M$.

\subsubsection{}
The product metric $g^U\oplus g^{[0,1)^k}$ on $N(k,U)$ extends
naturally to a
product metric $g^U\oplus g^{(-\infty,1)^k}$ on $\bar N(k,U)$.
Thus any admissible metric on $N(k,U)$ extends naturally to $\bar
N(k,U)$. More general,
if $g^M$ is an admissible Riemannian metric on $M$, then it extends
naturally to $\bar M$. If $M$ is compact, then
$\bar M$ becomes a complete Riemannian manifold.

\subsubsection{}

Let $M$ be a manifold with corners and $g^M$ be an admissible
Riemannian metric.
Let $\cV$ be a Dirac bundle over $M$.
\begin{ddd}
\index{Dirac bundle!admissible}We say that $\cV$ is admissible, if for each chart
$f:N(U,k)\rightarrow M$
there exists a neighborhood of the singular stratum $U\times \{0\}$
on which $f^*\cV$ is locally of product type and satisfies
$R^{f^*V}(X,Y)=0$ for all $X,Y\in T[0,1)^k$.
\end{ddd}

\subsubsection{}
If $f:N(U,k)\rightarrow M$ is a local chart of $M$, then by
$\bar f:\bar N(U,k)\rightarrow \bar M$ we denote its canonical
extension to the extension.

\begin{lem}
\index{extension!of a Dirac bundle}\index{Dirac bundle!extension}Assume that $\cV$ is admissible.
Then there exists a Dirac bundle $\bar \cV$ over $\bar M$ which
extends $\cV$, and which is for all charts $\bar f:\bar
N(U,k)\rightarrow \bar M$
locally of product type  and satisfies
$R^{\bar f^*\bar V}(X,Y)=0$ for all $X,Y\in T(-\infty,0)^k$ over
$\bar f(U\times (-\infty,0)^k)$.
The extension $\bar \cV$ is unique up to a canonical isomorphism.
\end{lem}
\proof
Let $f:N(U,k)\rightarrow M$ be a local chart of $M$. After reducing $U$
 there exists an $\epsilon>0$ such that $f^*\cV_{|U\times (0,\epsilon)^k}$ 
satisfies  the curvature assumptions.
By Lemma \ref{loglo} there exists a Dirac bundle $\cW$ on $U$ and an
isomorphism $\Phi:f^*\cV_{|U\times (0,\epsilon)^k}\rightarrow \cW*(0,\epsilon)^k$.
We define the extension of $f^*\cV$ to the extension $\bar N(U,k)$ by glueing
$\cW*(-\infty,\epsilon)^k$ with $f^*\cV$ over $U\times (0,\epsilon)^k$
using the isomorphism $\Phi$. This provides a Dirac bundle $\bar\cV_f$ over
$\bar f(\bar N(U,k))$.
Assume that $N(U^\prime,k^\prime)$ is another local chart
with $f^\prime(N(U^\prime,k^\prime))\subset f(N(U,k))$.  
 The same construction applied to this chart provides
the bundle $\bar \cV_{f^\prime}$ over $\bar f^\prime(\bar
N(U^\prime,k^\prime))$. We have
a canonical isomorphism $\bar\cV_{f|f^\prime(N(U^\prime,k^\prime))}\cong \cV_{|f^\prime(N(U^\prime,k^\prime))}\cong \bar\cV_{f^\prime|f^\prime(N(U^\prime,k^\prime))}$.
We extend this isomorphism to an isomorphism $\bar \cV_{f|\bar f^\prime((\bar
N(U^\prime,k^\prime)))}\cong \bar \cV_{f^\prime}$ using parallel
transport along paths with constant first coordinate (i.e. along path
moving in the
normal directions).
Using these isomorphisms we want to glue the locally defined extensions
$\bar \cV_f$ of the
Dirac bundle $\cV$. In fact, if
$f^{\prime\prime}:N(U^{\prime\prime},k^{\prime\prime})\rightarrow M
$ is a third chart with image contained in $N(U^\prime,k^\prime)$,
then must check a cocycle relation over $\bar f^{\prime\prime}(\bar
N(U^{\prime\prime},k^{\prime\prime}))$. But this relation obviously
holds true over
$f^{\prime\prime}(N(U^{\prime\prime},k^{\prime\prime}))$
and therefore also over the extension.

Assume, that $\cU$ is another Dirac bundle which satisfies the
assumptions of an extension of $\cV$ to $\bar M$. Then we can extend
the canonical isomorphism $\cU_{|M}\cong \cV\cong \bar \cV_{|M}$ to an
isomorphism $\cU\cong \bar \cV$ again using parallel transport in the
normal directions.
\hB

\subsubsection{}\label{zzuwwedd}

\begin{ddd}
\index{geometric!manifold}\index{$\cM_{geom}$}\index{$(M,g^M,\orient,\cV)$}A geometric manifold $\cM_{geom}$ is a tuple $(M,g^M,\orient,\cV)$, where
\begin{enumerate}
\item $M$ is a manifold with corners and fixed admissible face decomposition,
\item $g^M$ is an admissible Riemannian metric,
\item $\orient$ is an orientation of $M$,
\item $\cV$ is an admissible Dirac bundle on $M$.
\end{enumerate}  
\end{ddd}
Let $M_i$, $g^{M_i}$, $i=0,1$, be Riemannian manifolds.
A diffeomorphism $f:M_0\rightarrow M_1$ is an isometry if $f^*g^{M_1}=g^{M_0}$.
\index{isomorphism!of geometric manifolds}Let $\cM_{geom, i}$ be geometric manifolds. An isomorphism
$(f,f_{\cV}):\cM_{geom, 0}\stackrel{\sim}{\rightarrow} \cM_{geom, 1}$ is an
isometry  $f$ of the underlying  Riemannian manifolds which
preserves
the admissible face decompositions and orientations together with an
isomorphism of Dirac bundles $f_\cV:f^*\cV_1 \stackrel{\sim}{\rightarrow} \cV_0$.

If $\cM_{geom}$ is a geometric manifold and
$i\in I_k(M)$, $k\le 1$, then we will  define an
isomorphism class of 
geometric manifolds $\partial_i\cM_{geom}$.
If $k=1$, then  $\partial_i\cM_{geom}$ is called a boundary face.
In fact, the construction will be more rigid. If $k=0$, then
$\partial_i\cM_{geom}$ has a canonical model. If $k=1$, then
we represent  $\partial_i\cM_{geom}$ by distinguished models which are uniquely
determined up to a canonical isomorphism. Furthermore, if we fix a
distinguished model $(N,g^N,\orient_N,\cW)$, then
we have a distinguished embedding $T:\bar N\times (-\infty,0)\rightarrow
\bar M$ and a canonical product structure $\Pi:T^*\bar \cV\rightarrow
\bar \cW*(-\infty,0)$. This rigidity will be crucial
later where we want to lift additional structures (e.g. a taming) from
$\partial_i\cM_{geom}$ to $\cM_{geom}$ (see also \ref{mottuw4}).

First we consider the case $k=0$, i.e.  
\index{$\partial_i\cM_{geom}$} $i\in I_0(M)$. Then $\partial_i\cM_{geom}$ has a canonical model
which has as underlying manifold with corners and admissible face
decomposition the canonical model of $\partial_iM$, and the other
geometric structures are obtained by restriction.

\index{distinguished!model}If $i\in I_1(M)$, then a distinguished model of
$\partial_i\cM_{geom}$ is defined as follows. We take any model $(N,f)$ of
$\partial_i M$ which is uniquely defined up to a canonical isomorphism.
We equip $N$ with the  orientation $\orient_N$ and Riemannian metric
$g^N$ induced by $f$. These structures are preserved by the canonical
isomorphisms of the models. The model comes with an embedding
$T:\bar N\times (-\infty,0)\rightarrow \bar M$ which is an orientation
preserving isometry. The pull-back $T^*\bar \cV$ satisfies the assumptions
of Lemma \ref{loglo}. Recall that we have fixed once and for all a
spinor bundle $S(\R)$ and hence, by restriction, $S((-\infty,0))$.
Therefore, we have a Dirac bundle $\bar \cW$  over
\index{$\Pi$}\index{Dirac bundle!product structure}$\bar N$ (we set $\cW:=\bar \cW_{|N}$) which is well-defined up to a canonical isomorphism.
We further get a product
structure $\Pi:T^*\bar\cV\rightarrow \bar \cW*(-\infty,0)$.  
We let $(N,g^N,\orient_N,\cW)$ be the distinguished model of
$\partial_i\cM_{geom}$. 

If $(N^\prime,f^\prime)$ is another model of $\partial_iM$, and
$g:N^\prime\rightarrow N$ induces the canonical isomorphism,
i.e. $f\circ g=f^\prime$, then we obtain the canonical isomorphism
$G:g^*\cW\rightarrow\cW^\prime$ as follows. In order to choose a specific
representative of the reductions $\cW$ and $\cW^\prime$ we take the
base point $-1\in (-\infty,0)$. Let us abbreviate
$\dots//((-\infty,0),-1,S((-\infty,0)))$ by $//(-\infty,0)$ for the moment.
Then  $G:g^*\cW\rightarrow \cW^\prime$
is given by 
\begin{eqnarray*}g^*\cW&=&g^*(T^{*}\bar \cV//(-\infty,0))\\& 
\stackrel{g^\sharp}{\rightarrow} &(g\times \id_{(-\infty,0)})^*T^*
\bar \cV//(-\infty,0)\\&\stackrel{T^\prime=T\circ (g\times
  \id_{(-\infty,0)})}{\cong}&T^{\prime *} \bar
\cV//(-\infty,0)\\&=&\cW^\prime \ .
\end{eqnarray*}
Let $\Pi^\prime:T^{\prime*}\bar \cV\rightarrow \cW^\prime*(-\infty,0)$
be the product structure associated to the second model. Then we have
the following compatibility of product structures.
\begin{lem}\label{ffhu}
The composition
\begin{eqnarray*}
T^{\prime *}\bar \cV&\stackrel{T^\prime=T\circ (g\times
  \id_{(-\infty,0)})}{\rightarrow}&(g\times \id_{(-\infty,0)})^*T^{*}\bar \cV\\
&\stackrel{(g\times \id_{(-\infty,0)})^*\Pi}{\rightarrow}&(g\times
\id_{(-\infty,0)})^*(\cW*(-\infty,0))\\
&\stackrel{g_\sharp}{\rightarrow}&g^*\cW*(-\infty,0)\\
&\stackrel{G*(-\infty,0)}{\rightarrow}& \cW^\prime*(-\infty,0)
\end{eqnarray*}
coincides with the product structure $\Pi^\prime$.
\end{lem}
\proof
This equality can be checked in a straight forward manner.
\hB

\subsubsection{}

\begin{ddd}\label{bounfgg}
For $k\le 1$ and $i\in I_k(M)$ the isomorphism class of geometric manifolds
$\partial_i\cM_{geom}$  is the class
represented by a distinguished model as constructed above, where two
distinguished models are related by the canonical isomorphism also
constructed above.
\end{ddd}

\subsubsection{}
Note that we can not define $\partial_j\cM_{geom}$ for $j\in I_k(M)$, $k\ge
2$ in a similar manner since e.g. there is no canonical orientation of
$\partial_jM$ (see Lemma \ref{hiphop}).

\subsubsection{}
\begin{ddd}
\index{opposite!geometric manifold}
\index{$\cM_{geom}^{op}$}The opposite $\cM^{op}_{geom}$ of the geometric manifold $\cM_{geom}=(M,g^M,\orient,\cV)$ is defined by
$\cM_{geom}^{op}:=(M,g^M,-\orient,\cV^{op})$.
\end{ddd}
If $(f,F):\cM_{1,geom}\rightarrow \cM_{2,geom}$ is an isomorphism of
geometric manifolds, then it induces an isomorphism
$(f,F)^{op}:\cM_{1,geom}^{op}\rightarrow \cM_{2,geom}^{op}$ in a
natural way.

\subsubsection{}

Let $\cM_{geom}$ be a geometric manifold $i\in I_1(M)$, $k\in
I_1(\partial_iM)$ and $j\in I_1(M)$ be adjacent to $i$ with respect
to $k$. Then we also have $k\in I_1(\partial_jM)$.

\begin{lem}\label{hiphop}
We have
$$\partial_k\partial_i\cM_{geom}=
(\partial_k\partial_j\cM_{geom})^{op}\ .$$
\index{canonical isomorphism of corner models}In fact, if we choose distinguished models $\cN_{ki,geom},\cN_{kj,geom}$
of both classes, then this equality is implemented by a canonical
isomorphism $(g,L):\cN_{ki,geom}^{op}\rightarrow \cN_{kj,geom}$.
The canonical isomorphism
$(g^\prime,L^\prime):\cN_{kj,geom}^{op}\rightarrow \cN_{ki,geom}$ is
the inverse of the opposite of $(g,L)$.
\end{lem}
\proof
Let $\cN_{i,geom}:=(N_i,g^{N_i},\orient^{N_i},\cW_i)$ and
$\cN_{j,geom}:=(N_j,g^{N_j},\orient^{N_j},\cW_j)$
be distinguished models of $\partial_i\cM_{geom}$ and
$\partial_j \cM_{geom}$. Let $T_i:\bar N_i\times (-\infty,0)\rightarrow
\bar M$ and $T_j:\bar N_j\times (-\infty,0)\rightarrow \bar M$ be the
corresponding canonical inclusions with images $U_i$ and $U_j$, and let 
$\Pi_i:T_i^*\bar\cV\rightarrow \bar \cW_i*(-\infty,0)$ and
$\Pi_j:T_j^*\bar\cV\rightarrow \bar \cW_j*(-\infty,0)$ be the
corresponding canonical product structures. 

Let $\cN_{ki,geom}:=(N_{ki},g^{N_{ki}},\orient^{N_{ki}},\cW_{ki})$ and
$\cN_{kj,geom}:=(N_{kj},g^{N_{kj}},\orient^{N_{kj}},\cW_{kj})$
be distinguished models of $\partial_k \cN_{i,geom}$ and
$\partial_k \cN_{j,geom}$. Then we have the
corresponding canonical inclusions $T_{ki}:\bar N_{ki}\times (-\infty,0)\rightarrow
\bar N_i$ and $T_{kj}:\bar N_{kj}\times (-\infty,0)\rightarrow
\bar N_j$ with images $U_{ki}$ and $U_{kj}$, and the canonical product structures 
$\Pi_{ki}:T_{ki}^*\bar\cW_i\rightarrow \bar \cW_{ki}*(-\infty,0)$ and
$\Pi_{kj}:T_{kj}^*\bar\cV\rightarrow \bar \cW_{kj}*(-\infty,0)$. 
In order to work with specific Dirac bundles we always choose the point
$-1\in (-\infty,0)$ in order to define the reductions.

Let $f_{ki}:N_{ki}\rightarrow M$ and $f_{kj}:N_{kj}\rightarrow M$ be
the inclusions such that the pairs $(N_{ki},f_{ki})$ and
$(N_{kj},f_{kj})$ are models of $\partial_k M$.
We therefore have a canonical diffeomorphism $g:N_{ki}\rightarrow
N_{kj}$ such that $f_{ki}=f_{kj}\circ g$. Note that
$g$ is an orientation reversing isometry.
It induces a canonical orientation reversing isometry
$G:T_{kj}\circ (\bar g\times \id_{(-\infty,0)}) \circ
T_{ki}^{-1}:U_{ki}\rightarrow U_{kj}$.
We first construct a canonical extension of $G$ to an isomorphism
$(G,\Gamma):\cU_{ki,geom}^{op}\rightarrow \cU_{kj,geom}$ of geometric manifolds
which are
defined by restriction of structures to $U_{ki}$ and $U_{kj}$,
respectively.
Thus we must construct an isomorphism of Dirac bundles
$\Gamma:G^*\cW_{j|U_{kj}}\rightarrow \cW^{op}_{i|U_{ki}}$.

We choose the distinguished embedding $S:\bar N_{ki}\times
(-\infty,0)^2\rightarrow \bar M$ with image $U_k$ such that
$S(a,r,s)=T_i(T_{ki}(a,r),s)$.
We write $\Cl^1:=\Cl(T_{-1}((-\infty,0)))$ and $S:=S((-\infty,0))_{-1}$.
Furthermore, we define homomorphisms $\rho_i,\rho_j:\Cl^1\rightarrow
\End(\cV_{|U_k})$ by $\rho_j(\xi):=c(dS(\partial_r))$ and
$\rho_i(\xi):=c(dS(\partial_s))$, where $\xi$ is the generator of $\Cl^1$.

Let $a\in \bar N_{ki}$ and $r\in (-\infty,0)$.
The fiber of $G^*\cW_{j}$ at $T_{ki}(a,r)$ is canonically isomorphic to
$\Hom_{\Cl^1,\rho_j}(S,V_{S(a,-1,r)})$ if $M$ is even-dimensional,
and to 
$\Hom_{\Cl^1,\rho_j}(S\oplus S^{op},V_{S(a,-1,r)})$ if $M$ is
odd-dimensional. The fiber of $\cW_i$ at $T_{ki}(a,r)$ is canonically
isomorphic to $\Hom_{\Cl^1,\rho_i}(S,V_{S(a,r,-1)})$ if $M$ is even-dimensional,
and to 
$\Hom_{\Cl^1,\rho_i}(S\oplus S^{op},V_{S(a,r,-1)})$ if $M$ is
odd-dimensional.
We must define a canonical isomorphism
$$\Gamma(a,r):\Hom_{\Cl^1,\rho_j}(S,V_{S(a,-1,r)})\rightarrow 
\Hom_{\Cl^1,\rho_i}(S,V_{S(a,r,-1)})$$ if $M$ is even-dimensional, and
$$\Gamma(a,r):\Hom_{\Cl^1,\rho_j}(S\oplus S^{op},V_{S(a,-1,r)})\rightarrow 
\Hom_{\Cl^1,\rho_i}(S\oplus S^{op},V_{S(a,r,-1)})\ ,$$
if $M$ is odd-dimensional.

We define the unitary isomorphism $U\in \End(\cV_{|U_k})$ by
$U:=\frac{1}{\sqrt{2}}(\rho_i(\xi)-\rho_j(\xi))$.
One can check that $U\rho_iU^{*}=-\rho_j$,  $U\rho_jU^*=-\rho_i$,
$UzU^*=-z$, and $U c(dS(X))U^*=-c(dS(X))$ for all $X\in T\bar N_{ki}$.
Furthermore, let $\gamma_{(a,r)}$ be the path $\gamma_{(a,r)}:[0,1]\rightarrow
U_k$ given by $\gamma_{(a,r)} (t)=S(a,(-1,r)(1-t)+(r,-1)t)$.
By $\|(\gamma)$ we denote the parallel transport in $\cV$ along the
path $\gamma$.
Then we define for $\phi\in G^*\cW_{j,T_{ki}(a,r)}$
$$\Gamma(a,r)(\phi):=U\circ \|(\gamma_{(a,r)})\circ \phi\ .$$
        
It is now easy to check, that $\Gamma$ is an isomorphism of
Dirac bundles as required. Note that once we have chosen the
models for $\cN_{i,geom}$ and $\cN_{j,geom}$ this construction of the
isomorphism is completely canonical. If we interchange the roles of
$i$ and $j$, then the same construction gives an isomorphism
$(G^\prime,\Gamma^\prime):\cU_{kj,geom}^{op}\rightarrow \cU_{ki,geom}$.
Its opposite $(G^\prime,\Gamma^\prime)^{op}:\cU_{kj,geom}\rightarrow
\cU^{op}_{ki,geom}$ turns out to be the inverse of $(G,\Gamma)$. 

In the next step we extend  $g:N_{ki}\rightarrow N_{kj}$ to an
isomorphism $(g,L):\cN_{ki,geom}^{op}\rightarrow \cN_{kj,geom}$
in a canonical manner.
We must define an isomorphism of Dirac bundles $L:g^*
\cW_{kj}\rightarrow \cW_{ki}^{op}$.

We consider the actions $\rho_{ki}:\Cl^1\rightarrow \End(\cW_{i|U_i})$,
respectively $\rho_{kj}:\Cl^1\rightarrow \End(\cW_{j|U_j})$ given by
$\rho_{ki}(\xi):=c_{\cW_i}(dT_i(\partial_r))$, respectively $\rho_{kj}(\xi):=c_{\cW_j}(dT_j(\partial_r))$.
Let $a\in N_{ki}$. 
 Then the fiber of
$g^*\cW_{kj}$ at $a$ is given by
$\Hom_{\Cl^1,\rho_{kj}}(S\oplus S^{op},\cW_{j,T_j(g(a),-1)})$ if $M$ is
even-dimensional, and by 
$\Hom_{\Cl^1,\rho_{kj}}(S,\cW_{j,T_j(g(a),-1)})$ if $M$ is odd-dimensional.
The fiber of $\cW_{ki}$ at $a$ is given by 
$\Hom_{\Cl^1,\rho_{ki}}(S\oplus S^{op},\cW_{i,T_i(a,-1)})$ if $M$ is
even-dimensional, and by 
$\Hom_{\Cl^1,\rho_{ki}}(S,\cW_{j,T_i(a,-1)})$ if $M$ is odd-dimensional.
Note that $\cW_{j,T_j(g(a),-1)})=(G^*\cW_j)_{T_i(a,-1)}$.
Therefore we can define $L:g^*\cW_{kj}\rightarrow \cW_{ki}^{op}$ by
$L(\phi):=\Gamma(a,-1)\circ \phi$, $\phi\in (g^*\cW_{kj})_a$.

One can check, that $L$ is in fact an isomorphism of Dirac bundles.
Its construction is completely canonical once we have chosen the
distinguished 
models. If we interchange the roles of $i$ and $j$, then the same
construction provides an isomorphism
$(g^\prime,L^\prime):\cN_{kj,geom}^{op}\rightarrow \cN_{ki,geom}$.
Its opposite is the inverse of $(g,L)$.
\hB

\subsubsection{}

Let $(M,g^M,\orient)$ be a manifold with corners, admissible face
decomposition, admissible Riemannian metric and orientation.
Assume in addition that $M$ has a spin structure. The spinor bundle
$\cS(M)$ is then an admissible Dirac bundle.
Therefore the data $(M,g^M,\orient)$ together with the spin structure
determine a geometric manifold $\cM_{geom}:= (M,g^M,\orient,\cS(M))$.

Each atom of faces of codimension one of $M$
aquires an induced spin structure (see \ref{spininduced}). 
Let $i\in I_1(M)$ and $(N,f)$ be a model of $\partial_iM$
with induced metric $g^N$, orientation $\orient_N$, and admissible
face decomposition. It  also has an
induced spin structure. We thus obtain a geometric manifold
$(N,g^N,\orient_N,\cS(N))$.
Our construction of the boundary faces of a geometric manifold  is set up
such that
this geometric manifold is a model of 
$\partial_i \cM_{geom}$.  

\subsubsection{}\label{spininduced}
In order to fix conventions we recall the definition of the induced
spin structure in the set-up started in  \ref{setupor}.
Let $V\rightarrow X$ be a real oriented $n$-dimensional vector bundle.
A metric $g^V$ on $V$ induces a reduction of the structure group
of $\Fr(V)^+$ to $SO(n)$ which is represented by the
bundle of oriented orthonormal frames $\Fr(V,g^V)^+\subset \Fr(V)^+$.
A spin structure is then a reduction of the structure group
of $\Fr(V,g^V)^+$ to $Spin(n)$, i.e. an isomorphism class of pairs
$(Q,\Psi)$, where $Q\rightarrow X$ is a $Spin(n)$-principal bundle and
$\Psi:Q\times_{Spin(n)}SO(n)\stackrel{\sim}{\rightarrow}
\Fr(V,g^V)^+$. We can consider $\Psi$ as a two-fold covering
$\Psi:Q\rightarrow \Fr(V,g^V)^+$.

Assume now that $V$ fits into an exact sequence
$$0\rightarrow W\rightarrow V\rightarrow X\times \R\rightarrow 0\ .$$
Let $s:X\times \R\rightarrow V$ be a split which is normalized such
that it is an isometry with the orthogonal complement of $W$. Let $g^W$ be the restriction of
$g^V$ to $W$. Then we have the reduction of structure groups of
$\Fr(V,g^V)^+$ to $SO(n-1)$ represented by $\Fr(W,g^W)^+\subset
\Fr(V,g^V)$,
where the inclusion maps the frame $(w_1,\dots,w_{n-1})$ of $W$ to
the frame $(w_1,\dots,w_{n-1},s(1))$ of $V$.
We define the reduction of structure groups of $Q$ to $Spin(n-1)$ as
the pre-image  $R:=\Psi^{-1}(\Fr(W,g^W)^+)\subset Q$. The covering
$\Psi_{|R}:R\rightarrow \Fr(W,g^W)^+$ represents the induced spin
structure of $W$.

In the case of a Riemannian spin manifold with corners
we consider the canonical map $i:\partial M\rightarrow M$.
Then we have the sequence
$$0\rightarrow T\partial M\rightarrow i^*TM\rightarrow N\rightarrow 0\
.$$
We choose the section of $i^*TM$ given by the inward pointing normal
vector field. It provides the trivialization $N\cong \partial M\times
\R$
and the split $s:N\rightarrow i^*TM$.
The construction above provides the induced spin structure on
$T\partial M$.

\subsubsection{}

Let $\cN_{geom}=(N,g^N,\orient,\cW)$ be a geometric manifold and $H$
be an oriented Riemannian spin manifold with corners. 
\begin{ddd}\label{pps}
\index{$\cN_{geom}*H$}We define the geometric manifold
 $\cN_{geom}*H$ 
such that the underlying manifold with corners is
$M\times H$ with induced admissible face decomposition, product metric
and orientation, and the Dirac bundle is given by $\cW*H$.
\end{ddd}

\subsection{Pre-taming}\label{sectan}

\subsubsection{}

We consider a geometric manifold $\cM_{geom}=(M,g^M,\orient,\cV)$
such that the underlying manifold $M$ is compact. Let $\bar M$ be its
extension. With the canonical extension of the Riemannian metric $\bar
M$ is a complete Riemannian manifold. The Dirac bundle
$\cV$ has a canonical  extension $\bar \cV$ to $\bar M$.
Let $D(\cM_{geom}):=D(\bar \cV)$ be the associated Dirac operator.
It  is  essentially selfadjoint as an unbounded operator on the Hilbert space $L^2(\bar M,\bar V)$
with the domain $\dom (D(\cM_{geom}))= C_c^\infty(\bar M, \bar V)$ of smooth compactly
supported sections. The first order Sobolev space $H^1(\bar M,\bar
V)$ is the closure of $\dom(D(\cM_{geom}))$ with respect to the norm
$$\|\psi\|^2_{H^1}:=\|\psi\|^2_{L^2}+\|D(\cM_{geom})\psi\|_{L^2}^2\ .$$
In general $D(\cM_{geom})$ is neither invertible
nor Fredholm as an operator from $H^1(\bar M,\bar V)$ to $L^2(\bar
M,\bar V)$.

\subsubsection{}

A taming $\cM_t$ of the underlying geometric manifold $\cM_{geom}$ 
is given by the choice of smoothing operators on all faces of $M$, which when lifted to $\bar M$ provide a certain perturbation
$D(\cM_t)$ of $D(\cM_{geom})$ such that $D(\cM_t)$ is invertible.
A similar construction in the context of boundary value problems has been introduced
by Melrose and Piazza \cite{melrosepiazza97}.
We now describe the notion of taming in detail.

\subsubsection{}\label{tamam34}
Let $k\ge 0$ and $j\in I_k(M)$.
We consider a model $(N,f)$ of the face $\partial_j M$ and a
distinguished embedding $T:\bar N\times (-\infty,1)^k\rightarrow \bar
M$. Let $U_j$ denote the cylinder over the face $j$ in $\bar M$.
 Note that $N$ has an induced metric and orientation.
We choose a Dirac bundle $\cW$ over $N$ such that there is a product structure
$\Pi:T^* \cV\rightarrow \cW*(-\infty,1)^k$.
Assume that $R$ is an operator on $C^\infty(\bar N,W)$ with smooth
compactly supported integral kernel. We further assume that $R$ is
selfadjoint, and that it is odd if $N$ is even-dimensional. 
Then we can form the operator 
$\Pi^{-1}\circ \bL^{U_j}_{\bar N}(R)\circ \Pi$ on $C^\infty(U_j,\bar V_{|U_j})$.

\begin{ddd}\label{tamtam1234}
\index{pre-taming}\index{$\cM_t$}Any operator on $C^\infty(U_j,\bar V_{|U_j})$ which arrises by this
construction
is called a pre-taming of the face $j$.
\end{ddd}

\begin{ddd}
A pre-taming of $M$ is a collection of operators $\{Q_j|j\in I_*(M)\}$
such that $Q_j$ is a pre-taming of $j$.
\end{ddd}

Our notation for a geometric manifold with distinguished pre-taming is
$\cM_t$.

\subsubsection{}

Let $\cM_{geom}$ and $\cM^\prime_{geom}$ be geometric manifolds
and $(f,F):\cM_{geom}\rightarrow\cM^\prime_{geom}$ be an isomorphism.
Let $\cM_t$ and $\cM_t^\prime$ be pre-tamings of the underlying
geometric manifolds given by $\{Q_j|j\in I_*(M)\}$ and
$\{Q_j^\prime|j\in I_*(M^\prime)\}$.
Then we have a bijection $j\mapsto j^\prime$ from
$I_*(M)$ to $I_*(M^\prime)$.
\begin{ddd}
\index{isomorphism!of pre-tamed manifolds}We say that $(f,F)$ is an isomorphism of pre-tamed manifolds, if
$F^*Q_{j^\prime}=Q_j$ for all $j\in I_*(M)$.
\end{ddd}



\subsubsection{}
If $k\le 1$ and $i\in I_k(M)$, then we will define $\partial_i \cM_t$.
First of all, this is a well-defined isomorphism class. Again we have
a more rigid structure. Namely, if we choose a distinguished model for
$\partial_i \cM_{geom}$, then it has an induced pre-taming.

In the case $k=0$ the class  $\partial_i \cM_t$ has a canonical representative.
Let us therefore discuss the more complicated  case $k=1$ in detail.
We take a distinguished model $\cN_{geom}=(N,g^N,\orient,\cW)$ of
$\partial_i\cM_{geom}$. It comes with a canonical embedding $T:\bar
N\times (-\infty,0)\rightarrow \bar M$ and a canonical product
structure $\Pi: T^*\bar\cV\rightarrow\cW*(-\infty,0)$.
We construct the induced  pre-taming of $\cN_{geom}$ as follows.

Let $l\ge 0$ and $j\in I_l(N)$. Let $U_{ij}\subset \bar N$ be the
cylinder over the face $j$ in $\bar N$. We further define the subset
$U_i(j):=T(U_{ij}\times (-\infty,0))\subset U_i$ of the cylinder over
$i$ in $\bar M$. It coincides with the cylinder over the face $j$ of
$M$, i.e. the set where $Q_j$ is defined.
By  Lemma \ref{ssignn} there is now a unique pre-taming
$R_j$ of the face $j\in I_{l}(N)$ (an operator on $C^\infty(U_{ij},\bar \cW)$)  such that
$$\Pi^{-1}\circ \bL^{U_i(j)}_{U_{ij}}(R_j)\circ \Pi=Q_j\ .$$
\begin{ddd}\label{bounftt}
\index{$\partial_i\cM_t$}The isomorphism class $\partial_i \cM_t$ is represented by the
pre-taming of $\cN_{geom}$ given by $\{R_j|j\in I_*(N)\}$ as
constructed above.
\end{ddd}

\subsubsection{}

Note that the Dirac bundles of $\cM_{geom}$ and $\cM_{geom}^{op}$
coincide as vector bundles so that we can compare operators on
sections of these bundles.
Let $\{Q_j|j\in I_*(M)\}$ be a taming $\cM_t$.
\begin{ddd}
\index{$\cM_t^{op}$} \index{opposite!pre-tamed manifold}The opposite $\cM_t^{op}$ of $\cM_t$ is 
given by $\{-Q_j|j\in I_*(M)\}$.
 \end{ddd}

\subsubsection{}

We have the following generalization of Lemma \ref{hiphop}
to the tamed case.  
\begin{lem}\label{hiphop1}
We have
$\partial_k\partial_i\cM_t\cong (\partial_k\partial_j\cM_t)^{op}$.
In fact, if we choose distinguished models $\cN_{ki,geom}$ and
$\cN_{kj,geom}$ as in the proof of Lemma \ref{hiphop}
with their induced tamings,
then this isomorphism
implemented by the canonical
isomorphism $(g,L):\cN_{ki,geom}^{op}\rightarrow \cN_{kj,geom}$
constructed there.
\end{lem}
\proof
It suffices to show that
the isomorphism  $(g,L):\cN_{ki,geom}^{op}\rightarrow \cN_{kj,geom}$
induces an isomorphism of pre-tamed manifolds.
This is a straight-forward but tedious calculation. \hB

\subsection{Taming}

\subsubsection{}

A pre-taming gives rise to a perturbation $D(\cM_t)$ of the Dirac operator $D(\cM_{geom})$. The condition which characterizes tamings among pre-tamings 
will be formulated in terms of the spectral theory of this operator.
In this subsection all geometric manifolds are compact.

\subsubsection{}
Let $\cM_t$ be a pre-taming of $\cM_{geom}$
given by $\{Q_j|j\in
I_*(M)\}$. 
In order to write out
analytical details in a sufficiently readable way we use the following
simplification of the notation. Namely, we will often write
$Q_j=\bL_{\overline{\partial_j M}}^{U_j}(R_j)$.
It is then implicitly understood that the symbol $R_j$ has a
meaning as an operator only after choosing a model
$(N_j,f_j)$ of $\partial_jM$, the distinguished embedding $T_j$, a
Dirac bundle $\cW_j$, and a
product structure
$\Pi_j$.

Another simplification of notation is the following. If $j\in I_k(M)$,
then we often write $\partial_j \bar M\times
(-\infty,0)^k\subset \bar M$ for the cylinder over $j$ in $\bar M$.
 Explicitly this means that we choose
(without mentioning)
a model $(N_j,f_j)$ of $\partial_j M$ and a distinguished embedding
$T_j:N\times (-\infty,0)^k\rightarrow \bar M$ in order to parameterize
this cylinder. 

In order to write out certain constructions  
we often we identify (again without mentioning) the bundle
$\bar \cV_{|U_j}$ with a bundle of the form $\cW_j*(-\infty,0)^k$ using a
product structure.

\subsubsection{}
We choose compatible cut-off functions $\rho_j$ for all faces
$j$ of $M$.
Fix once and for all a function $\rho\in C^\infty(\R)$ be such that $\rho(r)=0$ for $r\ge -1/4$ and $\rho(r)=1$ for
$r\le -1$. It gives rise to cut-off functions in the following way.

If $M$ is a manifold with corners and $j\in I_1(M)$, then we choose a
model
$(N,f)$ of $\partial_jM$. We then  have a unique distinguished
embedding 
$T:\bar N\times (-\infty,0)\rightarrow \bar M$, and we denote its
image, the cylinder over the face $j$,  by $U_j$.  
We define the function
$\rho_j\in C^\infty(\bar M)$ such that it is supported on the 
cylinder $U_j$ and satisfies $T^*\rho_j(x,r)=\rho(r)$.

This construction is compatible with restriction in the following
sense. Let $i\in I_1(\partial_jM)$
and $j^\prime\in I_1(M)$ be the adjacent face (Def. \ref{conjface})
to $j$ with respect to $i$.
Let $\bar f:\overline{\partial_j M}\rightarrow \bar M$ be the natural embedding.
Then we have
$\rho_i=\bar f^*\rho_{j^\prime}$. 

For $j\in I_k(M)$, $k\ge
2$, we set $\rho_j:=\prod_{i\in I_1(M), j<i} \rho_i$. 
For $i\in I_0(M)$ we let
$\rho_i$ be the characteristic function of the corresponding face $\partial_iM$.

\subsubsection{}

Since   $M$ is compact, the face $\partial_jM$ is compact for all $j\in
I_*(M)$. Since the operator $R_j$ has a compactly supported smooth
integral kernel and $Q_j=\bL_{\overline{\partial_j M}}^{U_j}(R_j)$ is
local in the normal directions the operator $\rho_jQ_j$
acts on $C^\infty(\bar M,\bar \cV)$ and  preserves smooth sections of compact
support. It extends to a selfadjoint bounded operator on the Hilbert
space $L^2(\bar M,\bar V)$. It also acts as bounded operator on the
Sobolev space $H^1(\bar M,\bar V)$.

\subsubsection{}
We define  the operator
$$D(\cM_t):=D(\cM_{geom}) + \sum_{k\ge 0}\sum_{j\in I_k(M)} \rho_j Q_j$$ acting on the domain
$\dom (D(\cM_t)):= C_c^\infty(\bar M, \bar{V})$. It is a bounded
perturbation of $D(\cM_{geom})$ which is also
essentially selfadjoint. 
By the spectrum $\spec(D(\cM_t))\subset \C$
we mean the spectrum of its unique selfadjoint extension
$\overline{D(\cM_t)}$. We say that $\cD(\cM_t)$ is invertible, if
$0\not\in \spec(\cD(\cM_t))$.

\index{essential spectrum}\index{$\spec_{ess}$}Let $\spec_{ess}(D(\cM_t))\subseteq
\spec(D(\cM_t))$ denote the essential spectrum. A point
\index{Weyl sequence}$\lambda\in\C$
belongs to the essential spectrum iff there exists a Weyl sequence,
i.e. a sequence $(\phi_i)$, $\phi_i\in C^\infty_c(\bar M,\bar V)$,
with $\|\phi_i\|_{L^2}=1$, $\phi_i\to 0$ weakly, and
$\|D(\cM_t)\phi_i-\lambda\phi_i\|_{L^2}\to 0$.

\subsubsection{}
The closure of $D(\cM_t)$
acts as a bounded operator
$\overline{\cD(\cM_t)}:H^1(\bar M,\bar V)\rightarrow L^2(\bar M,\bar V)$.
We say that $\cD(\cM_t)$ is a Fredholm operator if this bounded
operator is a Fredholm operator.
By spectral theory, $D(\cM_t)$ is Fredholm iff $0\not\in \spec_{ess}(\cD(\cM_t))$.
\begin{prop}\label{tamm}
\index{Fredholm operator}The operator $D(\cM_t)$ is Fredholm iff
$D(\partial_i\cM_t)$ is invertible 
for all $i\in I_1(M)$.
\end{prop}
\proof
Assume that $D(\partial_i\cM_t)$ is invertible 
for all $i\in I_1(M)$.
We will then construct left and right  parameterizes $R_l$ and $R_r$ for $\overline{\cD(\cM_t)}$.
In order to check compactness of remainder terms we employ the
following variant of Rellich's Lemma valid on complete Riemannian
manifolds: If $f\in L^\infty(\bar M)$ satisfies
$\lim_{r\to \infty} \sup_{x\in \bar M\setminus {}_r\bar M}|f(x)|=0$, 
 then multiplication by $f$ is a compact operator
$f:H^1(\bar M,\bar V)\rightarrow L^2(\bar M,\bar V)$.

We choose $T_0>0$ sufficiently large such that the integral kernel of
$R_j$ is supported in ${}_{T_0}\overline{\partial_jM}\times
{}_{T_0}\overline{\partial_jM}$ for all $j\in I_*(M)$.
If $\phi\in C^\infty_c(\bar M,\bar V)$ is supported on
$\overline{\partial_iM}\times (-\infty,-T_0)$, then so is $D(\cM_t)\phi$.
Therefore we can restrict $D(\cM_t)$ to $\overline{\partial_iM}\times (-\infty,-T_0)$.
This restriction extends to a $\R$-invariant operator 
on $\overline{\partial_i M}\times \R$ which is isomorphic to 
$$D_i=c(\partial_s)\partial_s+\bL_{\overline{\partial_i
    M}}^{{\overline{\partial_i M}}\times \R}(D(\partial_i\cM_{t}))\ ,$$
where $s$ is the normal variable.
We have
$$D_i^2=-\partial^2_s + \bL_{\overline{\partial_i
    M}}^{{\overline{\partial_i M}}\times \R}(D(\partial_i\cM_{t}))^2=-\partial^2_s + L_{\overline{\partial_i
    M}}^{{\overline{\partial_i M}}\times
  \R}(D(\partial_i\cM_{t})^2\otimes 1)$$
(omit  "$\otimes 1$" in the case where $\partial_i M$ is even-dimensional).
By Fourier transformation in the $s$-variable this operator is
unitary equivalent to the sum of  commuting  non-negative operators
$\sigma^2$ and $L_{\overline{\partial_i
    M}}^{{\overline{\partial_i M}}\times
  \R}(D(\partial_i\cM_{t})^2\otimes 1)$, where $\sigma\in\R$ is the
conjugated variable.
We conclude that $\inf
\spec(D_i^2)=\inf\spec_{ess}(D_i^2)=\inf\spec(D(\partial_i\cM_{t})^2)>0$.
Thus $D_i$ is invertible. Let $R_i$ denote its inverse which is
considered as a bounded operator $R_i:L^2\rightarrow H^1$.

\index{parametrix}The left and right parameterizes $R_l$ and $R_r$ of $D(\cM_{t})$ are built from the operators $R_i$
and an interior parametrix by the usual glueing construction. The main
point is to choose the cut-off functions in a careful manner.
\newcommand{\grad}{{\tt grad}}

Let $q\in\nat$ be larger than the highest codimension of a corner of
$M$. Then we consider the intersection $Q_q:=S^{q-1}\cap (-\infty,0]^q$.
The group $\Sigma^{q}$ acts on $Q_q$ by permutation of coordinates.
We choose a smooth partition of unity 
$(\chi_1,\dots,\chi_q)$ on $Q_q$ with the following properties.
\begin{enumerate}
\item There exists $\epsilon>0$ such that $\chi_i(r_1,\dots,r_q)=0$ if $|r_i|\le \epsilon$.
\item $\sigma^*\chi_i=\chi_{\sigma^{-1}(i)}$ for $\sigma\in \Sigma^q$.
\end{enumerate}
Note that $\chi_i(r_1,\dots,r_q)=1$ if $|r_i|>1-\frac12 \epsilon$.
Note that for $h\le q$ we have $Q_{h}=Q_q\cap
\{r_{h+1}=0,\dots,r_{q}=0\}$.
By restriction we obtain a similar partition of unity
$(\chi_1,\dots,\chi_h)$ on $Q_h$.

We now define functions $\phi_i$ on $\bar M$ for all $i\in
I_1(M)$ as follows. Let $T_1:=T_0+1$.  We first fix these functions on
$\bar M\setminus {}_{T_1}\bar M$. 
We can decompose $\bar M\setminus \inter({}_{T_1}\bar M )$ as a union of closed
subsets
$V_j:={}_{T_1}\overline{\partial_j M}\times (-\infty,-T_1]^k$, $j\in I_k(M)$,
which meet along common boundaries.
Consider now $i\in I_1(M)$. Let $k>0$ and $j\in I_k(M)$. If 
$V_j\cap \overline{\partial_iM}\times (-\infty,-T_1)=\emptyset$, then
we set $\phi_{i|V_j}=0$.
Otherwise there exists a unique $h\in \{1,\dots,k\}$ such that the
normal variable $r_h$ on $V_j={}_{T_1}\overline{\partial_j M}\times (-\infty,-T_1]^k$
coincides with the normal variable to $\partial_iM$. In this case we
define the function $\phi_i$ on $V_j$ by 
$\phi_i(n,r_1,\dots,r_k)=\chi_h(\frac{r_1}{r},\dots,\frac{r_k}{r})$,
where
$(n,r_1,\dots,r_k)\in{}_{T_1}\overline{\partial_j M} \times (-\infty,-R_1]^k$ and
$r=\sqrt{r_1^2+\dots+r_k^2}$.
The functions $\phi_i$ thus defined are continuous and almost
everywhere differentiable with 
piecewise continuous derivatives. 
We extend $\phi_i$ as a continuous functions with piecewise
continuous derivatives supported on $U_i\setminus {}_{T_0}\bar M$.
Note that $|d\phi_i|\in L^\infty(\bar M)$ satisfies
$\lim_{r\to \infty} \sup_{x\in \bar M\setminus {}_r\bar M}|d\phi_i(x)|=0$.
We set $\phi_0:=1-\sum_{i\in I_1(M)}\phi_i$. Then
 $(\phi_0,\{\phi_i\}_{i\in I_1(M)})$ is a partition of unity subordinate to the covering
$(\inter({}_{T_1}\bar M),\{U_i\setminus {}_{T_0}\bar M\}_{i\in
  I_1(M)})$ of $\bar M$.



In a similar manner (by choosing appropriate functions on $Q_q$) we construct almost everywhere differentiable cut-off functions $\psi_l$,
 $l\in
\{0\}\cup I_1(M)$, such that $\psi_l\phi_l=\phi_l$ and 
$\supp(\psi_0)\subset {}_{T_1}\bar M$, respectively $\supp(\psi_i)\subset
U_i\setminus {}_{T_0}\bar M$ for all $i\in I_{1}(M)$, and such that
$|d\psi_i|\in L^\infty(\bar M)$,  
$\lim_{r\to \infty} \sup_{x\in \bar M\setminus {}_r\bar
  M}|d\psi_i(x)|=0$.
Note that multiplication by $\psi_i$ or $\phi_i$ acts as a bounded
operator on $H^1(\bar M,\bar V)$. 

Let $R_0:L^2_c(\bar M,\bar V)\rightarrow H^1_{loc}(\bar M,\bar V)$ be a local parametrix of $D(\cM_{geom})$ constructed
e.g. using pseudo-differential calculus. Then we have  proper
smoothing operators $S_r:=D(\cM_{geom})R_0-1$ and $S_l:=R_0D(\cM_{geom})-1$.
We define $R_r:L^2(\bar
M,\bar V)\rightarrow H^1(\bar M,\bar V)$ by 
$$R_r:=\psi_0R_0\phi_0+\sum_{i\in I_1(M)}\psi_i R_i \phi_i$$
(making the obvious identifications so that we can consider the
terms  $\psi_i R_i \phi_i$ as continuous linear maps
from $L^2(\bar M,\bar V)$ to $H^1(\bar M,\bar V)$).
Then we get (writing $D:=D(\cM_t)$ for the moment)
\begin{eqnarray*}
DR_r&=&D\psi_0R_0\phi_0+\sum_{i\in I_1(M)}D\psi_i R_i \phi_i\\
&=&\psi_0 DR_0 \phi_0 + \sum_{i\in I_1(M)}\psi_i D_i R_i \phi_i+
[D,\psi_0]R_0\phi_0 + \sum_{i\in I_1(M)} [D,\psi_i] R_i \phi_i\\
&=&  \psi_0  \phi_0 + \sum_{i\in I_1(M)}\psi_i
\phi_i+\psi_0S_l\phi_0+\psi_0(D-D(\cM_{geom}))R_0\phi_0\\&&+c(\grad(\psi_0))R_0\phi_0
+ [(D-D(\cM_{geom})),\psi_0]R_0\phi_0\\&&+ \sum_{i\in I_1(M)} c(\grad(\psi_i))
R_i \phi_i+ \sum_{i\in I_1(M)} [(D-D(\cM_{geom})),\psi_i]
R_i \phi_i\\
&=&1+\psi_0S_l\phi_0+\psi_0(D-D(\cM_{geom}))R_0\phi_0\\&&+c(\grad(\psi_0))R_0\phi_0
+ [(D-D(\cM_{geom})),\psi_0]R_0\phi_0\\&&+ \sum_{i\in I_1(M)} c(\grad(\psi_i))
R_i \phi_i+ \sum_{i\in I_1(M)} [(D-D(\cM_{geom})),\psi_i]
R_i \phi_i\ .\\
\end{eqnarray*}
We now argue that all remainder terms are compact.
First of all the term $\psi_0S_l\phi_0$ has a compactly supported
integral kernel and is therefore compact.
The terms involving gradients of the cut-off functions are compact
by the variant of Rellich's Lemma mentioned above.
The term $\psi_0(D-D(\cM_{geom}))R_0\phi_0$ factors through a continuous
map from $L^2( \bar M,\bar \cV)$
to $H^1_c(\bar M,\bar \cV)$ and induces therefore a compact operator.
The term $[(D-D(\cM_{geom})),\psi_0]R_0\phi_0$ is compact, since
$(D-D(\cM_{geom}))$ is bounded on $L^2$ and $H^1$, and
$\psi_0$ has compact support.
An inspection of the definitions shows that
$[(D-D(\cM_{geom})),\psi_i]$ is compactly supported for all $i\in I_1(M)$.
Fix $i\in I_1(M)$ and $j\in I_k(M)$.
Then distribution kernel of the  term $ \rho_j Q_j$ is supported on $({}_{T_0}\overline{\partial_j M}\times
(-\infty,0)^k)\times( {}_{T_0}\overline{\partial_j M}\times
(-\infty,0)^k)$. Furthermore, it is local with respect to the second variable.
The restriction $\psi_{i|{}_{T_0}\overline{\partial_j M}\times
  \{r_1,\dots,r_k\}}$ is independent of the coordinate in 
${}_{T_0}\overline{\partial_j M}$ as long as $r_i\ge T_1$ for at least
one $i\in \{1,\dots,k\}$. Therefore, the distribution kernel of 
$[(D-D(\cM_{geom})),\psi_i]$ is supported in
${}_{T_1}\bar M\times {}_{T_1}\bar M$.
It now follows that also the remaining terms involving
$[(D-D(\cM_{geom})),\psi_i]$ are compact operators.

As a right parametrix we can take the adjoint of $R_r$.

Assume now that $D(\partial_i\cM_t)$ is not invertible for some $i\in I_1(M)$.
Then $0\in\spec_{ess}(D_i)$ and we can construct a Weyl sequence
$(\Psi_k)$ for $0$ such that $\supp(\Psi_k)\subset
\overline{\partial_i M}\times (-\infty,-T_0]\subset U_i$ for all
$k\in\nat$. We can consider $(\Psi_k)$ as a Weyl sequence for $0$ of
the operator
$D(\cM_t)$. Thus $0\in\spec_{ess}(D(\cM_t))$, and hence
$D(\cM_t)$ is not Fredholm.
\hB
\subsubsection{}
\begin{ddd}
\index{taming}A pre-taming $\cM_t$ of the underlying geometric manifold $\cM_{geom}$ is called a taming
if the operator $D(\cM_t)$ is invertible.
\end{ddd}

It follows from Proposition  \ref{tamm} that if $\cM_t$ is a taming of $\cM_{geom}$, then the induced pre-taming $\partial_i \cM_t$
of $\partial_i \cM_{geom}$ is a taming for all $i\in I_k(M)$, $k\le 1$.

\subsection{Obstructions against taming}\label{obst}

\subsubsection{}
Consider a compact geometric manifold $\cM_{geom}$. One can now naturally ask
if there exists a taming of $\cM_{geom}$. It turns out that there are
obstructions against the existence of tamings. Before we turn to the
general case let us discuss some instructive examples.
In the present section all geometric manifolds are compact.

\subsubsection{}
Assume first that $M$ is closed. Then $D(\cM_{geom})$ is a Fredholm
operator. For a geometric manifold $\cN$ let $\cN^{ev}$ denote the even-dimensional 
part.
\begin{ddd}\label{infg}
\index{$\ind(\cM_{geom})$}\index{index!of a closed geometric manifold}We define the function
$$\ind(\cM_{geom}):I_0(M)\rightarrow \Z$$ such that
$\ind(\cM_{geom})(i)\in\Z$ is the index of the Fredholm operator
$$\overline{D(\cM_{geom})}^+:H^1(\partial_iM^{ev},V^+)\rightarrow
L^2(\partial_iM^{ev},V^-)\ .$$
\end{ddd}
\begin{lem}
If $M$ is closed, then $\cM_{geom}$ admits a taming iff
$\ind(\cM_{geom})=0$.
\end{lem}
\proof
This is a special case of Lemma \ref{btt} proved below. \hB

\subsubsection{}
Note that we must require the vanishing of the index on every face
$i\in I_0(M)$ separately, since the smoothing operators $W_i$ 
are not allowed to mix faces. Let e.g. $\cN_{geom}$ be an
even-dimensional closed connected geometric manifold. Then $I_0(N)$
consists of one element $c$. Assume that $\ind(\cN_{geom})(c)\not=0$.
Let $\cM_{geom}:=\cN_{geom}\sqcup \cN_{geom}^{op}$. Then
$I_0(M)=\{c,c^{op}\}$. We have
$\ind(\cM^{op}_{geom})(c^{op})=-\ind(\cM_{geom})(c)$.
The geometric manifold $\cM_{geom}$ does not admit a taming.
But we can consider its reduction $\cM_{geom}^{red}$ which has one
face of codimension zero $(c,c^{op})$. In this case
$\ind(\cM_{geom}^{red})((c,c^{op}))=\ind(\cN_{geom})(c)+\ind(\cN^{op}_{geom})(c^{op})=0$,
so that $\cM_{geom}^{red}$ admits a taming.

\subsubsection{}
Next we consider the case of a manifold with boundary.
Assume that the underlying manifolds with corners of $\cM_{geom}$ 
is a connected odd-dimensional manifold with boundary.
If $\cM_{geom}$ admits a taming, then first of all
$\partial_i\cM_{geom}$ must admit a taming for all $i\in I_1(M)$.
Assume that $I_1(M)$ consists of one element $b$. Then
$\partial_b\cM_{geom}$ is zero-bordant. Therefore the index of the
corresponding Dirac operator vanishes,
i.e. $\ind(\partial_b\cM_{geom})=0$. It follows that
$\partial_b\cM_{geom}$
admits a taming $\partial_b\cM_t$.
Let $I_0(M):=\{c\}$. For the moment we consider the pre-taming
$\cM_t$ of $\cM_{geom}$ which is induced by the smoothing operators
$W_b$ defining the taming $\partial_b\cM_t$ and $W_c:=0$.
By Proposition \ref{tamm} the operator $\cD(\cM_{t})$ is Fredholm.
It is now easy to see (special case of Lemma \ref{btt}) that we can
choose a possibly different $W_c$ such that
$(W_b,W_c)$ induces a taming $\cM_{t}$.

\subsubsection{}
If we consider $\cM_{geom}$ with a different face decomposition
where $I_1(M)$ has more than one element, then it can happen that
$\cM_{geom}$ does not admit a taming since already $\partial_i
\cM_{geom}$ does not admit a taming for some $i\in I_1(M)$.

\subsubsection{}
 The simplest example for this effect is  the unit interval $I=[0,1]$
with the standard metric and orientation, and with the spinor bundle
$\cS(I)$. Let us first consider the atomic face decomposition and
\index{$\cI_{geom}$}denote the corresponding geometric manifold by $\cI_{geom}$.
Then we have $I_1(I)=\{0,1\}$. The Dirac bundle $\cW_i$ of $\partial_i\cI_{geom}$
is one-dimensional. Therefore it does not admit any non-trivial odd
operator. In particular, $D(\partial_i\cI_{t})=0$ for any pre-taming.
It follows that none of the  pre-tamings $\cI_t$ can  be a taming.
\newcommand{\cJ}{{\cal J}}

\subsubsection{}\index{$\cJ_{geom}$}
Let $\cJ_{geom}$ be the interval as above, but now with a different
face decomposition such that
$I_1(J)=\{(0,1)\}$.  
We choose identifications of the spaces of sections of $\cW_i$ with
$\C$ for $i=0,1$.
This induces an identification of the space of sections of the Dirac
bundle $\cW$ of $\partial_{(0,1)}\cJ_{geom}$ with $\C\oplus\C$
as a $\Z/2\Z$-graded vector space.
We set
$$W_{(0,1)}:=\left(\begin{array}{cc}0&1\\1&0\end{array}\right)\in
\End(C^\infty(\partial_{(0,1)}J,\cW))\cong \Mat(2,\C)\ .$$
For the moment we take $W_c:=0$ for $c\in I_0(J)$.
These two operators together induce a pre-taming of $\cJ_{t}$.
One can check, that
$\spec_{ess}(D(\cJ_t))=(-\infty,-1]\cup [1,\infty)$.
Thus after perturbing $W_c$, if necessary, the pair
$(W_c,W_{(0,1)})$ defines a taming $\cJ_t$.

\subsubsection{}
We now consider the question of the existence of a taming in general.
Let $\cM_{geom}$ be a geometric manifold. 
We try to construct a taming of $\cM_{geom}$, i.e. the operators
$W_j$,  inductively by decreasing codimension of $\partial_jM$. 
In
each step we encounter obstructions which we will investigate by
methods of index theory. At some places we will employ the special case
of the local index theorem for families with corners \ref{famind} for the family
over a point. Note that the proofs of the index theorem and related
results are completely independent of the present section.

\subsubsection{}
\begin{ddd} A pre-taming 
\index{boundary pre-taming}$\cM_t$ is called a boundary pre-taming if $W_i=0$ for all $i\in
I_0(M)$. The notation for a boundary pre-tamed manifold is\index{$\cM_{bt}$} $\cM_{bt}$.
\end{ddd}

\begin{ddd} A boundary pre-taming is called a boundary taming iff one of the following equivalent conditions is satisfied: 
\begin{enumerate}
\item
$D(\cM_{bt})$ is Fredholm. 
\item
$\partial_i\cM_t$ is a taming for all $i\in I_1(M)$.
\item
$D(\partial_i \cM_t)$ is invertible for all $i\in I_1(M)$.
\end{enumerate}
\end{ddd}
In fact, the equivalence of 1. and 3. is Proposition \ref{tamm}.
The equivalence of 2. and 3. is the definition of a taming.

\subsubsection{}
If the underlying manifold of $\cM_{geom}$ is closed, then
this manifold is boundary-tamed for trivial reasons. 
We can extend the definition of the index function \ref{infg}
to general boundary-tamed manifolds $\cM_{bt}$.
\begin{ddd}\label{infg1}
\index{index!of a boundary tamed manifold}We define the function
$$\ind(\cM_{bt}):I_0(M)\rightarrow \Z$$ such that
$\ind(\cM_{bt})(i)\in\Z$ is the index of the Fredholm operator
$$\overline{D(\cM_{bt})}^+:H^1(\overline{\partial_i M}^{ev},\bar V^+)\rightarrow
L^2(\overline{\partial_iM}^{ev},\bar V^-)\ .$$
\end{ddd}

\subsubsection{}
Let $\cM_{bt}$ be a boundary-tamed manifold. An extension of the
boundary-taming to a taming is a taming which is obtained from the
boundary taming by changing the operators $W_i$, $i\in I_0(M)$, while
keeping fixed the others. The main step in the inductive construction
of a taming is the extension of a boundary taming to a taming.

\begin{lem}\label{btt}
\index{extension!of a boundary taming to a taming}
The boundary taming $\cM_{bt}$ can be extended to a taming if and only
if $\ind(\cM_{bt})=0$.
\end{lem}
\proof 
We can construct the operators $W_i$ for each face $i\in I_0(M)$ separately.
Therefore we can assume that $\cM_{geom}$ is reduced. Let $I_0(M)=\{c\}$.
Furthermore, we can assume that $\partial_cM$ is either
even-dimensional or odd-dimensional, since again we can construct the
operator $W_c$ on the even- and odd-dimensional parts separately.

Assume that the boundary taming $\cM_{bt}$ can be extended to a taming
$\cM_t$.
If $M$ is even-dimensional, then
$D(\cM_t)^+$ is a compact perturbation of the Fredholm operator
$D(\cM_{bt})^+$ which is in addition invertible. Thus
$\ind(\cM_{bt})(c)=0$.
If $M$ is odd-dimensional, then we have $\ind(\cM_{bt})(c)=0$ by definition.

We now show the converse. Assume that $\ind(\cM_{bt})(c)=0$.
First we consider the case that $M$ is odd-dimensional.
Let $P$ be the orthogonal projection onto $\ker(\overline{D(\cM_{bt})})$.
Then $P$ is a finite-dimensional smoothing operator and $D(\cM_{bt})+P$ is invertible.
If $\chi\in C_c^\infty(\bar M)$ is a cut-off function, then we consider 
the smoothing operator $\tilde P:=\chi P\chi$ with compactly supported
integral kernel.  
Since $P$ is compact, we can make the operator norm $\|\tilde P-P\|$ as
small as we want if we choose $\chi$ such that $\chi_{|{}_r\bar M}=1$
for sufficiently large $r$. In fact, if we choose a sequence of such functions $\chi_i$ with $r_i\to \infty$, then the sequence  of multiplication operators $\chi_i$ converges strongly to the identity.
Since $P$ is compact, it follows that $\chi_i P$ converges to $P$ in the norm topology.
Therefore $\chi_i P \chi_i=(\chi_i P)(P\chi_i)$ also converges in norm to $P$.

If $\|\tilde P-P\|$ is sufficiently small, then $D(\cM_{bt})+\tilde P$ is invertible, and we can set $W_c:=\tilde P$.

Now we consider the case that $M$ is even-dimensional.
Since  $\ind(\overline{D(\cM_{bt})^+})=0$ we can find a unitary operator $U:\ker(\overline{D(\cM_{bt})}^+)\stackrel{\sim}{\rightarrow}
\ker(\overline{D(\cM_{bt})}^-)$. We then define  $P:=U+U^*$ on
$\ker(\overline{D(\cM_{bt})})$ and extend it by zero to the orthogonal complement
of the kernel. Then $P$ is odd with respect to the
$\Z/2\Z$-grading, and $D(\cM_{bt})+P$ is invertible. As in the
odd-dimensional case we construct
a compactly supported perturbation $\tilde P$ such that
$D(\cM_{bt})+\tilde P$ is invertible. We then set $W_c:=\tilde P$.
\hB

In some situations the Dirac bundle might have additional symmetries like real or quarternionic structures.
For example, in the case of the spin Dirac operator in dimensions $1,2 (8)$ this leads to index invariants in $\Z/2\Z$. Our tamings need not be compatible with these additional symmetries. In particular, the $\Z/2\Z$-valued index invariants do not obstruct the extension of a boundary taming to a taming. Consider e.g. the closed manifold $S^1$ with the non-bounding spin structure. The associated Dirac operator has a one-dimensional kernel and therefore non-trivial $\Z/2\Z$-valued index.  Nevertheless it admits a taming by the construction above.

\subsubsection{}
Let $\cM_{bt}$ be a boundary-tamed manifold. Then we may have
different extensions of the boundary taming to a taming. Homotopy
classes of extensions are distinguished by a spectral flow invariant.
 Recall that the spectral flow $\Sf((F_u)_u)$ of a continuous family of selfadjoint Fredholm
operators $(F_u)_{u\in [0,1]}$ such that $F_0$ and $F_1$ are invertible
is the net number
of eigenvalues of $F_u$ which cross zero from the positive to the
negative side as $u$ tends from $0$ to $1$.
Let $\cM_t$ and $\cM_t^\prime$ be two extensions of $\cM_{bt}$ to tamings.
Then we  consider the family of Fredholm operators
$D_u:=(1-u)D(\cM_{t})+uD(\cM_t^\prime)$ parameterized by the interval
$[0,1]$.
\begin{ddd}\label{sflowdef}
\index{spectral flow}\index{$\Sf$}We define the function
$$\Sf(\cM_t,\cM_t^\prime):I_0(M)\rightarrow \Z$$
such that
$\Sf(\cM_t^\prime,\cM_t)(i)$ is the spectral flow of the
family of Fredholm operators
$(D_{i,u})_{u\in [0,1]}$ on $L^2(\overline{\partial_i M}^{odd},\bar V)$ induced by
$D_u$. If $\partial_iM^{odd}=\emptyset$, then we set $\Sf(\cM_t,\cM_t^\prime):=0$.
\end{ddd}
Note that $D(\cM_t^{op})\cong -D(\cM_t)$.
We have 
$$\Sf(\cM_t,\cM_t^\prime)=-\Sf(\cM_t^{op},(\cM^\prime_t)^{op})=-\Sf(\cM_t^\prime,\cM_t)\ .$$

\subsubsection{} 
Let $\cM_{bt}$ be a boundary-tamed manifold and $\cM_t$ be an
extension of the boundary taming to a taming.
\begin{lem}\label{spg}

Let $i\in I_0(M)$ be such that $\partial_iM^{odd}$ is not empty.
Then for any $n\in\Z$ there exists an extension $\cM_t^\prime$ of
$\cM_{bt}$
to a taming which differs from $\cM_t$ only on the face $i$ such that 
$$\Sf(\cM_t^\prime,\cM_t)(i)=n\ .$$
\end{lem}
\proof
We can assume that $M$ is reduced and odd-dimensional.
The case $n=0$ is trivial. Therefore we can assume that $n\not=0$.
Furthermore, without loss of generality we can assume that $n>0$,
since otherwise we consider the opposite case.

Let $R$ be the orthogonal projection onto an $n$-dimensional space
spanned by eigenfunctions to positive eigenvalues of $D(\cM_{t})$.
Take $C=1+\|RD(\cM_t)\|$.
The spectral flow of the family $D_u:=D(\cM_t)-uCR$ is given by
$\Sf((D_u)_{u})=n$. 
Given $c>0$ let  $\tilde R=\chi R\chi$ be a compactly supported approximation such that $\|R-\tilde R\|\le c$
(see the proof of Lemma \ref{btt}).
If $c$ is sufficiently small, then we can take $W_i^\prime:=-C\tilde R + W_i$.
\hB

\subsubsection{}
Let $\cM_{bt}$ be a boundary-tamed manifold. Then 
$\ind(\cM_{bt})$ depends on the boundary taming.
If we allow for changes of the boundary taming on higher-codimensional
faces, then it  seems to be very complicated to understand how
$\ind(\cM_{bt})$ depends on the data.
 However, its dependence on the taming of the boundary faces is very
explicit. Let $\cM_{bt}$ and $\cM_{bt}^\prime$ two boundary tamings
which differ only on boundary faces, i.e.
for each $i\in I_1(M)$ the tamings
$\partial_i\cM_{bt}$ and $\partial_i\cM_{bt}^\prime$ are extensions of
the same boundary taming 
(the precise meaning of this assertion is that if we fix a distinguished
model $\cN_{geom}$ of $\partial_i\cM_{geom}$, then the boundary tamings
$\cN_{bt}$ and $\cN_{bt}^\prime$ induced by
the boundary tamings of $\cM_{bt}$ and $\cM_{bt}^{\prime}$ coincide).

\begin{lem}
For $i\in I_0(M)$ we have
$$\ind(\cM_{bt}^\prime)(i)-\ind(\cM_{bt})(i)=\sum_{j\in
  I_0(\partial \partial_i M)}\Sf(
\partial\partial_i\cM_{bt}^\prime,\partial \partial_i\cM_{bt})(j)\ .$$
\end{lem}
\proof 
First we reduce to the case that the face decomposition of $M$ is
reduced. If $M$ is odd-dimensional, then the assertion is trivial.
Therefore let us now assume that $M$ is even-dimensional.

We employ the  special case of the index theorem \ref{famind} where
$B$ is a point.
 It is then  a generalization to manifolds with corners  of the index theorem
of Atiyah, Patodi and Singer.
If $\cN_{bt}$ is a boundary-tamed manifold, then we have the
$\eta$-invariant $\eta^0(\cN_{bt})\in \R$ which is a spectral
invariant of $D(\cN_{bt})$ and generalizes the $\eta$-invariant of
Atiyah, Patodi and Singer. If $\cN_t$ and $\cN^\prime_t$ are two
extensions of the boundary taming to a taming,
then by Lemma \ref{diieta} 
$$\eta^0(\cN_t^\prime)-\eta^0(\cN_t)=\sum_{j\in I_0(N)}\Sf(\cN_t^\prime,\cN_t)(j)\ .$$
The index theorem for boundary-tamed manifolds applied to $\cM_{bt}$
yields
$$\ind(\cM_{bt})(i)=\Omega^0(\partial_i\cM_{geom})+ \eta^0(\partial\partial_i \cM_{bt})\ ,$$
where $\Omega^0(\partial_i\cM_{geom})\in\R $ only depends on the underlying
geometric manifold.
Thus 
\begin{eqnarray*}
\ind(\cM_{bt}^\prime)(i)-\ind(\cM_{bt})(i)&=&
\eta^0(\partial\partial_i \cM_{bt}^\prime)
-\eta^0(\partial\partial_i \cM_{bt})\\&=& 
\sum_{j\in
  I_0(\partial\partial_i M)}\Sf(
\partial \partial_i\cM_{bt}^\prime,\partial \partial_i\cM_{bt})(j)\ .
\end{eqnarray*}
\hB


 
\subsubsection{}
Let $\cM_{geom}$ be a geometric manifold.
The set of faces $I_*(M)$ is a partially ordered set where $i\le k$
if $i$ is a face of $k$. In addition, this partially ordered set is
graded, i.e. $I_*(M)=\sqcup_{k\ge 0} I_k(M)$. A special property of
this partially ordered graded set is the presence of adjacent pairs.
Assume that $k\ge 1$,  $i\in I_k(M)$, $l\in I_{k+1}(M)$ and $l<i$. Then
there is a unique $i^\prime\in I_k(M)$ which is adjacent to $i$
with respect to $l$, or equivalently, $l<i^\prime$. Note that $i$ is
\index{face!complex}\index{$\Face(M)$}then adjacent to $i^\prime$ with respect to $l$. It is the presence
of adjacent pairs which allows for the construction of the face
complex $\Face(M)$.
\newcommand{\bi}{\mathbf{i}} 
\newcommand{\bl}{\mathbf{l}} 
\newcommand{\bk}{\mathbf{k}} 
\newcommand{\bj}{\mathbf{j}} 
\newcommand{\bI}{\mathbf{I}} 

By an oriented face of $M$ we mean a pair $\bi =(i,\orient)$, where $\orient$ is an
orientation of $\partial_iM$. We let $\bi^{op}:=(i,-\orient)$ denote the
face with the opposite orientation. Let $\tilde \bI_*(M)$ be the set of
oriented faces of $M$. Note that faces may have several connected
components. So if $j$ is a boundary of $i$, then in general $\bj$ is not
induced from $\bi$ or $\bi^{op}$. 

Therefore we define a subset
$\bI(M)\subset \tilde{\bI}(M)$ of compatibly oriented faces
inductively as follows.
The subset $\bI_0(M)$ consists of those oriented faces and their
opposites which are the
underlying oriented manifolds of $\partial_i\cM_{geom}$, $i\in
I_0(M)$. For $k\ge 0$ the elements of $\bI_{k+1}(M)$ are those orientations
of faces with are induced as boundaries from the oriented faces
in $\bI_{k}(M)$.

We introduce the function
$\kappa:\bI_{*}(M)\times\bI_{*}(M)\rightarrow \Z$
which is the characteristic function of the relation ``is oriented
boundary of''. Thus $\kappa(\bj,\bi)=0$ unless
$\bi\in \bI_{k}(M)$ and $\bj\in \bI_{k+1}(M)$ for
some $k\ge 0$ and $j<i$, and the
\index{$\kappa(\bj,\bi)$}orientation of $\bj$ is induced from the orientation of $\bi$.
In the latter case $\kappa(\bj,\bi)=1$.
Note that the group $\Z/2\Z$ acts on $\bI_*(M)$ by
reversing orientations. This action preserves the grading and the
function $\kappa$.

\index{$\tilde \delta$}\index{$\widetilde{\Face}(M)$}Let $\widetilde{\Face}(M)$ be the $\Z$-graded free abelian group
generated by $\bI_*(M)$. We define the operator 
$\tilde \delta :\widetilde{\Face}(M)\rightarrow \widetilde{\Face}(M)$
of degree $-1$
by
$$\tilde \delta (\bj):=\sum_{\bi\in \bI_{*}(M),j<i}
\kappa(\bj,\bi) \bi\ .$$  
The group $\Z/2\Z$ acts on $\widetilde{\Face}(M)$ preserving degree
and $\tilde \delta$ such that the non-trivial element of $\Z/2\Z$
sends $\bi$ to $-\bi^{op}$.
\begin{ddd}
The face complex of $\cM_{geom}$ is defined as the graded abelian
group of
coinvariants 
$$\Face(M):=\widetilde{\Face}(M)_{\Z/2\Z}\ .$$
We let\index{$\delta$}
$$\delta:\Face(M)\rightarrow \Face(M)$$
be the operator of degree $-1$ induced by $\tilde \delta$.
\end{ddd}
Let $[\bi]$  denote the class in $\Face(M)$ represented by $\bi$. Then we
have the relation $[\bi^{op}]=-[\bi]$. 
We now check that $\delta^2=0$. In fact, let $\bj\in \bI_k(M)$. 
Then we have
$$
\delta^2 [\bj]=\sum_{\bk,\bi,j<i<k}
\kappa(\bj,\bi)\kappa(\bi,\bk)[\bk]\ .
$$
There are exactly two faces $i$ and $i^\prime$ which contribute to the
coefficient at $[\bk]$. In fact, the faces
$i$ and $i^\prime$ of $k$ are adjacent  with respect to $j$.
Let $\bi$ and $\bi^\prime$ be such that
$\kappa(\bj,\bi)\not=0$ and $\kappa(\bj,\bi^\prime)\not=0$.
Let $\bk$ be the orientation induced from $\bi$.
Then $\bk^{op}$ is induced from $\bi^\prime$ so that
$$\sum_{\bl,\bi,j<i<k, [\bl]=\pm [\bk]}
\kappa(\bj,\bi)\kappa(\bi,\bl)[\bl]=[\bk]+[\bk^{op}]=0\ .$$

\subsubsection{}
The obstructions against the existence of a taming of $\cM_{geom}$
will be homology classes of $\Face(M)$. The homology $H_*(\Face(M))$ depends on
the face decomposition of $M$ and may have torsion. We demonstrate this by
the following examples.

\subsubsection{}
The face complex of the geometric manifold $\cI_{geom}$ (the interval with atomic
face decomposition) is the complex
$$
0\rightarrow
 \Z[\mathbf{0}]\oplus
\Z[\mathbf{1}]\rightarrow \Z[\mathbf{c}]  \rightarrow 0\ ,$$
where $\mathbf{c},\mathbf{1}$ and $\mathbf{0}$ have the induced
orientations.
The differential is given by
$\delta[\mathbf{1}]=[\mathbf{c}]$
and $\delta[\mathbf{0}]=[\mathbf{c}]$.
The cohomology of this complex is $\Z$ concentrated in degree one and
generated by class of $[\mathbf{1}]-[\mathbf{0}]$.

\subsubsection{}
The face complex of the geometric manifold $\cJ_{geom}$ (the interval
with one boundary face $(0,1)$) is the complex
$$
0\rightarrow \Z[\mathbf{(0,1)}] 
\rightarrow  \Z[\mathbf{c}]\rightarrow 0\ ,$$
where $\mathbf{c},\mathbf{(0,1)}$ have the induced
orientations.
The differential is given by
$\delta[\mathbf{0,1}]=[\mathbf{c}]$.
The cohomology of this complex is trivial.
\subsubsection{}
In order to demonstrate the possibility of torsion in $H_*(\Face(M))$
we present an example found by Th. Schick.
We construct a geometric manifold $\cQ^4_{geom}$ with underlying manifold with corners
given by the unit cube
$Q^4:=[0,1]^4$. The group $\Z/2\Z$ acts by reflection with
respect to the center. We consider the face decomposition such that
a face is a pair of atoms of faces related by this action.
One checks that this decomposition is admissible.
With its standard orientation, metric, and spin-structure we obtain the
geometric manifold $\cQ^4_{geom}$.

We now compute $H_*(\Face(\cQ^4))$. 
The boundary of $Q^4$ is $\Z/2\Z$-equivariantly homeomorphic to
$S^3$. The atomic face decomposition of $\partial Q^4$ provides a $\Z/2\Z$-equivariant
cell decomposition of $S^3$. The face complex can be identified with
the associated cochain complex (up to a reversion and a shift
of the grading) augmented by the integration map. Therefore,
$\Face(\cQ^4_{geom})$ is related in this way to the augmented cochain complex of the induced
cell decomposition of $S^3/(\Z/2\Z)\cong \R P^3$.
It follows that 
$$H_k(\Face(\cQ^4_{geom}))\cong \tilde H^{k-1}(\R P^3,\Z)\cong \left\{\begin{array}{cc}
0&k=0,1,2\\\Z/2\Z&k=3\\
\Z&k=4
\end{array}\right. \ .$$

\subsubsection{}
We now discuss the inductive construction of a taming of $\cM_{geom}$
by decreasing codimension of faces.
Let $k\in\nat\cup\{0\}$ and assume that we are given operators $Q_j$ for all
$j\in I_l(M)$, $l>k$, 
such that they induce boundary tamings $\partial_m \cM_{bt}$, for all
$m\in I_k(M)$.
The precise meaning of this assertion is the following.
Let $(N_m,f)$ be a model of $\partial_m M$,
$T:N\times (-\infty,0)^k\rightarrow \bar M$ be a distinguished
embedding,
$\cW$ be a Dirac bundle over $N$ and $\Pi:T^* \cV\rightarrow \bar
\cW*(-\infty,0)^k$ be a product structure. We have an induced
structure of a geometric manifold $\cN_{m,geom}$.
Let $U_{l}$ be the
cylinder over $l\in I_*(N_m)$ in $\bar N_m$ and define
$U(l):=T(U_{l}\times (-\infty,0)^k)$.
Then for all $l\in I_h(N_m)\subseteq I_{h+k}(M)$, $h>0$, we define operators
$R_l$ such that
$$\Pi^{-1}\circ \bL_{U_{l}}^{U(l)}(R_l)\circ \Pi=Q_{l|U(l)}\ .$$
Our assumption is now that $\{R_l|l\in I_{>0}(N_m)\}$ is a boundary
taming $\cN_{m,bt}$ of $\cN_{m,geom}$.

We call $\{Q_j|j\in I_{>k}(M)\}$ a taming in codimension $k+1$.
Our notation for a geometric manifold equipped with this kind of data
is $\cM_{t,k+1}$.
Therefore a boundary taming is a taming in codimension one.


If $\cM_{t,k+1}$ is a taming of $\cM_{geom}$ in codimension $k+1$, then we have a
natural notion of its extension to a taming in codimension $k$. 
To give such an extension is equivalent to give extensions of the
boundary tamings $\cN_{m,bt}$ to tamings for all $m\in I_k(M)$.

\subsubsection{}
The task is now to extend the given taming $\cM_{t,k+1}$ in codimension $k+1$
to a taming $\cM_{t,k}$ in codimension $k$.
Equivalently, we must define extensions of the boundary tamings
$\cN_{m,bt}$ to tamings  for all $m\in I_k(M)$. The obstruction is given by Lemma \ref{btt}. 
\newcommand{\bm}{\mathbf{m}}
We let $\bm:=(m,\orient_{N_m})$ and define the chain
$$C(\cM_{t,k+1}):=\sum_{m\in I_k(M)} C_m [\bm]\in\Face_k(\cM_{geom})\ ,$$ where $C_m:=\ind(
\cN_{m,bt})(m)$ (note that $I_0(N_m)=\{m\}$ in a canonical manner).
Then an extension of the taming $\cM_{t,k+1}$ to a taming in
codimension $k$ exists iff $C(\cM_{t,k+1})=0$.

\subsubsection{}
If $C(\cM_{t,k+1})\not=0$, then there is no extension of $\cM_{t,k+1}$. In order to
proceed in this case we must revise the taming $\cM_{t,k+1}$.
A common feature with obstruction theory is the fact, that
the dependence of the obstruction chain on the choices made in the
previous step is under control, while the dependence on the other
parts of the data is complicated. We will see in Lemma \ref{faccoh} that $C(\cM_{t,k+1})$ is in fact
closed, and that by altering the extension $\cM_{t,k+1}$ of $\cM_{t,k+2}$
 it can be
changed to every element it its cohomology class. Consequently, the
taming $\cM_{t,k+2}$ can be extended twice  iff
 the cohomology class of $C(\cM_{t,k+1})$ vanishes.

\subsubsection{}
We fix the taming $\cM_{t,k+2}$ in codimension $k+2$.
Let $\cM_{t,k+1}$ be an extension of $\cM_{t,k+2}$.
\begin{lem}\label{faccoh}
The chain $C(\cM_{t,k+1})$ is closed. Its cohomology class
$$[C(\cM_{t,k+1})]\in H_k(\Face(\cM_{geom}))$$ only depends on 
$\cM_{t,k+2}$. 
For all $C^\prime\in [C(\cM_{t,k+1})]$ there exists an extension
$\cM^\prime_{t,k+1}$ of $\cM_{t,k+2}$ such that
$C^\prime=C(\cM^\prime_{t,k+1})$.
\end{lem}
\proof
Note that $\Face(\cM_{geom})$ is a complex of free $\Z$-modules. Therefore we can
consider it as a sub-complex of $\Face(\cM_{geom})\otimes_\Z\R$. 
We employ the index theorem Theorem \ref{etaprop} in order to
express $C(\cM_{t,k+1})$. We  will use the corresponding notation introduced above.
For  $j\in I_{k+1}(M)$ we have geometric manifolds $\cN_{j,geom}= 
(N_j,g^{N_j},\orient_{N_j},\cW_j)$, where $\orient_{N_j}$
is induced by the embedding $T_j:\bar N_j\times
(-\infty,0)^{k+1}\rightarrow \bar M$.
The taming $\cM_{k+1,t}$ in codimension $k+1$ induces tamings $\cN_{j,t}$. We set
$\bj:=(j,\orient_{N_j})$
and $\bm:=(m,\orient_{N_m})$ for $j\in I_{k+1}(M)$ and $m\in I_k(M)$.
Let us consider the chains with real coefficients 
\begin{eqnarray*}\eta&:=&\sum_{j\in I_{k+1}(M)} \eta^0(
  \cN_{j,t}) [\bj]\in\Face(\cM_{geom})\otimes_\Z \R\\
\Omega&:=&\sum_{m\in I_k(M)}  \Omega^0(\cN_{m,geom})
[\bm]\in\Face(\cM_{geom})\otimes_\Z \R\ .
\end{eqnarray*}
Then by Theorem \ref{etaprop} we have the relation  $$C(\cM_{t,k+1})=\Omega+\delta
\eta\ .$$
We compute
\begin{eqnarray*}
\delta \Omega 
&=& \sum_{m\in I_k(M)} \sum_{\bi,i\in I_{k-1}(M)}
\Omega^0(\cN_{m,geom})  \kappa(\bm,\bi) [\bi]\\
&=& \sum_{i\in I_{k-1} (M)}
\left[\sum_{\epsilon=\emptyset,op}\sum_{m\in I_{k}(M)}
  \kappa(\bm,\bi^\epsilon)  \Omega^0(\cN_{m,geom})
  [\bi^\epsilon]\right]\\
&=&0\ .
\end{eqnarray*}
In fact, for $i\in I_{k-1}(M)$  the term
$\left[\sum_{\epsilon=\emptyset,op}\sum_{m\in I_{k}(M)} (-1)^\epsilon
  \kappa(\bm,\bi^\epsilon)  \Omega^0(\cN_{m,geom})
  \right]$ is the integral of the closed local index
form  over the
boundary of $\partial_i M$.  
It follows that $\delta C(\cM_{t,k+1})=\delta \Omega =0$.

Let $\cM^\prime_{t,k+1}$ be another extension of $\cM_{t,k+2}$.
Let $\cN_{j,t}^\prime$ be the induced taming
 and
$\eta^\prime$ be the corresponding $\eta$-chain.
Then by Lemma \ref{diieta} we have
$$\eta^\prime-\eta
=\sum_{j\in I_{k+1}(M)} \Sf(\cN_{j,t}^\prime,\cN_{j,t}) [\bj]\ .$$
This  difference is an integral chain thus belongs to $\Face_{k+1}(\cM_{geom})$.
We conclude that
$$C(\cM_{t,k+1}^\prime)-C(\cM_{t,k+1})=\delta(\eta^\prime-\eta)\ ,$$
and thus $[C(\cM_{t,k+1}^\prime)]=[C(\cM_{t,k+1})]$.

Let now $C^\prime\in [C(\cM_{t,k+1})]$. Then there is a chain
$D\in \Face_{k+1}(\cM_{geom})$ such that
$C^\prime=C(\cM_{t,k+1})+\delta D$.
We write $D=\sum_{j\in I_{k+1}(M)} d_j [\bj]$.
By Lemma \ref{spg} we can alter for each $j\in I_{k+1}(M)$ the extension $\cN_{j,t}$ of
the boundary taming $\cN_{j,bt}$ (which is fixed by $\cM_{t,k+2}$)
to a taming $\cN_{j,k+1}^\prime$ such that
$$\Sf(\cN_{j,t},\cN_{j,t}^\prime)=-d_j\ .$$
This data gives $\cM^\prime_{t,k+1}$.
We compute
\begin{eqnarray*}
C(\cM_{t,k+1}^\prime)&=&C(\cM_{t,k+1})+\delta(\eta^\prime-\eta)\\
&=&C(\cM_{t,k+1})-\delta\sum_{j\in I_{k+1}(M)} \Sf(\cN_{j,t},\cN^\prime_{j,t}) [\bj]\\
&=&C(\cM_{t,k+1})+\delta D\\
&=&C^\prime 
\end{eqnarray*}
\hB





\subsubsection{}
While $\Face(\cM_{geom})\rightarrow \Face(\cM_{geom})\otimes_\Z\R=:\Face(\cM_{geom},\R)$ is
an embedding, the induced map $x\mapsto x_\R$ on homology is not injective in general
because of the presence of torsion.

We have seen that
the cohomology class $[C(\cM_{t,k+1})]\in
H_k(\Face(\cM_{geom}))$ only depends on $\cM_{t,k+2}$.
We expect that a description of this dependence in general is
complicated. But by surprise its image
$[C(\cM_{t,k+1})]_\R\in H_k(\Face(\cM_{geom})\otimes_\Z\R)$
turns out to be an invariant of the differential-topological structure
underlying $\cM_{geom}$. 

\begin{lem}
The class $[C(\cM_{t,k+1})]_\R\in H_k(\Face(M,\R))$  
only depends on the differential-topological structure underlying
$\cM_{geom}$.
\end{lem}
\proof
We consider the dual cochain complex $\widehat \Face(M)$
\index{$\widehat \Face(\cM_{geom})$}which is defined as the space of $\Z/2\Z$-invariants
$$\widehat \Face(\cM_{geom}):=\Hom_\Z(\widetilde\Face(\cM_{geom}),\Z)^{\Z/2\Z}\
.$$
The differential $d:\widehat \Face(M)\rightarrow \widehat \Face(M)$ is
induced by $\delta$. 
There is a natural pairing
$$\Face(\cM_{geom})\otimes_\Z \widehat\Face(\cM_{geom})\rightarrow
\Z$$ which induces a pairing
$$\langle .,.\rangle:H_k(\Face(\cM_{geom}))\otimes H^k(\widehat
\Face(M_{geom}))\rightarrow \Z\ .$$
Note that a class $[C]_\R\in H_k(\Face(\cM_{geom},\R))$ is uniquely determined by
the numbers $\langle  [C],u\rangle\in \Z$, $u\in H^k(\widehat
\Face(\cM_{geom}))$.

Let us choose for each $i\in I_*(M)$ some orientation $\bi$.
Then we have basis $(\bi,\bi^{op})_{i\in I_*(M)}$ of
$\widetilde\Face(\cM_{geom})$.
Let $(\hat \bi,\hat \bi^{op})_{i\in I_*(M)}$ be the dual basis of $\Hom_\Z(\widetilde\Face(\cM_{geom}),\Z)$.
Then a basis of $\widehat \Face(\cM_{geom})$ is given by
$(\hat \bi-\hat \bi^{op})_{i\in I_*(M)}$. This basis is dual to
$([\bi])_{i\in I_*(M)}$. 
For $j\in I_k(M)$ we have 
$d(\hat \bj-\hat \bj^{op})=\sum_{i\in I_{k+1}(M),i<j}\kappa(i,j)(\hat
\bi-\hat \bi^{op})$.

Let $u$ be represented by $U=\sum_{j\in I_k(M)} u_j (\hat \bj-\hat \bj^{op})$.
Using Theorem \ref{etaprop} and
$dU=0$ we obtain (notation of the proof of Lemma \ref{faccoh})
\begin{eqnarray*}
\langle  [C(\cM_{t,k+1})],u\rangle 
&=&<\Omega,U>+<\delta\eta,U>\\
&=&<\Omega,U>+<\eta,dU>\\
&=&<\Omega,U>\\
&=&\sum_{j\in I_k(M)} \Omega^0(\partial_j \cM_{geom})    u_j \ .
\end{eqnarray*}
The right-hand side of this equation is an integer which only depends on $\cM_{geom}$. Since any two geometries on the same underlying differential-topological structure can be joined by a path
we see  by continuity that $\langle  [C(\cM_{t,k+1})],u\rangle$ is independent of the geometry. \hB

\section{Geometric families}\label{uuzz66}

\subsection{Families of manifolds with corners}\label{laq}

\subsubsection{}
Let $B$ be a smooth manifold and $M$ be a manifold with corners.
\index{bundle of manifolds with corners}A locally trivial bundle of manifolds with corners over $B$ with fiber $M$
is a manifold with corners $E$ together with a morphism of manifolds
with corners  $\pi:E\rightarrow B$ and an atlas of local trivializations
$\Phi_U:\pi^{-1}(U)\stackrel{\sim}{\rightarrow} M\times U$ for suitable open subsets
$U\subseteq B$ such the transition maps $\Phi_{U,V}:\Phi_U\circ
\Phi_V^{-1} :M\times (U\cap V)\rightarrow M\times (U\cap V)$
are of the form $(m,v)\mapsto (\phi_{U,V}(v)(m),v)$, where $\phi_{U,V}(v):M\stackrel{\sim}{\rightarrow} M$
is an automorphism of manifolds with corners for all $v\in U\cap V$.  


\subsubsection{}
Let $F$ be an atom of faces of codimension $k$ of $E$.
It comes with an induced map $q:F\rightarrow B$. This is again a
locally trivial bundle of manifolds with corners.
In fact, let $f\in F$ and $U$ be a neighborhood of $q(f)$
which is the domain of a local trivialization $\phi:E_{|U}\stackrel{\cong}{\rightarrow}
M\times U$. Then there is a  collection of atoms of faces $A_\alpha$ of $M$ such that
$\phi_{|F_{|U}}:F_{|U}\stackrel{\cong}{\rightarrow} 
\sqcup_\alpha A_\alpha\times U$. This is an element of the atlas of local
trivializations of $q:F\rightarrow B$.

Let $T:F\times  [0,1)^k \rightarrow E$ be a distinguished embedding.
Then we have $\pi\circ T=q\circ \pr_F$. Thus, if we consider
$q\circ \pr_F:F\times [0,1)^k\rightarrow B$ as a bundle, then
$T$ is a local isomorphism of bundles.

\subsubsection{}
Let $\pi:E\rightarrow B$ be a locally trivial bundle of manifolds with
corners. 
Since $E$ is a manifold with corners
\index{$\bar E$}\index{${}_r\bar E$}we can form the extension $\bar E$ and the enlargements ${}_r\bar E$, $r\ge 0$. The same objects can be obtained by the corresponding fiber-wise constructions. In particular, $\bar E$ is a locally trivial bundle with fiber $\bar M$, and
${}_r \bar E$ is a locally trivial bundle of manifolds with corners over $B$ with fiber ${}_r\bar M$.

\subsubsection{}
We now start with the discussion of the additional geometric structures which
refine $\pi:E\rightarrow B$ to a geometric family.

First of all, an admissible face decomposition is an admissible face
decomposition of $E$. Note that it may happen, that the fiber  admits
an admissible face decomposition, but $E$ does not.
Consider e.g the bundle $E\rightarrow S^1$ with fiber $\Delta^2$ which
is obtained as the quotient $ \Delta^2\times \R/\Z$, where $\Z$ acts by
$(y,x)\mapsto (\sigma y,x+1)$, and where
$\sigma:\Delta^2\rightarrow\Delta^2$ is a non-trivial rotation. 
Then $E$ has one atom of faces of codimension one and two,
respectively. The boundary atom self-intersects in the corner.
 
\subsubsection{}
A geometric manifold comes with a Riemannian metric and an
orientation. In the case of families we have corresponding structures
on the vertical bundle.
The vertical bundle $T^v\pi$ is the sub-bundle of $TE$ given by 
$T^v\pi:=\ker(d\pi)$.\index{vertical bundle}
\index{$T^v\pi$}

\index{fibre-wise!orientation}\index{$\orient_\pi$} A fiber-wise orientation $\orient_\pi$ is an orientation  of the
vertical bundle $T^v\pi$ (see \ref{setupor}).
 
A vertical admissible Riemannian metric is a metric on $T^v\pi$ which induces
an admissible Riemannian metric (see \ref{udiwqdwqdq}) on each fiber by restriction.
Note that a locally trivial bundle of manifolds with corners admits a
vertical Riemannian metric. In fact, we can choose an admissible
metric on the model fiber $M$. This gives admissible metrics on the
restriction $E_{|U}=\pi^{-1}(U)$ for each domain of a local
trivialization chart $E_{|U}\cong M\times U$. These metrics can be
glued using a partition of unity over $B$. 
Here we use the following facts. If $\phi:M\to M$ is an automorphism of manifolds with corners (see \ref{hdiqwdqwd}) and $g$ is an admissible metric on $M$, then $\phi^* g$ is again admissible. Furthermore, if $\lambda\in [0,1]$ and $g_0,g_1$ are admissible metrics on $M$, then
$\lambda g_0+(1-\lambda)g_1$ is again admissible.

\subsubsection{}

\index{horizontal!distribution}
\index{$T^h\pi$}
A horizontal distribution $T^h\pi\subseteq TE$ is a complement to the
vertical bundle in $T^v\pi$. We require that the horizontal
distribution has a product structure near the singularities.
In the following we describe this condition in detail.
Let $i\in I_k^{atom}(E)$, $k\ge 1$. Let $(F,f)$ be a model of
$\partial_iM$, $q:F\rightarrow B$ be the induced bundle and
$T:F\times [0,1)^k\rightarrow E$ be a distinguished embedding.
If $T^hq\subseteq TF$ is a horizontal bundle of $q$, then it induces a
horizontal bundle $\pr_F^*(T^hq)$ of $q\circ\pr_F:F\times [0,1)^k\rightarrow B$ in the
natural way. 
\begin{ddd} We call a horizontal distribution  $T^h\pi$ admissible, if for all
$k\ge 1$, $i\in I_k^{atoms}(E)$, model   $(F,f)$ of $\partial_iE$ and
distinguished embedding $T:F\times [0,1)^k\rightarrow E$ there exists
a horizontal distribution $T^hq$ of $q:F\rightarrow B$ such that
$T^*T^h\pi=\pr_F^*(T^hq)$.
\end{ddd}

\newcommand{\tT}{{\tt T}}

A locally trivial bundle of manifolds with corners admits an
admissible horizontal distribution. In fact we can construct
one by pasting horizontal distributions on local trivializations
using a partition of unity.

First we represent a horizontal distribution $T^h\pi$ by its connection one
form which is the unique element of  $\omega\in C^\infty(E,\Hom(T^*E,
T^v\pi))$ such that $T^h\pi=\ker(\omega)$ and $\omega_{|T^v\pi}=\id_{T^v\pi}$.

On $U\times M$ we take  $T^h\pr_U:=TU$ and let $\omega_U$ be the
corresponding connection one-form.
If $E_{|U}\cong M\times U$ is a local trivialization, then we let
$\omega_U$ denote the corresponding connection one-form on $E_{|U}$.
We choose a locally finite cover $(U_\alpha)$ of $B$ of domains of
local trivialization charts.
We can paste the forms $\omega_{U_\alpha}$ using a subordinated partition of unity. 
The resulting connection one-form  form defines an admissible horizontal distribution.

\subsubsection{}\label{curintro}

The horizontal distribution $T^h\pi$ determines a decomposition
$TE=T^v\pi\oplus T^h\pi$ and therefore a projection
\index{horizontal!curvature}\index{$\tT$}$\pr_{T^v\pi}:TE\rightarrow T^v\pi$. The curvature tensor $\tT\in C^\infty(E,\Lambda^2
(T^h\pi)^*\otimes T^v\pi)$ of a horizontal distribution is given by
$$\tT(X,Y):=\pr_{T^v\pi} [X,Y]\ ,X,Y\in C^\infty(E,T^h\pi)$$
(one easily checks that this formula defines a tensor). 

\subsubsection{}\label{duiqdwqqwdq}

Finally we consider the notion of a family of admissible Dirac bundles.
Let $\pi:E\rightarrow B$ be equipped with an admissible vertical Riemannian metric and an admissible horizontal distribution $T^h\pi$. Then we have an induced connection $\nabla^{T^v\pi}$
on the vertical bundle (\cite[Prop. 10.2]{berlinegetzlervergne92}). 
A family of admissible Dirac bundles is given by
a tuple $\cV=(V,h^V,\nabla^V,c,z)$ if the fibers are even-dimensional, and by $\cV=(V,h^V,\nabla^V,c)$ in the case of odd-dimensional fibers.
Here $(V,h^V,\nabla^V)$ is a hermitian vector bundle with connection over $E$,
$z$ is a parallel $\Z/2\Z$-grading, and
$c\in \Hom(T^v\pi,\End(E))$ is a  parallel sections with respect to the connection induced by
$\nabla^{T^h\pi}, \nabla^V$  are such that for all $b\in B$ the restriction
$\cV_{|E_b}$ of $\cV$ to the fiber $E_b:=\pi^{-1}(\{b\})$ is an admissible  Dirac bundle

\subsubsection{}
\begin{ddd}\label{geo987}
\index{geometric!family}\index{$\cE_{geom}$} A geometric family $\cE_{geom}$ over $B$ is given by the following structures:
\begin{enumerate}
\item a manifold with corners $M$,
\item a locally trivial fiber bundle $\pi:E\rightarrow B$ with fiber $M$,
\item an admissible face decomposition of $E$,
\item an admissible vertical Riemannian metric $g^{T^v\pi}$,
\item an admissible horizontal distribution $T^h\pi$,
\item a family of admissible Dirac bundles $\cV$,
\item a fiber-wise orientation.
\end{enumerate} 
\end{ddd}

\subsubsection{}

\index{fibre-wise!spin structure}A fiber-wise spin structure is a spin structure of $T^v\pi$.
If we have a locally trivial fiber bundle $\pi:E\rightarrow B$ equipped with admissible vertical Riemannian metric and horizontal distribution, fiber-wise orientation and fiber-wise spin structure, then the fiber-wise spinor bundle
$S(\pi)$ admits the structure $\cS(\pi)$ of a family of admissible
Dirac bundles. More general examples can be constructed by twisting
the fiber-wise spinor bundle with auxiliary $\Z/2\Z$-graded hermitian vector bundles with connection. Locally over $E$ every family of admissible Dirac bundles is isomorphic to a twisted spinor bundle.

\subsubsection{}

Let now $\cE_{geom}$ be a geometric family.
Using the local trivializations of $E$ we can consider $\cE_{geom}$ as a family
of geometric manifolds $\cM_{geom,b}$, $b\in B$, with structures smoothly parameterized by $B$.
\index{fibre-wise!pre-taming}Thus it makes sense to speak of a fiber-wise pre-taming. 
In order to state the smoothness assumptions precisely let 
us introduce this notion locally. 
Assume that $E=M\times U$ is a trivial bundle of manifolds with
corners. Let $k\ge 0$ and $i\in I_k(M)$. Let $(N,f)$ be a model of $\partial_i
M$ with distinguished embedding $T:N\times (-\infty,0]^k\rightarrow
\bar M$. It induces a fiber-wise orientation and a family of admissible Riemannian metrics on
$N\times U$. We further get a complex vector bundle $W$ over $N\times
U$ such that the family of its restrictions to $N\times \{u\}$, $u\in U$, 
refines to a family of admissible Dirac bundles
$(\cW_u)_{u\in U}$ over
$N$ with product structures $\Pi_u:T^*\cV_{|\{u\}\times M}\rightarrow
\cW_{u}*(-\infty,0]^k$. 
Let $(R_u)_{u\in U}$ be a family of compactly supported smoothing operators 
on $C^\infty(\bar N,\bar \cW_u)$ which for each $u$ are as in 
\ref{tamam34}. We assume that this family is smooth in the sense that the
 integral kernel is a smooth section of $\pr_{1}^*W\otimes \pr_2^*
W^*$ over
$ N \times N\times U$. Then the family operators
$Q_u:=\Pi^{-1}_u\circ L_{\bar N_u}^{U_i}(R_u)\circ \Pi_u$ defines a
fiber-wise  pre-taming of the face $i$.


\subsubsection{}

Let $\cE_{geom}$ be a geometric family over $B$. 
\begin{ddd}
\index{$\cE_{geom}$}A pre-taming of $\cE_{geom}$ is given by (fiber-wise) pre-tamings of
all faces of $\cE_{geom}$.
 \end{ddd}

Our notation for a geometric family with distinguished pre-taming
is $\cE_t$.

\begin{ddd} 
The pre-taming of $\cE_{t}$ is called a taming if it induces a taming
of all the fibers $\cM_{geom,b}$, $b\in B$. 
\end{ddd}
In a similar way we define the notions of a boundary pre-taming and a
boundary taming of a geometric family.

\subsubsection{}

If $\cE_{\sharp}$, $\sharp\in \{geom,bt,t\}$, is a (decorated) family over $B$, and
$f:B^\prime\rightarrow B$ is a smooth map, then we can define the (decorated) family
$f^*\cE_{\sharp}$ over $B^\prime$ in a natural way. 

If $\cE_{i,\sharp}$, $i=1,2$, are two (decorated) families over $B$,
then there is a  natural notion of a fiber-wise sum $\cE_{1,\sharp}\sqcup_{B}
\cE_{2,\sharp}$.

To a family $\cE_{\sharp}$ we associate its reduction
$\cE_{\sharp,red}$.

\subsubsection{}

Let $\cE_{geom}$ be a geometric family  and $H$ be a Riemannian spin
manifold with corners. Then we can form the geometric family
$\cE_{geom}*H$ by applying the construction \ref{pps} fiber-wise.

\subsubsection{}

For all $i\in I_1(E)$ we can form the boundary $\partial_i
\cE_{\sharp}$. This is a well-defined isomorphism class of decorated
families. We leave it to the reader to write down the family
version of the constructions \ref{bounfgg} and \ref{bounftt}.

The boundary of a boundary tamed family is an isomorphism class of
tamed geometric families.
 
\subsubsection{}

\index{opposite!geometric family}\index{$\cE_\sharp^{op}$}We form  the opposite family $\cE_\sharp^{op}$ of
$\cE_\sharp$ by taking the opposite fiber-wise. 
Sometimes we will write $- \cE_\sharp:=\cE_\sharp^{op}$.

\subsubsection{}

\index{$\cE_{bt}$}\index{Dirac operator!of a boundary tamed family}Assume that $\cE_{bt}$ is a boundary tamed geometric family.
Then we have a family of (unbounded) selfadjoint Fredholm operators $(D(\cE_{bt,b}))_{b\in B}$.
If the fibers are even-dimensional, then using the grading $z$ we decompose
$$D(\cE_{bt,b})=\left(\begin{array}{cc}0&D(\cE_{bt,b})^-\\D(\cE_{bt,b})^+&0\end{array}\right)\ .$$

We have a decomposition 
 $D(\cE_{bt,b})=\oplus_{\epsilon\in \{0,1\},i\in
   I_0(E)}D(\cE_{bt,b})(\epsilon,i)$, where
$D(\cE_{bt,b})(\epsilon,i)$ is the restriction of $D(\cE_{bt,b})$
to the even- ($\epsilon=0$) or odd-dimensional ($\epsilon=1$)
part of the face $i$.
 If $B$ is compact, then
a family of Fredholm operators  parameterized by
$B$ has an index in the complex $K$-theory of $B$.
We refer to \ref{indele} for details,  particular for the explanation of the necessary continuity assumptions. The conventions in the odd-dimensional case are fixed in \ref{sshix}.
\begin{ddd}
\index{index!of a boundary tamed family}\index{$\ind(\cE_{bt})$}We define the index of the boundary tamed family 
as the function
$$\ind(\cE_{bt}):I_0(E)\rightarrow K(B)$$
such that
$$\ind(\cE_{bt})(i):=\ind(D(\cE_{bt,b})^+(0,i))\oplus
\ind(D(\cE_{bt,b})(1,i))\in K^0(B)\oplus K^1(B)\ .$$
\end{ddd}
Note that $\ind(\cE_{bt}^{op})=-\ind(\cE_{bt})$.

If the boundary tamed family $\cE_{bt}$ is reduced we will often write
$\ind(\cE_{bt}):=\ind(\cE_{bt})(*)$, where $*$ is the unique face
of codimension zero.

\subsubsection{}
The following Lemma generalizes Lemma \ref{btt} to families.
\begin{lem}\label{fambtt}
Assume that $\cE_{bt}$ is a boundary tamed family over a compact base
$B$. If the boundary taming can be extended to a taming, then
$\ind(\cE_{bt})=0$. For the converse we  assume in addition  that no $i\in I_0(E)$ contains a zero-dimensional  component. If $\ind(\cE_{bt})=0$, then
the boundary taming $\cE_{bt}$
can be extended to a taming. 
\end{lem}
\proof
It is clear that if the boundary taming can be extended
to a taming, then a compact perturbation of the Dirac operator becomes
invertible so that the index vanishes. Let us show the converse in
detail.

We can assume that $E$ is reduced and has either even- or
odd-dimensional fibers.
By \ref{jkl12}, \ref{se564} there exists a continuous family of compact operators
$(K_b)_{b\in B}$ which are self-adjoint (and odd in the
even-dimensional case) such that
$D(\cE_{bt,b})+K_b$ is invertible for all $b\in B$.

By an arbitrary small deformation of $K_b$ we can assume that
$b\mapsto K_b$ is a smooth family of compact operators. In fact
using the local trivializations of $E$ we trivialize the associated bundle of
Hilbert spaces $(L^2(\bar E_b,\bar V_{|E_b}))_{b\in B}$. Then locally we can consider
$(K_b)_{b\in B}$ as a family of compact operators on a fixed Hilbert space.
By the usual convolution technique locally we can find  smooth
families of compact operators
arbitrary close to $(K_b)_{b\in B}$. Gluing these local smoothings using 
a partition of unity on $B$ we find a smooth approximation
of $(K_b)_{b\in B}$.

Let us now assume that $(K_b)_{b\in B}$ depends smoothly on $b\in B$.
For $t>0$ we define the smooth family of smoothing operators
$$K_b(t):=\ee^{-tD(\cE_{bt,b})^2} K_b \ee^{-tD(\cE_{bt,b})^2}\
.$$  Note that $K_b(t)\stackrel{t\to 0}{\rightarrow} K_b$
uniformly on $B$. Thus we can fix a sufficiently small $t>0$ such that 
$D(\cE_{bt,b})+K_b$ is invertible for all $b\in B$.
Finally we use cut-off functions as in the proof of Lemma \ref{btt} in
order to produce a smooth family of compactly
supported smoothing 
operators $(\tilde K_b(t))_{b\in B}$ close to $(K_b(t))_{b\in B}$ which gives the extension
of the boundary taming to a taming. 
\hB 

\subsubsection{}
The additional assumption in Lemma \ref{fambtt} is necessary in
general since there may exist stably trivial but non-trivial vector
bundles. See \ref{glor1}.

\subsubsection{}
One could develop an obstruction theory for families in a similar manner as
in Subsection \ref{obst}.

\subsection{Some examples}\label{exex}

\subsubsection{}\label{ghghghgh11}

Let $\bV=(V,h^V,\nabla^V,z_V)$ be a $\Z/2\Z$-graded complex vector bundle over $B$ with hermitian metric and metric connection. Then we obtain a geometric family $\cE(\bV)_{geom}$ as follows.
The underlying fiber bundle of $\cE(\bV)_{geom}$ is $\pi:=\id_B:B\rightarrow B$.
The vertical bundle is trivial (zero-dimensional). Therefore $\cE(\bV)_{geom}$
has a canonical fiber-wise orientation, vertical Riemannian metric, and horizontal distribution. The family of Dirac bundles is
$\cV:=(V,h^V,\nabla^V,0,z_V)$.

Since the fibers of $\cE(\bV)_{geom}$ are closed (a point is a closed manifold) this family is boundary tamed. We have 
$$\ind(\cE(\bV)_{bt})=[V]:=[V^+]-[V^-]\in K^0(B)\ , $$ where $[V^\pm]$ are
the classes represented by $V^\pm$.

The family $(\cE(\bV)_{geom}\sqcup_B \cE(\bV)^{op}_{geom})_{red}$
admits a taming, while $\cE(\bV)_{geom}\sqcup_B \cE(\bV)^{op}_{geom}$
 only admits a taming if $V^+\cong V^-$.

\subsubsection{} \label{glor1}
We have the following consequence of Lemma \ref{fambtt}.
\begin{kor}
Let $\cE_{bt}$ be a boundary tamed family over a compact base $B$ such
that
$\ind(\cE_{bt})=0$.
Then there exists $n\ge 0$ such that the boundary taming of 
$(\cE_{bt}\sqcup_B \cE(\bW)_{geom})_{red}$ admits an extension to a taming.
Here $W:=B\times (\C^n\oplus \C^n)$ with grading $z_W:=\diag(1,-1)$.
\end{kor}

\subsubsection{}\label{eufiewfewfew}

We now consider the odd-dimensional analog of this construction.
Let $n\in \nat$ and $F:B\rightarrow U(n)$ be a smooth map.
Then we construct a geometric family $\cE(F,*)_{geom}$ as follows
(the argument $*$ shall indicate, the $\cE(F,*)_{geom}$ depends on
additional choices).

We consider the bundle $\pr_B:\R\times B\rightarrow B$
with the obvious orientation, spin-structure, horizontal distribution, and vertical metric.
Over this bundle we consider the hermitian vector bundle
$\tilde W:=\R\times B\times \C^n\rightarrow \R\times B$.
We define a fiber-wise free proper $\Z$-action on $\R\times B$ by
$\Z\times \R\times B\ni(k,t,b)\mapsto (t+k,b)\in \R\times B$.
This action preserves the orientation, spin-structure, horizontal
distribution, and the vertical metric. 
Furthermore, this  action lifts to $\tilde W$ such that
$\Z\times \R\times B\times \C^n\ni(k,t,b,v)\mapsto (t+k,b,F(b)^kv)\in \R\times B\times \C^n$.
We choose a $\Z$-invariant metric connection on $\tilde W$.
We thus obtain a $\Z$-equivariant geometric bundle $\tilde \bW$.
Then we form the twisted bundle $\tilde \cV:=\cS(\pr_B)\otimes \tilde \bW$.

If we identify $S^1\cong \R/\Z$, then the underlying fiber bundle of
$\cE(F,*)_{geom}$ is $S^1\times B\rightarrow B$ with induced structures.
The Dirac bundle is $\cV:=\tilde \cV/\Z$.

\index{$\cE(F,*)_{geom}$}The family $\cE(F,*)_{geom}$ has closed fibers and is thus boundary tamed.
We have
\begin{equation}\label{uuiwed44}
\ind(\cE(F,*)_{bt})=[F]\in K^1(B)\ ,
\end{equation}
where $[F]\in [B,U(\infty)]\cong K^1(B)$. 
To see this formula note the following facts (see also \ref{sq21}). Let $U:=\C^n\times S^1\times B\to S^1\times B$ be the trivial bundle, and set  $W:=\tilde W/\Z\to S^1\times B$. The class
$[W]-[U]\in K^0(S^1\times B)$ corresponds 
to $[F]\in K^1(B)$ under the suspension map
$$s:K^1(B)\cong  \tilde K^0(S^1\wedge B_+)\stackrel{S^1\times B\to S^1\wedge B_+}{\to} K^0(S^1\times B)\ .$$ The fibre-wise spin structure of the bundle $\pi:S^1\times B\to B$ gives a $K$-orientation of $\pi$ and therefore an integration  map $\pi_!:K^0(S^1\times B)\to K^1(B)$. Explicitly, for a vector bundle $W\to S^1\times B$ the class $\pi_![W]$
is represented by the family of fibre-wise Dirac operators along $\pi$ twisted by $W$. In other words,
$\ind(\cE(F,*)_{bt})=\pi_![W]$. The equation (\ref{uuiwed44}) now follows from
$\pi_!\circ s=\id_{K^*(B)}$ and $\pi_![U]=0$.

\subsubsection{}

Let us illustrate the definition of the index of a boundary tamed family in a simple example.
We consider the bundle $\pi:E:=[-1/2,1/2]\times S^1\to S^1$ with fibre $[-1/2,1/2]$. 
We combine the endpoints of the interval into one face of codimension
one, i.e.  $I_1(E):=\{c\}$. The bundle $E$  comes with the product metric and product horizontal distribution. It has a canonical fibre-wise spin structure, and we consider the Dirac bundle $\cS(\pi)$ given by the fibre-wise spinors. This data defines a geometric family $\cE_{geom}$.

 A boundary taming $\cE_{bt}$ is a taming of $\partial_c \cE_{geom}$. Let $\cF_{geom}$ be the geometric family over $S^1$ with underlying bundle $\id_{S^1}:S^1\to S^1$, canonical geometric structures and the trivial Dirac bundle $\cW:=S^1\times \C\to \C$.
Then we have an isomorphism $\partial_c\cE_{geom}\cong (\cF_{geom}\sqcup_{S^1}\cF_{geom}^{op})_{red}$.
A taming of this family is given by an invertible  selfadjoint section of $\End(\cW\oplus \cW^{op})^ {odd}$. For simplicity  we assume that it is unitary. Then   it is given
by a matrix $$\left(\begin{array}{cc}0&u^*\\u&0\end{array}\right)$$ for some smooth function  $u:S^1\to U(1)$. 

Let us now calculate $\ind(\cE_{bt})\in K^1(S^1)\cong \Z$ in terms of the function $u$. In fact, this index is the spectral flow of the loop of selfadjoint Fredholm operators
$(D(\cE_{bt,b}))_{b\in S^1}$.  Let us describe these operators explicitly.
The extension of $\cE_{geom,b}$ is the extension $\R$ of $[-1/2,1/2]$.
Let $\rho\in C^\infty(\R)$ be a cut-off function such that $\rho(x)=1$ for $x\ge 1$ and
$\rho(x)=0$ for $x<3/4$. Using the canonical trivialization of the Dirac bundle $\cS\cong \C\times E$   we can consider  $D(\cE_{bt,b})$ as an operator on functions on $\R$. We have
$$D(\cE_{bt,b})f(x)=i \frac{d}{dx}f(x) + \rho(x) \bar u(b) f(-x) + \rho(-x) u(b) f(x)\ .$$

We translate this non-local operator to a two-component local operator on $[0,\infty)$ by
setting $(f_1(x),f_2(x))=(f(x),f(-x))$ for $x\ge 0$.
In terms of $(f_1,f_2)$ we have
$$D(\cE_{bt,b})\left(\begin{array}{c}f_1\\f_2\end{array}\right)=  \left(\begin{array}{cc}i\frac{d}{dx}&\bar u(b)  \rho \\u(b)\rho&-i\frac{d}{dx}\end{array}\right) \left(\begin{array}{c}f_1\\f_2\end{array}\right)\ ,$$
and we must impose the transmission boundary \index{transmission boundary condition}conditions
$$f_1(0)=f_2(0)\ .$$
We can deform this family of operators to the family of constant coefficient operators
$$A_v=\left(\begin{array}{c}f_1\\f_2\end{array}\right)=  \left(\begin{array}{cc}i\frac{d}{dx}&\bar v  \\v&-i\frac{d}{dx}\end{array}\right) \left(\begin{array}{c}f_1\\f_2\end{array}\right)$$ where $v:=u(b)$
without changeing the spectral flow. By a simple calculation we see that
$A_v$ has a real eigenvalue if and only if $\Imm(v)<0$. In this case the unique real eigenvalue is $\Ree(v)$, and 
the corresponding $L^2$-eigenvector is given by 
$$\left(\begin{array}{c}f_1(x)\\f_2(x)\end{array}\right)=\ee^{\Imm(v) x} \left(
\begin{array}{c}1\\1 \end{array}\right)\ .$$
Let $\deg(u)\in \Z$ denote the winding number of the function $u:S^1\to U(1)$.
We get $\Sf(A_{u(b)})_{b\in S^1}=-\deg (u)$ and hence 
$$\ind(\cE_{bt})=-\deg(u)\ .$$

\subsection{Local index forms}

\subsubsection{}
If $B$ is a smooth manifold, then by \index{$\cA_B$}\index{sheaf of differential forms}$\cA_B$ we denote the $\Z$-graded sheaf of real smooth differential forms on $B$.


\subsubsection{}
Let $\cV$ denote the Dirac bundle of $\cE_{geom}$.
Locally on $E$ we can write $\cV$ as a twisted spinor bundle
$\cV:=\cS(\pi)\otimes \bW$,
where $\bW=(W,h^W,\nabla^W,z^W)$ is a $\Z/2\Z$-graded hermitian
vector bundle with connection which is called the twisting bundle.
We have
$W:=\Hom_{\Cl(T^v\pi)}(S(\pi),V)$ with induced structures if the fibers are even-dimensional, and
\index{Chern!form}\index{$\ch(\nabla^{W})$}$W:=\Hom_{\Cl(T^v\pi)}(S(\pi)\oplus S(\pi)^{op},V)$, if the fibers are odd-dimensional. The Chern form 
$$\ch(\nabla^{W}):=\tr_s \exp\left(-\frac{R^{\nabla^W}}{2\pi i}\right)\in \cA_E(E)$$ is globally
defined. 

\subsubsection{}

The fiber-wise orientation induces an integration map\index{integration over the fibre}
$\int_{E/B}:\cA_E(E)\rightarrow\cA_B(B)$ of degree $\dim(B)-\dim(E)$.
In order to fix signs we define this integration such that for a
trivial
bundle $M\times U\rightarrow U$ we have
$$\int_{M\times U/U} \alpha\times \beta =(\int_M \alpha)\beta\ ,$$
$\alpha\in \cA_M(M)$, $\beta\in \cA_U(U)$ and
$\alpha\times\beta=\pr_M^*\alpha\wedge \pr_U^*\beta$. 
Stoke's theorem takes the form
$$d \int_{E/B} \omega + \int_{\partial E/B} \omega
= (-1)^{\dim(E)-\dim(B)} \int_{E/B} d\omega 
 $$
 for $\omega\in \cA(E)_E$.

\subsubsection{}

The vertical bundle $T^v\pi$ of a geometric family $\cE_{geom}$ has a natural connection
\index{$\hA$-form}\index{$\hA(\nabla^{T^v\pi})$}$\nabla^{T^v\pi}$ (see \ref{duiqdwqqwdq}). We define the
$\hat A$-form
$$\hA(\nabla^{T^v\pi}):=\det {}^{1/2}\left(\frac{\frac{R^{\nabla^{T^v\pi}}}{4\pi }}{\sinh\left(\frac{R^{\nabla^{T^v\pi}}}{4\pi}\right)}\right)\in \cA_E(E)
$$

\begin{ddd}\label{eet0}
\index{local!index form}\index{$\Omega(\cE_{geom})$}We define the local index form 
$$\Omega(\cE_{geom}):=\int_{E/B}
\hA(\nabla^{T^v\pi})\ch(\nabla^{W})\in\cA_B(B)\ .$$
\end{ddd}

\subsubsection{}

\index{local!index theorem} The main result of local index theory for families is the following theorem.
We assume that $\cE_{geom}$ is reduced and has closed fibres.
Then we have $\ind(\cE_{bt})\in K(B)$. Let $\ch:K(B)\to H(B,\Q)$ be the Chern character.

\begin{theorem}\label{loinde}
The local index form  $\Omega(\cE_{geom})$ is a closed form.
Its de Rham cohomology class $[\Omega(\cE_{geom})]\in H_{dR}(B)$
is equal to the image of $\ch(\ind(\cE_{bt}))\in H(B,\Q)$ under the de Rham map
$\dR:H(B,\Q)\rightarrow H_{dR}(B)$.
\end{theorem}

\subsubsection{}

In the case of even-dimensional fibres a proof of this theorem by methods of local index theory was given by Bismut. For a detailed presentation we refer to the book \cite{berlinegetzlervergne92}.
The odd-dimensional case can be reduced to the even-dimensional case in a standard manner as follows.

\subsubsection{}\label{sq21}
We have a split exact sequence
$$0\rightarrow K^0(B)\stackrel{\pr_B^*}{\rightarrow} K^0(S^1\times B)\rightarrow
K^1(B)\rightarrow 0\ ,$$ where the split is given by
$i^*:K^0(S^1\times B)\rightarrow K^0(B)$, and where $i:B\rightarrow
S^1\times B$ is induced by the inclusion of a point $*\rightarrow
S^1$. \index{suspension isomorphism}The map $K^0(S^1\times B)\rightarrow
K^1(B)$ is induced by the inverse of the suspension isomorphism
$K^1(B)\cong \tilde K^0(S^1\wedge B)\cong \ker(i^*)$.

\subsubsection{} 
Assume that $B$ is connected.  
If $F:B\rightarrow U(n)$ represents $[F]\in K^1(B)$,
then we define the bundle $V_F$ over $S^1\times B$ as in \ref{eufiewfewfew}. Then $[V_F]-[S^1\times B\times
\C^n]\in \ker(i^*)\subset  K^0(S^1\times B)$ is  the element corresponding to $[F]$. 
We define the odd Chern\index{Chern! odd character}  character $\ch:K^1(B)\to H^{odd}(B,\Q)$  so that it becomes compatible with the suspension.
Therefore we set
$$\ch([F]):=\int_{S^1\times B/B} \ch([V_F])\ ,$$
where $\int_{S^1\times B/B}:H^{ev}(S^1\times B,\Q)\rightarrow H^{odd}(B,\Q)$ is the integration over the fibre in rational cohomology.

\subsubsection{}\label{switch231}

In order to switch from the even to the odd dimensional case we consider the following geometric family $\cS_{geom}$ over $S^1$. The underyling fibre bundle is $\pr_2:S^1\times S^1\to S^1$. As geometric structures we  choose the trivial connection and the fibrewise metric given by the standard metric on $S^1$ of volume one.  Furthermore, we choose the fibrewise orientation and spin structure induced from the standard orientation and non-bounding  spin structure of $S^1$. 

We furthermore choose a complex line bundle $L\to S^1\times S^1$ with a metric $h^L$ and connection $\nabla^L$  such that
$\int_{S^1\times S^1/S^1} c_1(\nabla^L)$ is the normalized volume form $\vol_{S^1}$ of $S^1$.
On the level of topology $c_1(L)\in H^2(S^1\times S^1,\Z)\cong \Z$ is the positively oriented generator. As the Clifford bundle of the geometric family $\cS_{geom}$ we take the fibrewise spinor bundle twisted by $\cL:=(L,h^L,\nabla^L)$.

The index $\ind(\cS_{geom})\in K^1(S^1)$ is a generator, and we can fix the sign conventions for odd index theory by declaring $\ind(\cS_{geom})=1$ under the supension isomorphism
$K^1(S^1)\cong \tilde K^0(S^0)$. Furthermore we have
$\Omega(\cS_{geom})=\vol_{S^1}$.

\subsubsection{}\label{reduodd}

Let now $\cE_{geom}$ be a reduced geometric family
with closed odd-dimensional fibers. Then we consider the even-dimensional geometric family
$$\cF_{geom}:=\pr_1^*\cS_{geom}\times_B \pr_2^*\cE_{geom}$$
over $S^1\times B$.
By construction we have
$$\ind(\cF_{geom})=\pr_1^*\ind(\cS_{geom})\cup \pr_2^*\ind(\cE_{geom})\in \ker(i^*)\subset K^0(S^1\times B)\ .$$

In fact, the construction  which associates $\cF_{geom}$ to $\cE_{geom}$ is another (different from \ref{eufiewfewfew}) explicit construction  of the suspension of the  class $\ind(\cE_{geom})\in K^1(B)$, in this case realized by a family of Dirac operators.

One  checks that $\Omega(\cF_{geom})=\pr_1^*\vol_{S^1}\wedge \pr_2^*\Omega(\cE_{geom})$ and hence  $\int_{S^1\times B/B} \Omega(\cF_{geom})=\Omega(\cE_{geom})$. Hence the even-dimensional case of the local index theorem implies the odd-dimensional case.

An alternative proof will be given in Subsection \ref{rrq}, where the general case of boundary tamed family is considered.

\subsubsection{}

If the fibers of $\cE_{geom}$ 
are not closed but are manifolds with corners, then the form
$\Omega(\cE_{geom})$ is not closed in general. 
\begin{lem}\label{omeg}
We have $d \Omega(\cE_{geom})+\Omega(\partial \cE_{geom})=0$.
\end{lem}
\proof
The assertion is a consequence of Stoke's formula. Let $m:=\dim(E)-\dim(B)$.
\begin{eqnarray*}
d\Omega(\cE_{geom})&=&(-1)^{m}\int_{E/B}
d ( \hA(\nabla^{T^v\pi})\ch(\nabla^{W}))- \int_{\partial E/B
  }(\hA(\nabla^{T^v\pi})\ch(\nabla^{W}))_{|\partial E}\\
 &=&-\int_{\partial E/B
  }\hA(\nabla^{T^v\pi_{|\partial E}})\ch(\nabla^{W_{|\partial E}})\\
&=&-
 \Omega(\partial \cE_{geom})\ .
\end{eqnarray*}
Here we have employed  the equation $\hA(\nabla^{T^v\pi})_{|\partial E}=\hA(\nabla^{T^v\pi_{|\partial E}})$ which follows from the product structure.
\hB
\subsubsection{}

Finally note the following simple consequences of the definition:
\begin{lem}
\begin{enumerate}
\item
If $U\subseteq B$ is open, then
$\Omega(\cE_{geom|U})=\Omega(\cE_{geom})_{|U}$.
\item
We have 
$\Omega(\cE_{geom}\sqcup_B \cE^\prime_{geom})=\Omega(\cE_{geom})+\Omega(\cE_{geom}^\prime)$.
\item
We have 
$\Omega(\cE^{op}_{geom})=-\Omega(\cE_{geom})$.
\end{enumerate}
\end{lem}

 \subsection{Eta forms}\label{ttrr4411}

\subsubsection{}

We consider a reduced geometric family $\cE_{geom}$  over some base $B$.
Let $\cE_t$ be a taming. 
In the present subsection we define the eta form
$\eta(\cE_t)\in\cA_B(B)$.
\subsubsection{}
The following relations will follow immediately from the definition:
\begin{lem} 
\begin{eqnarray*}
\eta(\cE_t^{op})&=&-\eta(\cE_t)\ ,\\
\eta(\cE_{t|U})&=&\eta(\cE_t)_{|U}\ , U\subseteq B\ , \\
\eta(\cE_t\sqcup_B \cE_t^\prime)&=&\eta(\cE_t)+\eta(\cE_t^\prime)\ .
\end{eqnarray*}
\end{lem}
\subsubsection{}
For $k\in\nat\cup\{0\}$ let $\eta^k(\cE_t)\in \cA_B^k(B)$ denote the degree $k$-component.
If $\cE_{bt}$ is a boundary taming of the underlying geometric family $\cE_{geom}$,
then we denote by $\ind_0(\cE_{bt})$ the locally constant $\Z$-valued function
$B\ni b\mapsto \ind(D(\cM_{bt,b}))\in\Z$.
The main result of the present subsection are the following relations:
\begin{theorem}\label{etaprop}
\begin{enumerate}
\item $\Omega^k(\cE_{geom})= 
  d\eta^{k-1}(\cE_t)-\eta^k(\partial \cE_{t}) $, $\quad k\ge 1$.
\item $\Omega^0(\cE_{geom})=-\eta^0(\partial \cE_{bt})+\ind_0(\cE_{bt})$.
\end{enumerate}
\end{theorem}

\subsubsection{}
We now develop the details as a generalization of the  constructions of \cite{bunkekoch98}, Sec. 3.
Let $\Gamma(\cE_{geom})$ denote the bundle of Hilbert spaces with fiber
$L^2(\bar E_b,\bar V_{|\bar E_b})$ over $b\in B$.
We fix once and for all a function
$\chi\in C^\infty(\R)$ such that $\chi(t)=0$ for $t\le 1$ and $\chi(t)=1$ for
$t\ge 2$.
Then we define the rescaled super-connections $A_t(\cE_t)$ on $\Gamma(\cE_{geom})$ 
as follows.

\index{$A_t(\cE_t)$}If the dimension of the fibers of $\cE$ is even, then
we set
$$A_t(\cE_t):=t \left[(1-\chi(t))D(\cE_{geom})+\chi(t) D(\cE_t)\right]+\nabla^{\Gamma(\cE_{geom})} +\frac{1}{4t} c(T)  \ .$$
Here $D(\cE_{geom})$ (resp.  $D(\cE_t)$) denote the family of Dirac operators associated to the family of geometric
(resp. tamed) manifolds $(\cE_{geom,b})_{b\in B}$ and $(\cE_{t,b})_{b\in B}$. 
The connection $\nabla^{\Gamma(\cE_{geom})}$ is defined in
\cite{berlinegetzlervergne92}, Prop. 9.13. The last term is Clifford
multiplication by the curvature tensor \ref{curintro}.

If the dimension of the fibers of $\cE_{geom}$ is odd, then we set
$$A_t(\cE_t):=t \sigma \left[(1-\chi(t))D(\cE_{geom})+\chi(t) D(\cE_t)\right]+\nabla^{\Gamma( \cE_{geom})} +\frac{1}{4t} \sigma c(T)
\ .$$
where $\sigma$ is the  generator of the Clifford algebra $\Cl^1$ satisfying
$\sigma^2=1$. The grading of the Clifford algebra $\Cl^1$ induces a grading of $\Gamma(\cE_{geom})\otimes \Cl^1$, and $A_t(\cE_t)$ is an odd operator in the odd-dimensional case, too.  
\index{Bismut super-connection}For $t\le 1$  the super-connection $A_t(\cE_t)$ is the usual rescaled Bismut super-connection.
For large $t$ it  differs from the Bismut super-connection by the terms coming from the taming.

\subsubsection{}
By $\rho_r$ we denote the characteristic function of $ {}_r\bar E$. It acts as multiplication operator
on $\Gamma(\cE_{geom})$. Furthermore let $\lambda(\bar E/B)$ denote the bundle of fibre-wise densities. 
Let $\Lambda(\bar E/B)$ denote the  vector bundle over $B$ with fiber 
$\Lambda(\bar E/B)_b=C^\infty(\bar E_b,\lambda(\bar E/B)_{|\bar E_b})$ over $b\in B$.
Assume that
$H\in \cA_B(B, \End(\Gamma(\cE_{geom})))$ has coefficients in the smoothing operators on 
$\Gamma(\cE_{geom})$. In the even-dimensional case we denote by 
 $\tr_s H\in \cA_B(B,\Lambda(\bar E/B))$ the 
local super trace of the integral kernel.
As a fibre-wise density the local super trace  can be integrated along the fibre.

By abuse of notation (the operator on the left-hand side may not be of trace class)
we define
$$\Tr_s \rho_r H:=\int_{{}_r\bar E/B} \tr_s H\in \cA_B(B)\ .$$

\index{$\Tr_s$}\index{$\Tr_s^\prime$} Furthermore, we set 
$$\Tr_s^\prime(H):=\lim_{r\to\infty} \Tr_s \rho_r H $$
provided that this limit exists.

In the odd-dimensional case if $H\in \cA_B(B, \End(\Gamma( \cE_{geom} ))\otimes \Cl^1)$ is of the form
$H=H_1+\sigma H_2$, and $H_i\in \cA_B(B, \End(\Gamma( \cE_{geom})))$
have coefficients in the smoothing operators, then we define
$\tr_s H:=\tr H_2$, $\Tr_s \rho_r H:=\int_{{}_r\bar E/B}\tr_s H$,  and
$$\Tr_s^\prime(H):=\lim_{r\to\infty} \Tr_s \rho_r H_2\ ,$$
provided that the limit exists.

Note that $\Tr_s\rho_r$ and $\Tr_s^\prime$ are not super traces in the ordinary sense since they do not vanish super commutators, in general. It is actually this trace defect of $\Tr^\prime_s$ which leads to the boundary contributions in the two formulas claimed in Theorem \ref{etaprop}.

Note that $\Tr_s^\prime$ corresponds to the $b$-trace in the $b$-calculus approach
of Melrose (see e.g. \cite{melrose93}). The material of the present subsection  has
a $b$-calculus analog, see \cite{melrosepiazza97} and \cite{mrp2}.

\subsubsection{}

\begin{lem}\label{exx} 
$$\Tr^\prime_s \ee^{-A_t( \cE)^2}\in \cA_B(B)$$
exists. Moreover, derivatives with respect to $t$ and $b\in B$ can be
interchanged with $\Tr^\prime_s$. 
\end{lem}
\proof
The main point is that $\tr_s  \ee^{-A_t( \cE)^2}$ and its derivatives are  rapidly decaying on $\bar E_b$
locally uniformly with respect to $b\in B$ and $t>0$.

Let $j\in I_1(E)$ and $\partial_j E$ be the corresponding  boundary face of $E$. 
For $s\le 0$ we have a half cylinder $U_{j,s}=\overline{\partial_j E}\times (-\infty,s]\subset \bar E$.
The restriction of $A_t(\cE_t)$ to $U_{j,s}$ for sufficiently small $s$ extends to a $\R$-invariant
super-connection $A_t(Z_j)$ on the cylinder  $Z_j:= \overline{\partial_j E}\times\R$.
We consider $U_{j,s}$ as a subset of $Z_j$ as well.

We consider the induced pre-tamed family $\cZ_{j,t}$.
The reflection at a point of $\R$ induces an isomorphism
$\cZ_{j,t}\cong \cZ_{j,t}^{op}$.
Hence $\tr_s \ee^{-A_t(Z_j)^2}=0$. The usual finite propagation speed comparison\footnote{This method was invented by \cite{cgt}. See \cite{bunke95} for a similar application in the context of Dirac operators} gives constants $c,C\in \R$, $c>0$, such that
\begin{eqnarray}
|\tr_s \ee^{-A_t(\cE_t)^2}(x)| &=&
|\tr_s \ee^{-A_t(\cE_t)^2}(x)-\tr_s \ee^{-A_t(Z_j)^2}(x)|\nonumber \\
&\le &C\ee^{-c\frac{s^2}{t^2}}\label{compa}\ ,
\end{eqnarray}
uniformly for all $x\in U_{j,s}$, $s<-1$, $t>0$,  
 and locally uniformly on $B$. This implies the existence of
$\Tr^\prime_s \ee^{-A_t^2(\cE)}$
locally uniformly with respect to the base $B$ and $t>0$.
Using Duhamel's principle in order to express the derivatives
of $\ee^{-A_t(\cE_t)^2}$ with respect to the base variable $b\in B$ or with respect to time $t$, 
and using a similar finite propagation speed comparison estimate
one shows that one can interchange derivatives with respect to $t$ or $b$
and $\Tr^\prime_s$. \hB

\subsubsection{}
\index{$\epsilon$}Let $\epsilon\in\Z/2\Z$ be the parity of the dimension of the fibers of $\cE$.

\begin{lem}\label{derr}
We have the following identity
\begin{eqnarray*}
\partial_t \Tr^\prime_s \ee^{-A_t(\cE_t)^2}
&=&  -\left\{ \begin{array}{cc} \frac{1 }{2\imath\sqrt{\pi}} &\epsilon=1\\
\frac{1}{\sqrt{\pi}} & \epsilon=0\end{array}\right\}
 \Tr^\prime_s \partial_t A_t(\partial \cE_t)\ee^{-A_t(\partial \cE_t)^2}\\
&& - d \Tr^\prime_s \partial_t A_t(\cE_t) \ee^{-A_t(\cE_t)^2}\ .\end{eqnarray*}
\end{lem}
\proof
We first show that
$$\partial_t \Tr^\prime_s \ee^{-A_t(\cE_t)^2}=-\Tr^\prime_s[ A_t(\cE_t),\partial_t A_t(\cE_t) \ee^{-A_t(\cE_t)^2}]\ .$$
We write out the details in the case $\epsilon=0$ and indicate the necessary modifications for the case $\epsilon=1$.

Using Duhamel's formula we get
\begin{eqnarray*}
\partial_t \Tr^\prime_s \ee^{-A_t(\cE_t)^2}&=&
-\Tr^\prime_s \int_0^1 \ee^{-s A_t(\cE_t)^2}\partial_t A_t(\cE_t)^2 \ee^{-(1-s)A_t(\cE_t)^2} ds\\
&=&-\lim_{r\to\infty} \int_0^1  \Tr_s \rho_r \ee^{-s A_t(\cE_t)^2}\partial_t 
A_t(\cE_t)^2 \ee^{-(1-s)A_t(\cE_t)^2} ds\\
&=&-\lim_{r\to\infty} \int_0^1 \Tr_s[\partial_t A_t(\cE_t),A_t(\cE_t)] \ee^{-s A_t(\cE_t)^2}\rho_r\ee^{-(1-s)A_t(\cE_t)^2} ds\\
&=&-\lim_{r\to\infty} \lim_{v\to\infty} \int_0^1 \Tr_s \rho_v 
[\partial_t A_t(\cE_t),A_t(\cE_t)] \ee^{-s A_t(\cE_t)^2}\rho_r\ee^{-(1-s)A_t(\cE_t)^2}ds\\
&=&-\lim_{v\to\infty} \lim_{r\to\infty}  \int_0^1 \Tr_s \rho_v [\partial_t A_t(\cE_t),A_t(\cE_t)] \ee^{-s A_t(\cE_t)^2}\rho_r\ee^{-(1-s)A_t(\cE_t)^2}ds\\
&=&-\lim_{v\to\infty}  \int_0^1 \Tr_s \rho_v [\partial_t A_t(\cE_t),A_t(\cE_t)] \ee^{- A_t(\cE_t)^2} ds \\
&=&-\Tr^\prime_s [A_t(\cE_t),\partial_t A_t(\cE_t)\ee^{-A_t(\cE_t)^2}]\ .\end{eqnarray*}
In order to justify that the limits $\lim_{r\to\infty}$ and $\lim_{v\to\infty}$
can be interchanged one can again use a comparison with model cylinders as in the proof of Lemma \ref{exx}.

We further compute
\begin{eqnarray*}
-\Tr^\prime_s [A_t(\cE_t),\partial_t A_t(\cE_t)\ee^{-A_t(\cE_t)^2}]&=&
-\Tr^\prime_s [\nabla^{\Gamma(\cE_{geom})} ,\partial_t A_t(\cE_t)\ee^{-A_t(\cE_t)^2}]\\&&-\Tr^\prime_s[t D(\cE_{geom})
,\partial_t A_t(\cE_t)\ee^{-A_t(\cE_t)^2}]\\
&=&
-d\Tr^\prime_s  \partial_t A_t(\cE_t)\ee^{-A_t(\cE_t)^2}  -\Tr^\prime_s[t D(\cE_{geom}),
\partial_t A_t(\cE_t)\ee^{-A_t(\cE_t)^2}]\ ,
\end{eqnarray*}
(where we must  replace $D(\cE_{geom})$ by $\sigma D(\cE_{geom})$ in the case $\epsilon=1$)
by checking that
$$\Tr_s^\prime \left[ 
\left(
A_t(\cE_t)-tD(\cE_{geom})-\nabla^{\Gamma(\cE_{geom})}
\right)
 , \partial_t A_t(\cE_t)\ee^{-A_t(\cE_t)^2}\right]=0\ .$$
By integration by parts we get
\begin{eqnarray*}
\lefteqn{-\Tr_s^\prime \left[t D(\cE_{geom}), \partial_t A_t(\cE_t)\ee^{-A_t(\cE_t)^2}\right]} \hspace{0.5cm}\\
&=&-\Tr_s^\prime tD(\cE_{geom}) \partial_t A_t(\cE_t)\ee^{-A_t(\cE_t)^2}
-\Tr_s^\prime  \partial_t A_t(\cE_t) \ee^{-A_t(\cE_t)^2}tD( \cE_{geom})\\
&=&-\lim_{r\to\infty} \int_{{}_r\bar E/B} \tr_s 
 tD( \cE_{geom}) \partial_t A_t(\cE_t)\ee^{-A_t(\cE_t)^2}
-\lim_{r\to\infty}     \int_{{}_r\bar E /B}
\tr_s  \partial_t A_t(\cE_t)\ee^{-A_t(\cE_t)^2} tD(\cE_{geom}) 
\\
&=&t \lim_{r\to\infty}  
 \int_{\partial ({}_r\bar E) /B} \tr_s c(\cN) \partial_t A_t(\cE_t)\ee^{-A_t(\cE_t)^2}\ ,
\end{eqnarray*}
\index{$\cN$} where $\cN$ denotes the inner unit normal field of $\partial ({}_r\bar E) $
(Here again we must  replace $D(\cE_{geom})$ and $c(\cN)$  by $\sigma D(\cE_{geom})$
and $\sigma c(\cN)$ in the case $\epsilon=1$).
Using the comparison with our model cylinder $Z=\cup_{j\in I_1(E)}Z_j$  
and $\R$-invariance of $A_t(Z)$  
we obtain
$$\lim_{r\to\infty}  
 \int_{\partial ({}_r\bar E) /B} \tr_s c(\cN) \partial_t A_t(\cE_t)\ee^{-A_t(\cE_t)^2}
=  
\int_{ \overline{ \partial E}\times\{0\} /B} \tr_s c(\cN) \partial_t A_t(Z)
\ee^{-A_t(Z)^2}\ ,
$$
where here $\cN$ denotes the unit vector field generating
the $\R$-action by translations on the cylinder $Z:=\overline{\partial E}\times \R$, and $\overline{\partial E}$ denotes the disjoint union
of the boundary faces  (and where we  replace $D(\cE_{geom})$ and $c(\cN)$  by $\sigma D(\cE_{geom})$
and $\sigma c(\cN)$ in the case $\epsilon=1$).

We identify the Dirac bundle over $Z$ with
$\partial \cV*\R$. Then we can write
$$A_t(Z)=tc(\cN)\cN +L_{\overline{\partial E}}^Z (A(\partial \cE_t))$$ in the case $\epsilon=0$, and
$$A_t(Z)=t\sigma c(\cN)\cN+ L_{\overline{\partial E}}^Z (A(\partial \cE_t)^{odd}) + \sigma  L_{\overline{\partial E}}^Z (A(\partial \cE_t)^{even}) $$
in the case $\epsilon=1$, where the superscripts $even,odd$ indicate the form degree.
By an easy computation using this explicit form
of the super-connection $A_t(Z)$ we get in the case $\epsilon=0$ 
$$\ee^{-A_t(Z)^2}(r,s)= \frac{\ee^{-(r-s)^2/4t}}{t\sqrt{4\pi}} 
L_{\overline{\partial E}}^{Z} (\ee^{-A_t(\partial \cE_t)^2})\ ,$$
where $r,s$ are coordinates in $\R$, and we consider
$\ee^{-A_t(Z)^2}(r,s)$ as an element of
$\cA_B(B,\End(\Gamma(\partial \cE_{geom})))$.

If $\epsilon=1$, then 
$$\ee^{-A_t(Z)^2}(r,s)= \frac{\ee^{-(r-s)^2/4t^2}}{t\sqrt{4\pi}} (
L_{\overline{\partial E}}^{Z} (\ee^{-A_t(\partial \cE_t)^2})^{even} +
\sigma L_{\overline{\partial E}}^{Z} (\ee^{-A_t(\partial \cE_t)^2})^{odd} )
\ .$$

If $\epsilon=0$, then we obtain that 
\begin{eqnarray*}
\tr_s c(\cN) \partial_t A_t(Z)
\ee^{-A_t(Z)^2}&=&\frac{1}{t\sqrt{4\pi}}\tr_s c(\cN)  L_{\overline{\partial E}}^{Z} (\partial_t A_t(\partial \cE_t) \ee^{-A_t(\partial \cE_t)^2})\\
&=&-\frac{2}{t\sqrt{4\pi}}\tr_s \partial_t A_t(\partial \cE_t) \ee^{-A_t(\partial \cE_t)^2}\ .
\end{eqnarray*}
If $\epsilon=1$, then we obtain that
\begin{eqnarray*}
 \tr_s \sigma c(\cN) \partial_t A_t(Z)
\ee^{-A_t(Z)^2}&=&
\frac{1}{t\sqrt{4\pi}}\tr_s \sigma c(\cN)  (
L_{\overline{\partial E}}^{Z} (\partial_t A_t(\partial \cE_t)\ee^{-A_t(\partial\cE_t)^2})^{odd} \\&&+
\sigma L_{\overline{\partial E}}^{Z} (\partial_t A_t(\partial \cE_t)\ee^{-A_t(\partial \cE_t)^2})^{even} )
\\
&=&\frac{1}{t\sqrt{4\pi}}\tr  c(\cN) L_{\overline{\partial E}}^{Z} (\partial_t A_t(\partial \cE_t)\ee^{-A_t(\partial \cE_t)^2})^{odd} \\
&=&\frac{\imath}{ t\sqrt{4\pi}} \tr_s \partial_t A_t(\partial \cE_t)  \ee^{-A_t(\partial \cE_t)^2}\ .
\end{eqnarray*}
The last equalities can be seen using the explicit
constructions of \ref{ccoa}.  The proof also shows that there is no extra contribution from faces of codimension $\ge 2$. 
\hB

\subsubsection{}

\begin{ddd}\label{eet1}
\index{$eta$ form}\index{$\eta^k(\cE_{t})$}We  define the eta forms by
\begin{eqnarray*}
\eta^{2k-1}(\cE_t)&:=&(2\pi\imath)^{-k} \int_0^\infty \Tr^\prime_s
\partial_t A_t(\cE_t) \ee^{-A_t(\cE_t)^2} dt \ , \quad \epsilon=0\\ 
\eta^{2k}(\cE_t)&:=&-(2\pi\imath)^{-k}  \pi^{-1/2}  \int_0^\infty \Tr^\prime_s 
\partial_t A_t(\cE_t) \ee^{-A_t(\cE_t)^2} dt\ ,\quad \epsilon=1\ .
\end{eqnarray*}
\end{ddd}
Note that $\eta^0$ is minus half of the eta invariant of Atiyah-Patodi-Singer \cite{atiyahpatodisinger75}. We choose this different convention in order to have a simple
scheme of signs and prefactors.
The higher eta forms where introduced by Bismut and Cheeger \cite{bismutcheeger89} 
in connection with the study of the adiabatic limit of eta invariants. They also appear in the
index theorem for families of APS-boundary value problems due to Bismut-Cheeger 
\cite{bismutcheeger90}, \cite{bismutcheeger901} and its extension to the ($b$-calculus version of the) boundary tamed case by Melrose and Piazza \cite{melrosepiazza97}.

\subsubsection{}

The standard  small time asymptotic expansion of the local super traces
of the heat kernel of the Bismut super-connections
and the estimate (\ref{compa}) show that these integrals converge at $t=0$.
In order to see that we have convergence at $t=\infty$ we
use the fact that the tamed Dirac operator is invertible.
We therefore have an estimate of the local super trace
by $C\ee^{-ct^2}$ which is uniform on $\bar E$ locally over $B$.
Combined with (\ref{compa}) (where we set $s:=t$) 
we obtain an estimate of the integrands
of the $\eta$-forms by $C\ee^{-ct}$ which is uniform for large $t$ and locally on $B$.

\subsubsection{}
We now finish the proof of Theorem \ref{etaprop}.
The first assertion follows from the local index theorem:
\begin{eqnarray*}
\frac{-1}{\sqrt{\pi} (2\pi\imath)^k}\lim_{t\to 0}[ \Tr_s^\prime \ee^{-A_t(\cE_t)^2}]_{2k+1}=
  \Omega^{2k+1}(\cE_{geom})&&\epsilon=1\\
\frac{1}{(2\pi\imath)^k}\lim_{t\to 0}[ \Tr_s^\prime \ee^{-A_t(\cE_t)^2}]_{2k}=
  \Omega^{2k}(\cE_{geom})&&\epsilon=0\ ,
\end{eqnarray*}
Lemma \ref{derr}, the definition of the $\eta$-forms,
and the estimate
$ |\Tr_s^\prime \ee^{-A_t(\cE_t)^2}|\le C\ee^{-ct^2}$
for large times. 
In order to show the second assertion of Theorem  \ref{etaprop}
we must modify the argument above. Since we now only have a boundary taming
the Dirac operator $D(\cE_{bt})$ is Fredholm and $0$ may be in the spectrum.
We have
$$\int_0^\infty \partial_t [\Tr_s^\prime \ee^{-A_t(\cE_t)^2}]_0 dt=
-\Omega^0(\cE_{geom})+\ind_0(\cE_{bt})\ .$$
Now the second assertion of Theorem \ref{etaprop} follows again
from Lemma \ref{derr}. \hB

\subsubsection{}
Assume that $B$ is a point. Then a geometric family is just a
geometric manifold.
Assume that $\cM_{bt}$ is odd-dimensional and boundary tamed.
Let $\cM_t$ and $\cM^\prime_t$ be two extensions of the boundary
taming to a taming. 
\begin{lem}\label{diieta}
We have 
$$\eta^0(\cM^\prime_t)-\eta^0(\cM_t)=\sum_{j\in
  I_0(M)}\Sf(\cM_t^\prime,\cM_t)(j)\ .$$
\end{lem}
\proof
We consider a path of pre-tamings which extend the boundary taming
and connects the two given tamings. Then we study the jumps of the
spectral flow and the eta invariants as eigenvalues cross zero.
If the eigenvalue crosses from the positive side to the negative then the spectral flow increases by one. The eta invariant also increases by one. In order to see this  one isolates the contribution of the small eigenvalue  to the eta invariant.
\hB

\subsection{An index theorem for boundary tamed families}\label{rrq}

\subsubsection{}\label{ufiwfwef}
In this subsection we show that Theorem \ref{etaprop} implies an
index theorem for a boundary tamed family.
We again assume that $\cE_{geom}$ is a reduced geometric family over a
compact manifold. For simplicity, we also assume that the fibre dimension is greater than zero.
Let $\cE_{bt}$ be a boundary taming. Assume that the dimension of the fibers of $\cE$ has parity $\epsilon\in\Z/2\Z$.
The family of Fredholm operators $D(\cE_{bt})$ gives rise to an element
$\ind(\cE_{bt})\in K^\epsilon(B)$. 
\begin{theorem}\label{famind}
The form
$\Omega(\cE_{geom}) + \eta(\partial \cE_{bt})$ is closed,
and its de Rham cohomology class $[\Omega(\cE_{geom})+ \eta(\partial \cE_{bt})]\in H_{dR}(B)$ represents the image of   
$\ch(\ind(\cE_{bt}))$ under the de Rham map
$\dR:H(B,\Q)\rightarrow H_{dR}(B)$. 
\end{theorem}
\proof
We first consider the case that $\epsilon=0$. 
Let $V=V^+\oplus V^-$, be a complex $\Z/2\Z$-graded vector bundle over $B$ such that
$[V]=-\ind(\cE_{bt})$ in $K^0(B)$. We choose a hermitian metric and 
a metric connection preserving the grading and thus obtain hermitian bundle with connection $\bV$. Since the fibers of the family $\cE(\bV)_{geom}$
(see Subsection \ref{exex}) are closed
the boundary taming of $\cE_{bt}$ induces a boundary taming 
$(\cE_{geom}\cup_B\cE(\bV)_{geom})_{red,bt}$. Moreover, we have by
construction $\ind((\cE_{geom}\cup_B\cE(\bV)_{geom})_{red,bt})=0$ so that
by Lemma \ref{fambtt} the 
boundary taming  admits an extension to a taming
$(\cE_{geom}\cup_B\cE(\bV)_{geom})_{red,t}$.

We compute using Theorem \ref{etaprop}
$$\Omega(\cE_{geom}\cup_B\cE(\bV)_{geom}) +\eta(\partial \cE_{bt})=d\eta((\cE_{geom}\cup_B\cE(\bV)_{geom})_{red,t})\ .$$
Since $\Omega(\cE(\bV)_{geom})$ is closed, 
and $\Omega(\cE_{geom}\cup_B\cE(\bV)_{geom})=\Omega(\cE_{geom}) +\Omega(\cE(\bV)_{geom})$  it follows that $\Omega(\cE_{geom}) +\eta(\partial \cE_{bt})$ is closed, too.
Moreover, in de Rham cohomology we have
\begin{eqnarray*}
[\Omega(\cE_{geom}) + \eta(\partial \cE_{bt})]&=&-[\Omega(\cE(\bV)_{geom})]\\
&=&-\dR(\ch([V]))\\& =&\dR(\ch(\ind(\cE_{bt})))\ .\end{eqnarray*}

Let now $\epsilon=1$.
Let $F:B\rightarrow U(n)$ be a smooth map such that
$[F]\in[B,U(\infty)]\cong K^1(B)$ represents 
$-\ind(\cE_{bt})\in K^1(B)$. 
Let $\cE(F,*)_{geom}$ be a geometric family associated with $F$ as introduced in Subsection \ref{exex}. Then we have 
$\ind(\cE(F,*)_{geom})= -\ind(\cE_{bt})$. The form
$\Omega(\cE(F,*)_{geom})$ is closed and $[\Omega(\cE(F,*)_{geom})]=-\dR(\ch(\ind(\cE_{bt})))$.
Now we argue as in the even-dimensional case 
\hB

\subsubsection{}
A similar result was previously shown by Melrose and Piazza \cite{melrosepiazza97} in the even-dimensional  and \cite{mrp2} in the odd-dimensional case.

\subsubsection{}
Let $\cE_{geom}$ be a geometric family over $B$ and $\cE_t$ and $\cE_t^\prime$
be two tamings. Then $\eta(\cE_t)-\eta(\cE_t^\prime)$ is a closed form.
\begin{kor}\label{etajump}
There exists $\psi\in K(B)$ such that
$[\eta(\cE_t)-\eta(\cE_t^\prime)]=\dR (\ch(\psi))$.
\end{kor}
\proof
Let $I$ be the unit interval with two boundary faces $0$ and $1$.
We consider the geometric family
$\cE_{geom}*I$. We identify
$\partial_0 (\cE_{geom}*I)\cong \cE_{geom}$ and
$\partial_1(\cE_{geom}*I)\cong \cE_{geom}^{op}$.
We let $\cE_t$ and $(\cE_t^\prime)^{op}$ induce tamings
of $\partial_0 (\cE_{geom}*I)$ and $\partial_1(\cE_{geom}*I)$.
The result is a boundary taming
$(\cE_{geom}*I)_{bt}$. Since $\Omega(\cE_{geom}*I)=0$ we get by
\ref{famind} that
$$[\eta(\partial
(\cE_{geom}*I)_{bt})]=[\eta(\cE_t)-\eta(\cE_t^\prime)]=\dR(\ch(\ind((\cE_{geom}*I)_{bt})))\
.$$
\hB

\newpage

\part{Analytic obstruction theory}

\section{The filtration of $K$-theory and related obstructions}\label{hhgg1}

\subsection{Fredholm operators classify $K$-theory}\label{hfefjefewkjfekjwfewjk78}

\subsubsection{}

\index{$\Fred$}\index{space!of Fredholm operators}Let $\Fred$ be the space of Fredholm operators on a separable Hilbert
space $\tH$ with the topology induced by the operator norm. It is well-known (see \cite{atiyah???})  that
it has the homotopy type of the classifying space of the 
$K$-theory functor $K^0$. \index{space!classifying $K$-theory}
For two spaces $X$ and $Y$ let $[X,Y]$\index{$[X,Y]$} denote the set of homotopy classes of
continuous maps from $X$ to $Y$. 
 If $X$ is compact, then the natural transformation $\Psi^0_X:[X,\Fred]\rightarrow K^0(X)$ is induced by the index bundle construction (see \ref{jkl12}) \index{$\Psi^0_X$}
(if we consider the definition of $K^0(X)$ in terms of $\Z_2$-graded vector bundles). 

\subsubsection{}
Let $\Fred^*$ denote the \index{$\Fred^*$}space of selfadjoint Fredholm operators on $\tH$.
It has three components. The component $\Fred^*_0$ is distinguished by the property that its elements have infinite positive as well as infinite negative spectrum.
This component has the homotopy type of the classifying space of the functor
$K^1$  (see \cite{atiyah???}). 
Recall from \ref{sq21} that we have a split exact sequence
\begin{equation}\label{dqwhdwqhjdqdjwq667}
0\rightarrow
K^0(X)\stackrel{\stackrel{i^*}{\leftarrow}}{\stackrel{\pr_{S^1}^*}{\rightarrow}}
K^1(S^1\times X)\rightarrow K^1(X)\rightarrow 0
\end{equation}
which identifies $K^1(X)$ with a direct summand of $K^0(S^1\times X)$.
 For compact $X$ the natural \index{$\Psi^1_X$}transformation $\Psi^1_X:[X,\Fred^*_0]\rightarrow K^1(X)$
(i.e. the transformation which associates to $F:X\rightarrow
\Fred^*_0$ a vector bundle over $S^1\times X$)
is more complicated to describe. We give a description in terms of Dirac operators in \ref{sshix}. 

%

\subsubsection{}\label{jkl12}

\index{$\Comp$}\index{space! of compact operators}
Let $\Comp$ denote the space of compact operators on $\tH$. 
Let $F:X\rightarrow \Fred$ represent $\Psi_X^0([F])\in K^0(X)$. Then
$\Psi_X^0([F])=0$ iff there exists a continuous map
$K:X\rightarrow \Comp$ such that $F+K$ is invertible, i.e.
for each $x\in X$ the operator $F(x)+K(x)$ has a bounded inverse.

In order to see this we argue as follows.
First we construct the index bundle using the method of
\cite{mfsss}. \newcommand{\Proj}{{\tt Proj}}
 Using the compactness of $X$ and the Fredholm property of $F$ we choose families of finite-dimensional projections $P,Q:X\to \Proj(\tH)$ such that
$(1-Q)F(1-P): (1-P)\tH \to (1-Q)\tH$ is a family of invertible operators. The index bundle is then represented by the difference of the ranges of $P$ and $Q$, i.e.
$\Psi^0_X([F])=[P\tH]-[Q\tH]\in K^0(X)$.

If there exists a family of compacts $K$ such that $F+K$ is invertible,
then $F+tK$, $t\in [0,1]$ is a homotopy from $F$ to $F+K$ of families of Fredholm operators.
Hence $\Psi^0_X([F])=\Psi^0_X([F+K])$. But $\Psi^0_X([F+K])=0$ since in this case we can take
$P=Q=0$.

We now consider the converse.
If $\Psi^0_X([F])=0$, then $P\tH$ and $Q\tH$ are stably isomorphic. We show that in this case  one can modify $P,Q$ such that $P\tH\cong Q\tH$.
 Assume that $P\tH$ and $Q\tH$ become isomorphic after adding a trivial $n$-dimensional vector bundle. By Kuiper's theorem an infinite-dimensional bundle of Hilbert spaces is trivial. Therefore
$(1-P)\tH$  is trivial, and we can find a trivial $n$-dimensional subbundle.
We let $P^\prime$ denote the projection onto this subbundle. Furthermore, we let
$Q^\prime$ denote the projection onto the trivial $n$-dimensional bundle $(1-Q)FP^\prime \tH$.
If we replace $P$ and $Q$ by $P+P^\prime$ and $Q+Q^\prime$, then we have  $P\tH\cong Q\tH$.

Finally  we choose an isomorphism of bundles  $U:P\tH\to Q\tH$.
Then $(1-Q)F(1-P)+QUP$ is a family of invertible operators, and
$K:=(1-Q)F(1-P)+QUP-F$ is compact.

\subsubsection{}

Let $\Comp^*\subset \Comp$ denote the subspace of selfadjoint compact operators.
Let  $F:X\rightarrow \Fred^*_0$ represent $\Psi_X^1([F])\in K^1(X)$.
Then $\Psi_X^1([F])=0$ iff there exists a continuous map
$K:X\rightarrow \Comp^*$ such that $F+K$ is invertible.

\newcommand{\tB}{{\tt B}}\newcommand{\tQ}{{\tt Q}}
For the sake of completeness we will again prove this statement.
Let us first assume that $\Psi_X^1([F])=0$. Then there exists a homotopy of maps $F_t:X\to \Fred_0^*$, $t\in [0,1]$ from  $F_0\equiv F$ to a constant map $F_1$ with value $E$, where $E$ is a base point of $\Fred^*_0$, a selfadjoint involution with infinite-dimensional $\pm 1$-eigenspaces.
Since $X$ is compact there exists $\epsilon>0$ such that
$\sigma_{ess}(|F_t(x)|)\cap [0,\epsilon]=\emptyset$ for all $(t,x)\in [0,1]\times X$.
We choose a positive function $g\in C([0,\infty))$ such that $g(s)=1$ for $s<\frac{\epsilon}{2}$ and
$g(s)=\frac{1}{s}$ for $s\ge \epsilon$. Since $F$ is bounded the operator $g(|F|)$ is invertible.
We define $G_t:=g(|F_t|)F_t$.
Note that $\sigma_{ess}(G_t(x))= \{1,-1\}$, in particular  $G_t(x)^2-1$ is compact for all $x\in X$.
Therefore $G:[0,1]\times X\to \tB$, $\tB:=B(\tH)$, induces a family
of involutions $\bar G:[0,1]\times X\to \tQ$ in the Calkin algebra $\tQ:=\tB/\Comp$.
We apply  \cite[Prop. 4.3.3]{black} to the family of projections
$\frac{1+\bar G_t}{2}\in C(X,\tQ)$ and find a family of unitaries
$\bar U_t\in U(C(X,\tQ))$ such that $\bar U_1=1$ and
$\bar U_t \bar G_t \bar U_t^{-1} =  \bar E$ for all $t\in [0,1]$.
In particular we have 
\begin{equation}\label{wefuwieufcw}
\bar U_0 \bar G_0 \bar U_0^{-1}= \bar E\ .
\end{equation}
We now apply \cite[Corollary 3.4.4]{black} to the surjective homomorphism of $C^*$-algebras  $$C(X,\tB)\to C(X,\tQ)\ .$$ Since $\bar U_0\in U(C(X,\tQ))$ is connected by a path to the identity  there exist a lift $U_0\in U(C(X,\tB))$.
Equation (\ref{wefuwieufcw}) implies  that  $\tilde K:=U_0^{-1} E U_0-G_0$ is a family of selfadjoint compact operators. Since
$U_0(G_0+\tilde K)U_0^{-1}=E$ is invertible, $G_0+\tilde K$ is invertible, too.
It follows that
$g(|F|)^{-1/2}(G_0+K)g(|F|)^{-1/2}=F+K$ with $K:=g(|F|)^{-1/2}\tilde K g(|F|)^{-1/2}$ is invertible, where $K$ is a family of selfadjoint compact operators.

For the converse assume that $F+K$ is invertible for some family of selfadjoint compact operators $K$. In a first step we connect $F$ and $F+K$ by the path
$F+tK$, $t\in [0,1]$ of maps $X\to \Fred_0^*$. In the next step we connect
$F+K$ with the involution $R:=\frac{F+K}{|F+K|}$ in $C(X,\tB)$ by the path 
$\frac{F+K}{|F+K|^s}$, $s\in [0,1]$ in $C(X,\Fred_0^*)$. Since $K^0(C(X,\tB))=0$ (observe that $C(X,\tB)$ is the stable multiplier algebra of $C(X)$ and use \cite[Prop. 12.2.1]{black}) the space of projections
$\Proj(C(X,\tB))$ is connected. Hence there exists  a homotopy of projections from $\frac{1+R}{2}$ to $\frac{1+E}{2}$. This homotopy induces a homotopy from $R$ to $E$ as maps $X\to \Fred_0^*$.
Therefore $F$ is homotopic to the constant map with value $E$.

\subsubsection{}
\index{$K^0(X,Y)$}\index{relative $K$-theory}Let now $(X,Y)$ be a pair of compact spaces.
Then we can represent $K^0(X,Y)$ as the set of homotopy classes
of pairs $(F,K)$, $F:X\rightarrow \Fred$, $K:Y\rightarrow \Comp$, such that
$F_{|Y}+K$ is invertible.

In a similar manner we represent $K^1(X,Y)$ as the set of homotopy classes of
pairs $(F,K)$, $F:X\rightarrow \Fred^*_0$, $K:Y\rightarrow \Comp^*$, such that
$F_{|Y}+K$ is invertible.

\subsection{The filtration of $K$-theory and the Atiyah-Hirzebruch spectral sequence}

\subsubsection{}

We describe a natural decreasing filtration   \index{$K_p^*(X)$}\index{Atiyah-Hirzebruch!filtration}
$$\dots\subseteq	 K_p^*(X)\subseteq K_{p-1}^*(X)\subseteq\dots\subseteq  K_0^*(X)=K^*(X)\ .$$
\begin{ddd}\label{fillies}
Let $p\in\nat\cup\{0\}$ and $\psi\in K^*(X)$.
We say that $\psi\in K^*_p(X)$ if  $f^*\psi=0$ for all $CW$-complexes $Y$ of dimension $<p$
and continuous maps $f:Y\rightarrow X$.
\end{ddd}
 This filtration has been studied by Atiyah-Hirzebruch
\cite{atiyahhirzebruch???} where it gave rise to the celebrated 
Atiyah-Hirzebruch spectral sequence relating $K$-theory with integral cohomology (see \ref{tzqwdwqd}).
\subsubsection{}

We have for all $p\ge 0$
\begin{eqnarray}
K_{2p+1}^0(X)&=&K^0_{2p+2}(X)\label{zueure}\\
K^1_{2p}(X)&=&K^1_{2p+1}(X)\label{zueure1} \ .
\end{eqnarray}

In order to see (\ref{zueure}) let $\psi\in K^0_{2p+1}(X)$. We must show that
 $\psi\in K^0_{2p+2}(X)$.
We consider a continuous map $f:Y\to X$, where $Y$ is a $CW$-complex of dimension $\le 2p+1$.
We must show that $f^*\psi=0$.
Let $Y^{2p}\subseteq Y$ be the $2p$-skeleton.
Then we have a push-out diagram
$$\xymatrix{\vee_{\alpha} S^{2p}\ar[d]\ar[r]&Y^{2p}\ar[d]\\\vee_{\alpha}D^{2p+1}\ar[r]&Y }\ .$$
We consider the following segment of the associated Mayer-Vietoris sequence:
$$\dots \prod_\alpha K^1(S^{2p})\to K^0(Y)\stackrel{\alpha}{\to} K^0(Y^{2p})\oplus \prod_{\alpha} K^0(D^{2p+1})\to \dots\ .$$
Since $\psi\in K^0_{2p+1}(X)$ we have  $\alpha(f^*\psi)=0$. Since in addition
$K^1(S^{2p})\cong 0$ we conclude that $f^*\psi=0$.

The proof of (\ref{zueure1}) is similar.

\subsubsection{}

The filtration is also compatible with the  ring structure on $K^*(X)$, i.e. we have $$K^*_p(X)K_q^*(X)\subseteq K^*_{p+q}(X)\ .$$
For a proof see \cite{atiyahhirzebruch???}.

\subsubsection{}\label{tzqwdwqd}

Assume that $X$ is a $CW$-complex with filtration $\emptyset\subseteq X^0\subseteq X^1\dots\subseteq \dots X$ by skeletons. The transition from $X^{r-1}$ to $X^{r}$ is given by a push-out diagram
$$\xymatrix{\vee_{\alpha\in I_r} S^{r-1}\ar[d]\ar[r]&X^{r-1}\ar[d]\\\vee_{\alpha\in I_r}D^{r}\ar[r]&X^r }\ ,$$
where $I_r$ is the set of $r$-cells of $X$. Following \cite{atiyahhirzebruch???}
 the $E_1$-page of the Atiyah-Hirzebruch spectral sequence is given by
$E_1^{p,q}:=K^{p+q}(X^{p},X^{p-1})$, and the differential
$d_1:E_1^{p,q}\to E_1^{p+1,q}$
is given by the boundary operator associated to the long exact sequence of the triple
$(X^{p+1},X^p,X^{p-1})$
{\small $$ K^{p+q}(X^{p+1},X^{p-1})\to K^{p+q}(X^{p},X^{p-1})\stackrel{d_1}{\to} K^{p+q+1}(X^{p+1},X^{p})\to K^{p+q+1}(X^{p+1},X^{p-1})$$}

By excision we have an isomorphism
$E_1^{p,q}=K^{p+q}(X^p,X^{p-1})\cong \prod_{\alpha\in I_p} K^{p+q}(D^p,S^{p-1})$.
We fix once and for all identifications (coming from Bott periodicity)
$K^{2r}(D^{2p},S^{2p-1})\cong \Z$ and $K^{2r+1}(D^{2p-1},S^{2p-2})\cong \Z$ for all $p\ge 1$, $r\in \Z$.
We adjust the signs such that equation (\ref{udiqdiqw}) holds.
If $X$ is a $CW$-complex, then by $C^*(X)$ we denote the cellular
cochain complex of $X$. Note that $H^*(X,\Z)\cong H^*(C^.(X))$. 
The indentifications above induce isomorphisms of complexes  
\begin{equation}\label{tzqwdwqd1}
E_1^{*,2k}\cong C^*(X)\ ,\quad E_1^{*,2k+1}\cong 0\ .
\end{equation}
In particular
$E_2^{p,2k}\cong H^{p}(X,\Z)$ and $E_2^{p,2k+1}\cong 0$.

\subsubsection{}

\index{$\Gr K(X)$}\index{Atiyah-Hirzebruch!spectral sequence}The Atiyah-Hirzebruch spectral sequence converges to $\Gr K(X)$, the associated graded group of the group $K(X)$ with the  filtration \ref{fillies}.
If $x\in K^r_p(X)$, then we can choose a lift
$\hat x\in K^r(X^{p},X^{p-1})$ of $x_{|X^{p}}\in K^r(X^p)$ under the restriction map
$K^{r}(X^{p},X^{p-1})\to K^r(X^p)$. The class $\hat x\in E_1^{p,r-p}\cong C^p(X)$ is called a $1$-symbol
of $x$.  It descents to the $E_\infty$-term where it represents the class $[x]\in \Gr^p K^r(X)$ corresponding to $x$.
The  $1$-symbol of $x$  is well-defined modulo the images of all higher differentials ending at this place $E_*^{p,r-p}$.  This non-uniqueness corresponds to the non-uniquenes of the choice of the lift above.

\subsubsection{}\label{uefiwefwefwef4444}

\index{$l$-symbol}The class in $E_l^{p,r-p}$ represented by a $1$-symbol of $x$ will be called an $l$-symbol of $x$.
In particular, a $2$-symbol of $x$ is a class  $z\in H^{p}(X,\Z)$.  By \cite{atiyahhirzebruch???},
2.5 (2. Corollary (ii)) and the fact that the Atiyah-Hirzebruch spectral sequence degenerates rationally we see that the image $z_\Q\in H^p(X,\Q)$ of $z$ in rational cohomology satisfies
\begin{equation}\label{udiqdiqw}
\ch_p(x)=z_\Q
\end{equation}\
(the odd case follows from the even by suspension, see also \ref{ifowee5}).

\subsection{Obstruction theory}

\subsubsection{}
\index{obstruction!theory}We describe the obstruction theory related to the non-trivial steps of the filtration of $K$-theory. In particular, we will connect the abstract algebraic topology constructions with the more concrete picture of $K$-theory (described in Subsection \ref{hfefjefewkjfekjwfewjk78}) involving Fredholm operators and invertible perturbations by compact operators.

Let $X$ be a finite $CW$-complex and 
$X^0\subseteq X^1\subseteq X^2\subseteq \dots$ its filtration by skeletons.
Let $\psi\in K^0_{2p}(X)$ be represented by 
$F:X\rightarrow \Fred$. Then $\psi_{|X^{2p-1}}=0$ so that there exists
a map $K:X^{2p-1}\rightarrow \Comp$ such that
$F_{|X^{2p-1}}+K$ is invertible.

Let $\chi:D^{2p}\rightarrow X^{2p}$ be the characteristic map of a $2p$-cell $E^{2p}$ of $X$. Then $(\chi^* F,\chi_{|\partial D^{2p}}^*K)$ represents
the element $c^{2p}(F,K)(E^{2p})\in K^0(D^{2p},\partial D^{2p})\cong \Z$.
It turns out that $c^{2p}(F,K)$ is a closed $2p$-cochain in the cochain complex
$C^*(X)$.  In fact, we can consider the pair $(F_{|X^{2p}},K)$ as a choice of a  lift of $\psi_{|X^{2p}}\in K^0(X^{2p})$ to a class $[F_{|X^{2p}},K]\in K^0(X^{2p},X^{2p-1})$.
The cochain $c^{2p}(F,K)$  is the corresponding  $1$-symbol 
under the identification $E_1^{2p,-2p}\cong C^{2p}(X)$ (see (\ref{tzqwdwqd1})).
Since the $1$-symbol is anihilated by $d_1$ we conclude that $c^{2p}(F,K)$ is closed.

\begin{ddd}By $\bo^{2p}(F,K)\in H^{2p}(X,\Z)$ we denote the cohomology class
\index{$\bo^{k}(F,K)$}represented by $c^{2p}(F,K)$.
\end{ddd}
In other words, $\bo^{2p}(F,K)$ is the $2$-symbol induced by the $1$-symbol  $c^{2p}(F,K)$ of $\psi$.

\subsubsection{}

Let $\psi\in K^1_{2p+1}(X)$ be represented by 
$F:X\rightarrow \Fred^*_0$. Then $\psi_{|X^{2p}}=0$ so that there exists
a map $K:X^{2p}\rightarrow \Comp^*$ such that
$F_{|X^{2p}}+K$ is invertible. 

Let $\chi:D^{2p+1}\rightarrow X^{2p+1}$ be the characteristic map of a $2p+1$-cell $E^{2p+1}$  of $X$. Then $(\chi^* F,\chi_{|\partial D^{2p}}^*K)$ represents
the element $c^{2p+1}(F,K)(E^{2p+1})\in K^1(D^{2p+1},\partial D^{2p+1})\cong \Z$.
It turns out that $c^{2p+1}(F,K)$ is a closed $2p+1$-cochain in the cochain complex
$C^*(X)$. As above $[F_{|X^{2p+1}},K]\in K^1(X^{2p+1},X^{2p})$ defines a $1$-symbol
of $\psi$ in $E_1^{2p+1,-2p}$, and $c^{2p+1}(F,K)$ corresponds to this one-symbol under the identification $C^{2p+1}(X)\cong E_1^{2p+1,-2p}$ (see (\ref{tzqwdwqd1})).

\begin{ddd}By $\bo^{2p+1}(F,K)\in H^{2p+1}(X,\Z)$ we denote the cohomology class
\index{$\bo^{k}(F,K)$} represented by $c^{2p+1}(F,K)$.
\end{ddd}
Again, the class $\bo^{2p+1}(F,K)$ is the $2$-symbol represented by the $1$-symbol $c^{2p+1}(F,K)$.

\subsubsection{}

The cocycle $c^*(F,K)$ and therefore the class $\bo^*(F,K)$ only depend on the homotopy class of the pair
$(F,K)$. In fact, homotopic pairs define the same lifts $[F_{|X^*},K]\in K(X^{*},X^{*-1})$.

Let $F_t$, $t\in [0,1]$ be a homotopy and $K_0$ for $F_0$ be given,
then $K_0$ extends to a family $K_t$, $t\in[0,1]$, accordingly. 

Given $F$, there may be various homotopy classes of maps $K$.
We define the set $\bo^*(F)$ as the set of classes
$\bo^*(F,K)$ for $K$ running over all families as above.
This set only depends on the homotopy class of $F$, i.e. only on  the
element $\psi=\Psi_X^*([F])\in K^*(X)$. 
\begin{ddd}\label{obset} We  will 
write $\bo^*(\psi)$ for this set. \index{$\bo^*(\psi)$}
\end{ddd}
Of course, $\bo^*(\psi)\subseteq H^*(X,\Z)$ corresponds exactly to the set of $2$-symbols of $\psi$.
An explicit description of this set would involve the knowledge of the images of the higher differentials of the Atiyah-Hirzebruch spectral sequence and is therefore difficult.

\subsubsection{}

The following Lemma justifies to call $\bo(\psi)$ an obstruction set.

\begin{lem}
If $\psi\in K^0_{2p}(X)$ and $0\in \bo^{2p}(\psi)$, then
$\psi\in K^0_{2p+2}(X)$. Similarly, if $\psi\in K^1_{2p+1}(X)$ and
$0\in\bo^{2p+1}(\psi)$, then $\psi\in K^1_{2p+3}(X)$.
\end{lem}
\proof
If $0\in \bo^{2p}(\psi)$, then by the interpretation of $\bo^{2p}(\psi)$ as the set of $2$-symbols of $\psi$ we have $[\psi]=0$ in  $\Gr^{2p} K^0(X)$. Therefore, $\psi\in K^0_{2p+1}(X)=K^0_{2p+2}(X)$, where we use (\ref{zueure}).
The odd-dimensional case is similar.

Alternatively one can argue directly as follows.
Let $\psi\in K_k(X)$ be represented by $[F]$ such that $0\in
\bo^k(\psi)$. By this assumption we can find $K$ defined on the $k-1$-skeleton
such that $\bo^k(F,K)=0$.

In a first step one shows that, after stabilization of $(F,K)$,  one can change the cycle $c^k(F,K)$
by arbitrary boundaries by altering $K$ on the $k-1$-skeleton
with fixed  restriction to the $k-2$-skeleton.

Let us discuss the $K^0$-case. Then $k$ is even.  We choose $n\in \nat$ sufficiently large such that
we have an isomorphism $\deg_{k-1}:[(D^{k-1},S^{k-2}),(U(n),1)]\stackrel{\sim}{\to}\Z$.
Let $z\in C^{k-1}(X)$. We can extend the constant map $1: X^{k-2}\to U(n)$ to a map
$\kappa:X^{k-1}\to U(n)$ such that $\deg_{k-1}(\kappa_{|e^{k-1}},)=z(e^{ k-1})$ for each cell $e^{k-1}$ of $X^{k-1}$. Now observe that
$c^k(F\oplus 0_{\C^n},K\oplus \kappa)=c^{k}(F,K) \pm d z$
(the sign depends on the identifications of relative $K$-groups of pairs $(D^{*},S^{*-1})$ and the $(k-1)$'th homotopy group of $U(n)$ with $\Z$ which we have fixed above).
The $K^1$-case is similar.

We see that if $0\in \bo^k(\psi)$, then we can find a representative $(F,K)$ of $\psi$ such that
$c^k(F,K)=0$. 
This $K$ extends over the $k$-skeleton. \hB

\subsection{Chern classes of the obstructions}\label{choob}

\subsubsection{}\label{fgr1}
\index{Chern! even classes}We explain the relation between the obstruction set $\bo^*(\psi)$ and Chern classes. Chern classes are natural transformations from the $K$-theory functor 
to the integral cohomology functor. In the present paper we write
$c_p:K^{[p]}(\dots)\rightarrow H^p(\dots,\Z)$ in order to simplify the notation.
In the standard notation $c_{2p}$ corresponds to $c_p$, and
$c_{2p+1}$ corresponds to $c^{odd}_p$. In order to define
the odd Chern classes \index{Chern!odd classes}we use the identification $K^1(X)\cong
\ker(i^*)\subset K^0(S^1\times X)$ (see (\ref{dqwhdwqhjdqdjwq667}) for notation). Then we have by definition for odd $p$
$$c_{p}(\psi):=\int_{S^1\times X/X} c_{p+1}(\tilde \psi)\ ,$$
where $\tilde \psi\in K^0(S^1\times X)$ corresponds to
$\psi\in K^1(X)$. \index{$c_p$}

\subsubsection{}

Now assume that $\psi\in K^0_{2p}(X)$ and $z\in\bo^{2p}(\psi)$.
Then we have by  \cite{kervaire59}, Lemma 1.1,  that
\begin{equation}\label{hccwcwec}
(-1)^{p-1}(p-1)! z=c_{2p}(\psi)\ .
\end{equation}
Actually, the cited result is up to sign. The sign is obtained from the compatibilty of (\ref{udiqdiqw}), (\ref{hccwcwec}), and
(\ref{d3d3ruiwedwedw}).

If $\psi\in K^1_{2p+1}(X)$, then $\tilde\psi\in K_{2p+2}^0(S^1\times X)$.
We have $\bo^{2p+1}(\psi)=\int_{S^1\times X/X} \bo^{2p+2}(\tilde \psi)$.
Thus, if $z\in \bo^{2p+1}(\psi)$, then $$ (-1)^pp! z = c_{2p+1}(\psi)\ . $$

\subsubsection{}\label{ifowee5}
 
If $\psi\in K_{2p}^0(X)$, and $z\in \bo^{2p}(\psi)$, then
\begin{equation}\label{uefiweefwef8}
z_\Q=\ch_{2p}(\psi)\ ,
\end{equation} where $z_\Q$ denotes the image of $z$ under
$H^{2p}(X,\Z)\rightarrow H^{2p}(X,\Q)$, and $\ch_{2p}(\psi)$ is the
degree $2p$-component of $\ch(\psi)$.
This follows from the fact that $c_{2k}(\psi)=0$ for all $k<p$
and 
\begin{equation}
\label{d3d3ruiwedwedw}\ch_{2p}=\frac{(-1)^{p-1}}{(p-1)!} c_{2p} + \mbox{\em polynomial in lower Chern
  classes}\ .
\end{equation}
Using the interpretation of $\bo^{2p}(\psi)$ as the set of $2$-symbols of $\psi$, equation (\ref{uefiweefwef8}) has been observed previously in \ref{uefiwefwefwef4444}.

Analogously, if $\psi\in K_{2p+1}^1(X)$ and $z\in \bo^{2p+1}(\psi)$, then
we have $$z_\Q=\ch_{2p+1}(\psi)\ .$$

\subsection{The \v{C}ech cohomology picture}\label{edele}

\subsubsection{}

Let $X$ be a topological manifold. Then it is homotopy equivalent to a
$CW$-complex, but not in a unique manner. 
A simplicial complex (which is a particularly nice $CW$-complex) which is homotopy equivalent to $X$ can be constructed as the geometric realization of the nerve of a good covering of $X$.
It is therefore more appropriate to describe the obstruction theory directly in terms of the open covering, i.e. in the framework of \v{C}ech cohomology.

\subsubsection{}\label{hhjdiuqwiudqw}
\index{$\cU$}
\index{good covering}
\index{$\bN$}
\index{nerve of a covering}
Let $\cU=\{U_l\}_{l\in L}$ be a good  covering of $X$ and $\bN$ be its  nerve.
$\bN$ is a simplicial set. A $p$-simplex $x\in\bN[p]$ is a map
$x:[p]\rightarrow L$ such that 
$U_x:=\cap_{i\in [p]} U_{x(i)}\not=\emptyset$, where $[p]:=\{0,1,\dots, p\}$.
The condition that the covering $\cU$ is good is that $U_x$ is
contractible for all simplices $x\in \bN$.

\subsubsection{}

For each monotone map $\partial:[p-1]\rightarrow [p]$ we have a map
$\partial^*:N[p]\rightarrow N[p-1]$ defined by
$\partial^* x:=x\circ \partial$.

The geometric realization $|\bN|$ of $\bN$ is the simplicial complex
\index{geometric!realization}
$$|\bN|:=\bigsqcup_{p\in\nat\cup\{0\}} \bigsqcup_{x\in \bN[p]} \Delta^p_x /\sim$$
with the equivalence relation generated by 
$u\sim v$ if $u\in \Delta^p_x$ and $v\in \Delta^{p-1}_y$,
$y=\partial^* x$, and $u=\partial_* v$, where
$\partial :[p-1]\rightarrow [p]$ is monotone, and
$\partial_*:\Delta^{p-1}\rightarrow \Delta^p$ is the embedding of the corresponding face. 

\subsubsection{}

\index{$\tilde X$}The space $|\bN|$ is homotopy equivalent to $X$, and an equivalence can be constructed as follows.
We consider the space
$$\tilde X:=\bigsqcup_{p\in\nat\cup\{0\}} \bigsqcup_{x\in \bN[p]} U_x\times \Delta^p_x/\sim\ .$$
Here the relation is generated by
$(a,u)\sim (b,v)$ if
$u\in \Delta^p_x$ and $v\in \Delta^{p-1}_y$,
$y=\partial^* x$, and $u=\partial_* v$,
$a=b$.
There are natural maps
$p_2:\tilde X\rightarrow |\bN|$, $p_2(a,u):=u$, and
$p_1:\tilde X\rightarrow X$, $p_1(a,u):=a$. Both maps have contractible fibers and are homotopy equivalences.

 \subsubsection{}\label{tfc21}
The spaces $|\bN|$ and $\tilde X$ have natural filtrations such that
\begin{eqnarray*}
|\bN|^q&:=&\bigsqcup_{q\ge p\in\nat\cup\{0\}} \bigsqcup_{x\in \bN[p]} \Delta^p_x /\sim\\
\tilde X^q&:=&\bigsqcup_{q\ge p\in\nat\cup\{0\}} \bigsqcup_{x\in \bN[p]}U_x\times  \Delta^p_x/\sim\ .\end{eqnarray*}
The map $p_2$ respects this filtration and
$(p_2)_{|\tilde X^q}:\tilde X^q\rightarrow |\bN|^q$ is a homotopy equivalence
for all $q\ge 0$. We choose a homotopy inverse $r:|\bN|\rightarrow 
\tilde X$ which is compatible with the filtrations. Then $p_1\circ r:|\bN|\rightarrow X$ is a homotopy equivalence. Finally, we define $r^q:=r_{|\:|\bN|^q}$.

\subsubsection{}\label{checkjdiisodic}

Let us fix our conventions concerning the \v{C}ech complex.
\index{\v{C}ech complex}Let $\cS$ be any sheaf of abelian groups over $X$. Then we define the 
\v{C}ech complex of $\cS$ associated to the covering $\cU$ (not necessarily good) by \index{$\check{C}^p(\cU,\cS)$}
$$\check{C}^p(\cU,\cS):=\prod_{x\in\bN[p]}\cS(U_x)\ .$$
The differential $\delta:\check{C}^{p-1}(\cU,\cS)\rightarrow \check{C}^p(\cU,\cS)$ is given by 
$$\delta \prod_{y\in \bN[p-1]}\phi_y:=\prod_{x\in \bN[p]} \sum_{j\in[p]}  (-1)^j (\phi_{\partial_j^* x})_{|U_x}\ ,$$
where $\partial_j:[p-1]\rightarrow [p]$ is the unique monotone map such that
$\image(\partial_j)=[p]\setminus\{j\}$.

\subsubsection{}

If $\cU^\prime$ is a refinement of $\cU$, then we have
a morphism of complexes $\check{C}(\cU,\cS)\rightarrow \check{C}(\cU^\prime,\cS)$.
We define
$$\check{C}(X,\cS):=\lim_{\longrightarrow} \check{C}(\cU,\cS)\ ,$$
where the limit is taken over the directed system of open coverings of $X$.
By $\check{H}(X,\cS)$ we denote \index{$\check{H}(X,\cS)$} the cohomology of $\check{C}(X,\cS)$.

\index{constant sheaf}\index{$\underline{G}_X$}If $G$ is an abelian group, then let $\underline{G}_X$ denote the constant sheaf on $X$ with value $G$.
If $\cU^\prime$ is a good refinement of a good covering $\cU$, then
$\check{C}(\cU,\underline{G}_X)\rightarrow \check{C}(\cU^\prime,\underline{G}_X)$
is a quasi-isomorphism. Since the directed subsystem of good coverings is cofinal in all coverings we have
$\check{H}(X,\underline{G}_X)\cong  H^*(\check{C}(\cU,\underline{G}_X))$
for all good coverings $\cU$.

\subsubsection{}\label{zt1234}

From now on we assume that $X$ is compact.
Consider $F:X\rightarrow \Fred$ and assume that 
$\psi:=\Psi^0_X([F])\in K^0_{2q}(X)$. Then $r_{2q-1}^* (p_1^*\psi)_{|\tilde X^{2q-1}}=0$. Therefore $(p_1^*\psi)_{|\tilde X^{2q-1}}=0$ so that we can find
$K:\tilde X^{2q-1}\rightarrow \Comp$ such that
$(p_1)^*_{|\tilde X^{2q-1}} F+K$ is invertible.
Let now $x\in \bN^{2q}$  and fix some $a\in U_x$.
Then $(p_1^*F_{|\{a\}\times \Delta^{2q}_x},K_{|\{a\}\times \partial \Delta^{2q}_x})$ represents an element $\check c^{2q}(F,K)(x)\in K^0(\Delta^{2q},\partial \Delta^{2q})\cong \Z$ which is independent of the choice of $a$.
It turns out that $\check{c}^{2q}(F,K)$ is a \vC ech cocycle 
in $\check{C}(\cU,\underline{\Z}_X)$. Let $\check{\bo}^{2q}(F,K)\in \check H^{2q}(X,\underline{\Z}_X)$ be its cohomology class. 
We further define the set $\check{\bo}^{2q}(\psi)\subset  \check H^{2q}(X,\underline{\Z}_X)$ of all classes $\check{\bo}^{2q}(F,K)$ for varying $K$. \index{$\check{\bo}^{2q}(F,K)$}

Under the natural
identification 
$H^*(|\bN|,\Z)\cong \check{H}^*(|\bN|,\underline{\Z}_{|\bN|})$
of cellular  and \vC ech cohomology we have
$$r^*\circ p_1^*\check{\bo}^{2q}(\psi)\cong \bo^{2q}(r^*\circ p_1^*\psi)\ .$$
This is our description of the obstruction set in the \vC ech cohomology picture. For $\psi\in K^1_{2q+1}(X)$ there is an analogous construction of
$\check{\bo}^{2q+1}(\psi)\subset \check H^{2q+1}(X,\underline{\Z}_X)$.

\subsection{The obstruction set and tamings}\label{indele}

\subsubsection{}

In this subsection we apply the construction of the obstruction set to the
$K$-theory classes which arise as the index of families of Dirac operators.
We represent the obstruction sets in terms of tamings of the family.
For this purpose it is useful to work with unbounded operators.

\subsubsection{}

Let $\tilde \Fred$ be \index{$\tilde \Fred$} the space of unbounded densely defined operators $D$ on $\tH$ such that $(D^*D+1)^{-1}$ and $(DD^*+1)^{-1}$ are compact. We equip
$\tilde \Fred$ with the smallest topology such that
$\tilde \Fred\ni   D\mapsto D(D^*D+1)^{-1/2}\in B(\tH)$,
and  $\tilde \Fred\ni   D\mapsto D^*(DD^*+1)^{-1/2}\in B(\tH)$
are
continuous w.r.t. the strict topology on $B(\tH)$, and
$\tilde \Fred\ni D\mapsto (D^*D+1)^{-1}\in B(\tH)$,
$\tilde \Fred\ni D\mapsto (DD^*+1)^{-1}\in B(\tH)$
are norm continuous.
The space $\tilde \Fred$ also has the homotopy type of the classifying space
of $K^0$. Similarly, the subspace $\tilde \Fred^*_0\subset \tilde \Fred$
of selfadjoint operators with infinite positive and negative spectrum
classifies $K^1$. In the obstruction theory above we can replace
$\Fred$ and $\Fred^*_0$ by $\tilde\Fred$ and $\tilde\Fred^*_0$.
For details we refer e.g. to \cite{bunkejoachimstolz02}.

\subsubsection{}

Let $B$ be a compact smooth manifold and $\cE_{geom}$ be a reduced geometric family over $B$ with closed fibers.
Invoking Kuiper's theorem about the contractibility of the unitary group of a Hilbert space the \index{$\Gamma(\cE_{geom})$} bundle of Hilbert spaces $\Gamma(\cE_{geom})$ with fiber
$L^2(\bar E_b,\bar V_{|\bar E_b})$ over $b\in B$ 
can be trivialized and identified with the  trivial bundle $B\times \tH$ in unique way up to homotopy. In the case of even-dimensional fibers we have a decomposition $\Gamma(\cE_{geom})=\Gamma(\cE_{geom})^+\oplus \Gamma(\cE_{geom})^-$
given by the $\Z/2\Z$-grading, and we identify both bundles separately with 
with $B\times \tH$.

\subsubsection{}\label{se564}
The family $D(\cE_{geom})^+$ (resp. $D(\cE_{geom})$) gives rise to a family of Fredholm operators
$D:B\rightarrow \tilde \Fred$ (resp. $D:B\rightarrow \tilde \Fred_0^*$). 
The homotopy class $[D]$ is well-defined independent of the choice of trivializations, and it represents $\ind(\cE_{geom})$.

Exactly if $\ind(\cE_{geom})=0$ we can find a family of compact
operators $K$ such that $D+K$ is invertible. By an approximation
argument (see e.g. the proof of Lemma \ref{fambtt}) we can assume that
$K$ is a smooth family of smoothing operators.

We now apply this reasoning in order to define the obstruction set
$\bo(\ind(\cE_{geom}))$. \index{$\bo(\ind(\cE_{geom}))$} \index{obstruction!set}

\subsubsection{}
Assume that $\ind(\cE_{geom})\in K_p^*(B)$.
Let $\cU$ be a good covering with nerve $\bN$, and let
$B\stackrel{p_1}{\leftarrow} \tilde B \stackrel{p_2}{\rightarrow} |N|$ the corresponding diagram of homotopy equivalent spaces.

Since $p_1^*\ind(\cE_{geom})_{|\tilde B^{p-1}}=0$ we can find a family of smoothing operators $K$ over $\tilde B^{p-1}$ such that
$D(p_1^*\cE_{geom})+K_x$ is invertible. We can assume that
$K$ is smooth over each piece $U_x\times \Delta_x^q\subset \tilde B$.

Translated back to $B$ we have the following.
For $q<p$ and $x\in\bN[q]$ we have a smooth family of fiber-wise smoothing operators $K_x$ on $p^*_{1,x}\Gamma(\cE_{geom})$, where $p_{1,x}:U_x\times \Delta^q_x\rightarrow U_x$ is the projection, such that $p^*_{1,x}D(\cE_{geom})+K$ is invertible.
If $y=\partial^*x$, then we have the compatibility 
$(K_y)_{|U_x\times \Delta^{q-1}_y}=(\id_{U_x}\times \partial_*)^* K_x$,
where $\partial_*:\Delta^{q-1}\rightarrow\Delta^q$ is the embedding corresponding to $\partial$.

\subsubsection{}\label{indf3}

\index{$\bK$}In terms of the family $\bK:=(K_x)_{x\in \bN[q],q<p}$ we can define the 
\index{$\check{c}^{q}(\cE_{geom},\bK)$}chain $\check{c}^{q}(\cE_{geom},\bK)\in\check{C}^p(\cU,\underline{\Z}_B)$ as follows.
Let $x\in \bN[p]$. Then we define $K_x$ on $U_x\times \partial \Delta_x^p$  
such that $(1\times \partial_*)^*K_x=(K_y)_{|U_x\times \Delta^{p-1}_y}$
for all monotone maps $\partial:[p-1]\rightarrow [p]$, where $y=\partial^*x$.
Because of the compatibility relations satisfied by $\bK$ we see that $K_x$ is well-defined. Furthermore, 
$p_{1,x}^*D(\cE_{geom})_{|U_x\times \partial \Delta_x^p}+K_x$ is invertible.
Therefore, after choosing some $a\in U_x$, we can define 
\begin{equation}\label{coooa}\check{c}^{p}(\cE_{geom},\bK)(x):=(p_{1,x}^*D(\cE_{geom})_{|\{a\}\times  \Delta_x^p},(K_x)_{|\{a\}\times \partial \Delta^p_x})\in K^*(\Delta^p,\partial\Delta^p)\cong \Z\ .\end{equation}
This chain is closed, independent of $a\in U_x$, and it represents the class
$\bo^p(D,K)\in\check{H}^p(B,\underline{\Z}_B)$.

We conclude :

\begin{prop}\label{zzuduwedw}
Assume that $\ind(\cE_{geom})\in K_p^*(B)$.
The obstruction set
$$\check{\bo}^p(\ind(\cE_{geom}))\subseteq \check{H}^p(B,\underline{\Z}_B)$$
is given by the set of classes represented by chains of the form
$\check{c}^{p}(\cE_{geom},\bK)$ for varying good coverings $\cU$
and choices of families $\bK$ as above.
\end{prop}

\subsubsection{}\label{sshix}

In order to fix signs we describe the natural transformation
$\Psi^1_B:[B,\Fred_0^*]\rightarrow K^1(B)$ in terms of Dirac operators.
Here we assume that $B$ is a compact manifold.

Consider the Kuenneth formula
$K^0(S^1\times B)\cong K^0(S^1)\otimes K^0(B)\oplus K^1(S^1)\otimes K^1(B)$.
\index{$\theta$}The second summand is $\ker(i^*)$ where $i$ is as in \ref{sq21}. Let $\theta\in K^1(S^1)$ be the generator which corresponds to
$1\in \tilde K^0(S^0)\cong \Z$. Then the identification
$K^1(B)\stackrel{\sim}{\to} \ker(	i^*)$ is given by the left cup product with $\theta$.
Note that $\ind(\cS_{geom})=\theta$, where $\cS_{geom}$ was defined in \ref{switch231}.
 
Let $\cE_{geom}$ be a geometric family with odd-dimensional fibers over $B$.
It gives rise to a class $[D]\in [B,\tilde\Fred_0^*]$.
We form the geometric family $\cF_{geom}:=\pr_1^*\cS_{geom}\times_B \pr_2^*\cE_{geom}$ with even-dimensional fibers
over $S^1\times B$. Here $\pr_1:S^1\times B\to S^1$ and $\pr_2:S^1\times B\to B$ denote the projections on the factors. Whatever definition of the index of a family of Dirac operators in the odd-dimensional case one uses, it should be multiplicative so that
$$\ind(\cF_{geom})=\pr_1^*\ind(\cS_{geom})\cup \pr_2^*\ind(\cE_{geom})=\pr_1^*\theta\cup \pr_2^*\ind(\cE_{geom})$$
(compare \ref{reduodd}). 
We can thus define $\ind(\cE_{geom})$ by this formula, and then we set
$$\Psi^1_B([D]):=\ind(\cE_{geom})\ .$$

Note that every element of $K^1(B)$ can be represented as the index of  a family of Dirac operators.
In one definition of $K^1$ we have
$K^1(B)\cong [B,\lim_{\stackrel{N\to\infty}{\to}} U(N)]$.
In particular, every element of $K^1(B)$ can be represented by a map $F:B\to U(N)$ for some sufficiently large $N$. Then we have $\ind(\cE(F,*)_{geom})=[F]$ as explained in 
\ref{eufiewfewfew}. For example, if  $F:S^1\to U(1)\cong S^1$ is the identity,
then we have $\cS_{geom}\cong \cE(F,*)_{geom}$ for an appropriate choice of the geometry on the right-hand side.

\section{Localization over the base}\label{d555}

\subsection{Chains and bordism of chains}

\subsubsection{}\label{langg}
In order to simplify the language let us agree about the following.
Let $\cF_{geom}$ and $\cE_{geom}$ be geometric families.
Let $i\in I_{1}(\cE_{geom})$. By an isomorphism
$\cF_{geom}\cong \partial_i\cE_{geom}$ we understand the
identification of $\cF_{geom}$ with a distinguished model of
$\partial_i\cE_{geom}$. Therefore we have the canonical
embedding $T:\bar F\times (-\infty,0)\rightarrow \bar E$
and the product structure $\Pi:T^*\bar \cV\rightarrow \bar \cW*(-\infty,1)$,
where $\cV$ and $\cW$ denote the families of  Dirac bundles of
$\cE_{geom}$ and $\cF_{geom}$.
We will apply the same convention to pre-tamed families.

If $j\in I_k(\cF_{geom})$, then we can consider $j\in
I_{k+1}(\cE_{geom})$. A pre-taming of the face $j$ of $\cF_{geom}$
induces a pre-taming of the face $j$ of $\cE_{geom}$.

\subsubsection{}

Let $i,j\in I_1(\cE_{geom})$ be adjacent with respect to $k\in
I_2(\cE_{geom})$. Assume that we have geometric families
$\cF_{i,geom}$ and $\cG_{k,geom}$.
Assume further that we have isomorphisms
$\partial_i\cE_{geom}\cong \cF_{i,geom}$,
$\partial_j\cE_{geom}\cong \cF_{j,geom}$,
and
$\partial_k\cF_{i,geom}\cong \cG_{k,geom}$,
$\partial_k\cF_{j,geom}\cong \cG_{k,geom}^{op}$.
\begin{ddd}
\index{corner condition} We say that the tuple of these three isomorphisms satisfy the corner condition
if the induced isomorphism
$$\partial_k\cF_{j,geom}\cong  \cG_{k,geom}^{op}\cong
\partial_k\cF_{i,geom}^{op}$$
is up to sign the canonical isomorphism of canonical models constructed in Lemma
\ref{hiphop}.
\end{ddd}
By Lemma \ref{hiphop1} this has the following consequence.

\subsubsection{}\label{conscor}
Let $h\ge 0$ and $l\in I_h(\cG_{k,geom})$. Then we have
$l\in I_{h+1}(\cF_{i,geom})$, $l\in I_{h+1}(\cF_{j,geom})$,
and $l\in I_{h+2}(\cE_{geom})$.
Assume that we have a pre-taming of the face
$\partial_l\cG_{k,geom}$. It induces pre-tamings of
$\partial_l \cF_{i,geom}$ and $\partial_l \cF_{j,geom}$.
These two pre-tamings induce two pre-tamings of $\partial_l\cE_{geom}$.
If we assume the corner condition, then these two pre-tamings 
of $\partial_l\cE_{geom}$ are equal.

 \subsubsection{}

Let $\Delta$ be the category whose objects are the ordered sets $[n]:=\{0,1,\dots,n\}$ and whose morphisms are the monotone maps. A simplical set is a functor $S:\Delta^{op}\to \Sets$, i.e an object of the category of functors $\Sets^{\Delta^{op}}$. The nerve associated to an open covering of a space  was defined in explicit terms in \ref{hhjdiuqwiudqw}. It is  an example of a simplicial set.

\newcommand{\Mf}{{manifolds}}
\newcommand{\bM}{{\mathbf{M}}}
In similar manner we define a simplicial manifold as a functor
\index{simplicial manifold}$\bM:\Delta^{op}\to \Mf$, where $\Mf$ denotes the category of smooth manifolds and smooth maps. In explicit terms a simplicial manifold is given by a sequence of manifolds $\bM[p]$, $p=0,1,\dots$ and smooth maps $\partial_f:\bM[p]\to \bM[q]$ for all monotone map
$f:[q]\to [p]$ satisfying the relations $\partial_{f\circ g}=\partial_g\circ \partial_f$. 
In particular, for $j=0,\dots p+1$, we let
$j:[p]\to [p+1]$ be the monotone map whose image misses $j\in [p]$. It induces the map
$\partial_j:\bM[p+1]\to \bM[p]$.

\subsubsection{}\label{zuedqdqwdwq}

If $\{U_\alpha\}_{\alpha\in L}$ is an open covering of a manifold $B$ with associated nerve $\bN$, then
we can define the simplicial manifold
$\bM$ such that
$\bM[p]=\sqcup_{x\in \bN[p]} U_x$ (see \ref{hhjdiuqwiudqw} for notation) and
$\partial_f:\bM[p]\to \bM[q]$ is induced by the collection of incusions
$U_x\to U_{x\circ f}$, $x\in \bN[p]$, where
$f:[q]\to [p]$ is monotone.

\subsubsection{}

 Let $\bM$ be a simplicial manifold. Consider $k\in \nat$.
For $j\le k$ let
$I_j([k])$ be the set of $j$-element subsets of $[k]$. Note that
$[k]\cong I_1([k])$.
\begin{ddd}
\index{geometric $k$-chain!over simplicial manifold}A geometric $k$-chain over $\bM$ consists of 
\begin{enumerate}
\item geometric families $\cE_{geom}[p]$ over $\bM[p]$ for $p=0,1,\dots,k$,
\item identifications $I_j(\cE_{geom}[p])\cong I_j[p]$ for all $p=0,1,\dots,k$ and $j=0,\dots p$,
\item isomorphisms of geometric families $(-1)^p\partial_j \cE_{geom}[p]\cong  (-1)^j\partial_j^* \cE_{geom}[p-1]$ for $1\le p\le k$, $j\in [p]$ such that all possible corner conditions are satisfied.
\end{enumerate}
\end{ddd}
In order to illustrate the corner conditions consider $0\le i\le j\le p\le k-2$.
Let $k:[p]\to [p+2]$ be the monotone map whose image misses $i$ and $j$.
Then we have $\partial_j\circ \partial_i=\partial_i\circ \partial_{j-1}=\partial_k:[p]\to [p+2]$.
The datum of the geometric $k$-chain gives a chain of isomorphisms
\begin{eqnarray*}
\partial_i\partial_j\cE_{geom}[p+2]&\cong& (-1)^{p+j}\partial_i \partial_j^* \cE_{geom}[p+1]\\
&\cong& (-1)^{1+j+i}\partial_i^* \partial_j^* \cE_{geom}[p]\\
&\cong&(-1)^{1+j+i}\partial_k^*  \cE_{geom}[p]\\
&\cong&(-1)^{1+j+i}\partial_{j-1}^*\partial_i^*  \cE_{geom}[p]\\
&\cong&(-1)^{p+j}\partial_{j-1}^*\partial_i  \cE_{geom}[p+1]\\
&\cong&-\partial_{j-1}  \partial_{i} \cE_{geom}[p+2]\ .
\end{eqnarray*}
The corner condition requires that 
this isomorphism
$$\partial_i\partial_j\cE_{geom}[p+2]\cong -\partial_{j-1}  \partial_{i} \cE_{geom}[p+2]$$
is the canonical one from Lemma \ref{hiphop}.

\subsubsection{}

As an illustration of the notion of a geometric $k$-chain we consider the following example.
Let $G$ be a Lie group which acts smoothly on a manifold $B$. Then we
define a simplicial manifold $\bB(G,B)$ by the bar construction as follows.
We set $$\bB(G,B)[p]:=\underbrace{G\times \dots\times G}_{p\:factors }\times B\ .$$
For $j\in [p+1]$ we let
$\partial_j:\bB(G,B)[p+1]\to \bB(G,B)[p]$ be given by
\begin{eqnarray*}
\partial_j(g_1,\dots,g_{p+1},m)&:=&(g_1,\dots,g_{j}g_{j+1},g_{p+1},m)\ , 
j=1,\dots,p\\\partial_0(g_1,\dots,g_{p+1},m)&:=&(g_2,\dots,g_{p+1},m)\\
\partial_{p+1}(g_1,\dots,g_{p+1},m)&:=&(g_1,\dots,g_{p+1}m)\ .
\end{eqnarray*}
Let now $\cE_{geom}$ be a $G$-equivariant geometric family over $B$ with closed fibres.
Then we can define a geometric $k$-chain for all $k\ge 0$ as follows.
We define
$$\cE_{geom}[p]:=(-1)^{c(p)}\pr_B^*\cE_{geom}*\cDelta_{geom}^p\ ,$$ where
$\pr_B:\bB(G,B)[p]\to B$ is the projection, $\cDelta_{geom}^p$ is the geometric $p$-simplex
\ref{symy}, and $c(p):=\sum_{i=0}^p i=\frac{p(p+1)}{2}$.
We have $I_j(\cE_{geom})\cong I_j(\cDelta_{gem})\cong I_j[p]$.
For $j=0,\dots,p-1$ the identification 
$$(-1)^p\partial_j \cE_{geom}[p]\cong  (-1)^j\partial_j^* \cE_{geom}[p-1]$$
 is given by
\begin{eqnarray*}
 (-1)^p\partial_j \cE_{geom}[p]&\cong& (-1)^{p+c(p)}\pr_B^* \cE_{geom}*\partial_j\cDelta^p
\\&\cong& (-1)^{j+c(p-1)} \pr_B^* \cE_{geom}*\cDelta^{p-1}\\
&\cong&(-1)^j \partial_j^*\cE_{geom}[p-1]\ .
\end{eqnarray*}
 For $j=p$ we use the action of $G$.
Let $m:=(g_1,\dots,g_p,b)\in \bB(G,B)[p]$.
On the one hand the fibre of
$(-1)^p\partial_p \cE_{geom}[p]$ over $m$ is
$(-1)^{c(p)}\cE_{geom,b}*\cDelta^{p-1}$. On the other hand the fibre of
$(-1)^p\partial_p^*\cE_{geom}[p]$ over $m$ is
$(-1)^{p+c(p-1)}\cE_{geom,g_pb}*\cDelta^{p-1}$.
The $G$-action on $\cE_{geom}$ gives an isomorphism
$g_p:\cE_{geom,b}\stackrel{\sim}{\to} \cE_{geom,g_pb}$.
This induces (note that $c(p-1)+p=c(p)$) 
$$(-1)^p\partial_p \cE_{geom}[p]\cong (-1)^p\partial_p^* \cE_{geom}[p-1]\ .$$

This example will play an importand role in equivariant generalizations of the theory, but it is not the main example of the present paper.

\subsubsection{}

Let $B$ be a smooth manifold.
We consider an open covering 
$\cU=\{U_\alpha\}_{\alpha\in L}$ of $B$ with associated nerve $\bN$
and simplicial manifold $\bM$ as in \ref{zuedqdqwdwq}.
 \begin{ddd}\label{chd}
\index{$Z$}\index{$(Z^0,\dots, Z^k)$}\index{geometric $k$-chain!over manifold $B$}A (geometric) $k$-chain over $B$ (w.r.t. $\cU$) is
a geometric $k$-chain over $\bM$. In explicit terms it consists of
\begin{enumerate}
\item 
a $k+1$-tuple  $Z:=(Z^0,\dots, Z^k)$, where 
$Z^p$ associates to each $x\in \bN[p]$ a reduced  geometric family
$Z^p(x)$ over $U_x$,
\item identifications,
$I_{j}(Z^p(x))\cong I_j([p])$ for all $1\le j\le p$,
\item isomorphisms $$(-1)^p\partial_j Z^p(x)\cong (-1)^j
  Z^{p-1}(\partial^*_j x)_{|U_x}$$ for all  $1\le p\le k$,  $x\in
  \bN[p]$, and $j\in [p]$, such that
\item these isomorphisms satisfy all possible corner conditions.
\end{enumerate}
\end{ddd}
Note that the definition of $\partial_j Z^p(x)$ for $j\in [p]$
involves the identification $I_{1}(Z^p(x))\cong [p]$.

In Definition \ref{uiuuiiufwefwe} we associate to a geometric family $\cE_{geom}$ over $B$ geometric $k$-chains over $B$ (w.r.t $\cU$) for all $k\ge 0 $ which we call geometric $k$-resolutions.

\subsubsection{}

There is a natural notion of an isomorphism of $k$-chains (w.r.t. $\cU$). 
The set of isomorphism classes of $k$-chains (w.r.t. $\cU$) is denoted
by $\tilde G^k_{\cU}(B)$ . It forms an abelian  semigroup with respect
to disjoint union over $B$ followed by reduction of the face
decompositions.

This reduction in detail means the following. 
Let $Z$ and $Z^\prime$ be $k$-chains (w.r.t. $\cU)$, $0\le l\le
k$ and $x\in \bN[l]$. Then we first form  
$Y^l(x):=Z^l(x)\sqcup_{U_x} Z^{\prime l}(x)$.
Then $I_j(Y^l(x))=I_j(Z^l(x))\sqcup I_j(Z^{\prime l}(x))\cong I_j([l])\sqcup
I_j([l])$. We define a new admissible face decomposition $Y^l(x)_r$ such
that $I_j(Y^l(x)_r)=I_j([l])$,
where the face corresponding to
$A\in I_j([l])$ is the set of atoms of faces which belong
the faces of $Y^l(x)$ corresponding to the two copies of $A$ in
 $I_j(l)\sqcup
I_j(l)$.
We define
$(Z+Z^\prime)^l(x):=Y^l(x)_r$.

\subsubsection{}

If $f:\bM^\prime\to \cM$ is a morphism of simplicial manifolds, and $\cE$ is a geometric $k$-chain
over $\bM$, then we can define its pull-back $f^*\cE$, a geometric $k$-chain over $\bM^\prime$, in a natural way.

Let $\cU^\prime=\{U^\prime_\beta\}_{\beta\in L^\prime}$ , $L^\prime\rightarrow
L$, be  a refinement of $\cU$.
Then we have a natural morphism between the simplicial manifolds $\bM^\prime\to \bM$ associated to the coverings as in \ref{zuedqdqwdwq}. It induces a pull-back of geometric $k$-chains. 
In this way we get  a homomorphism $\tilde G^k_{\cU}(B)\rightarrow \tilde G^k_{\cU^\prime}(B)$.
\index{$\tilde G^k(B)$}\index{$\tilde G^k_{\cU}(B)$}By $\tilde G^k(B)$ we denote the abelian  semigroup
$$\tilde G^k(B):=\lim_{\longrightarrow} \tilde G^k_{\cU}(B)\ ,$$
where the limit is taken over the directed system of open coverings of $B$.
The elements $\tilde z\in \tilde G^k(B)$ will be called $k$-chains.

\subsubsection{}

Let $\bM$ be a simplicial manifold. Then there is a natural notion of a zero-bordism of a geometric $k$-chain over $\bM$. We leave it to the interested reader to write out the details of the definition in this generality. For the purpose of the present paper we only need the special case where $\bM$ is associated to an open covering of $B$. 

Let us present the details in this case.
We consider a $k$-chain $\tilde z\in \tilde G^k(B)$.
\begin{ddd}\label{bbo1}
\index{zero bordism of a geometric $k$-chain}A zero bordism of $\tilde z$ is given by
\begin{enumerate}
\item
an open covering $\cU$ such that $\tilde z$ is represented by a $k$-chain (w.r.t. $\cU$) $Z\in \tilde G^k_{\cU}(B)$,
\index{$W$}\index{$(W^0,\dots, W^{k})$}\item a $k$-tuple  $W:=(W^0,\dots, W^{k})$, where 
$W^p$ associates to each  $x\in \bN[p]$ a reduced geometric family $W^p(x)$ over $U_x$, 
\item for all  $0\le p\le k$ and $x\in \bN[p]$ an identification of $I_1(W^p(x))\cong [p]\cup\{*\}$,
\item isomorphisms $$(-1)^p\partial_* W^p(x) \cong  Z^p(x)$$ for all
  $0\le p\le k$ and   $x\in \bN[p]$,  and 
\item isomorphisms 
$$(-1)^p\partial_j W^p(x)\cong (-1)^j  W^{p-1}(\partial_j^*
x)_{|U_x}$$ for all  $1\le p\le k$,  $x\in \bN[p]$, and $j\in[p]$
isomorphisms such that
\item these isomorphisms satisfy all possible corner conditions. 
 \end{enumerate} 
\end{ddd}

\subsubsection{}
Let us again comment on the corner condition.
Fix e.g. $1\le p\le k$ and $j\in [p]$. Let $x\in \bN[p]$.
Then the data of $Z$ and $W$ give isomorphisms
$$(-1)^p\partial_j \partial_* W^p(x)\cong \partial_j Z^p(x)\cong
(-1)^{j+p} Z^{p-1}(\partial^*_j x)_{|U_x}$$
and
$$-(-1)^{p} \partial_*\partial_j W^p(x)\cong -(-1)^j\partial_*
W^{p-1}(\partial_j^*x)_{|U_x}\cong (-1)^{j+p} Z^{p-1}(\partial^*_j
x)_{|U_x}\ ,$$
and thus
$$\partial_j \partial_* W^p(x)\cong \partial_*\partial_j
W^p(x)^{op}\ .$$
The corner condition requires that this is the canonical isomorphism
of Lemma \ref{hiphop}.

\subsubsection{}
We say that $\tilde z$ is zero-bordant if it admits a zero bordism. 
The set  $\tilde G^k_0(B)\subseteq \tilde G^k(B)$ of zero bordant $k$-chains
forms a sub-semigroup. In fact, a zero bordism of $\tilde z+\tilde
z^\prime$
is given by a sum of zero bordisms of $\tilde z$ and $\tilde z^\prime$
which
is  defined in a similar manner as a sum of $k$-chains.

\begin{ddd}
\index{$G^k(B)$}We define $G^k(B):= \tilde G^k(B)/\tilde G^k_0(B)$.
\end{ddd}

\subsubsection{}
\begin{lem}\label{bord}
$G^k(B)$ is an abelian group.
\end{lem}
\proof
Let $\tilde z\in \tilde G^k(B)$. We claim that
\index{opposite!geometric $k$-chain}$\tilde z^{op}$ is the inverse of $\tilde z$.
Let $\tilde z$ be represented by the $k$-chain $Z$ (w.r.t $\cU$).
Then $\tilde z^{op}$ is represented by the $k$-chain
$Z^{op}:=(Z^{0,op},\dots,Z^{k,op})$ (w.r.t $\cU$) given by
$Z^{p,op}(x):=Z^p(x)^{op}$.

We consider the unit interval $I$ as a Riemannian spin manifold with
one boundary face $o$.
We define the zero bordism $W:=(W^0,\dots,W^k)$ of $\tilde z+\tilde z^{op}$ by
$W^p(x):=Z^p(x)*I$.  Let $*\in I_1(W^p(x))$ be the boundary face $Z^p(x)\times \partial_o I$.
The remaining boundary faces  are $\partial_j Z^p(x)\times I$, $j\in [p]$.
We fix the natural identifications $(-1)^p\partial_* W^p(x) = Z^p(x) * \partial_o I\cong (Z +Z^{op})^{p}(x)$,
$(-1)^p\partial_j W^p(x)=(-1)^p \partial_j Z^p(x) * I\cong (-1)^j  Z^{p-1}(\partial_j^* x)_{|U_x}* I =(-1)^jW^{p-1}(\partial_j^* x)_{|U_x}$.
Thus $W$ is a zero bordism of $\tilde z+\tilde z^{op}$.
\hB 
We will often write $-\tilde z$ for $\tilde z^{op}$.

\subsubsection{}
\index{tamed $k$-chain}If we replace geometric families by  tamed families in
the definition of chains \ref{chd} then we obtain the notion of tamed $k$-chains.
Of course we require that the isomorphisms now respect the tamings.

In a similar manner we get the notion of a tamed zero-bordism $W$ of a tamed $k$-chain $Z$.
In the definition  \ref{bbo1}  we now require that the families $W^p(x)$  are tamed for $p<k$, and
that $W^k(x)$ is boundary tamed. Again, all isomorphisms in this definition must now respect the tamings.

By  $\tilde G^k_{0,t}(B)\subseteq \tilde G^k_t(B)$ we denote the corresponding semigroups of tamed and  of zero bordant tamed $k$-chains. 
\begin{ddd}
\index{$G^k_t(B)$}\index{$\tilde G^k_t(B)$}We define the semigroup
$G^k_t(B):= \tilde G^k_t(B)/\tilde G^k_{0,t}(B)$.
\end{ddd}
We will see later in Lemma  \ref{gtamegroup} that $G^k_t(B)$ is a group, too.

\subsubsection{}
Note that the correspondences $B\mapsto G^k(B)$ and $B\mapsto G^k_t(B)$
are contravariant functors on the category of smooth manifolds with values in
(semi)groups. On morphisms these functors are given by pull-back.
We leave it to the reader to write out the details.

\subsubsection{}
There is a natural commutative diagram of homomorphisms
$$\begin{array}{ccc}
\tilde G^k_{0,t}(B)&\rightarrow&\tilde G^k_{0}(B) \\
\downarrow &&\downarrow\\
\tilde G^k_{t}(B)&\rightarrow&\tilde G^k(B) 
\end{array}\ ,
$$
which induces a homomorphism
$G^k_t(B)\rightarrow G^k(B)$.
Furthermore, we have natural homomorphisms
$G^k(B)\rightarrow G^{k-1}(B)$, $G^{k}_t\rightarrow G^{k-1}_t(B)$ such that $$\begin{array}{ccc}
 G^k_{t}(B)&\rightarrow&G^{k-1}_t(B) \\
\downarrow &&\downarrow\\
G^k(B)&\rightarrow&G^{k-1}(B) 
\end{array}
$$
commutes.

\subsection{Obstruction theory: taming of chains}\label{yy}

\subsubsection{}
Let $0\le p\le k$. Then there is a natural forgetful map
\index{$\cF$}\index{forgetful map}$\cF:\tilde G^k(B)\rightarrow \tilde G^p(B)$. Furthermore, there is a forgetful map
$\cF:\tilde G^p_t(B)\rightarrow  \tilde G^p(B)$
(we use the symbol $\cF$ to denote various forgetful maps).

\subsubsection{}
Next we introduce the concept of partially tamed chains.
\index{$(Z,Z_t)$}\index{$\tZ$}\index{$\tilde F^k_{p,\cU}(B)$}Let $\cU$ be an open covering of $B$.
We consider pairs $\tZ:=(Z,Z_t)$, where $Z$ is a geometric $k$-chain (w.r.t. $\cU$),
$Z_t$ is a tamed $p$-chain (w.r.t. $\cU$), and $\cF(Z_t) \cong   \cF(Z)$
as geometric $p$-chains. This isomorphism is part of the structure.
By $\tilde F^k_{p,\cU}(B)$ we denote the set of isomorphism classes of
such pairs which we call partially tamed chains (w.r.t $\cU$).
It is again a semigroup under the operation of disjoint sum over $B$ followed by reduction of the face decomposition.
If $\cU^\prime$ is a refinement of $\cU$, then we have a homomorphism
$\tilde F^k_{p,\cU}(B)\rightarrow \tilde F^k_{p,\cU^\prime}(B)$.
We define the semigroup \index{$\tilde F^k_{p}(B)$}
$$\tilde F^k_{p}(B):=\lim_{\longrightarrow} \tilde F^k_{p,\cU}(B)\ ,$$
where the limit is taken over the system of open coverings of $B$.

Furthermore, 
we define $\tilde F^k_{-1}(B):= \tilde G^k(B)$ and identify $\tilde F^k_k(B)= \tilde G^k_t(B)$.
For $p\le q\le k$ there is a forgetful map
$\cF^q_p:\tilde F^k_q(B)\rightarrow \tilde F^k_p(B)$.
The elements of $\tilde F^k_{p}(B)$ will be called partially tamed chains.
\index{partially tamed chains}
\subsubsection{}

On $\tilde F^k_p(B)$ we define the following notion of a zero bordism.
A zero bordism of a partially tamed chain  $z\in \tilde F^k_p(B)$ is
given by a covering $\cU$ of $B$, a partially tamed chain
$\tZ:=(Z,Z_t)$ (w.r.t $\cU$) representing
$z$,   and a pair $\tW:=(W,W_{t})$
such that $W$ is a zero bordism of $Z$, 
$W_{t}$ is a zero bordism of $Z_t$, and
$\cF (W)\cong \cF(W_{t})$ in a way which is compatible with the
isomorphism $\cF(Z)\cong  \cF(Z_t)$.

\subsubsection{}
Let $\tilde F^k_p(B)_0\subseteq \tilde F^k_p(B)$ be the sub-semigroup of
partially tamed chains which are zero bordant.

\begin{ddd}
\index{$F^k_p(B)$}We  define $F^k_p(B):=\tilde F^k_p(B)/\tilde F^k_p(B)_0$.
\end{ddd}

Note that there is natural homomorphism $\cF:F^k_p(B)\rightarrow G^p_t(B)$. 
We further define $F^k_{-1}(B):= G^k(B)$ and identify $F^k_k(B)= G^k_t(B)$.
For $p\le q\le k$ we have a homomorphism $\cF^q_p:F^k_q(B)\rightarrow  F^k_p(B)$.
We will see in Lemma \ref{gtamegroup} that the semigroups $F^k_p(B)$ are in fact groups.

\subsubsection{}

Let $p\le k-1$ and a partially tamed chain  $z\in \tilde F^k_{p-1}(B)$ be given. In the present subsection we study the question
under which conditions there exists a partially tamed chain $z^\prime\in \tilde F^k_{p}(B)$ such that 
$\cF^p_{p-1}(z^\prime)=z$.

\subsubsection{}

Let $\cU$ be an open covering of $B$ such that $z$ is represented
by a partially tamed chain $\tZ:=(Z,Z_t)\in\tilde F^k_{p-1,\cU}(B)$
(w.r.t $\cU$).
If $x\in \bN[p]$ and $j\in [p]$, 
then $(-1)^p\partial_j Z^{p}(x)\cong (-1)^j \cF (Z^{p-1}_t(\partial_j^* x)_{|U_x})$, so that we obtain a
boundary taming $Z^p_{bt}(x)$ of the underlying geometric family
\index{$\bz$}$Z^p(x)$. In fact, it is a consequence of the corner conditions that
this boundary taming is well-defined (see \ref{conscor}).
We consider the chain (see \ref{checkjdiisodic} for conventions) 
\index{index!of a partially tamed chain}\index{$\ind(\tZ)$}$\ind(\tZ)\in \check{C}^p(\cU,\underline{\Z}_B)$
given by \begin{equation}\label{ggg45673}\ind(\tZ):=\prod_{x\in \bN[p]} \ind_0(Z^p_{bt}(x)) \ . \end{equation}

\subsubsection{}
Let $\tZ\in\tilde F^k_{p-1,\cU}(B)$ represent $z\in \tilde
F^k_{p-1}(B)$
and $\bz\in F^k_{p-1}(B)$. Assume that $0\le p\le k-1$.
\begin{lem}\label{olem}
\begin{enumerate}
\item We have  
$\delta  \ind(\tZ)=0$.
\index{$o^p(\bz)$}\item The class $$o^p(\bz) :=[\ind(\tZ)]\in
\check{H}^{p}(B,\underline{\Z}_B)$$ only depends on the class 
$\bz\in F^k_{p-1}(B)$.
\item
The cohomology class
$o^p(\bz)$  only depends on the image $\cF^{p-1}_{p-2}(\bz)$, where 
$\cF^{p-1}_{p-2}:F^k_{p-1}(B)\rightarrow F^k_{p-2}(B)$.
\item If $o^p(\bz)=0$, then after refining the covering and altering the taming of
$Z^{p-1}_t$ on codimension zero faces we can find   a pair $\tZ^\prime=(Z,Z^\prime_t)\in \tilde F^k_{p,\cU^\prime}(B)$
representing 
$ z^\prime \in \tilde F^k_{p}(B)$ such that 
$\cF^p_{p-1}(z^\prime)=z$.
\item If $o^p(\bz)=0$, then there exists $\bz^\prime\in F^k_p(B)$ such
  that $\cF^p_{p-2}(\bz^\prime)=\cF^{p-1}_{p-2}(\bz)$.
 \end{enumerate}
\end{lem}
\proof
It is clear that 4. implies 5.
The Assertions 1.,2. and 3. of the Lemma are  only non trivial if the
dimension of the fibers of $Z^p$ is even.
During the following proof this will be a standing assumption.
The proof of Assertion 4. in the case of odd-dimensional fibers can be done using parts of the
arguments for Assertion 4. in the even-dimensional case.

Recall that $\cA_B^q$ denotes the sheaf of $q$-forms on $B$.
For $r,q\in\nat\cup\{0\}$ we define the following chains (see \ref{checkjdiisodic}, \ref{eet0} and \ref{eet1} for notation):
\begin{eqnarray*}
\Omega^q(Z^r)&:=&\prod_{x\in \bN[r]}\Omega^q(Z^r(x))\in \check{C}^r(\cU,\cA_B^q)\ , \quad r\le k \\
\eta^q(Z_t^{r})&:=&\prod_{x\in \bN[r]}\eta^q(Z^{r}_{t}(x))\in \check{C}^{r}(\cU,\cA_B^q) \ ,\quad r\le p-1  \ .
\end{eqnarray*}
We consider $\underline{\Z}_B$ as a sub-sheaf of $\cA_B^0$.
Then we compute using Theorem \ref{etaprop}, (2.), and Definition \ref{chd}, 2. that
$\ind(\tZ)=\Omega^0(Z^p)+(-1)^{p}\delta \eta^0(Z_t^{p-1})$.
It follows that $\delta \ind(\tZ)=\delta \Omega^0(Z^p)$.
By Lemma \ref{omeg} and  Definition \ref{chd}, 2. we have the general relation
$$\delta \Omega^q(Z^p)= (-1)^{p} d \Omega^{q-1}(Z^{p+1})\ .$$
It is here where we use that $p\le k-1$.
For $q=0$ we get in particular $\delta \Omega^0(Z^p)=0$.
This shows 1.

In order to show Assertion 2.
it suffices to show that $[\ind(\tZ)]=0$ if
$\tZ$ admits a zero bordism
$\tW=(W,W_t)$.
By Lemma \ref{omeg} and Definition \ref{bbo1}, 4. and 5.,  we have the general relation
$$\delta \Omega^q(W^{p-1})+  \Omega^q(Z^{p})= (-1)^{p+1} d \Omega^{q-1}(W^{p})\ .$$
For $q=0$ we obtain 
$-\delta \Omega^0(W^{p-1})=\Omega^0(Z^{p})$.
Let $\partial W_{bt}^{p-1}$ denote the object which associates to $y\in\bN[p-1]$ the boundary tamed family
$\partial W^{p-1}_{bt}(y)$. 
This boundary taming is again well-defined because of the corner conditions.
We define 
\begin{eqnarray*}
\eta^0(\partial W_{bt}^{p-1})&=&\prod_{y\in\bN[p-1]} \eta^0(\partial W^{p-1}_{bt}(y)) \in \check{C}^{p-1}(\cU,\cA^0_B)\\
\ind(W^{p-1}_{bt})&:=&\prod_{y\in \bN[p-1]} \ind_0(W^{p-1}_{bt}(y))\in \check{C}^{p-1}(\cU,\underline{\Z}_B)\ .
\end{eqnarray*}
Since we have 
$(-1)^{p-1} \eta^0(\partial W_{bt}^{p-1}) = \eta^0(Z_t^{p-1}) + \delta \eta^0(W_{t}^{p-2})$
we see that $\delta\eta^0(Z_t^{p-1})= (-1)^{p-1} \delta \eta^0(\partial W_{bt}^{p-1})$.
It follows that 
\begin{eqnarray*}
\ind(\tZ)&=&\Omega^0(Z^p)+(-1)^p\delta \eta^0(Z_t^{p-1})\\
&=&-\delta\left(\Omega^0(W^{p-1}) +  \delta \eta^0(\partial W_{bt}^{p-1})\right)\\
&=&-\delta \ind(W^{p-1}_{bt})\ .
\end{eqnarray*}
We now see that $[\ind(\tZ)]=0$.
We have thus shown Assertion 2.

Let $Z^{\prime,p-1}_t$ be an alteration of
the taming of $Z^{p-1}_t$ on the codimension zero faces.
Let $\tZ^\prime$ denote the corresponding pair.
We define the chain
\begin{equation}\label{relsf}c(Z^{\prime,p-1}_t,Z^{p-1}_t):=\prod_{y\in \bN[p-1]} -\Sf(Z^{\prime,p-1}_t(y),Z^{p-1}_t(y))\in  \check{C}^{p-1}(\cU,\underline{\Z}_B)\ ,\end{equation}
where $ \Sf(Z^{\prime,p-1}_t(y),Z^{p-1}_t(y))$ is the locally constant integer valued function
on $U_y$ given by the fiber-wise spectral flow as introduced in
Definition \ref{sflowdef}.
Then we have (Lemma \ref{diieta})
$$\eta^0(Z_t^{\prime,p-1})-\eta^0(Z_t^{p-1})=-c(Z^{\prime,p-1}_t,Z^{p-1}_t)\ .$$
We see that
$$\ind(\tZ^\prime)-\ind(\tZ)=(-1)^{p+1} \delta c(Z^{\prime,p-1}_t,Z^{p-1}_t)\ .$$
It follows that
$[\ind(\tZ)]=[\ind(\tZ^\prime)]$.
This finishes the proof of Assertion 3.

Assume now that $o^p(\bz)=0$.
Let $c\in \check{C}^{p-1}(\cU,\underline{\Z}_B)$ such that $\delta c= \ind(\tZ)$.
We assume that $\cU$ is a good covering, i.e., that all intersections $U_x$, $x\in \bN(\cU)^q$,
$q\in\nat\cup\{0\}$, are contractible. Then we can write
$c=\prod_{y\in \bN[p-1]} c_y$, where $c_y\in\Z$.
If $b\in U_y$, then we can find by Lemma \ref {spg} an alteration of the taming
of the fiber $Z^{\prime,p-1}_{t}(y)_b$ in codimension zero such that
$-\Sf(Z^{\prime,p-1}_{t}(y)_b,Z^{p-1}_{t}(y)_b)=(-1)^pc_y$.
By continuity this holds true on a neighborhood of $b$. Thus, after a good refinement of the good covering
we can find an alteration $Z^{\prime,p-1}_t$ of the taming of $Z^{p-1}_t$ in codimension zero
such that $c=(-1)^p c(Z^{\prime,p-1}_t,Z^{p-1}_t)$.
Then we have $\ind(\tZ^\prime)=\ind(\tZ)+(-1)^{p+1} \delta c(Z^{\prime,p-1}_t,Z^{p-1}_t)=\ind(\tZ)-\delta c=0$.

Now assume that $\ind(\tZ)=0$. If $x\in \bN[p]$, then we have
$\ind(Z^p_{bt}(x))=0$. If $b\in U_x$, then we can extend the boundary taming
of the fiber $Z^p_{bt}(x)_{b}$ to a taming $Z^p_{t}(x)_{b}$. Again, by continuity, we obtain
an extension of the boundary taming to a taming over a neighborhood of $b$. Thus after refining the covering we obtain 
a partially tamed $\tZ^\prime$ representing 
$z^\prime\in \tilde F^k_{p}(B)$ such that 
$\cF^p_{p-1}(z^\prime)=z$.
This finishes the proof of Assertion 4. \hB

\subsubsection{}

In Theorem \ref{compara} we will show the relation between the set of obstruction classes $\bo^{p}(z)$ for $z$ represented by 
$(Z,Z_t)$ (with varying $Z_t$)
and the obstruction set $\bo^{p}(\ind(\cE_{geom}))$ introduced in \ref{obset} in the special case that $Z$ comes from
a geometric $k+1$-resolution (see Definition \ref{uiuuiiufwefwe}) of a geometric family $\cE_{geom}$ over $B$ with closed fibres such that
$\ind(\cE_{geom})\in K_p(B)$.

\subsection{Obstruction theory: taming of zero bordisms}
 
\subsubsection{}
Let $\tZ^\prime:=(Z,Z^\prime_t)\in \tilde F^k_{p,\cU}(B)$ be a partially
tamed chain (w.r.t. $\cU$) representing 
$z^\prime \in \tilde F^k_p(B)$.
\index{$\tW$}\index{$(W,W_t)$}Let $\tZ:=(Z,Z_t):=\cF(\tZ^\prime)\in\tilde F^k_{p-1,\cU}(B) $ represent $z:=\cF(z^\prime)\in \tilde F^k_{p-1}(B)$.
Assume that $\tW:=(W,W_t)$ is a zero bordism of $\tZ$.
In the present subsection we study the question under which conditions
we can extend the taming of $W_t$ to $W_t^\prime$ such that $\tW^\prime:=(W,W^\prime_t)$ is a zero
bordism of $\tZ^\prime$.
\subsubsection{}

Note that for $y\in\bN[p-1]$ we have a boundary tamed family $W_{bt}^{p-1}(y)$.
Assume that we can extend the boundary taming to a taming $W_{t}^{p-1}(y)$ for
all $y\in \bN[p-1]$.
Since for $x\in\bN[p]$ we have 
$(-1)^{p}\partial_* W^{p}(x)=\cF(Z^p_t)$ and
$(-1)^{p}\partial_j W^{p}(x)=(-1)^j \cF(W^{p-1}_t(\partial_j^* x))_{|U_x}$, $j\in[p]$,
we would get a boundary taming $W^p_{bt}(x)$. In this way we can define $W_t^\prime$.
 
\subsubsection{}
We consider the chain 
\index{$\ind(W^{p-1}_{bt})$}$$\ind(W^{p-1}_{bt}):=\prod_{x\in \bN[p-1]} \ind_0(W^{p-1}_{bt}(x))\in \check{C}^{p-1}(\cU,\underline{\Z}_B)\ .$$
\begin{lem}\label{ww}
\begin{enumerate}
\item We have $\delta \ind(W^{p-1}_{bt})=0$.
\index{$p(\tZ^\prime,\tW)$}\item The class $p(\tZ^\prime,\tW):=[\ind(W^{p-1}_{bt})]\in
\check{H}^{p-1}(B,\underline{\Z}_B)$
is independent of the choice of the taming $W_t^{p-2}$ in codimension zero.
\item  If $p(\tZ^\prime,\tW)=0$, then after refining the covering and altering the taming
of $W_t^{p-2}$ in codimension zero we can extend the boundary taming of $W_{bt}^{p-1}$
to a taming $W_t^{p-1}$ so that the resulting pair $\tW^\prime$ is a zero bordism of $\tZ^\prime$.
\end{enumerate}
\end{lem}
\proof
Assertions 1. and 2. are only nontrivial if the dimension of the fiber of $W^{p-1}$ is even.
This will be the standing assumption in the following proof. Assertion 3. in the
odd-dimensional case can be proved using a part of the arguments for the even-dimensional case.

We have
$$\ind(W^{p-1}_{bt})=\Omega^0(W^{p-1}) + \eta^0(\partial
W^{p-1}_{bt})\ .$$  Since
$$\eta^0(\partial W^{p-1}_{bt})=(-1)^{p-1} \delta \eta^0(W^{p-2}_t) + (-1)^{p-1}\eta^0(Z^{p-1}_t)$$
we have 
$$\delta \ind(W^{p-1}_{bt})=\delta \Omega^0(W^{p-1})+ (-1)^{p-1}\delta \eta^0(Z^{p-1}_t)\ .$$
Furthermore,
$$\delta \eta^0(Z^{p-1}_t)=(-1)^p \eta^0(\partial Z^{p}_{t})$$
and 
$$\delta \Omega^0(W^{p-1})+\Omega^0(Z^p)=(-1)^p \Omega^0(\partial W^p)=0$$ 
so that
\begin{eqnarray*}
\delta \ind(W^{p-1}_{bt} )
&=&
-\Omega^0(Z^p) - \eta^0(\partial Z^{p}_{t}) \\
&=&-\ind(Z^p_t)\\
&=&0\ .
\end{eqnarray*}
This proves Assertion 1.

If we change the taming $W^{p-2}_t$ to $W_t^{\prime,p-2}$ in codimension zero, then we have
\begin{eqnarray*}
\ind(W^{\prime,p-1}_{bt})-\ind(W^{p-1}_{bt}) &=&(-1)^{p-1}\delta(\eta^0(W_t^{\prime,p-2})-\eta^0(W^{p-2}_t))\\
&=&(-1)^{p}\delta c(W_t^{\prime,p-2},W^{p-2}_t)
\end{eqnarray*}
(see (\ref{relsf}) for a definition of $c(W_t^{\prime,p-2},W^{p-2}_t)$).
It follows that $p(\tZ,\tW^\prime)=p(\tZ,\tW)$.
This shows Assertion 2.

Assume now that $p(\tZ,\tW)=0$.
Let  $c\in  \check{C}^{k-2}(\cU,\underline{\Z}_B)$ be such that
$\delta c = \ind(W^{p-1}_{bt})$.
As in the proof of Lemma \ref{olem} we find (after refinement of the covering) an alteration $W_t^{\prime,p-2}$
of the taming of $W^{p-2}_t$ in codimension zero such that $c(W_t^{\prime,p-2},W^{p-2}_t)=(-1)^{p-1}c$.
Then $\ind(W^{\prime,p-1}_{bt})=0$. Now (after further refinement of the covering) we find an extension $W^{\prime,p-1}_{t}$ of the boundary taming $W^{\prime,p-1}_{bt}$ to a taming. This provides $\tW^\prime$ as required.
We thus have shown Assertion 3. \hB

\subsection{$F^k_p(B)$ and $G^k_t(B)$ are groups}\label{grus}

\subsubsection{}\label{someis}

\index{$\check{\bC}(B,\cS)$}If $\cS=(\cS^q,d)$ is a complex of sheaves on $B$, then by 
$\check{\bC}(B,\cS)$ we denote the \index{total complex of a double complex} total complex of the double complex
\index{$\bd$}$(\check{C}^p(B,\cS^q),d,\delta)$ with differential $\bd c^{p,q}=(-1)^p d c^{p,q}-\delta c^{p,q}$
for  $c^{p,q}\in \check{C}^p(B,\cS^q)$. By $\check{\bH}(B,\cS)$ we denote the cohomology of $\check{\bC}(B,\cS)$\index{hyper-cohomology}\index{$\check{\bH}(B,\cS)$}
which is usually called the hyper-cohomology of $\cS$.

\subsubsection{}\label{tzt6565s}
We apply this construction to the complex $\cA_B$ of differential forms on $B$.
The embedding $\underline{\R}_B\hookrightarrow \cA_B$ 
(here we consider $\underline{\R}_B$ as a complex of sheaves)
is a quasi-isomorphism and thus
induces
an isomorphism $\check{H}(B,\underline{\R}_B) \stackrel{\sim}{\rightarrow} \check{\bH}(B,\cA_B)$.
The embedding $\cA_B(B)\hookrightarrow \check{C}^0(B,\cA_B)$ induces an isomorphism
$H_{dR}(B)\stackrel{\sim}{\rightarrow} \check{\bH}(B,\cA_B)$ since 
the sheaves $\cA_B^*$ are soft.
The composition of the first with the inverse of the second
listed isomorphism gives the de Rham isomorphism
$\check{H}(B,\underline{\R}_B) \stackrel{\sim}{\rightarrow}H_{dR}(B)$
\footnote{This fixes in particular the signs}.

\subsubsection{}

If $z\in \check{H}(B,\underline{\Z}_B)$, then let $z_\R\in \check{H}(B,\underline{\R}_B)$ denote
its image under the natural homomorphism $\check{H}(B,\underline{\Z}_B)\rightarrow \check{H}(B,\underline{\R}_B)$.
The class $z$ is a torsion class iff $z_\R=0$. The condition $z_\R=0$ is equivalent to the condition $\bbz=0$, where 
$\bbz\in  \check{\bH}(B,\cA_B)$ is the image of $z$ under $\check{H}(B,\underline{\Z}_B)\rightarrow \check{H}(B,\underline{\R}_B)\rightarrow  \check{\bH}(B,\cA_B)$.
Let $z$ be represented by a \v{C}ech cocycle $c\in \check{C}(B,\underline{\Z}_B)$.
Let $\bc\in \check{C}(B,\cA^0)\subseteq \check{\bC}(B,\cA_B)$ be the image of $c$.
If there is a chain $\bb\in\check{\bC}(B,\cA_B)$ such that $\bd\bb=\bc$, then $z$ is a torsion class.
We employ this sort of argument in the proofs of the following two lemmas.

\subsubsection{}

Let $k-1\ge p$.
\begin{lem}\label{tore}
If $\bz\in\ker(\cF: F^k_{p-1}(B)\rightarrow G^k(B))$, then
the obstruction $o^p(\bz)\in \check{H}^{p}(B,\underline{\Z}_B)$
is a torsion class.
\end{lem}
\proof
Let $\bz$ be represented by a partially tamed chain $\tZ=(Z,Z_t)\in
\tilde F^k_{p-1,\cU}(B)$ (w.r.t.
  $\cU$),  and let $W$ be a zero bordism of $Z$.
The class $o^p(\bz)$ is represented by the \v{C}ech cocycle
$\ind(\tZ)\in \check{C}^p(\cU,\underline{\Z}_B)$. Let $\bc \in   \check{\bC}^p(\cU,\cA_B)$ be the image of $\ind(\tZ)$.
We define the following chains:
\begin{eqnarray*}
\bbeta&:=&(\eta^{p-1}(Z^0_t),\dots, \eta^{0}(Z^{p-1}_t)) \in \check{\bC}^{p-1}(\cU,\cA_B)\\
\bOmega(W)&:=&((-1)^0\Omega^{p-1}(W^0),(-1)^1 \Omega^{p-2}(W^1),\dots,(-1)^{p-1}\Omega^{0}(W^{p-1}))\in \check{\bC}^{p-1}(\cU,\cA_B)\\
\bOmega(Z)&:=&((-1)^0\Omega^{p}(Z^0),(-1)^1\Omega^{p-1}(Z^1,)\dots,(-1)^p\Omega^{0}(Z^{p}))\in\check{\bC}^{p}(\cU,\cA_B)\ .
 \end{eqnarray*}
It follows from Lemma \ref{omeg} that $\bd \bOmega(W)+\bOmega(Z)=0$.
Furthermore, by Theorem \ref{etaprop}
$$\bd \bbeta=\bOmega(Z)-(-1)^p\Omega^0(Z^p)- \delta \eta^0(Z^{p-1}_t)
=\bOmega(Z)-(-1)^p\bc\ .$$
We conclude that $\bc=(-1)^ {p+1}\bd( \bOmega(W)+ \bbeta)$.\hB

\subsubsection{}

\begin{lem}\label{constor}
If $\tilde z^{\prime}\in \tilde G^{k+1}(B)$ is such that 
$\tilde z:=\cF(\tilde z^{\prime})\in \tilde G_0^{k}(B)$, then there exists a number $N\in\nat$ 
and $\tilde u\in \tilde G^k_{t}(B)$ such that $\cF(\tilde u)= N \tilde z$.  
\end{lem}
\proof
We consider 
$v_{-1}:=\tilde z^{\prime}\in F_{-1}^{k+1}(B)$.
We define inductively lifts $v_p\in \tilde F^{k+1}_p(B)$ of
$N_p v_{-1}$ for suitable $N_p\in\nat$.

Assume that we have already found a lift
$v_{p-1}\in \tilde  F^{k+1}_{p-1}$(B) of $N_{p-1}v_{-1}$ and $p\le k$. Then we have
$\bv_{p-1}\in \ker(\cF: F^{k+1}_{p-1}(B)\rightarrow G^k(B))$.
Therefore by Lemma  \ref{tore} the class $o^p(\bv_{p-1})$ is a torsion class and
$Lo^p(\bv_{p-1})=0$ for some $L\in\nat$.
We define $N_p:=N_{p-1}L$.
By Lemma \ref{olem} we now find a lift $v_p\in \tilde F^{k+1}_p(B)$ of $Lv_{p-1}$.
Eventually we obtain an element $v_k\in \tilde F^{k+1}_k(B)$ which lifts
$N_kv_{-1}$. 

We set $N:=N_k$ and let $\tilde u\in\tilde G^k_t(B)$ be given by
$\cF(v_k)$.\hB

\subsubsection{}

Let $\tilde z_t\in \tilde G^k_{t}(B)$. We form $\tilde y_t=\tilde z_t + \tilde z_t^{op}$.
Let $\tilde z_t$ be represented by the tamed chain $Z_t\in\tilde G^k_{t,\cU}(B)$ (w.r.t. $\cU$),  
and let  $Y_t:=Z_t+ Z_t^{op}$ be the corresponding representative of $\tilde y_t$.
Let $W$ be the cylinder introduced in the proof of Lemma \ref{bord}
which gives the zero bordism of $Y:=\cF(Y_t)$.

Let $p\le k-1$, and let $y \in\tilde F^k_{p-1}(B)$
be represented by the partially tamed chain  $\tY:=(Y,\cF(Y_t))$
(w.r.t $\cU$). 
Assume that $W$ admits a tamed lift $W_t$ such that $\tW:=(W,W_t)$
is a zero bordism of $\tY$.
\begin{lem}\label{tto}
The obstruction $p(\tY,\tW)\in \check{H}^{p-1}(B,\underline{\Z}_B)$ (defined in Lemma \ref{ww})
is a torsion class.
\end{lem}
\proof
The class $p(\tY,\tW)\in \check{H}^{p-1}(B,\underline{\Z}_B)$
is represented by the  \v{C}ech cocycle
$\ind(W^{p-1}_{bt})\in \check{C}^{p-1}(\cU,\underline{\Z}_B)$.
Let $\bc \in   \check{\bC}^{p-1}(\cU,\cA_B)$ be the image of $\ind(W^{p-1}_{bt})$.
We define the chain
$$\bbeta:=(\eta^{p-2}(W^0_t),\dots,\eta^{0}(W^{p-2}_t))\in \check{\bC}^{p-2}(\cU,\cA_B)\ .$$
By Theorem \ref{etaprop}, the fact that the local index form for a cylinder vanishes,  
and that $\eta^q(Y^{p-2-q}_t)=0$ (since $Y^{p-2-q}_t\cong Y^{p-2-q,op}_t$) we obtain
$\bd \bbeta=-\delta \eta^0(W_t^{p-2})$.
It follows
\begin{eqnarray*}
\bc&=&\Omega^0(W^{p-1})+\eta^0(\partial W^{p-1}_{bt})\\
&=&(-1)^{p-1}\delta \eta^0(W^{p-2}_t)\\
&=&  (-1)^p\bd \bbeta\ .
\end{eqnarray*}
 \hB

\subsubsection{}

\begin{lem}\label{gtamegroup}
\begin{enumerate}
\item
Let $\tilde z_t\in \tilde G^k_{t}(B)$ and $\tilde y_t=\tilde z_t + \tilde z_t^{op}$.
Then there exists $N\in\nat$ such that $N\tilde y_t\in \tilde G^k_{0,t}(B)$.
\item $G^k_{t}(B)$ is a group.
\item $F^k_p(B)$ is a group.
\end{enumerate}
\end{lem}
\proof
Assertion 1. is a consequence (by a similar argument as for Lemma  \ref{constor}) of Lemma \ref{ww} and of Lemma \ref{tto}.

Assertion 2. follows from the first. In fact,
the inverse of the class  $z_t\in G^k_{t}(B)$
represented by $\tilde z_t$
is given by the class represented by $(N-1)\tilde z_t+N \tilde z_t^{op}$.

Let $\bz\in  F^k_p(B)$ be represented by $\tZ=(Z, Z_t)\in \tilde F^k_p(B)$.
Then we form $\tY:=(Y,Y_t)$,
where $Y_t:=Z_t + Z_t^{op}$ and $Y:=Z + Z^{op}$.
It again follows from Lemma \ref{tto} that $N \tY \in \tilde F^k_p(B)_0$
for a suitable $N\in\nat$. Thus the element
$(N-1) \bz + N \bz^{op}\in F^k_p(B)$
is the inverse of $\bz$.
This proves Assertion 3.
\hB

\section{Resolutions}\label{kj65}

\subsection{The $n$-simplex} \label{symy}

\subsubsection{}

In Lemma \ref{symy1} we have constructed a manifold with corners
version of the $n$-simplex $\Delta^n$ such that $\Sigma^{n+1}$ acts by
automorphisms. 

The group $\Sigma^{n+1}$ acts naturally on $[n]$ and thus on the set
$I_{k}([n])$ of $k$-element subsets of $[n]$.
For all $n\ge 0$ we fix a $\Sigma^{n+1}$-equivariant identification $I_k(\Delta^n)\cong
I_{k}([n])$. In view of the construction of $\Delta^n$ we can for each $n\ge 1$ fix
an identification of $\Delta^{n-1}$ with the underlying manifold of a canonical model of
$\partial_0\Delta^n$. 
We fix the orientations of $\Delta^n$ in the unique way such that
$\Delta^{n-1}\cong \partial_0\Delta^n$ is orientation preserving
and such that the point $\Delta^0$ is positively oriented.

For $j\in [n]$ let
$\sigma_j\in \Delta^{n+1}$ be the permutation
$(0,1,\dots ,n)\mapsto (1,\dots, j,0,j+1,\dots n)$. 
By composing the identification $\Delta^{n-1}\cong \partial_0\Delta^n$
with $\sigma_j$ we obtain an identification
$\partial_j \Delta^n\cong (-1)^j \Delta^{n-1}$ for all $j\in [n]$.

For all $n\in\nat$ we equip $\Delta^n$ with a $\Sigma^{n+1}$-invariant
admissible  Riemannian metric such that
$\Delta^{n-1}\cong \partial_0\Delta^n$ is an isometry.
\index{$\cDelta^n_{geom}$}Note that $\Delta^n$ has a unique spin structure.
The collection of this data will be the geometric manifold
$\cDelta^n_{geom}$.

Let $j\in [n]$ and $T_j:\Delta^{n-1}\times (0,1)\rightarrow \Delta^n$
be the corresponding distinguished  embedding.
We require that the product structure $\Pi_j:T_j^*\cS(\Delta^n)\rightarrow
\cS(\Delta^{n-1})*(0,1)$ preserves the real or quaternionic structures.
Then it is unique up to sign. 
 
We thus have fixed for all  $j\in[n]$ an isomorphism
$\partial_j\cDelta^n_{geom}\cong (-1)^j\cDelta^{n-1}_{geom}$
in the sense of \ref{langg}.
The corner condition is automatic.

\subsubsection{}

Note that the admissible Riemannian metric on $\cDelta^n_{geom}$
is not flat. Consider e.g. the parallel transport in $T\Delta^2$
along the closed curve given by the boundary of $\Delta^2$ which
produces a $\pi/2$-rotation. There is no reason that the form
$\hA(\nabla^{T\Delta^{4n}})_{4n}$ vanishes locally.

\subsubsection{}
Nevertheless we can show that its integral vanishes.
\begin{lem}\label{ha0}
For $n\in\nat$ we have
$$\int_{\Delta^{4n}} \hA(\nabla^{T\Delta^{4n}})_{4n}=0\ .$$
\end{lem}
\proof
Let $\sigma\in\Sigma^{4n+1}$ be an odd permutation. 
It acts by isometries on $\Delta^{4n}$ and therefore $$\sigma^*
\hA(\nabla^{T\Delta^{4n}})_{4n}= \hA(\nabla^{T\Delta^{4n}})_{4n}\ .$$
Since it changes the orientation we have 
$$\int_{\Delta^{4n}} \sigma^*
\hA(\nabla^{T\Delta^{4n}})_{4n}=-\int_{\Delta^{4n}}  
\hA(\nabla^{T\Delta^{4n}})_{4n}\ .$$
Thus $\int_{\Delta^{4n}} \hA(\nabla^{T\Delta^{4n}})_{4n}=0$.
\hB

An alternative argument\footnote{suggested by the referee} goes as follows. The construction of 
the manifold with corners $\Delta^{4n}$ as a subspace of $\R^{4n+1}$ in \ref{symy1} gives a picewise smooth metric
which is locally isometric to $S^k\times \R^{4n-k}$. For this metric we have $\hA(\nabla^{T\Delta^{4n}})_{4n}\equiv 0$.
We now deform this metric to a smooth one. Since $\hA(\nabla^{T\Delta^{4n}})_{4n}$ vanishes near the boundary for every admissible metric we see that 
$\int_{\Delta^{4n}} 
\hA(\nabla^{T\Delta^{4n}})_{4n}$
is deformation invariant. So this  integral also vanishes for the smooth admissible metric. 

\subsubsection{}

Let $\cE_{geom}$  be a geometric family  with closed fibers over a base $B$.
We consider the covering $\cU_0$ of $B$ consisting of one open set $U_o:=B$.
For all $p\in\nat\cup\{0\}$  the nerve $\bN[p]$ contains a single simplex
\index{$c(p)$}$o_p$. Let $c(p):=\sum_{i=0}^p i=\frac{p(p+1)}{2}$.

\begin{ddd}\label{uiuuiiufwefwe}
For $k\in\nat\cup\{0\}$ we define the chain
$\tilde z^k(\cE_{geom})\in \tilde G^k(B)$ to be the element which is represented by the chain
$Z(\cE_{geom})=(Z^0(\cE_{geom}),\dots,Z^k(\cE_{geom}))$ (w.r.t. $\cU_0$) such that  
$Z^p(\cE_{geom})(o_p)=(-1)^{c(p)}\cE_{geom}*\cDelta^p_{geom}$.
Then we have $$I_k(Z^p(\cE_{geom})(o_p))\cong
I_k(\cDelta^p_{geom})\cong I_k([p])\ ,$$ and
the isomorphisms
$$(-1)^p \partial_j Z^p(\cE_{geom})(o_p)\cong (-1)^j
Z^{p-1}(\cE_{geom})(o_p)$$ are induced by the identification
$\partial_j\cDelta^p_{geom}\cong (-1)^j\cDelta^{p-1}_{geom}$.

\index{resolution!geometric}The chain $\tilde z^k(\cE_{geom})$ is called the  geometric $k$-resolution of the geometric family
$\cE_{geom}$. Furthermore, by $\bz^k(\cE_{geom})\in G^k(B)$ we denote
the class
of $\tilde z^k(\cE_{geom})$. 
\end{ddd}

\subsection{The index form and the obstruction class $o$}
 
\subsubsection{}

In this subsection we study the relation between the obstruction against 
lifting the geometric $k$-resolution $\tilde z^k(\cE_{geom})\in \tilde G^k(B)$ 
to an element $z\in \tilde F^k_p(B)$ such that $\cF(z)=\tilde z^k(\cE_{geom})$
and the Chern character of the index of $\cE_{geom}$.
In a later Subsection \ref{intmean} we refine this relation
to the integral level. 

\subsubsection{}
\begin{lem}\label{high0}
For $1\le p$ we have $\Omega^k(\cE_{geom}*\cDelta^p_{geom})=0$.
 \end{lem}
\proof
Let $\pr:E\times \Delta^p \rightarrow E$, $\pi:E\rightarrow B$,  and $q:E\times \Delta^p\rightarrow B$ be the projections. 
Using Lemma \ref{ha0} we compute  
\begin{eqnarray*}
\Omega^{k}(\cE_{geom}*\cDelta^p_{geom})&=&\left[
 \int_{ (E\times \Delta^p)/B}
\hA(\nabla^{T^vq})
\ch(\nabla^{\pr^* W})
\right]_{k}\\
&=&
\left[\int_{E/B} 
\hA(\nabla^{T^v\pi})
\ch(\nabla^{W}) \right]_{k}  
\int_{\Delta^p} \hA(\nabla^{T\Delta^p})\\
&=&0\ .
\end{eqnarray*}
Here $(W,h^W,\nabla^W,z_W)$ denotes the (locally defined) twisting bundle of $\cE_{geom}$.
\hB

\subsubsection{}
 
Let $\tilde z^{k}(\cE_{geom})\in \tilde G^k(B)$ be the geometric $k$-resolution of the geometric family $\cE_{geom}$.
\begin{ddd}\label{tkres}
\index{resolution!tamed}A tamed $k$-resolution of $\cE_{geom}$ is a chain 
$\tilde z^{k}(\cE_{geom})_t\in \tilde G^k_t(B)$ such that $\cF(\tilde z^{k}(\cE_{geom})_t)=\tilde z^{k}(\cE_{geom})$.
\end{ddd}

\subsubsection{}
Let $0\le p\le k-1$ and assume that $\tilde z^{p}(\cE_{geom})_t$ is a
tamed $p$-resolution of $\cE_{geom}$. Thus there exists a covering
$\cU$ and tamed $p$-chain $Z_t$ (w.r.t. $\cU$) such that the geometric
$p$-chain $\cF(Z_t)$ (forget taming)
is equal to $\cF(Z)$ (reduce length), where $Z$ is the  geometric $k$-chain (w.r.t. $\cU$) induced by the
geometric $k$-resolution of $\cE_{geom}$ .  
We let $z\in \tilde F^k_p(B)$ be represented by the partially tamed
chain
$\tZ=(Z,Z_t) \in \tilde F^k_{p,\cU}(B)$ (w.r.t $\cU$).
 Let $\bz\in F^k_p(B)$ denote the class represented by $z$.
By $o^{p+1}(\bz)_{dR}\in H_{dR}^{p+1}(B)$
we denote the image of $o^{p+1} (\bz)$
under the natural homomorphism (see \ref{tzt6565s}) $$\check{H}^{p+1}(B,\underline{\Z}_B)\rightarrow \check{H}^{p+1}(B,\underline{\R}_B)\stackrel{\sim}{\rightarrow}
H_{dR}^{p+1}(B)\ .$$
\subsubsection{}
 \begin{prop}\label{rt}
We have  $(-1)^{p+1}o^{p+1}(\bz)_{dR}=\dR(\ch_{p+1}(\ind(\cE_{geom})))$.
\end{prop}
\proof
The form $\Omega^{p+1}(\cE_{geom})\in\cA_B^{p+1}(B)$ is closed.
If  $[\Omega^{p+1}(\cE_{geom})]\in H_{dR}^{p+1}(B)$ denotes the corresponding de Rham cohomology class,
then we have by the local index theorem \ref{loinde}
$$\dR(\ch_{p+1}(\ind(\cE_{geom})))=[\Omega^{p+1}(\cE_{geom})]\ .$$
Let $\bo\in\check{\bH}^{p+1}(B,\cA_B)$ be the image of $o^{p+1}(\bz)$
under  the natural  homomorphism (see \ref{someis})
 $$\check{H}^{p+1}(B,\underline{\Z}_B)\rightarrow 
\check{H}^{p+1}(B,\underline{\R}_B)\stackrel{\sim}{\rightarrow}
\check{\bH}^{p+1}(B,\cA_B)\ .$$
Furthermore, let
$[\bOmega]\in \check{\bH}^{p+1}(B,\cA_B)$ be the image of the class 
$[\Omega^{p+1}(\cE_{geom})]$ under the natural isomorphism (see again \ref{someis})
$$ H_{dR}^{p+1}(B)\stackrel{\sim}{\rightarrow}
\check{\bH}^{p+1}(B,\cA_B)\ .$$
It suffices to show that $(-1)^{p+1}\bo=[\bOmega]$.

The class $[\bOmega]$ is represented by the chain
$$\bOmega:=\prod_{x\in\bN[0]}\Omega^{p+1}(Z^0(x))\in \check{C}^0(\cU,\cA_B^{p+1})\subset \check{\bC}^{p+1}(\cU,\cA_B)\ .$$
The class $\bo$ is represented by the chain
$$\ind(\tZ):=\prod_{x\in\bN[p+1]} \ind_0(Z^{p+1}_{bt}(x))\in  \check{C}^{p+1}(\cU,\cA_B^{0})\subset \check{\bC}^{p+1}(\cU,\cA_B)\ .$$
We now define the chain 
$$\bbeta:=(\eta^{p}(Z^0_t),\dots, \eta^{0}(Z^{p}_t)) \in \check{\bC}^{p}(\cU,\cA_B)\ .$$
By Theorem \ref{etaprop} and Lemma \ref{high0}
we obtain 
\begin{equation}\label{careful3}\bd \eta =\bOmega + (-1)^p\ind(\tZ)\ .\end{equation}
Here are some details of the verification of this formula with emphasis on signs.
By Theorem \ref{etaprop} we have
\begin{eqnarray*}
d\eta^{p-q}(Z^q_t)&=&\Omega^{p-q+1}(Z^q_t)+\eta^{p-q+1}(\partial Z^q_t)\\
&=&\Omega^{p-q+1}(Z^q_t)+(-1)^q \delta \eta(Z^{q-1}_t)\ ,
\end{eqnarray*}
i.e. $$(-1)^q d\eta^{p-q}(Z^q_t) - \delta \eta(Z^{q-1}_t)=(-1)^q \Omega^{p-q+1}(Z^q_t)\ .$$
The right-hand side vanishes for all $q>0$, and for $q=0$ we get
$$d\eta^p(Z^0_t) = \Omega^{p+1}(Z^0_t)\ .$$
At the other end we have
$$-\delta \eta^0(Z^p_t)=-(-1)^{p+1}\eta^0(\partial Z^{p+1}_{bt})=(-1)^p \ind(Z^{p+1}_{bt})\ .$$

This proves $(-1)^{p+1}\bo=[\bOmega]$. \hB 

\subsubsection{}

Let $H^*_{dR}(B,\Z)$ denote the image of $\dR:H^*(B,\Z)\rightarrow H^*_{dR}(B)$. It is the lattice of classes with integral periods.

\begin{kor}\label{zcoho}
If  the geometric family $\cE_{geom}$ admits a tamed $k$-resolution,
then we have $\dR(\ch_{k+1}(\ind(\cE_{geom}))) \in H^{k+1}_{dR}(B,\Z)$. 
\end{kor}

\subsection{Classification of tamings - finiteness}

\subsubsection{}

We will say that $B$ is finite if it 
\index{$\ch^\R$}is homotopy equivalent to a finite $CW$-complex. If $B$ is finite, then the Chern character induces
an isomorphism
$$\ch^\R:K(B)_\R\stackrel{\sim}{\rightarrow} H_{dR}(B)\ ,$$
where $K(B)_\R=K(B)\otimes_\Z \R$.

\subsubsection{}
Let $p,k\in\nat\cup\{0\}$, $p< k$. 
Let\index{$R^p_k(B)$} $R^p_k(B)\subseteq \check{H}^{p}(B,\underline{\Z}_B)$ be the subset of
elements which can be written in the form $p(\tZ,\tW)$, where
$\tZ\in \tilde F^k_{p+1,\cU}(B)$ for some
covering $\cU$ of $B$, and $\tW$ is a zero bordism of $\cF(\tZ)\in \tilde F^k_{p,\cU}(B)$.
The set $R^p_k(B)$ is a subgroup since it is closed under the sum and $-p(\tZ,\tW)=p(\tZ^{op},\tW^{op})$.

\subsubsection{}
\index{$O^p_k(B)$}Let $O^p_k(B)\subseteq \check{H}^{p}(B,\underline{\Z}_B)$ 
denote the subgroup $o^p( F^k_{p-1}(B))$.
 
\subsubsection{}
Finally, let  \index{$\tilde O^p_k(B)$}$\tilde O^p_k(B)$ be the subgroup of $O^p_k(B)$ of elements of the form $o^p(\bz)$, where
$\bz\in F^k_{p-1}(B)$ is represented by  a partially tamed chain 
$\tZ=(Z(\cE_{geom}),Z(\cE_{geom})_t)$ (w.r.t some covering) such that
$Z(\cE_{geom})$ represents the geometric $k$-resolution of  a geometric family 
$\cE_{geom}$.

\subsubsection{}
\begin{lem}\label{finit1}
\begin{enumerate}
\item  We have $O^p_k(B)\subseteq R_k^p(B)$.
\item  If $B$ is finite, then the quotient $R_k^p(B)/O^p_k(B)$ is
  finite.
\item If $B$ is finite, then the quotient $R_k^p(B)/\tilde O^p_k(B)$ is finite.
\end{enumerate}
\end{lem}
\proof
Let $\cU$ be some covering of $B$ and $\tZ\in \tilde F^k_{p-1,\cU}(B)$.
Then we can consider $\tZ$ as a zero bordism of the empty family $\emptyset$.
Comparing the definitions of the obstructions $p$ and $o$ (see \ref{ww} and \ref{olem}) we see that
$p^{p-1}(\emptyset,\tZ)=o^{p-1}(\bz)$, where $\bz$ is the class
represented by $\tZ$.
This shows Assertion 1.

Assertion 2. follows from 3.
since by definition $\tilde O^p_k(B)\subseteq O^p_k(B)$.

We now show 3.
Let $R_k^p(B)_{dR}$, $\tilde O^p_k(B)_{dR}$  denote the images of $R_k^p(B)$, $\tilde O^p_k(B)$
in  $H_{dR}^{p}(B)$.
Since $B$ is finite the order of the torsion subgroup
of $\check{H}(B,\underline{\Z}_B)$ is finite.
It therefore suffices to show that $R_k^p(B)_{dR}/\tilde O^p_k(B)_{dR}$ is finite.
This assertion is a consequence of the following stronger assertion:
The quotient $H_{dR}^p(B,\Z)/\tilde O^p_k(B)_{dR}$ is finite.

Let $T^p\subseteq K(B)$ be the subgroup of elements $v$ such that
$\ch_q(v)=0$ for all $q < p$. Note that $K_p(B)\subseteq T^p$, but $T^p$ also contains all torsion elements of $K(B)$ so that this inclusion is in general proper.
Let $v=v^0+v^1$ be the decomposition of $v$ such that $v^i\in K^i(B)$.
Let $v^0\in T^p$ be equal to the class $[V]$ of the graded vector bundle $V$.
 We choose a hermitian metric and a metric connection on $V$ and thus
 obtain 
the geometric bundle $\bV$.
Let $\cE(\bV)_{geom}$ be the corresponding geometric family (see \ref{ghghghgh11}).
Then the form $\Omega(\cE(\bV)_{geom})$ represents $\dR(\ch(v^0))$.
Let $v^1\in T^p$ be represented by a map $F:B\rightarrow U(n)$ for some $n\in\nat$.
Let $\cE(F,*)_{geom}$ be an associated geometric family (see  \ref{eufiewfewfew}).
Again, the form $\Omega(\cE(F,*)_{geom})$ represents $\dR(\ch(v^1))$.
We define $\cF_{geom}:=\cE(\bV)_{geom}  + \cE(F,*)_{geom}$.

Assume that $q<p-1$ and that we have a tamed $q$-resolution $\bz\in F^k_{q}(B)$ of $\cF_{geom}$.
 Then by Lemma \ref{rt} the class 
$(-1)^{q+1}o^{q+1}(\bz)_{dR}$ 
is represented by $\Omega^{q+1}(\cF_{geom})$, and it is trivial,
since  $v\in T^p$.
We see that $o^{q+1}(\bz)$
is a torsion class.

Let $N\in \nat$ be the order of the torsion subgroup
of $H(B,\Z)$. Then
$o^{q+1}(N\bz )=0$.
Thus $N\cF_{geom}$ admits a tamed $q+1$-resolution.  

After finite induction we conclude that 
$N^p\cF_{geom}$ admits a tamed $p-1$-resolution.
Moreover, $\dR(\ch_p(N^pv))=(-1)^po^p(N^p\bz)_{dR}$.
We see that $\dR(\ch_p(N^p T^p))\subseteq \tilde O^p_k(B)_{dR}$.
In order to show that  $H_{dR}^p(B,\Z)/\tilde O^p_k(B)_{dR}$ is finite it therefore
suffices to see that $H_{dR}^p(B,\Z)/\dR(\ch_p(N^pT^p))$ is finite.
Indeed, $H_{dR}^p(B,\Z)/\dR(\ch_p(N^pT^p))$ is finite since $\ch^\R_p(T^p_\R)
=H_{dR}^p(B)$.
\hB

\subsubsection{}
\begin{lem}\label{ffiinn}
Assume that $B$ is  finite. 
Let $p\le k-1$.
Then $$\ker(\cF:F^k_{p+1}(B)\rightarrow F^k_{p}(B))$$ is a finite group.
\end{lem}
\proof
We fix a finite set of pairs $(\tZ_x^\prime,\tW_x)$, $x\in  R^p_k(B)/O^p_k(B)$, where 
$\tZ^\prime_x\in \tilde F^k_{p+1,\cU_x}(B)$, $\tW$ is a zero bordism of $\tZ_x:=\cF(\tZ_x^\prime)\in \tilde F^k_{p,\cU_x}(B)$,
and $p(\tZ^\prime_x,\tW_x)\in R^p_k(B)$ represents $x\in  R^p_k(B)/O^p_k(B)$.
 
Let $\tZ^\prime\in \tilde F^k_{p+1,\cU}(B)$ represent some element
in $\ker(\cF:F^k_{p+1}(B)\rightarrow F^k_{p}(B))$  such that 
$\tZ:=\cF(\tZ^\prime)\in \tilde F^k_{p,\cU}(B)$ admits a zero bordism
$\tW$.
If  $\tY \in \tilde F^k_{p-1,\cU}$ represents $\by\in F^k_{p-1}(B)$,
then $\tW_1:=\tW+\tY$ is again a zero bordism of $\tZ$, and 
we have $p^p(\tZ^\prime,\tW_1)=p^p(\tZ^\prime,\tW)+o^p(\by)$.
 
Let now  $p^p(\tZ^\prime,\tW)$ represent the class $x\in  R^p_k(B)/O^p_k(B)$.
Let $\tY\in  \tilde F^k_{p-1,\cU}(B) $ represent some element $\by\in F^k_{p-1}(B)$ such that
$o^p(\by) = p^p(\tZ^\prime+\tZ^{\prime,op}_x,\tW+\tW_x^{op})$
(after refining $\cU$ if necessary).
Then $p^p(\tZ^\prime+\tZ^{\prime,op}_x,\tW+\tW_x^{op}+\tY^{op})=0$,
and a modification of the taming of $\tW+\tW_x^{op}+\tY^{op}$
in codimension zero
admits an extension which is a zero bordism of
$\tZ^\prime+\tZ^{\prime,op}_x$ (again after refining $\cU$ if necessary).

We see that every element of $\ker(\cF)$ can be represented by some
$\tZ_x^\prime$, $x\in R^p_k(B)/O^p_k(B)$,  so that $\sharp(\ker(\cF:F^k_{p+1}(B)\rightarrow F^k_{p}(B)))\le \sharp(R^p_k(B)/O^p_k(B))$.
\hB 

\subsubsection{}

Lemma \ref{finit1} has the following consequence.
\begin{kor}\label{kor1}
If $B$ is  finite and $k\in\nat\cup\{0\}$, then
the group  $\ker(\cF^k_{-1}:F^{k+1}_k(B)\rightarrow G^{k+1}(B))$
is finite.  
\end{kor}
   

\subsection{The filtration of $K$-theory and tamed resolutions}\label{intmean}

\subsubsection{}
Let $\cE_{geom}$ be a geometric family over a base manifold $B$.
\begin{ddd}
\index{$R^k(\cE_{geom})$}\index{resolution!set of resolutions}Let  $R^k(\cE_{geom})\subseteq F^{k+1}_k(B)$ be the (possibly empty) set represented by tamed $k$-resolutions of $\cE_{geom}$.
\end{ddd}
If $B$ is finite, then by Corollary \ref{kor1} the set $R^k(\cE_{geom})$ is finite. If $\bz\in R^k(\cE_{geom})$, then we have a class
$o^{k+1}(\bz)\in \check{H}^{k+1}(B,\underline{\Z}_B)$ which is the obstruction against prolonging the resolution.

\subsubsection{}
Recall from Subsection \ref{edele} the following notation: 
$K^*_{k+1}(B)$ is the $k+1$`th step of the filtration of $K$-theory, and if $\psi\in K_{k+1}^*(B)$, then $\check{\bo}^{k+1}(\psi)\subseteq
\check{H}^{k+1}(B,\underline{\Z}_B)$ denotes the obstruction set.
\begin{theorem}\label{compara}
\begin{enumerate}
\item
The family $\cE_{geom}$ admits a tamed $k$-resolution if and only if we have 
$\ind(\cE_{geom})\in K_{k+1}^*(B)$.
\item If $\ind(\cE_{geom})\in K_{k+1}^*(B)$, then
we have  the equality of sets $$o^{k+1}(R^k(\cE_{geom}))=(-1)^{c(k+1)}\check{\bo}^{k+1}(\ind(\cE_{geom}))\ ,$$
where $c(p):=\frac{p(p+1)}{2}$.
\end{enumerate}
\end{theorem}
The proof of this theorem almost occupies the remainder of the present subsection.
\subsubsection{}
Recall that for $j\in [n]$ we have fixed isomorphisms $\partial_j
\cDelta^n_{geom}\cong (-1)^j \cDelta^{n-1}_{geom}$ in the sense of \ref{langg}.
In particular, we
have fixed the embedding of the boundary faces
$\partial_{j*}:\Delta^{n-1}\rightarrow \Delta^n$.

\subsubsection{}  
Let $\cU$ be some covering of $B$ by open subsets.
\begin{ddd}\label{chdalter}
\index{resolution!local}A local $k$-resolution of $\cE_{geom}$ (w.r.t. $\cU$) consists of
\index{$X$}\index{$(X^0,\dots, X^k)$}a $k+1$-tuple  $X:=(X^0,\dots, X^k)$, where for $0\le p\le k$ the object
$X^p$ associates to each $x\in \bN[p]$ a tamed lift
$X(x)$ of the geometric family $\pr^*_x\cE_{geom}$, where
$\pr_x^*:U_x\times \Delta^p\rightarrow U_x$ is the projection onto the
first factor.
We require that for all $j\in [p]$ we have
$(\id_{U_x}\times \partial_{j*})^* X^p(x) =  X^{p-1}(\partial_{j}^* x)_{|U_x\times \Delta^{p-1}}$ under the canonical identification of the underlying geometric families.
\end{ddd}

\subsubsection{} 

To give a local $k$-resolution of $\cE_{geom}$ is equivalent to give the datum 
$\bK=(K_x)_{x\in\bN[p],p\le k}$ as considered in \ref{indf3}.   It was used 
 to define
the obstruction class $\check{o}^{k+1}(\cE_{geom},\bK)\in \check{H}^{k+1}(B,\Z)$ in the realm of families of Dirac operators. This obstruction class is represented by a cocycle
$u^{k+1}(X):=\check{c}^{k+1}(\cE_{geom},\bK)\in \check{C}^{k+1}(\cU,\underline{\Z}_B)$ defined in (\ref{coooa}).
 By  Proposition \ref{zzuduwedw} the set $\check{\bo}^{k+1}(\ind(\cE_{geom}))$ is the
set of cohomology classes 
$[u^{k+1}(X)]$, where $X$ runs over all local $k$-resolutions of
$\cE_{geom}$ (for varying $\cU$).

\subsubsection{}
\index{bordism}The idea of the proof is to introduce a relation $\sim$ (called bordism)  between tamed $k$-resolutions $\tZ$ and local $k$-resolutions $X$ such that
we have the following assertions.
\begin{enumerate}
\item For every tamed $k$-resolution $\tZ$ there exists a local  $k$-resolution $X$
such that $\tZ\sim X $.
\item For every local $k$-resolution $X$ there exists tamed $k$-resolution $\tZ$
such that $\tZ\sim X$.
\item If $\tZ\sim X$, then
$(-1)^{c(k+1)}u^{k+1}(X)=\ind(\tZ)$ (see \ref{ggg45673}).
\end{enumerate}
In view of the definition of $o^{k+1}(R^k(\cE_{geom}))$ as the set of classes
$[\ind(\tZ)]$ for tamed $k$-resolutions $\tZ$ of $\cE_{geom}$
it is clear that these three statements imply the theorem.

For clarity let us point out the difference between a tamed resolution and a local tamed resolution.  Let $x\in \bN[p]$ an element of the nerve.
The corresponding constituent of a tamed resolution 
$\tZ$ is the geometric  $(\cE_{geom})_{|U_x}*\cDelta^p$
with choice of a taming.

The constituent of  a local tamed resolution is the choice of a taming of the geometric family $\pr_x^*\cE_{geom}$ over $U_x\times \Delta^p$, where $\pr_x:U_x\times \Delta^p\to U_x$ is the projection.

In the first case the simplex becomes a factor of the fibre, while in the second case it is considered as a factor of the base.

\subsubsection{}

In order to define the notion of bordism we must extend the notion of taming 
to a certain class of perturbations of Dirac operators.

Let $\cM_{geom}$ be a closed reduced geometric manifold with Dirac bundle $\cV$, and let
$N$ be a compact manifold with corners equipped with a Riemannian metric $g^N$, orientation,
spin structure, and a reduced face decomposition. On the one hand we can form the geometric manifold $\cM_{geom}*N$.
The underlying manifold with corners is 
$M\times N$. 
On the other hand we can consider the trivial geometric family
$\cE_{geom}:=\cM_{geom}\times N$ over $N$. Its total space has the same
underlying manifold with corners.
 
\subsubsection{}\label{dqed}

Let $Q:N\rightarrow \End(C^\infty(M,V))$ be a smooth family of smoothing operators which provide a pre-taming $\cM_{t}$.
We assume that we can extend $Q$ to a smooth family $\bar Q$ defined
on $\bar N$ such that it is constant in the normal directions on the
cylinders $U_j$ over the faces 
 $j\in I_k(N)$, $k\in\nat$.
In this situation we introduce the notation 
$\cM_t*N$. The associated 
perturbed Dirac operator
$D(\cM_{t}* N)$ is defined by
$$D(\cM_{t}* N):=D(\cM_{geom}* N)+ \bL_{M}^{M\times \bar N}(
\bar Q)\ .$$
Thus $D(\cM_{t} *  N)$ is a bounded perturbation of
$D(\cM_{geom}* N)$, where the perturbation is local
w.r.t. the $\bar N$-variable and smoothing in the $M$-direction.

\subsubsection{}

\index{adiabatic limit}We call the situation where we replace the metric $g^N$ by $\epsilon^{-2}g^N$ for sufficiently small $\epsilon>0$ the adiabatic limit.
\begin{prop}\label{uduiduied}
\begin{enumerate}
\item If $Q$ defines a taming of $\cE_{geom}$, then in the adiabatic limit  the operator
$D(\cM_{t}* N)$ becomes invertible.
\item If $Q$ defines a taming of $(\cE_{geom})_{|\partial N}$, then in the adiabatic limit 
$D(\cM_{t}*N)$ is a Fredholm operator. 
\item Assume that $k+\dim(M)$ is even. Under the assumption of 2. and if $(N,\partial N)$ is
homotopy equivalent to $(D^{k},\partial D^k)$,
the integer $\ind(D(\cM_{t} * N))$ coincides with
the element $\ind(D(\cM_{geom}),Q_{|\partial N})\in K^{*}(N,\partial N)$ under the isomorphism $K^{*}(N,\partial N)\cong K^*(D^k,\partial D^k)\cong \Z$.
Here $*:=[k]\in\Z_2$.
\end{enumerate}
\end{prop}
\proof
In order to prove assertion 1. we consider  the square 
$D(\cM_{t}* N)^2$. To write out a formula for this operator completely we would have to introduce more notation like a lift $L_{\bar N}^{M\times \bar N}$. But in order to see the assertion it suffices to use the rough structure of $D(\cM_{t}* N)^2$.
For simplicity we will refrain from using the symbols $L_{M}^{M\times \bar N}$ etc. completely. We use the product structure
$M\times N$ in order to let $D(\cM_t)$ and $D(N)$ (the Dirac operator of the spin manifold $(N,g^N)$ on the spinor bundle) act on the sections of the Dirac bundle on $\cM_{geom}*N$ (in the precise formula we would have go through the $2\times 2$-matrix yoga as in \ref{ioioewf}).
The rough structure of the square after rescaling $g^N\mapsto \epsilon^{-2}g^N$ is
\begin{eqnarray*}
D(\cM_{t}* N)^2&=&
(D(\cM_{geom}) + Q +  \epsilon D(N))^2\\
&=&(D(\cM_{geom})+Q)^2+\epsilon^2  D(N)^2 +\epsilon \{ D(\cM_{geom}) + Q , D(N)\}\\
&=&D(\cM_t)^2+\epsilon^2  D(N)^2 +\epsilon \{Q , D(N)\}\ ,
\end{eqnarray*}
where we have used $\{D(\cM_{geom}),D(N)\}=0$.
The anticommutator $\{Q , D(N)\}$ is bounded since the derivatives of $D(N)$ are used to differentiate $Q$ in the $N$-directions. 
Since $D(\cM_t)^2$ is positive and $\epsilon^2  D(N)^2$ is non-negative we see that
$D(\cM_{t}* N)^2$ becomes positive for sufficiently small $\epsilon>0$.

For Assertion 2. we use the Assertion 1. on the cylinders over the boundary faces and apply the parametrix construction as
in the proof of Lemma \ref{tamm}.

Assertion 3. is just one analytic way
to invert the suspension map and to provide the isomorphism $K^{*}(N,\partial N)\cong K^*(D^k,\partial D^k)\cong \Z$. In terms of Kasparovs $KK$-theory (we use the unbounded picture \cite{baju}) one could argue as follows.

We consider the $C^*$-algebra $C(N,\partial N)$ of continuous functions on $N$ vanishing on the boundary. Let $V$ denote the Dirac bundle of $\cM_{geom}$.
We define the Hilbert $C(N,\partial N)$-modul
$E:=C(N,\partial N)\otimes L^2(\bar M;\bar V)$ (completion of the tensor product is understood implicitly).
The family of operators $D(\cM_t)$ parametrized by $N$ 
can be considered as an unbounded operator on $E$ and gives rise to an unbounded
Kasparov module $\{D(\cM_t)\}$ over the pair of $C^*$-algebras $(\C,C(N,\partial N))$ (for simplicity of presentation we are going to supress the distinction between the even and odd cases). Let $[D(\cM_t)]\in KK(\C,C(N,\partial N))$ be the class represented by $\{D(\cM_t)\}$. It corresponds to
$\ind(D(\cM_{geom},Q_{|\partial N}))\in K(N,\partial N)$ under the identification
$KK(\C,C(N,\partial N)) \cong K(N,\partial N)$.

The Dirac operator
$D(N)$ also gives rise to a Kasparov module $\{D(N)\}$ over $(C(N,\partial N),\C)$ with underlying Hilbert space $L^2(N,S(N))$ and action of  $C(N,\partial N)$ by multiplications. The class $[D(N)]\in KK^*(C(N,\partial N),\C)\cong K_*(N,\partial N)$ is the $K$-homology fundamental class of $(N,\partial N)$ associated to the $K$-orientation of $N$ determined by the choice of the $spin$-structure.

On the one hand, the  operator $D(\cM_t*N)$ represents the Kasparov product
$$[D(\cM_t)]\otimes_{C(N,\partial N)} [D(N)]\in KK(\C,\C)\cong K(\C)\cong \Z\ .$$
On the other hand, if $(N,\partial N)\cong (D^{k},D^{k-1})$, then the inverse of the suspension isomorphism
$\Z\to K(N,\partial N)$ is given by the pairing with the $K$-homology fundamental class of $(N,\partial N)$, i.e. also by the Kasparov product by $[D(N)]$. This shows Assertion 3. \hB

\subsubsection{}\label{eed4}

Now let $Q$ define a taming $\cM_t$.
Let $I$ be the unit interval with two boundary faces.
We consider the geometric manifold $\cM_{geom}*(N\times I)$.
 Let $t\in I$ be the coordinate, and let $\rho\in C^\infty(\R)$ be such that
$\rho(t)= 0$ for $t<1/4$ and $\rho(t)=1$ for $t>3/4$. We define the family of smoothing operators $R:N\times I\to \End(C^\infty(M,V))$ by $R(n,t):=\rho(t)Q(n)$ which goes into the definition of $D(\cM_{t} *(N\times I))$.

\index{special!taming}\index{$\cW_{st}$}We now introduce the notion of a special taming of the operator
$D(\cM_{t} *(N\times I))$.
After rescaling the metric $g^N$ we can assume by Proposition \ref{uduiduied},1. that 
the reductions  of this operator
to the extensions of the faces of $\cM_{geom}* (N\times \{1\})$ are
already invertible.

A special pre-taming is now a pre-taming $\{W_i|i\in I_*(\cM_{geom} * (N\times I))\}$
of $\cM_{geom} * (N\times I)$ such that $W_i=0$ 
for all faces $i$ which are contained in the face $M\times N\times \{1\}$ (this vanishing condition makes the pre-taming "special"). 
\index{special!pre-taming}
\index{special!boundary pre-taming}It is a special taming (subscript $(.)_{st}$), if 
$$D((\cM_{t}* (N\times I))_{st}):=D(\cM_{t}*(N\times
I))+\sum_i \rho_i W_i$$
is invertible.
A special boundary pre-taming is a special pre-taming where
$W_i=0$ also for the codimension zero face.
It is a special boundary taming (denote by the subscript $(.)_{sbt}$) if all boundary reductions 
of $D((\cM_{t} *(N\times I))_{st})$
are invertible. In this case $D((\cM_{t}* (N\times I)_{sbt})$ is a Fredholm operator. This can be seen by repeating the proof of Proposition \ref{tamm} 
with the corresponding modifications in order to include the presence of $R$.

\subsubsection{}\label{fuweufiewiufifoioi4}

All these notions can be extended to families parameterized by some auxiliary space. We can extend the theory of $\eta$-forms and Theorem \ref{etaprop}
to families operators of the form $D(\cM_{t}*N)$ with essentially the same proofs. 
If we form the rescaled super-connection, then we deal with the terms
coming from $Q$ in the same way as with the terms coming from usual tamings, i.e. we
insert a cut-off function, which switches off these terms for small
scaling parameters.

\subsubsection{}

We can now introduce the notion of a bordism between a tamed $k$-resolution and a local $k$-resolution.

Let $\tZ:=(Z,Z_t)$ and $X:=(X^0,\dots,X^k)$ be a tamed and a local $k$-resolution represented
with respect to a covering $\cU$. For $x\in \bN[p]$ and $b\in U_x$ 
the fiber $(-1)^{c(p)}(\cE_{geom})_b$, $\Delta^p$, and the restriction of
$(-1)^{c(p)}K_x$ to $\{b\}\times\Delta^p$ play the roles of $\cM_{geom}$, $N$, and $Q$
above. After a homotopy of $\bK$ we can assume that each $K_x$ extends
smoothly to the extension $U_x\times \overline{\Delta^p}$ so that this extension is independent of the normal variables of the cylinders over the faces of $\Delta^p$. We now consider the cylinder
$W:=(W^0,\dots W^k)$ with $W^p(x)=Z^p(x)* I\cong (-1)^{c(p)}
(\cE_{geom})_{|U_x}*(\Delta^p\times I)$.
Then a bordism between $\tZ$ and $X$ will be a given by a special
pre-taming $W_{st}:=(W_{st}^0,\dots,W^k_{st})$
such that for each $x\in \bN[p]$ and $b\in U_x$ 
$(W_{st}(x))_b$ is a special taming  of $(-1)^{c(p)}(\cE_{t})_{b}*( \Delta^p\times I)$
in the sense above, and 
 the restriction of the special taming to the faces of $Z^p(x)\times \{0\}$
is the taming $Z^p_t(x)$.

\begin{ddd}
We say that $\tZ$ and $X$ are bordant and write $\tZ\sim X$ iff
the cylinder $W$ admits a special taming as described above.
\end{ddd} 

\subsubsection{}

\begin{lem}\label{dee1}
If the tamed $k$-resolution $\tZ$ and the local $k$-resolution $X$ are bordant,
then we have $\ind(\tZ)=(-1)^{c(k+1)}u^{k+1}(X)$.
\end{lem}
\proof
We assume without loss of generality that $Z^{k+1}$ has even-dimensional fibers. Otherwise, both cocycles vanish. Consider $x\in\bN[k+1]$.
Then all faces of $W^{k+1}(x)$ are 
(specially) boundary tamed.  
The sum over the boundary components of $W^{k+1}(x)$ of the indices of the associated Fredholm operators vanishes. In order to see this we apply Theorem  \ref{etaprop}, 2., in its generalization to the present case.
If we sum up over all faces, then the contribution of the $\eta$-invariants cancels out. The sum of the contributions of the local index forms over the faces vanishes by Stokes theorem.

Note that all faces of $W^{k+1}(x)$ are in fact specially tamed with the exception of  
$Z^{k+1}(x)\times \{i\}$, $i=0,1$.
The index of the Fredholm operator associated to this face for $i=0$ is
equal to $\ind(\tZ)(x)$, while the index of the Fredholm operator associated to this face for $i=1$ is equal to $-(-1)^{c(k+1)}u^{k+1}(X)(x)$.
This implies the claim. 
\hB

\subsubsection{}

\begin{lem}\label{locexis}
Given a tamed $k$-resolution $\tZ$, there is a local $k$-resolution
$X$, which is bordant to $\tZ$. Vice versa, for any
local $k$-resolution $X$ there exists a bordant tamed
$k$-resolution $\tZ$.
\end{lem}
\proof
Given $\tZ$
we construct the local $k$-resolution $X$ inductively.
Assume that we already have constructed a local $k-1$-resolution
$X^\prime$ such that $\tZ^\prime$ is bordant to $X^\prime$, where $\tZ^\prime$ is the
tamed $k-1$-resolution induced by $Z$. We must define
the family of operators $K_x$ for $x\in\bN[k]$.

Let $x\in \bN[k]$. Then $K_x$ is already defined on $U_x\times \partial \Delta^{k+1}$ by the compatibility conditions.
Assume first that the fibers of $Z^k$ are even-dimensional.
All boundary  faces of $W^k(x)$ except the 
face $\{1\}\times Z^k(x)$ are already specially tamed. 
It follows by Lemma \ref{dee1} that
$u^k(X^\prime)=0$. Therefore, we can extend
$K_x$ to all of $U_x\times  \Delta^{k+1}$.
Since $W^k$ has odd-dimensional fibers we can now extend the special taming to $W^k$.

If the fibers of $Z^k$ are odd-dimensional, then we can always extend
$K_x$.
Each choice of such an extension completes a special boundary taming $W^k(x)_{sbt}$.
There is a unique choice up to homotopy such that $\ind(W^k(x)_{sbt})=0$. 
Taking this choice we can now extend the special taming of the cylinder.

In a similar manner we construct $\tZ$ given $X$.
\hB

This finishes the proof of the theorem \hB

\subsubsection{}
 
It now follows from Theorem  \ref{compara} and Subsection \ref{choob}
that we have the following equalities of sets.
\begin{kor}\label{kerv}
Under the natural isomorphism $H^*(B,\Z)\cong \check{H}^*(B,\underline{\Z}_B)$ we have 
$$
- (k-1)! \bo^{2k}(\ind(\cE_{geom})) \cong  \{c_{2k}(\ind(\cE_{geom}))\}$$
if $\ind(\cE_{geom})\in K^0_{2k}(B)$, and
$$ 
-k!  \bo^{2k+1}(\ind(\cE_{geom}))\cong \{c_{2k+1}(\ind(\cE_{geom}))\}$$
if $\ind(\cE_{geom})\in K^1_{2k+1}(B)$.
\end{kor}

\newpage

\part{Deligne cohomology valued index theory}\label{dd212}

\section{Deligne cohomology valued index theory}\label{klo98}

\subsection{Review of Deligne cohomology and Cheeger-Simons differential characters}\label{delche}

\subsubsection{}
\index{$\cK(k,R)_B$}If $R\subset\R$ is some  discrete subring and $k\in\nat\cup\{0\}$, then we let 
$\cK(k,R)_B$ be the complex of sheaves
$$0\rightarrow \underline{R}_B \stackrel{i}{\rightarrow}
\cA^0_B\stackrel{d}{\rightarrow} \dots \stackrel{d}{\rightarrow}
\cA_B^k\rightarrow 0\ .$$
Here the constant sheaf $ \underline{R}_B$ sits in degree $-1$.
We prefer his convention because of the fact that the sheaf $\cA_B^p$ of forms of degree $p$ sits in degree $p$. This helps to avoid confusing shifts at some places later, but unfortunately induces the degree-shift in Definition \ref{delcoh}.
 
\subsubsection{}
Recall from Subsection \ref{grus}, that if $\cK$ is a complex of sheaves, then $\check{\bH}(B,\cK)$ denotes
the hyper-cohomology of $\cK$. 

\begin{ddd}\label{delcoh}
\index{Deligne cohomology!$R$-valued}\index{$H_{Del}^k(B,R)$}For $k\in\nat\cup\{0\}$ the degree-$k$ Deligne cohomology of $B$ is defined by
$$H_{Del}^k(B,R):=\check{\bH}^{k-1}(B,\cK(k-1,R)_B)\ .$$
\end{ddd}
If $R=\Z$, then we write $H_{Del}^k(B):=H_{Del}^k(B,\Z)$.
\index{Deligne cohomology!integral}
\index{$H_{Del}^k(B)$}
\subsubsection{}  We refer to the book of Brylinski, \cite{brylinski93} for
an introduction to Deligne cohomology. What we define here is usually called smooth
Deligne cohomology as opposed to its algebraic geometric counterpart.
Our definition can be compared  with  \cite{brylinski93}, Def. 1.5.1. Note that
we work with the group $\Z$ instead of $\Z(p):=(2\pi\imath)^p\Z$,
and our differential forms are real valued as opposed to complex
valued forms in the reference.

The picture of  Deligne cohomology which is presented here is very
incomplete. So we do not touch the ring structure, integration over
fibers,
and the various applications and geometric interpretations.

\subsubsection{}
If $f:B^\prime\rightarrow B$ is a smooth map, then the differential of
$f$ induces a morphisms of sheaves of complexes
$f^*\cK(k-1,R)_B\rightarrow \cK(k-1,R)_{B^\prime}$.
This transformation gives functorial a pull-back
$f^*:H^k_{Del}(B,R)\rightarrow H^k_{Del}(B^\prime,R)$.

\subsubsection{}

If $x\in H_{Del}^k(B,R)$ is represented by the chain $\bc\in \check{\bC}^{k-1}(B,\cK(k-1,R)_B)$,
$\bc=(c^{p,q})_{p+q=k-1}$, $c^{p,q}\in \check{C}^p(B,\cK(k-1,R)_B^q)$, then
$\delta d c^{0,k-1}= 0$. Therefore there is a closed form $R^{\bc}\in \cA_B^k(B)$ which restricts
to $dc^{0,k-1}$. It only depends on the class $x$ and not on the
representative $\bc$. Thus we can write $R^x:=R^\bc$.

\begin{ddd}\label{delcoh2}
\index{curvature homomorphism}\index{$R^x$}The curvature homomorphism 
$R:H_{Del}^k(B,R)\rightarrow \cA_B^{k}(B)$
associates to $x\in H_{Del}^k(B,R)$ the closed form
$R^x\in\cA^{k}_B(B,d=0)$. 
\end{ddd}

\subsubsection{}

Let $x$ be represented by $\bc$ as above. Then
$c^{k,-1}\in \check{C}^{k}(B,\underline{R}_B)$ is a cocycle.
Its class $[c^{k,-1}]\in \check{H}^{k}(B,\underline{R}_B)$ only depends on $x$.

\begin{ddd}
\index{$\bv$}We define the homomorphism $\DD:H^k_{Del}(B,R)\rightarrow \check{H}^{k}(B,\underline{R}_B)$
such that $\DD(x):=[c^{k,-1}]$ if $x$ is represented by $\bc=(c^{p,q})_{p+q=k-1}$.
\end{ddd}

\subsubsection{}
From the definitions immediately follows (compare \ref{someis} and \ref{tzt6565s}): 
\begin{kor}
\begin{enumerate}
\item For $x\in H^k_{Del}(B,R)$ we have
$$\dR(\DD(x))=[R^x]\in H^{k}_{dR}(B)$$
\item For $x\in H^k_{Del}(B,R)$ we have 
$[R^x]\in H^k_{dR}(B,R)$
\end{enumerate}
\end{kor}
Here $H^k_{dR}(B,R)$ denotes the subgroup of classes in $H^k_{dR}(B)$
with periods in $R$.

\subsubsection{}
\begin{ddd}\label{pkgr}
\index{$P^k(B,R)$}We define the group
$$P^k(B,R):=\check{H}^k(B,\underline{R}_B)\times_{H^k_{dR}(B)} \cA_B^k(B,d=0)\ .$$
\end{ddd}
By $\cA^k_B(B,d=0,R)$ we denote the space of closed $k$-forms with
periods in $R$, i.e. the image of $\pr_2:P^k(B,R)\rightarrow \cA_B^k(B,d=0)$.
It is instructive to note the following natural exact sequences:
(compare \cite{brylinski93}, Thm. 1.5.3)
\begin{equation}\label{sq12}0\rightarrow \check{H}^{k-1}(B,\underline{\R/R}_B)\rightarrow H^k_{Del}(B,R)
\stackrel{R}{\rightarrow} \cA_B^k(B,d=0,R)\rightarrow 0\end{equation}
 \begin{equation}\label{sq121}0\rightarrow \cA_B^{k-1}(B)/\cA^{k-1}_{B}(B,d=0,R) \stackrel{a}{\rightarrow} H^k_{Del}(B,R)\stackrel{\DD}{\rightarrow} \check{H}^k(B,\underline{R}_B)\rightarrow 0\end{equation}
and
\begin{equation} \label{edrr}  0\rightarrow \check{H}^{k-1}(B,\underline{\R}_B)/ \check{H}^{k-1}(B,\underline{R}_B)\rightarrow H^k_{Del}(B,R)
\stackrel{(\DD,R)}{\rightarrow} P^k(B,R)\rightarrow 0
            \ . \end{equation}
The map $a:\cA_B^{k-1}(B)/\cA^{k-1}_{B}(B,d=0,R) \stackrel{a}{\rightarrow} H^k_{Del}(B,R)$ is induced by the obvious inclusion of complexes of sheaves
$$\cA^{k-1}_B\to \cK(k-1,R)_{B}\ ,$$
where we  view $\cA^{k-1}_B$ as a complex with only one non-trivial entry in degree $k-1$. This inclusion induces a map of hyper cohomology groups
$\cA_B^{k-1}(B)\to H^k_{Del}(B,R)$, and one checks that it factors over
the quotient $\cA_B^{k-1}(B)/\cA^{k-1}_{B}(B,d=0,R) $.

\subsubsection{}

We now give the definition of the group of Cheeger-Simons differential characters $\hat H^{k}(B,\R/R)$ which eventually turns out to be isomorphic to $H^{k+1}_{Del}(B,R)$
(see Cheeger and Simons \cite{cheegersimons83}
and  \cite{brylinski93}, Sec. 1.5, for an introduction).
\index{$Z_k(B)$}\index{$C_k(B)$}\index{smooth singular chains}Let $Z_k(B)\subset C_k(B)$ denote the group of smooth singular $k$-chains in $B$ and its subgroup of cycles.
\begin{ddd}\index{differential character}
\index{Cheeger-Simons}The group of Cheeger-Simons differential characters
\index{$\hat H^{k}(B,\R/R)$}$\hat H^{k}(B,\R/R)$ is defined by
\begin{eqnarray*}\lefteqn{\hat H^{k}(B,\R/R)}&&\\
&:=&\{\phi\in\Hom(Z_k(B),\R/R)\:|\: \exists \omega\in\cA_B^{k+1}(B)\:\forall c\in C_{k+1}(B) | \phi(\partial c)=[ \int_{c}\omega]_{\R/R}\}\ .\end{eqnarray*}
\end{ddd}

\subsubsection{}\label{holo87}
There is a natural isomorphism (see e.g. Gajer, \cite{gajer99} in the
case $R=\Z$)
$$H: H^{k+1}_{Del}(B,R)\stackrel{\sim}{\rightarrow} \hat H^{k}(B,\R/R) .$$
One way to define $H$ is as follows. 

Let $x\in H^{k+1}_{Del}(B,R)$.
We must construct the character $H(x):Z_k(B)\rightarrow \R/R$.
Let $z\in Z_k(B)$.
Then we can find an open neighborhood $i:U\rightarrow B$ of the trace of $z$ which is
homotopy equivalent to a $k$-dimensional $CW$-complex.
Here $i$ denotes the inclusion. We have $\DD(i^*x)=i^*\DD(x)=0$ so
that  by (\ref{sq121}) there exists a class $ [\omega]\in
\cA_B^{k}(U)/\cA^{k}_{B}(U,d=0,R)$
such that $a([\omega])=i^*x$.
Then $$H(x)(z):=[\int_z \omega]_{\R/R}\ .$$

\index{$H$}\index{holonomy}
We call $H(x)(z)\in \R/R$ the holonomy of $x$ along $z$.

\subsubsection{}
\index{homotopy formula}Let $x\in H^k_{Del}(I\times B,R)$ and $f_i:B\rightarrow I\times B$
be induced by the inclusion of the endpoints of the interval.
Then we have the following homotopy formula.
\begin{lem}\label{homotopiedel} We have
$$f_1^*x-f_0^*x=(-1)^ka(\int_{I\times B/B} R^x)\ .$$
In particular, if $R^x=0$, then
$$f_1^*x=f_0^*x\ .$$ 
\end{lem}
\proof
We will show that both sides of the equation  have the same holonomy. 
Let $z\in Z_{k-1}(B)$ be a smooth cycle.
Then we form the chain  $I\times z\in Z_{k}(I\times B)$.
We have $\partial (I\times z)=(-1)^{k-1}(0\times z-1\times z)$.
Therefore,
\begin{eqnarray*}
H( f_1^*x-f_0^*x)(z)&=&(-1)^k H(x)(\partial (I\times z))\\
&=&(-1)^k[\int_{I\times z} R^x]_{\R/R} \\
&=&H((-1)^ka(\int_{I\times B/B} R^x))(z)\ .
\end{eqnarray*}
 \hB

\subsection{The lift of obstruction set to  Deligne cohomology}

\subsubsection{}

We fix $k\in\nat\cup\{0\}$ and let $m\in \nat\cup\{0\}$ be defined by $k=2m$ or $k=2m-1$.
Let $\cE_{geom}$ be a geometric family
such that $\ind(\cE_{geom})\in K_{k}^*(B)$.
Then we have the obstruction set (see Definition \ref{obset} and \ref{zt1234})
$\check{\bo}^k(\ind(\cE_{geom}))\subseteq \check{H}^k(B,\underline{\Z}_B)$ and
the local index form $\Omega^k(\cE_{geom})\in \cA_{B}^k(B,d=0)$.
If $k>0$ and $x\in \check{\bo}^k(\ind(\cE_{geom}))$, then from the identity
$$(-1)^{m-1}
(m-1)!\ch_{k}(\ind(\cE_{geom}))=c_k^\R(\ind(\cE_{geom}))+\mbox{\em
  \small Polynomial
  in lower Chern classes}$$ 
and Lemma \ref{kerv} we conclude
$$\dR(x)=(-1)^m[\Omega^k(\cE_{geom})]\in H_{dR}^k(B)\ .$$
In fact, this also holds for $k=0$.
Therefore, the pair $((-1)^m x,\Omega^k(\cE_{geom}))$ defines an element
of the group $P^k(B)$ defined in \ref{pkgr}.

\subsubsection{}

In view of the exact sequence (\ref{edrr}) we may ask for a natural lift
of the pair $((-1)^m x,\Omega^k(\cE_{geom}))\in P^k(B)$ to
an element of $H^k_{Del}(B)$.
In the present subsection we define a natural set $\ind_{Del}^k(\cE_{geom},x)\subset H^k_{Del}(B)$ of such lifts. In Proposition \ref{denominators}
we obtain more information about the nature of this set.
 
\subsubsection{}
\begin{ddd}
\index{$\tilde S^{k-1}(\cE_{geom},x)$}\index{$\tilde S^{k-1}(\cE_{geom})$}For any geometric family $\cE_{geom}$ such that $\ind(\cE_{geom})\in K_{k}^*(B)$ and $x\in \check{\bo}^k(\ind(\cE_{geom}))$ we define
the set $\tilde S^{k-1}(\cE_{geom},x)\subseteq \tilde F^{k}_{k-1}(B)$
of those partially tamed resolutions $z\in \tilde F^{k}_{k-1}(B)$ of $\cE_{geom}$ 
which satisfy $o^k(\bz)=(-1)^{m}x$, where $\bz$ is the class represented by $z$. 
We furthermore define
$$\tilde S^{k-1}(\cE_{geom}):=\bigcup_{x\in\check{\bo}^k(\ind(\cE_{geom})) } \tilde S^{k-1}(\cE_{geom},x)\ .$$
\end{ddd}
It follows from Theorem \ref{compara}, that $\tilde S^{k-1}(\cE_{geom},x)$ is
not empty. Moreover,
$\tilde S^{k-1}(\cE_{geom})$ maps
onto $R^k(\cE_{geom})$ under the quotient map
$\tilde  F^{k}_{k-1}(B)\rightarrow  F^{k}_{k-1}(B)\ .$

\subsubsection{}\label{c824}

Let $z\in \tilde S^{k-1}(\cE_{geom},x)$
be represented by $\tZ=(Z,Z_t)\in F^{k}_{k-1,\cU}(B)$.
Let us write a chain $\bc\in \check{\bC}^{k-1}(B,\cK(k-1,\Z)_B)$ in the form $\bc=(c^{k-1},\dots,c^0,c^{-1})$, where $c^p\in \check{C}^{k-p-1}(B,\cA^p_B)$
for $p\ge 0$ and $c^{-1}\in \check{C}^{k}(B,\underline{\Z}_B)$.
\begin{ddd} We define \index{$\bbeta(z)$}
$$\bbeta(z)=(\eta^{k-1}(Z^0_t),\dots, 
\eta^0(Z^{k-1}_t),\ind(Z^{k}_{bt}))\in
\check{\bC}^{k-1}(B,\cK(k-1,\Z)_B)\ .$$
\end{ddd}
Note that this definition also makes sense for $k=0$.
It follows from Theorem \ref{etaprop} and Lemma \ref{high0}
that
\begin{eqnarray}
\bd \bbeta(z)&=&0\nonumber\\
R^{[\bbeta(z)]}&=&\Omega^k(\cE_{geom})\nonumber\\
\DD([\bbeta(z)])&=&(-1)^{m}x \label{pu9}\ .
\end{eqnarray}

\begin{ddd}
For $z\in \tilde S^{k-1}(\cE_{geom},x)$ we define
$\del(z):=[\bbeta(z)]\in H^k_{Del}(B)$.\index{$\del(z)$}
\end{ddd}
We thus have associated to each tamed $k-1$-resolution $z\in \tilde F^k_{k-1}(B)$
a Deligne cohomology class $\del(z)\in H^k_{Del}(B)$. This class will depend on more
than just the class $\bz\in F^k_{k-1}(B)$
(see \ref{borddep}).

\subsubsection{}
In order to define an invariant of $\cE_{geom}$ or of the pair $(\cE_{geom},x)$ we
introduce the following sets.


\begin{ddd}\label{delcoh1}
If $\cE_{geom}$ is a geometric family such that $\ind(\cE_{geom})\in K_{k}^*(B)$, then we define for $x\in \check{\bo}^k(\ind(\cE_{geom}))$
\index{$\ind_{Del}^k(\cE_{geom},x)$}\index{$\ind_{Del}^k(\cE_{geom})$}$$\ind_{Del}^k(\cE_{geom},x):=\{\del(z) \:|\: z\in \tilde S^{k-1}(\cE_{geom},x)\}\subseteq H^k_{Del}(B)\ .$$
Furthermore, we set
$$\ind_{Del}^k(\cE_{geom}):=\bigcup_{x\in\check{\bo}^k(\ind(\cE_{geom})) } \ind_{Del}^k(\cE_{geom},x)\ .$$
\end{ddd}

\subsubsection{}
This set is natural in the following sense.
Let $f:B^\prime \rightarrow B$ be a smooth map.
Then we have an induced map
$f^\sharp:\tilde S^{k-1}(\cE_{geom},x)\rightarrow \tilde S^{k-1}(f^*\cE_{geom},f^*x)$.
If $f^*_{Del}:H^k_{Del}(B)\rightarrow H^k_{Del}(B^\prime)$ denotes the induced map
in Deligne cohomology, then the following relation immediately follows from the definitions.
\begin{kor}\label{natur}
If $z\in \tilde S^{k-1}(\cE_{geom},x)$, then $f^*_{Del}( \del(z)) = \del(f^\sharp z)$.
In particular, $$f^*_{Del}(\ind_{Del}^k(\cE_{geom},x))\subseteq \ind_{Del}^k(f^*\cE_{geom},f^*x)\ .$$
\end{kor}

\subsubsection{}\label{borddep}
The set $ \tilde S^{k-1}(\cE_{geom},x)$ is huge. Therefore it is a
natural question how big the set $\ind_{Del}^k(\cE_{geom},x)$ is.
The relation that we have considered so far on $\tilde F^{k}_{k-1}(B)$
is bordism. In fact, if $B$ is finite, then the image of  $\tilde
S^{k-1}(\cE_{geom},x)$ in $F^{k}_{k-1}(B)$ is finite by Corollary
\ref{kor1}. 
Unfortunately $\del(z)$ is not a bordism invariant of $z$.

\subsubsection{}
Assume that $z\in \tilde S^{k-1}(\cE_{geom},x)$
is bordant to $z^\prime \in \tilde S^{k-1}(\cE_{geom})$.
By Lemma \ref{olem}, 3. we have $z^\prime\in  \tilde S^{k-1}(\cE_{geom},x)$, too.
Let $\tZ,\tZ^\prime\in \tilde F^k_{k-1,\cU}(B)$ be representatives of
$z,z^\prime$,   and let $\tW$ be a zero bordism (see Definition \ref{bbo1}) of
$\tZ+(\tZ^\prime)^{op}$.
Note that by Theorem \ref{etaprop}
\begin{eqnarray*}
\delta \ind(W_{bt}^{k-1})&=&\delta \Omega^0(W^{k-1})+\delta
\eta^0(\partial W^{k-1}_{bt})\\
&=&\delta \Omega^0(W^{k-1}) + (-1)^{k-1}\delta \eta^0(Z^{k-1}_t)-
(-1)^{k-1}\delta \eta^0(Z^{\prime,k-1}_t)\\&&+(-1)^k \delta\delta
\eta^0(W^{k-2}_t)\\
&=&\delta \Omega^0(W^{k-1}) -  \eta^0(\partial Z^{k}_t)+
\eta^0(\partial Z^{\prime,k}_t)\\&&
-\Omega^0(Z^k)+\Omega^0(Z^{\prime,k})\\
&=&\delta \Omega^0(W^{k-1}) - 
(\ind(Z_{bt}^{k})-\ind(Z_{bt}^{\prime,k}))\end{eqnarray*}
In particular, $\delta \Omega^0(W^{k-1})$ is integral.
We define 
\begin{eqnarray*}\bOmega(\tW)&:=&((-1)^0\Omega^{k-1}(W^0),\dots,(-1)^{k-1}\Omega^0(W^{k-1}),-\delta \Omega^0(W^{k-1}))\in\check{\bC}^{k-1}(B,\cK(k-1,\Z)_B)\\
\bkappa(\tW)&:=&(\eta^{k-2}(W^{0}_t),\dots,\eta^0(W^{k-2}_t),\ind(W^{k-1}_{bt}))\in \check{\bC}^{k-2}(B,\cK(k-1,\Z)_B)\ .\end{eqnarray*}

We use  Theorem \ref{etaprop} in order to show
$$\eta(z)-\eta(z^\prime)=\bd\kappa(\tW)-\bOmega(\tW)\ .$$  
The calculations going into the verification of this formula  are similar to the calculations for (\ref{careful3}).
Going over to the cohomology classes this formula implies the following result.
\begin{lem}
If $z,z^\prime\in \tilde S^{k-1}(\cE_{geom},x)$
are bordant and $\tW$ is a zero bordism of $z-z^\prime$, then we have
 
\begin{equation}\label{myx}\del(z)-\del(z^\prime)=-[\bOmega(\tW)]\in
  H^k_{Del}(B)\ .\end{equation}
\end{lem}
The right-hand side is non-zero, in general.
Thus $\del: \tilde S^{k-1}(\cE_{geom},x)\rightarrow H^k_{Del}(B)$ does not factor over bordism classes.
In Subsection \ref{next21} we investigate this defect in more detail.

\subsection{Deligne cohomology classes for local resolutions}

\subsubsection{}

Let $\cE_{geom}$ be a geometric family over some base $B$. 
Let $\tZ\in\tilde F^{k}_{\cU, k-1}(B)$
be a tamed lift of the geometric $k$-resolution of $\cE_{geom}$ representing $z\in \tilde \cF^k_{k-1}(B)$.
By Lemma \ref{locexis} there is a local $k-1$-resolution  $X$ which is bordant to $\tZ$.  Let $\bK=(K_x)_x$ be the corresponding family of smoothing operators.
We define the cochain
\index{$\eta(X)$}$$\eta(X):=((-1)^{c(0)}\eta^{k-1}(X^0),\dots,(-1)^{c(k-1)}\eta^0(X^{k-1}), (-1)^{c(k)}u(X))\in \check{\bC}^{k-1}(\cU,\cK(k-1,\Z)_B)\ ,$$
where for $x\in\bN[p]$  the form 
$\bbeta^{k-1-p}(X^{p}(x))\in \cA_B^{k-1-p}(U_x)$ is the eta form of the 'tamed' family $X^{p}(x)$ over $U_x$.
Here $X^p(x)=(\cE_{geom})_{|U_x}*\Delta^p$ is the underlying geometric family which is 'tamed' by the family $K_x$, i.e.
$\eta^{k-1-p}(X^{p}(x))$ is the eta form of the family of invertible
perturbed Dirac operators
$D((\cE_{t})_{|U_x}*\Delta^p)$ defined in 
\ref{dqed}. 
 
\subsubsection{}

\begin{prop}\label{locald}
The cochain $\bbeta(X)$ satisfies
\begin{eqnarray}
 \bd\bbeta(X)&=&0\label{de567}\\
R^{[\bbeta(X)]}&=&\Omega^k(\cE_{geom})\label{gffduffs8}\\
\DD([\bd\bbeta(X)])&=&(-1)^{c(k)} [u(X)]\label{euwfuiwefwefiwefewf}\ .
\end{eqnarray}
  We further have
$\del(z)=[\eta(X)]\in H^{k}_{Del}(B)$.
\end{prop}
\proof
Equations  (\ref{de567}) and (\ref{gffduffs8}) follow from the generalization of Theorem \ref{etaprop} to families of operators of the form $D((\cE_{t})_{|U_x}*\Delta^p)$ (compare \ref{fuweufiewiufifoioi4}).
For (\ref{euwfuiwefwefiwefewf}) we  must observe that for $x\in \bN[k]$ we have 
\begin{equation}
\label{eq44} 
  \delta \eta^0(X^{k-1})(x)=  \ind(\pr_x^* \cE_{geom}, ((\pr_x^*
  \cE_{geom})_{|\partial \Delta^{k} \times U_x})_t)\ .
\end{equation}
We choose a point $b\in U_x$. \newcommand{\ori}{{\tt or}}
By Theorem \ref{etaprop} and $\Omega^0(\pr^*_x\cE_{geom}*\Delta^k)=0$ we get  
$$\delta \eta^0(X^{k-1})(x)=\ind((\cE_{geom})_b*\Delta^k)_{bt}\ .$$
In the proof of Proposition \ref{uduiduied}, Assertion 3, we have explained that
$$\ind((\cE_{geom})_b*\Delta^k)_{bt}=<\ori^{K}_{\Delta^k,\partial \Delta^k},\psi>\ ,$$ where $\ori^{K}_{\Delta^k,\partial \Delta^k}\in K_k(\Delta^k,\partial \Delta^k)$ is the $K$-orientation and 
$$\psi:=\ind(\pr_x^* \cE_{geom}, ((\pr_x^* \cE_{geom})_{|U_x\times \partial\Delta^k})_t)\in K^*(\Delta^k,\partial \Delta^{k})$$
is the families index of a family of operators on $\Delta^k$ which is invertible over $\partial \Delta^k$.
Equation  (\ref{eq44}) now follows.

We now show the second assertion.
Let $\tW$ be the bordism (i.e. the specially tamed cylinder, see \ref{eed4}) between $\tZ$ and $X$. Then we define
$$\kappa(\tW):=(\eta^{k-2} (W^0_{st}),\dots,\eta^0(W_{st}^{k-2}),0)\in \check{\bC}^{k-2}(\cU,\cK(k-1,\Z)_B)\ .$$
By Theorem \ref{etaprop} we get
$\bd \kappa(\tW) = \bbeta(z)-\eta(X)$.
\hB
 
\subsubsection{}

Proposition \ref{locald} provides an alternative definition of
the set $\ind^k_{Del}(\cE_{geom},x)$.
\begin{kor}\label{allesloc}
The set $\ind^k_{Del}(\cE_{geom},x)\subseteq H^k_{Del}(B)$ is given
as the set of classes $[\eta(X)]$, where
$X$ runs over all local $k-1$-resolutions such that $[u(X)]=x$.
\end{kor}

\subsubsection{}
Recall that $\eta(X)$ is defined in the adiabatic limit, 
what so far means that the metric of the simplices  
$\Delta^p$ is scaled by $\epsilon^{-2}$ and $\epsilon$ is sufficiently
small, but non-zero (this was the convention adopted in \ref{eed4}).

It seems to be very likely that the variant of $\eta$-forms used here
behave as usual with respect to adiabatic limits.
We expect that
$$\lim_{\epsilon\to 0}\eta^{k-1-p}(X^p(x))= \int_{U_x\times
  \Delta_p/U_x} \eta^{k-1}((\pr^*\cE_{geom})_t)\ ,$$
where
$(\pr_x^*\cE_{geom})_t$ is the lifted family over $U_x\times\Delta^p$ tamed by
$K_x$. This would yield the connection to Lott's approach
\cite{lott01}. A proof would involve more involved analysis.
Since in the present paper we want to focus on the geometric and
algebraic aspects of the theory and keep analysis as easy as possible
we refrain from giving a proof of this conjecture here.

\subsection{Canonical lifts of Chern classes}\label{next21}

\subsubsection{}
In general the set $\ind^k_{Del}(\cE_{geom})$ may have several
distinct elements. In the present subsection we describe how these
elements are related.

Let $k,m\in\nat\cup\{0\}$ be such that $k=2m$ or $k=2m-1$.
Let $\cE_{geom}$ be a geometric family over $B$ such that
$\ind(\cE_{geom})\in K^*_k(B)$. Then the set
$\ind^k_{Del}(\cE_{geom})$ is not empty. 
If $x\in \ind^k_{Del}(\cE_{geom})$ and $k>0$, then by Corollary \ref{kerv} and (\ref{pu9})
we have $$(-1)^{m+1}(m-1)!\DD(x)=c_k(\ind(\cE_{geom}))\ .$$ 

We will see that $(-1)^{m+1}(m-1)!x\in H^k_{Del}(B)$ is in fact determined by
$\cE_{geom}$ and thus a canonical lift of $c_k(\ind(\cE_{geom}))$
to Deligne cohomology.

\subsubsection{}
In order to treat the case $m=0$ in a uniform manner we set $(-1)!:=1$.
\begin{prop}\label{denominators}
The set
$$(m-1)!\ind^k_{Del}(\cE_{geom})$$ contains exactly one 
element. 
\end{prop}
\proof
The case $k=0$ is obvious. We assume that $k\ge 1$.
Let $u_0,u_1\in \ind^k_{Del}(\cE_{geom})$. We show that
$(m-1)!(u_0-u_1)=0$.

In fact, we will show that
$(m-1)! (H(u_0)-H(u_1))=0$, where $H:H^k_{Del}(B)\rightarrow \hat H^{k-1}(B,\R/\Z)$ 
is the holonomy map from Deligne cohomology to Cheeger-Simons differential characters (see Subsection \ref{holo87}).

Let $z\in Z_{k-1}(B)$ be a smooth cycle. Let $j:U\rightarrow B$ be the
inclusion of an open neighborhood of the trace of $z$ which 
has the homotopy type of a $k-1$-dimensional $CW$-complex.
It suffices to show that
$(m-1)! (H(j^* u_0)(z)-H(j^*u_1)(z))=0$.

We have $j^*u_i\in \ind^k_{Del}(j^*\cE_{geom})$ by Corollary \ref{natur}.
By Corollary \ref{allesloc} we can find local $k-1$-resolutions $X_0$
and $X_1$ (w.r.t some open covering $\cU$ of $U$) of $j^*\cE_{geom}$ such that
$[\bbeta(X_i)]=u_i$.  

Let us assume that $\cU$ is a good covering.
We consider $r:\tilde U\rightarrow U$ and its restriction
$r^{k-1}:\tilde U^{k-1}\rightarrow U$
in the notation of \ref{tfc21}. Let $\cF_{geom}:=r^*\cE_{geom}$
and $\cF_{k,geom}:=(\cF_{geom})_{|\tilde U^{k-1}}$.
 The local $k-1$-resolutions $X_i$ are given by  families
$(\bK_{i,x})_{\dim(x)\le k-1}$, $i=0,1$, which can be interpreted as
tamings $\cF_{k,t,i}$ of $\cF_{k,geom}$.
The inclusion
$\tilde U^{k-1}\hookrightarrow \tilde U$ is a homotopy equivalence.
Therefore the tamings $\cF_{k,t,i}$ extend to 
tamings $\cF_{t,i}$ of $\cF_{geom}$.

Since $r:\tilde
U\rightarrow U$ is a homotopy equivalence,
there exists tamings
$(\cE_{geom})_{|U,t,i}$ such that
$\cF_{t,i}^\prime:=r^* (\cE_{geom})_{|U,t,i}$ are homotopic to
$\cF_{t,i}$.
The homotopy can be considered as a taming
$(\tilde \pr^*\cF)_{t,i}$ of $\tilde \pr^*\cF_{geom}$, where
$\tilde \pr:I\times \tilde U\rightarrow \tilde U$.

We now consider the covering  $\tilde \cU$ of $I\times U$
which is given by the sets $\tilde A:=I\times A$, $A\in \cU$.
Then the nerves of $\cU$ and $\tilde \cU$ coincide naturally and we
can identify $\widetilde{I\times U}\cong I\times \tilde U$,
where the left-hand side is again the construction of \ref{tfc21}
now applied to the covering $\tilde \cU$ of $I\times U$.

We can identify the restriction $((\tilde \pr^*\cF)_{t,i})_{\widetilde{I\times
  U}^{k-1}}$ with a local $k-1$-resolution $\tilde X_i$ of
$\pr^*j^*\cE_{geom}$, where $\pr:I\times U\rightarrow U$.
If $f_l:U\rightarrow I\times U$, $l=0,1$, are induced by the embeddings of the
endpoints of the interval $I$, then we have by construction
$f_0^* \tilde X_i=X_i$, while $f_1^*\tilde X_i$ is
the local $k-1$-resolution $X_i^\prime$ induced by
$(\cE_{geom})_{|U,t,i}$.

We form the classes $h_i=[\bbeta(\tilde X_i)]\in H^{k-1}_{Del}(I\times
U)$. Note that $R^{h_i}=\pr^*j^*\Omega^k(\cE_{geom})$.
In particular $\int_{I\times U/U} R^{h_i}=0$.
By the homotopy formula \ref{homotopiedel} we obtain
$$u_i=[\bbeta(X_i)]=[\bbeta(X_i^\prime)]\ .$$

We now employ the following fact.
Let $\cG_{geom}$ be a geometric family with closed fibers over some
base $B$ and let $Q$ be a family of smoothing operators which induce a
taming $\cG_{t}$.
Let $N$ be an oriented Riemannian spin manifold with corners giving
rise to $\cN_{geom}$.
Then $\cG_t*N$ is already tamed (in the adiabatic limit).
Therefore $\eta(\cG_t*N)$ is well-defined.

\begin{lem}
If $\cN_{geom}\cong \cN^{op}_{geom}$,
then $\eta(\cG_t*N)=0$.
\end{lem}
\proof
We have an induced isomorphism
$\cG_t*N\cong (\cG_t*N)^{op}$ so that
$\eta(\cG_t*N)=-\eta(\cG_t*N)$. \hB
 
This Lemma can be applied in the case $\cN_{geom}=\cDelta^p_{geom}$,
$p\ge 1$. We get
$$\bbeta(X_i^\prime)=(\eta^{k-1}((\cE_{geom})_{|U,t,i}),0,\dots,0)\in\check{\bC}^{k-1}(\cU,\cK(k-1,\Z)_U)\
,$$ i.e.
we have
$$a[\eta^{k-1}((\cE_{geom})_{|U,t,i})]=u_i\ .$$
We now see that 
$$H(u_i)(z)= [ \int_{z}\eta^{k-1} ((\cE_{geom})_{|U,t,i})]_{\R/\Z}\ .$$
 
By Corollary \ref{etajump} there exists $\psi\in K(U)$ such that
$$\dR(\ch_{k-1}(\psi))=[\eta^{k-1}((\cE_{geom})_{|U,t,1})-\eta^{k-1}((\cE_{geom})_{|U,t,0})]\ .$$
The rational cohomology class
$(m-1)! \ch_{k-1}(\psi)$ has integral periods.
We conclude that
$$(m-1)!
\int_{z}(\eta^{k-1}((\cE_{geom})_{|U,t,1})-\eta^{k-1}((\cE_{geom})_{|U,t,0}))
\in\Z\ .$$
Thus for all $z\in Z^{k-1}(B)$ we have
$(m-1)! H(u_1-u_0)(z)=0$.
This implies the Proposition.
\hB

\subsubsection{}
For later use we separate a consequence of the proof of the
Proposition. Let $k\ge 0$ and $\cE_{geom}$ be a geometric family over
$B$ such that $\ind(\cE_{geom})\in K_k(B)$.
Let $u\in \ind_{Del}^k(\cE_{geom})$ and $j:U\hookrightarrow B$ be the
inclusion of  an open subset
which is homotopy equivalent to a $CW$-complex of dimension $\le k-1$.
\begin{kor}
There exists a taming $(j^*\cE_{geom})_t$ such that
$$j^*u=a[\eta^{k-1}((j^*\cE_{geom})_t)]\ .$$
\end{kor}

\subsubsection{}

Let $k,m\ge 1$, $2m=k$ or $2m-1=k$, and $\cE_{geom}$ be a geometric family over
$B$ such that $\ind(\cE_{geom})\in K_k(B)$.
\begin{ddd}\label{chern3}
\index{$\hat c_k(\cE_{geom})$}\index{Chern! class in Deligne cohomology}
We define $\hat c_k(\cE_{geom})\in H^k_{Del}B)$ to be the unique element
of $(-1)^{m-1}(m-1)!\ind_{Del}^k(\cE_{geom})$.
\end{ddd}
We then have
\begin{eqnarray*}
\DD(\hat c_k(\cE_{geom}))&=&c_k(\ind(\cE_{geom}))\\
R^{\hat c_k(\cE_{geom})}&=&(-1)^{m-1} (m-1)! \Omega^k(\cE_{geom})\\
f^*\hat c_k(\cE_{geom})&=&\hat c_k(f^*\cE_{geom})\quad
\ ,\:\:f:B^\prime\rightarrow B
\end{eqnarray*}
Thus we can consider $\hat c_k(\cE_{geom})$ as the natural lift of the
Chern class $c_k(\ind(\cE_{geom}))$ to integral Deligne cohomology.

\subsubsection{}\label{holcomp}
We can compute the holonomy of $\hat c_k(\cE_{geom})$ as follows.
Let $z\in Z^{k-1}(B)$. We can assume that the trace of $z$ is contained
in an open subset $U\subseteq B$ which is homotopy equivalent to a
$CW$-complex of dimension at most $k-1$. Then we can choose a taming
$(\cE_{geom})_{|U,t}$ and get
$$H(\hat c_k(\cE_{geom}))(z)=[(-1)^{m-1}(m-1)!\int_z \eta^{k-1}((\cE_{geom})_{|U,t})]_{\R/\Z}\
.$$
 
\subsubsection{}

If $\bW$ is a geometric vector bundle over $B$, then
in \cite{cheegersimons83}  natural lifts of the Chern classes
$\hat c_{2m}(\bW)\in H_{Del}^{2k}(B)$ were defined for all $m\ge 1$.
The bundle  $\bW$ gives rise to the geometric family $\cE(\bW)_{geom}$ (see \ref{exex}).

If $[W]\in K^0_{2m}(B)$, then we also have defined
\index{$\hat c_{2m}(\bW)$}$\hat c_{2m}(\cE(\bW)_{geom})$.
In general we have
$$\hat c_{2m}(\bW)\not=\hat c_{2m}(\cE(\bW)_{geom})\ $$
already on the level of curvatures.
In fact, $R^{\hat c_{2l}(\bW)}$ is the Chern-Weyl representative of
the Chern class while $(-1)^{m-1}(m-1)!\Omega^{2m}(\cE_{geom})$ is the Chern-Weyl
representative of the
component
of the Chern character. These in general  do not coincide on the level of forms
since $R^{\hat c_{2l}(\bW)}$ must not vanish identically for $0\le l<
m$. But they do on the level of cohomology by the assumption that
$\ind(\cE_{geom})\in K_{2m}(B)$ since then  the forms
$R^{\hat c_{2l}(\bW)}$ are exact for $0\le l<
m$.

In a forthcoming paper \cite{bunkeschick031} we will introduce a model of smooth
$K$-theory\footnote{See the forthcoming paper \cite{bunkeschick031} for
  a formal definition of the concept of a smooth extension of a generalized cohomology theory. Roughly speaking smooth $K$-theory is related
  to $K$-theory in very much the same way  Deligne cohomology is related to
  integral cohomology}
based on geometric families. The observations above are closely
related to the problem of lifting the Chern classes to
natural transformations from smooth $K$-theory to smooth integral
cohomology
(i.e. Deligne cohomology).

In this direction note that Hopkins-Singer \cite{hs} propose models for smooth
$K$-theory, smooth integral cohomology, and smooth versions of other generalized cohomology theories. In their framework lifts of Chern classes can be constructed. It appears the interesting question  if a geometric family  gives rise to a natural smooth $K$-theory class in the picture of \cite{hs}, or equivalenty, if here is a natural transformation (which is then automatically an isomorphism)
between the smooth $K$-theories of \cite{hs} and \cite{bunkeschick031}.

\section{Examples}\label{fc34}

\subsection{The case $k=0$ : Index }\label{sma}

\subsubsection{}
If $k=0$, then $m=0$. We have a natural isomorphism $H^0_{Del}(B)\cong
\check{H}^0(B,\underline{\Z}_B)$.
In fact even the underlying complexes a equal.

\subsubsection{}

If $\cE_{geom}$ is a reduced geometric family, then we have
$\ind(\cE_{geom})\in K^*_0(B)=K^*(B)$. 
The set $\ind_{Del}^0(\cE)$ contains exactly  one element.
By definition, this element is represented by the locally constant
function
$$\ind_0(\cE_{geom}):B\rightarrow \Z\ .$$
Note that only the even-dimensional part of $\cE_{geom}$ contributes.

\subsection{The case $k=1$ : $\eta$-invariants}

\subsubsection{}

We have a natural isomorphism
$$H^1_{Del}(B)\cong C^\infty(B,\R/\Z)\ .$$
This isomorphism is given as follows.
Let $x\in H^1_{Del}(B)$. Then the corresponding $\R/\Z$-valued
function maps $b\in B$ to $H(x)(b)\in \R/\Z$, where we interpret the
point $b$ as a cycle
$b\in Z^0(B)$.

Let us spell out this isomorphism in terms of \v{C}ech cocycles, too.
Let $x\in H^1_{Del}(B)$ be represented by a cocycle (see \ref{c824} for the notation)
$\bc=(c^0,c^{-1})$ with respect to some open covering $(U_i)_{i\in L}$.
Then $c^{-1}$ is a collection of  continuous $\Z$-valued functions on the double intersections, and $c^0$ is a collection of  smooth $\R$-valued functions on the $U_i$ themselves.
The relation $-c^{-1}-\delta c^0$ implies that
$\exp(2\pi i c^0)$ is equal to the restriction of a global function, namely of
$H(x):B\to \R/\Z$.

\subsubsection{}

Let us assume that $\cE_{geom}$ is odd-dimensional.
Then $\ind(\cE_{geom})\in K^1_1(B)=K^1(B)$.
The set $\ind^1_{Del}(\cE_{geom})$ contains exactly one element,
namely $-\hat c_1(\cE_{geom})$.
In the following we describe the corresponding $\R/\Z$-valued function
in familiar terms.

\subsubsection{}

Let $b\in B$. 
We choose some  taming of the fiber $(\cE_{geom})_{b,t}$.
Then $$H(-\hat c_1(\cE_{geom}))(b)=[\eta^0(\cE_{t,b})]_{\R/\Z}\ .$$

Note that we also can define $\eta^0(D(\cE_{geom})_b)$
if this operator is not invertible.
Then we have
$$[\eta^0(\cE_{t,b})]_{\R/\Z}=[\eta^0(D(\cE_{geom})_b)+\frac{1}{2}\dim\ker(D(\cE_{geom})_b)]_{\R/\Z}\
.$$
We thus get:
\begin{lem}
If $\cE_{geom}$ is an odd-dimensional family, then
the holonomy of the  unique element of $\ind^1_{Del}(\cE_{geom})$
is given by
$$B\ni b\mapsto
[\eta^0(D(\cE_{geom})_b)+\frac{1}{2}\dim\ker(D(\cE_{geom})_b)]_{\R/\Z}\in\R/\Z\
.$$
\end{lem}

\subsubsection{}

Note that $$\ker(\bv:H^1_{Del}(B)\to \check{H}^1(B,\Z))= a(\cA^0_B(B)/\Z)\cong C^\infty(B,\R)/\Z$$
(do not confuse with $C^\infty(B,\R/\Z)$).
Assume that for some reason $\dim\ker(D(\cE_{geom})):B\to \Z$ is locally constant. Then $c_1(\ind(\cE_{geom}))=0$ and the function
\begin{equation}\label{dnhj238u}B\ni b\mapsto
\eta^0(D(\cE_{geom})_b)+\frac{1}{2}\dim\ker(D(\cE_{geom})_b)\in \R
\end{equation}
represents $\hat c_1(\cE_{geom})$ under the identification above.
This function contains more information than its class $\hat c_1(\cE_{geom})$ modulo $\Z$, and this additional information is of importance in some topological applications.
The condition that $\ker(D(\cE_{geom}))$ is constant can be satisfied for geometric or topological reasons.
Let us mention two cases.
 
 Let $\cE_{geom}$ be given by  family of Riemannian spin manifolds of positive scalar curvature with
underlying Clifford
bundle of the form $\cV:=\cS\otimes \bF$, where $\cS$ is the vertical spinor bundle and $\bF$ is a flat bundle.  Then $\ker(D(\cE_{geom}))=0$, hence in particular constant.
The function (\ref{dnhj238u}) for  such families plays an important role in the theory of positive scalar curvature metrics.
%

In the second case we let the Clifford bundle of $\cE_{geom}$ be (locally) of the form $\cV:=\cS\otimes \cS^{ungr}\otimes \bF$ (where  $\cS^{ungr}$ is the spinor bundle with trivial grading). Then $D(\cE_{geom})$ is the vertical signature operator. Its kernel is isomorphic to the cohomology of the fibre with coefficients in the local system given by $\bF$. These cohomology groups are  independent of geometric structures and their dimensions are locally constant on $B$. In this case the function (\ref{dnhj238u}) plays a role in Chern-Simons theory.

 
\subsection{The case $k=2$ : Determinant bundles}

\subsubsection{} 

Let $\Line(B)$ denote the group of isomorphism classes
of hermitian line bundles with connection on $B$, where the group operation is given
by the tensor product. Let $L\in \Line(B)$.
If $\gamma:S^1\rightarrow B$ is a smooth map,
then $L$ determines an element $\hol(L)(\gamma)\in \R/\Z$ such that
$\exp(2\pi i \hol(L)(\gamma))\in U(1)$ is the holonomy of the parallel
transport in $L$ along $\gamma$.
Note that an element $L\in \Line(B)$ is uniquely
determined by the values $\hol(L)(\gamma)$ on all maps
$\gamma:S^1\rightarrow B$.

There is a natural isomorphism $H^2_{Del}(B)\cong \Line(B)$
(see \cite{brylinski93}, Thm. 2.2.11). Let us describe this
isomorphism in more detail.

Let $x\in H^2_{Del}(B)$. Then we have the corresponding differential
character
$H(x):Z_1(B)\rightarrow B$. Let us fix once and for all a
triangulation of $S^1$.  Then we can consider $\gamma:S^1\rightarrow B$
as a smooth cycle $\gamma\in Z_1(B)$. 
The identification
$H^2_{Del}(B)\cong \Line(B)$
now asserts that there is a unique $L\in \Line(B)$ such that
$H(x)(\gamma)=\hol(L)(\gamma)$ for all $\gamma:S^1\rightarrow B$.
Note that this correspondence is independent of the choice of the
triangulation of $S^1$.

If one prefers \v{C}ech cocycles with respect to some covering $(U_i)_{i\in I}$, the transition from a class $x\in  H^2_{Del}(B)$ represented by a cocycle $\bc(c^1,c^0,c^{-1})$ to a hermitean  line bundle $L$ with connection goes as follows.
The object $c^0$ is a collection of smooth $\R$-valued functions on
the double intersections. Since $c^{-1}$ is $\Z$-valued on triple intersections the relation
$\delta c^0-c^{-1}=0$ shows that 
$\delta \exp(2\pi i c^0)=0$. The $U(1)$-valued cocycle
$\exp(2\pi i c^0)$ can be interpreted as the collection of transition  functions of a hermitean line bundle $L$
with the choice of local unit length sections $s_i$ of $L_{|U_i}$, $i\in I$, i.e.
$$\exp(2\pi i c^0((i,j))) (s_{i})_{|U_i\cap U_j}=(s_j)_{|U_i\cap U_j}\ ,\quad (i,j)\in \bN[1]\ .$$
The equation
$d c^0-\delta c^1=0$ can be rewritten as
$2\pi i  \delta c^1= d\log c^0$.
In other words, the collection of one forms
$2\pi i c^1(i)$ on the $U_i$, $i\in I$, has the interpretation of a collection of connection one forms which locally represent a connection $\nabla^L$ with $$\nabla^L s_i= 2\pi i c^1(i) s_i\ .$$

\subsubsection{}

Let us assume that $\cE_{geom}$ is a reduced even-dimensional bundle
such that $\ind(\cE_{geom})\in K_2^0(B)$. This means that
$\ind_0(\cE_{geom})$ vanishes identically.
Under this assumption the set $\ind^2_{Del}(\cE_{geom})$ contains a
unique element $\hat c_2(\cE_{geom})$.

\subsubsection{}
\index{determinant line bundle}An even-dimensional geometric family gives rise to a determinant
line bundle. The underlying complex line bundle is the determinant line bundle of the family of Dirac operators
$D(\cE_{geom})$. It comes with the Quillen metric and the Bismut-Freed
connection. By $\det(\cE_{geom})\in \Line(B)$ \index{$\Line(B)$}we denote its isomorphism class.

\subsubsection{}
\begin{prop}\label{gghba}
The hermitian line bundle with connection 
\index{$\det(\cE_{geom})$}which corresponds to $\hat c_2(\cE_{geom})$
is the isomorphism class of the determinant line bundle
$\det(\cE_{geom})\in \Line(B)$ of
$\cE_{geom}$.
\end{prop}
\proof
Let $\gamma:S^1\rightarrow B$ be a smooth map.
We will compare the holonomy $H(\hat c_2(\cE_{geom}))(\gamma)$
with the holonomy $\hol(\det(\cE_{geom}))(\gamma)$.
Let $j:U\rightarrow B$ be a neighborhood of the image
of $\gamma$ which has the homotopy type of a one-dimensional
$CW$-complex. We choose a taming $(j^*\cE_{geom})_t$ which is given by
a family of smoothing operators $Q$.
Then by \ref{holcomp} we have
$$H(\hat c_2(\cE_{geom}))(\gamma)=[\int_{S^1} \gamma^*
\eta^1((j^*\cE_{geom})_t)]_{\R/\Z}\ .$$

In order to compute $\hol(\det(\cE_{geom}))$ we use the results of Bismut-Freed \cite{bismutfreed861}, \cite{bismutfreed862}.
We equip $S^1$ with a Riemannian metric $g_\epsilon^{S^1}:=\epsilon^{-1} g^{S^1}$.
The pull-back bundle $\pi:\gamma^*\cE_{geom}$ comes equipped
with a horizontal distribution, a vertical metric, and a fiber-wise orientation.
Therefore the total space $\gamma^*E$ 
has an induced Riemannian metric.
We define the Dirac operator $D_\epsilon$ on the total space as
follows.
Let $\Gamma:\gamma^*E\rightarrow E$ be the canonical map over $\gamma$.
Locally on $\gamma^*E$ after fixing a spin structure of the vertical bundle
we can define a twisting bundle $\bW$ such that
$\Gamma^*\cV=\cS(T^v\pi)\otimes \bW$.
The choice of the local  vertical spin structure together with the spin structure of $S^1$ induces a local spin structure of $\gamma^*E$. 
Thus we can define locally the Dirac bundle
$\cS(T\gamma^*E) \otimes \bW$ on the total space.
One checks that this bundle is independent of the choices and therefore globally defined. We let $D_\epsilon$ be the Dirac operator associated to this bundle.

We consider $D_\epsilon(s)=sD_\epsilon+s\chi(s)Q$. Then we can define
$\tau(\epsilon,\chi)\in \R/\Z$ as in Section 4.4 of \cite{bunke011}
\footnote{Here we use  the identification $\exp(2\pi i\dots):
  \R/\Z\rightarrow U(1)$.}.
By \cite{bunke011},  Lemma 4.4, we have
$\tau(\epsilon,\chi)=\tau(\epsilon)$, where 
$$\tau(\epsilon)=\tau(\epsilon,0)=[\frac{\eta(D_\epsilon)+\dim\ker(D_\epsilon)}{2}]_{\R/\Z}\ \footnote{Here we use the usual conventions for the $\eta$-invariant. The connection is given by $\eta(D)=-2\eta^0(D)$.}
.$$
By Bismut-Freed \cite{bismutfreed861}, \cite{bismutfreed862},
$$\hol(\det(\cE_{geom}))(\gamma) = \lim_{\epsilon\to 0} \tau(\epsilon)\ .$$
As in \cite{bunke011}, Section 4.6 one can show that 
$$\lim_{\epsilon\to 0} \tau(\epsilon,\chi)=
[\int_{S^1}\gamma^*\eta^1((j^*\cE_{geom})_t)]_{\R/\Z}\ .$$
Combining these equalities we get
$$H(\hat c_2(\cE_{geom}))(\gamma)=\hol(\det(\cE_{geom}))(\gamma) \ .$$
\hB

\subsubsection{}

Note that there is a natural way to define
$\hat c_{2}(\cE_{geom})$ even in the case where
$\ind_0(\cE_{geom})\not\equiv 0$.
Let $\bV$ be  a trivial geometric bundle over $B$ such that
$\ind_0[\bV]=-
\ind_0(\cE_{geom})$.
Then $\ind_0(\cE_{geom}+\cE(\bV)_{geom})=0$ so that we can define
$\hat c_2(\cE_{geom}+\cE(\bV)_{geom})$.

\begin{lem}
The class 
$\hat c_2(\cE_{geom}+\cE(\bV)_{geom})$
does not depend on the choice of $\bV$. 
\end{lem}
\proof
This follows from Proposition \ref{gghba},
$$\det(\cE_{geom}+\cE(\bV)_{geom})=\det(\cE_{geom})\otimes
\det(\cE(\bV)_{geom})\ ,$$  and $\det(\cE(\bV)_{geom})=1$. \hB

\subsubsection{}

For a further investigation of properties of the determinant bundle, glueing formulas, and the  relation with  gerbes we refer to \cite{bunkepark}.

\subsection{The case $k=3$ : Index gerbes}

\subsubsection{}\label{gerde}
\index{gerbes}\index{$\Gerbe(B)$}Let $\Gerbe(B)$ denote the group of isomorphism classes of geometric
gerbes, i.e. of gerbes with connective structure and curving.
We refer to \cite{lott01} and \cite{hitchin99} for precise
definitions. Let $G\in\Gerbe(B)$. If $f:\Sigma\rightarrow B$ is a
smooth map from an oriented closed surface $\Sigma$ to $B$, then we
have a holonomy $\hol(G)(f)\in\R/\Z$. Note that as in the case of 
line bundles the isomorphism class $G\in\Gerbe(B)$ is completely
determined
by the values $\hol(G)(f)$ for all such $f:\Sigma\rightarrow B$.

\subsubsection{}
There is a natural isomorphism $H^3_{Del}(B)\cong \Gerbe(B)$ 
(see \cite{brylinski93}) which is
as in the case of line bundles given by the holonomy.
Let $x\in H^3_{Del}(B)$. Then the identification asserts that there is
a unique $G\in\Gerbe(B)$ such that
$H(x)(z_f)=\hol(x)(f)$ for all $f:\Sigma\rightarrow B$.
Here we choose some triangulation of $\Sigma$ and let $z_f\in Z_2(B)$ be
the corresponding cycle.

\subsubsection{}

\index{index gerbe}\index{$\gerbe(\cE)$}We assume that $\cE_{geom}$ has odd-dimensional fibers.
Then we have a canonical element
$\gerbe(\cE)\in  \Gerbe(B)$, the index gerbe of $\cE_{geom}$ which was
constructed by Lott \cite{lott01} in analogy to the determinant line
bundle. 

\subsubsection{}
Let us assume that $\ind(\cE_{geom})\in K_3^1(B)$.
This is equivalent to the vanishing of the spectral flow 
$c_1(\ind(\cE_{geom}))\in \check{H}^1(B,\underline{\Z}_B)$.
Under this assumption the set
$\ind^3_{Del}(\cE_{geom})$ contains a unique element $\hat
c_3(\cE_{geom})\in H^3_{Del}(B)$.

\subsubsection{}
\begin{prop}\label{c3ge}
The class $\hat
c_3(\cE_{geom})\in H^3_{Del}(B)$ corresponds to the index gerbe
$\gerbe(\cE_{geom})$ under the isomorphism described in \ref{gerde}
\end{prop}
\proof
Let $f:\Sigma\rightarrow B$ be a smooth map from a triangulated
oriented closed surface to $B$, and let $z_f\in Z_2(B)$ be the corresponding
cycle. We will show that
$H(\hat
c_3(\cE_{geom}))(z_f)=\hol(\gerbe(\cE_{geom}))(f)$.

Let $j:U\rightarrow B$ the inclusion of an open neighborhood of the
image of $f$ which has the homotopy type of a $2$-dimensional
$CW$-complex. We choose some taming $(j^*\cE_{geom})_t$ given by a
family of smoothing operators $Q$.
Then by \ref{holcomp} we have
$$H(\hat
c_3(\cE_{geom}))(z_f)=[\int_{\Sigma} f^*
\eta^2((j^*\cE_{geom})_t)]_{\R/\Z}\ .$$

 Using Lemma 4.6 of \cite{bunke011}, the method of the proof of
 \cite{bunke011}, Lemma 4.1, and the notation
 of that paper\footnote{We again identify $\R/\Z$ with
   $U(1)$. Furthermore, here we use the notation $\hol$ for the gerbe
   holonomy which was denoted by  $H$ in the reference}
 we get
\begin{eqnarray*}
[\int_\Sigma  f^* \eta^2((j^*\cE_{geom})_t)]_{\R/\Z}&=&
[ \int_{\Sigma } f^*\eta^2(j^*\cE_{geom},Q)]_{\R/\Z}\\
&=&\hol(\gerbe (f^*j^*\cE_{geom},f^*Q))(\id_{\Sigma})\\
&=&\hol(\gerbe( f^*(j^*\cE_{geom})))(\id_{\Sigma})\\
&=&\hol(\gerbe(\cE_{geom}))(f)\ .
\end{eqnarray*}
\hB

\subsubsection{}
Again we can define a natural element $\hat c_3(\cE_{geom})\in H^3_{Del}(B)$
even if $c_1(\ind(\cE_{geom}))\not=0$.
Let  $\cF_{geom}$ be a geometric family over $S^1$ with odd-dimensional fibers such that $\ind(\cF_{geom})=1\in K^1(S^1)\cong \Z$.
One could e.g. represent the generator $1\in K^1(S^1)$ by the identity map $\id:S^1\rightarrow S^1=U(1)$ and then take
$\cF_{geom}:=\cE(\id,*)_{geom}$ as introduced in Subsection \ref{exex}.
Furthermore, let $f:B\rightarrow S^1$ be the classifying map of 
$-c_1(\ind(\cE_{geom}))$,
i.e. $f^*c_1(\ind(\cF_{geom}))=-c_1(\ind(\cE_{geom}))$.
Then $\ind(\cE_{geom}+f^*\cF_{geom})\in K_3^1(B)$.

\begin{lem}
The class $\hat c_3(\cE_{geom}+f^*\cF_{geom})$
does not depend on the choice of $f$ or $\cF_{geom}$.
\end{lem}
\proof
Note that $H^3_{Del}(S^1)=0$. Therefore $\gerbe(\cF_{geom})=0$. We have
$$\gerbe(\cE_{geom})=\gerbe(\cE_{geom})+f^*\gerbe(\cF_{geom})=\gerbe(\cE_{geom}+f^*\cF_{geom})\
.$$
This implies the assertion in view of Proposition \ref{c3ge}. \hB

\subsection{Computations for $S^1$-bundles}

\subsubsection{}

Let $\pi:E\rightarrow B$ be an $U(1)$-principal bundle over a connected base $B$.  
We fix a basis vector of the Lie algebra $u(1)$. Then the 
vertical bundle $T^v\pi$  is trivialized by the corresponding  fundamental vector field.
We choose the orientation in which this vector field becomes positive.
 
\subsubsection{}
Recall that the circle $S^1$ admits two spin structures. One of them (the trivial one)
extends to a spin structure of the two dimensional  disc. If we take the metric of volume $2
\pi$, then the spectrum of the Dirac operator on $S^1$ with respect to this spin structure 
is $\frac{1}{2}+\Z$. In particular, it is invertible.

\subsubsection{}
We want to choose a spin structure on $T^v\pi$ which restricts to the trivial
spin structure on each fiber.
\begin{lem}
The vertical bundle $T^v\pi$ of an $U(1)$-principal bundle $\pi:E\rightarrow B$
admits a spin structure which restricts to the trivial spin structure
on the fibers iff the reduction of $c_1(E)$ modulo $2$ vanishes.
\end{lem}
\proof
Since $T^v\pi$ is trivial, it admits a spin structure
which restricts to the nontrivial spin structure on the fibers.
If $E\times SO(1)\cong P_{SO(1)}\rightarrow E$ is the $SO(1)$-principal bundle associated
to the trivial oriented bundle $T^v\pi$, then this spin structure is just the two-fold 
non-connected covering $P_{Spin(1)}\rightarrow P_{SO(1)}$.
We take this spin
structure as a base point so that the set of all spin structures
of $T^v\pi$ is in bijection with $H^1(E,\Z_2)$.

The Leray-Serre spectral sequence gives the exact Gysin sequence
$$0\rightarrow H^1(B,\Z)\rightarrow H^1(E,\Z) \stackrel{r}{\rightarrow} H^0(B,\Z)\stackrel{d_2}{\rightarrow}\ H^2(B,\Z)\ .$$
It is well-known, that $d_2(1)=-c_1(E)$
(see e.g. Borel-Hirzebruch \cite{borelhirzebruch59}, Thm. 29.4),
where $1\in H^0(B,\Z)\cong \Z$ is the generator.
The map $r: H^1(E,\Z)\rightarrow H^0(B,\Z)\cong \Z \cong H^1(F,\Z)$
is restriction to the fiber.
After reduction modulo two we obtain
$$0\rightarrow H^1(B,\Z/2\Z)\rightarrow H^1(E,\Z/2\Z)\stackrel{[r]}{\rightarrow} H^0(B,\Z/2\Z)\stackrel{[c_1(E)]}{\rightarrow}  H^2(B,\Z/2\Z)\ .$$
A spin structure of $T^v\pi$ corresponding to $x\in H^1(E,\Z/2\Z)$ restricts
to the trivial spin structure on the fibers iff $[r](x)\not=0$.
Since $H^0(B,\Z/2\Z)=\Z/2\Z$ the range of $[r]$ is non-trivial exactly if
$[c_1(E)]=0$. \hB

\subsubsection{}
From now on we assume that $c_1(E)$ is even, and that $T^v\pi$ is equipped 
with a spin structure which restricts to the nontrivial spin structure on the fibers.
We take an $U(1)$-invariant vertical metric such that the volume of the fibers is $2\pi$.
Furthermore, we choose an $U(1)$-connection $\omega$ which induces a horizontal distribution
$T^h\pi$. We consider the Dirac bundle bundle $\cS(T^v\pi)$.
Thus we have defined a geometric family $\cE_{geom}$.

\subsubsection{}
Note that the family of Dirac operators $D(\cE_{geom})$ is invertible.
Thus the geometric $k$-resolution is automatically a tamed $k-1$-resolution
with taming given by zero operators on the faces.
We let $z_{k-1}\in \tilde S^{k-1}(\cE_{geom})$ be the corresponding elements.

\subsubsection{}

Let $B_k$ denote the Bernoulli numbers which are defined by the generating series
$$\frac{\ee^x}{\ee^x-1}=\sum_{k=-1}^\infty B_{k+1} \frac{x^k}{(k+1)!}\ .$$
Furthermore, let $c_1(\omega):=\frac{-1}{2\pi\imath} F^\omega\in \cA_B^2(B)$ be the
first Chern form of $E$.

\begin{prop}
We have for $m\in\nat$
$$\del(z_{2m})=a[\frac{B_{m+1}}{(m+1)!} c_1(\omega)^{m}]\in H_{Del}^{2m+1}(B)\ .$$
In particular, $R^{\del(z_{2m})}=0$ and $\DD(\del(z_{2m}))=0$.
\end{prop}
\proof
Note that 
$\del(z_{2m})=a[\eta^{2m}(\cE_{t})]$.
The computation of the eta form by Goette \cite{goette00}, Lemma 3.4,
and Remark 3.5, gives
$$\eta^{2m}(\cE_{t})= \frac{B_{m+1}}{(m+1)!} c_1(\omega)^{m}\ .$$
\hB

\subsubsection{}
Let us specialize to the case $m=1$. In this case
$\del(z_{2})$ corresponds to $\gerbe(\cE_{geom})$. Note that $B_2=1/6$.
Therefore, we obtain the following computation of the index gerbe:
$$\del(z_2)=a [\frac{1}{12} c_1(\omega)]\ .$$
Since $c_1(E)$ is even, we see on the one hand  that $6\:\gerbe(\cE_{geom})=0$.
On the other hand, there exists nontrivial index gerbes.
E.g. take $B=\C P^1$ and let $E\rightarrow B$ be the square of the Hopf bundle.
Then $c_1(E)=2$, and $\del(z_2)\cong [1/6]_{\R/\Z}$ under
the isomorphism $H_{Del}^3(\C P^1)\cong \R/\Z$.


\printindex

\bibliographystyle{plain}

\end{document}